\magnification=1000
\documentstyle{amsppt}
\parindent 20 pt
\baselineskip 20pt
\define \sgn {\text {sgn}}
\define \supp{\text {supp}}
\define \vdelt {\overline \Delta}
\define \ndelt {\underline \Delta}
\define \p {\partial}
\define \bul{\bullet}

\define \hK {$\hat K$}
\define \Fr {\bold {Fr}}
\define \Di {\Cal D}
\define \Rm {\Bbb R^m}
\define \lh {\underline h}
\define \Deltax0{\Delta_{ \bold x ^0}}
\define \Om {\Omega}
\define \eps {\epsilon}
\define \mes{\roman {mes}}

\define \a{\alpha}
\define \be{\beta}
\define \Dl{\Delta}
\define \dl{\delta}
\define \g{\gamma}
\define \G{\Gamma}
\define \lm{\lambda}
\define \Lm{\Lambda}
\define \k{\kappa}
\define \om{\omega}
\define \r{\rho}
\define \s{\sigma}

\define \th{\theta}

\define \z{\zeta}

\define \iy{\infty}

\define \sm{\setminus}

\define \ri{\rightarrow}

\define \sbt{\subset}
\define \spt{\supset}

\define \edm{\enddemo}
\define \ep{\endproclaim}

\define \1{^{-1}}
\define \2{^{-2}}

\define \BC{\Bbb C}

\define \BG{\Bbb G}

\define \BN{\Bbb N}

\define \BR{\Bbb R}
\define \BT{\Bbb T}
\define \BZ{\Bbb Z}

\define \CG{\Cal G}
\define \CH{\Cal H}

\define \CM{\Cal M}

\define \var{\operatorname{Var\ }}

\define \liml{\lim\limits}
\define \limsupl{\limsup\limits}
\define \liminfl{\liminf\limits}
\define \intl{\int\limits}

\define \suml{\sum\limits}

\def\eps{\epsilon}

\define \clos{\text {clos }}
\define \grad {\operatorname {grad}}

\topmatter
\pageno=-1
\heading{\bf Growth Theory of Subharmonic Functions}\endheading
 \vskip .10in
\centerline{Vladimir Azarin}
\vskip .15in
\NoBlackBoxes
\abstract
In this course of lectures we give an account of the growth theory of subharmonic
functions, which is directed towards its applications to entire functions of one and
several complex variables.
\endabstract

\toc

\head{1}    Preface \page{1} \endhead

\head{2}    Auxiliary information. Subharmonic functions \page{3} \endhead
\subhead{2.1} Semicontinuous functions \page{3} \endsubhead
\subhead{2.2} Measures and integrals \page{8} \endsubhead
\subhead{2.3} Distributions \page{16} \endsubhead
\subhead{2.4} Harmonic functions \page{25} \endsubhead
\subhead{2.5} Potentials and capacities \page{34} \endsubhead
\subhead{2.6} Subharmonic functions \page{44} \endsubhead
\subhead{2.7} Sequences of subharmonic functions \page{54} \endsubhead
\subhead{2.8} Scale of growth. Growth characteristics of subharmonic functions \page{63} \endsubhead
\subhead{2.9} Representation theorem of subharmonic functions in ${\BR}^m$ \page{78} \endsubhead

\head{3}    Asymptotic behavior of subharmonic functions of finite order \page{93} \endhead
\subhead{3.1} Limit sets \page{93} \endsubhead
\subhead{3.2} Indicators \page{108} \endsubhead
\subhead{3.3} Densities  \page{124} \endsubhead

\head{4}    Structure of the limit sets \page{134} \endhead
\subhead{4.1} Dynamical systems \page{134} \endsubhead
\subhead{4.2} Subharmonic functions with prescribed limit set \page{152} \endsubhead
\subhead{4.3} Further properties of the limit sets \page{168} \endsubhead
\subhead{4.4} Subharmonic curves. Curves with prescribed limit sets \page{186} \endsubhead

\head{5}    Applications to entire  functions \page{189} \endhead
\subhead{5.1} Growth characteristics of entire and meromorphic functions \page{189} \endsubhead
\subhead{5.2} $\Di '$ -topology and topology of exceptional sets \page{191} \endsubhead
\subhead{5.3} Asymptotic approximation of subharmonic functions \page{197} \endsubhead
\subhead{5.4} Lower indicator of A.A.Gol$'$dberg. Description of lower indicator. Description
of pair: indicator -- lower indicator \page{204} \endsubhead
\subhead{5.5} Asymptotic extremal problems. Semiadditive integral \page{218} \endsubhead
\subhead{5.6} Entire functions of completely regular growth. Levin-Pfluger Theorem.
Balashov's theory \page{223}  \endsubhead
\subhead{5.7} General characteristics of growth of entire functions \page{227} \endsubhead
\subhead{5.8} A generalization of the Valiron-Titchmarsh Theorem \page{251} \endsubhead

\head{6}    Application to the completeness of exponential systems in a convex
domain and the multiplicator problem \page{256} \endhead
\subhead{6.1} Problem  of multiplicator \page{256} \endsubhead
\subhead{6.2} A generalization of $\r$-trigonometric convexity \page{281} \endsubhead
\subhead{6.3} Completeness of exponential systems in a convex domain \page{288} \endsubhead

\head{}    Notation \page{313} \endhead

\head{}    List of terms \page{315}  \endhead

\head{}    References \page{321} \endhead

\endtoc

\endtopmatter

\newpage

\pageno=1

\centerline {\bf 1.Preface}

This book aims to convert the noble art of constructing an entire
function with prescribed asymptotic behavior to a handicraft. This
is the aim of every Theory.

For this you should only construct the limit set that describes
the asymptotic behavior of the entire function, i.e., you should
consider the set $U[\r,\s]$ of subharmonic functions (that is,
$\{\text{$v$ is  subharmonic} :v(re^{i\phi })\leq \s r^\rho \}$)
and pick out the subset $U$ which characterizes its asymptotic
properties.

How to do it? The properties of limit sets are listed in \S3. All
the standard growth characteristics are expressed in terms of
limit sets in \S\S3.2,3.3,5.7. Examples of construction are to be found in
\S\S 5.4-6.3.

So you can use this book as a reference book for construction of
entire functions.

Of course, you need some terms. All the terms are listed on page
315-320. That is it.

If you want to study the theory, I recommend to solve the
Exercises that are in the text. Most of them are trivial. However,
I recommend to do all of them by the moment that they appear
trivial to you.

A few words about the history of this book. It arose from a course
of lectures that I gave in the Kharkov University in 1977. After
some time, under the pressure and with active help of
Prof.I.V.Ostrovskii, a rotaprint edition (Edition of KhGU)of this
course appeared: the first part in 1978, the second one in 1982.
Math.\ Review did not notice this fact.

From that time lots of new and important results were obtained. Part
of them were presented in Chapter 3 of the review \cite{GLO}.

In 1994, when I started to work in the Bar-Ilan University and
obtained a personal computer, my first wish was to study typing on
it in English. This was the first impulse for translating this
course to English (there are no more than 5 copies of this book in
the world, I believe, one of them being mine). I continued this
project while working in Bar-Ilan (1994-2006) but there was no
much time for this. And now I have finished.

\demo {Acknowledgements} I am indebted to many people. I start from
Prof. I.V.Ostrovskii, who supported this idea for many years,  and
Prof. A.A.Goldberg, who stimulated my mathematical activity all my
life by his letters and conversations.

I am indebted to Prof's A.Eremenko and M.Sodin, who, not being my
``aspirants,'' solved a lot of problems connected to limit sets,
and also to Dr's V.Giner, L.Podoshev and E.Fainberg that worked
with me to develop the theory.

I am indebted to Prof's L.H\"ormander and R.Sigurdsson that have
sent me the preprints of their papers that were not yet published.
I am indebted to Prof.I.F.Krasichkov-Ternovskii, who explained me
many years ago the connection between the multiplicator problem
and completeness of the exponent system in a convex domain.

I am indebted to Prof's.\ M.I.Kadec and V.P.Fonf for proving
Theorem 4.1.5.2 which is rather far off my speciality.

I am indebted to my coauthors Prof's D.Drasin and P.Poggi-Corradini;
I have exploited the results of our joint paper in \S6.2.

Of course, I am indebted to my late teacher Prof.B.Ya.Levin, who
has taught me entire and subharmonic functions and gave me the
first problems in this area. Actually, the theory of limit sets is
a generalization of the theory of functions of completely regular
growth.

I am also indebted very much to my grandson Sasha Sodin, who
transformed ``my English'' into English. \edm

\newpage
\centerline {\bf 2.Auxiliary information.Subharmonic functions}

\centerline {\bf 2.1.Semicontinuous functions}

\subheading {2.1.1} Let $x \in \Bbb R ^m $ be a point in an $m$-dimensional
Euclidean space, $E$ a Borel set and $f(x)$ a function on $E$ such that
$f(x) \neq \infty$ .

Set
$$M(f,x,\varepsilon):=\sup \{f(x'):|x-x'|< \varepsilon  ,\ x' \in E \}
\tag 2.1.1.1$$
The function
$$f^* (x):=\lim \limits_{\varepsilon \rightarrow 0} M(f,x,\varepsilon)$$
is called the {\it upper semicontinuous regularization} of the function
$f(x)$.

In the case of a finite jump, the regularization "raises" the values
of the function. However, there is no influence on $f^* (x)$, if $f(x)$
tends to $-\infty$ "continuously".
\proclaim {Proposition 2.1.1.1(Regularization Properties)} The following
properties hold:
$$f(x) \leq f^* (x);\tag rg1$$
$$(\alpha f)^* (x)=\alpha f^* (x);\tag rg2 $$
$$(f^* )^* (x)=f^* (x);\tag rg3 $$
$$\align
&(f_1 +f_2)^* (x) \leq f_1 ^* (x) + f_2 ^* (x);\\
&(\max (f_1,f_2))^* (x) \leq \max (f_1 ^* ,f_2 ^* ) (x); \tag rg 4\\
&(\min (f_1,f_2)) ^* (x) \leq \min (f_1 ^* ,f_2 ^* ) (x);
\endalign$$
\endproclaim

These properties are obvious corollaries of the definition of $f^* (x)$.

{\bf Exercise 2.1.1.1} Prove them.

\subheading {2.1.2} The function $f(x)$ is called {\it upper semicontinuous}
at a point $x$ if $f^* (x)=f(x)$.

We denote the class of upper semicontinuous functions on $E$ by $C^+
(E)$.

The function $f(x)$ is called {\it lower semicontinuous} if $-f(x)$ is upper
semicontinuous (notation $f \in C^- (E)$ ).
 Examples of semicontinuous functions are given by
\proclaim {Proposition 2.1.2.1(Semicontinuity of Characteristic Functions of
Sets)}

Let $G \subset \Bbb R ^m$ be an open set.Then its characteristic function
$\chi _G$ is lower
semicontinuous  in $\Bbb R ^m$. Let $F$ be a closed set, then
$\chi_F$ is
upper semicontinuous .
\endproclaim

The proof is obvious.

{\bf Exercise 2.1.2.1} Prove this.

\proclaim {Proposition 2.1.2.2(Connection with Continuity)} If $f\in C^+ \cap
C^-$, then $f$ is continuous.
\endproclaim

The assertion follows from the equalities
$$f^* (x)=\limsup \limits_{\varepsilon \rightarrow 0}
\{f(x'):|x-x'|<\varepsilon\};
-(-f)^* (x)=\liminf \limits_{\varepsilon \rightarrow 0}
\{f(x'):|x-x'|<\varepsilon\}.$$
\proclaim {Proposition 2.1.2.3($C^+$ -Properties)} The following  holds
$$f \in C^+ (E) \Rightarrow \alpha f \in C^+ (E),\ for\ \alpha \geq 0
\tag $C^+$ 1$$
$$f_1,f_2 \in C^+ \Rightarrow f_1 +f_2 ,\max (f_1,f_2),  \min (f_1,f_2)\ \in
 C^+ .$$
\endproclaim

These properties follow from the properties of regularization (Prop.2.1.1.1).

{\bf Exercise 2.1.2.2} Prove them.

Let $G$ be an open set. Set $G_A:=\{x \in G:f(x)<A\}$.

\proclaim {Theorem 2.1.2.4 (First Criterion of Semicontinuity)}One has
$f \in C^+ $ if and only if $G_A$ is open for all $A\in \BR$.
\endproclaim
\demo {Proof}Let $f(x)=f^* (x),\ x\in G$. Then
$\{f(x)<A\} \Longrightarrow \{f^* (x)<A\} \Longrightarrow \{M(f,x,
\varepsilon )
<A\}$ for all sufficiently small $\varepsilon. $
Thus the neighborhood of $x\ \mathbreak \ V_{\varepsilon ,x}:=
 \{x':|x-x'|<\varepsilon\}$ is contained in $G_A$.

Conversely, since the set $G_A$ is open for $A=f(x_0)+ \delta$, we have
\linebreak
$f^* (x_0)\leq f(x_0)+\delta$ for any $\delta >0$, hence for $\delta=0$. With
property (rg1) of Prop.2.1.1.1 this gives $f^*(x_0)=f(x_0).\ \ $\qed
\enddemo
Let $F$ be a closed set. Set $F^A:=\{x \in F:f(x)\geq A \}$.
An obvious corollary of the previous theorem is
\proclaim {Corollary 2.1.2.5}One has
$f \in C^+$ if and only if $F^A$ is closed for all $A$.
\endproclaim

{\bf Exercise 2.1.2.3} Prove the corollary.

 We denote compacts by $K.$ Set $M(f,K)=\sup \{f(x):x\in K\}$.
\proclaim {Theorem 2.1.2.6(Weierstrass)}Let $K\subset \Bbb R^m$ be a
compact set and $f\in C^+(K).$

Then there exists $x_0 \in K$such that $f(x_0)=M(f,K),$
\endproclaim
i.e., $f$ attains its supremum on any compact set.
\demo {Proof} Set $K_n :=\{x \in K:f(x) \geq M(f,K)-1/n \}.$

 The $K_n $ are closed by Cor.2.1.2.5, nonempty by definition of $M(f,K).$
 Their intersection is nonempty and is equal to the
set
$$K_{\max}:=\{x\in K:f(x) \geq M(f,K)\}. $$

It means that there exists $x_0$ in $K$ such that $ f(x_0)\geq M(f,K)$.

The opposite inequality holds for any $x$ in $K$. \ \ \qed
\enddemo

{\bf Exercise 2.1.2.4} Why?

The following theorem shows that the functional $M(f,K)$ is continuous with
respect to monotonic convergence of semicontinuous functions.
\proclaim {Proposition 2.1.2.7(Continuity from the right of $M(f,K)$)}

Let  $f_n\in C^+ (K),\ f_n\downarrow f,\ n=1,2,3...\ \ .$

Then $M(f_n,K) \downarrow M(f,K).$
\endproclaim

\demo {Proof} It is clear that
$\lim \limits_{n\rightarrow \infty} M(f_n,K):=M$ exists.

Set $K_n :=\{x\in K: f_n (x)\geq M\}.$ The intersection of the closed
nonempty sets $K_n$ is nonempty and has the following form:
$\underset n\to{\bigcap} K_n =\{x:f(x)\geq M\} .$
So $M(f,K)\geq M.$

The opposite inequality is obvious.\qed
\enddemo

{\bf Exercise 2.1.2.5} Why?

In the same way one  prove
\proclaim {Proposition 2.1.2.8(Commutativity of $inf$ and $M(\cdot)$)}Let
$\{f_{\alpha} \in C_+ (K) , \alpha \in (0;\infty) \}$ be an
arbitrary
decreasing family of semicontinuous functions.Then $\inf\limits_\alpha M(f_\alpha ,K) = M(\inf\limits_\alpha f_\alpha ,K).$
\endproclaim

{\bf Exercise 2.1.2.6} Prove this proposition.

\proclaim {Theorem 2.1.2.9(Second Criterion of Semicontinuity)}
 $f\in C^+ (K)$ iff there exists a  sequence $f_n$ of continuous functions
such that $f_n \downarrow f.$
\endproclaim
\demo {Sufficiency} Let $ f_n \in C^+ (K),\  f_n\downarrow f$.Set
$K_n ^A :=\{x\in K:f_n (x) \geq A\}$. This is a sequence of non-empty closed sets.
If the set $K^A :=\{x:f(x)\geq A\}$ is nonempty, then $K^A$ is closed because
$\underset n\to{\bigcap} K_n ^A =K^A $. Hence $f\in C^+ (K)$ by Corollary 2.1.2.5.
\enddemo
\demo {Necessity} Set $f_n (x,y) :=f(y)-n|x-y|.$

This sequence of functions has the following properties:

a) it decreases monotonically in $n$ and
$$\lim\limits_{n\rightarrow \infty} f_n (x,y) =\cases
f(x),\ \ for\ x=y;\\ -\infty ,\ \ for\ x\neq y;\endcases$$

b) for any fixed $n$ the functions $f_n$ are continuous in $x$ uniformly
with respect to $y$, because
$|f_n (x,y)-f_n (x',y)|\leq n|x-x'|;$

c) $f_n$ are upper semicontinuous in $y$.

  Prop.2.1.2.7 and c) imply that
$\lim\limits_{n\rightarrow \infty} M_y (f_n (x,y),K) =
M_y (\lim\limits_{n\rightarrow \infty} f_n (x,y),K)$

b) implies that the functions
$f_n (x):=M_y(f_n (x,y),K)$ are continuous, and a) implies that they decrease
monotonically to $f(x)$ .\ \qed
\enddemo
 \subheading {2.1.3} We will consider a family of upper
semicontinuous functions:$\{f_t:\ t\in T\subset (0,\infty)\}.$
It is easy to prove
\proclaim{Proposition 2.1.3.1}$f_t \in C^+ \Longrightarrow
\underset{t\in T}\to{\inf} f_t (x) \in C^+.$
\endproclaim

{\bf Exercise 2.1.3.1} Prove this Proposition.

Set $f_T (x):=\underset{t\in T}\to{\sup} f_t (x) .$ The function $f_T$ is
not, generally speaking, upper semicontinuous even if $T$ is countable and
$f_t$ are continuous. It is not possible to replace $\underset T\to{sup}$ in
the definition of $f_T$ by $\underset {T_0}\to{sup} ,$ where $T_0$ is a
countable set .However, the following theorem holds:
\proclaim{Theorem 2.1.3.2(Choquet's Lemma)}There exists such a countable set
$T_0\subset T$, that
$$(\underset{T_0}\to{\sup} f_t)^* (x)=(\underset{T}\to{\sup} f_t)^* (x).$$
\endproclaim
\demo{Proof}Let $\{x_n\}$ be a countable set that is dense in ${\Bbb R}^
m$ and
${\varepsilon}_j \downarrow 0$.Then the balls
$$K_{n,j}:=\{x:|x-x_n |<\varepsilon _j \}$$
cover every point $x\in {\Bbb R}^m $ infinitely many times.

Renumbering we obtain a sequence $\{K_l:l\in \Bbb N\}.$
For any $l$ there exists, by definition of $\underset{K_l}\to{sup}$, such a
point $x_0\in K_l$ that
$$\underset{K_l}\to{\sup} f_T (x)\leq f_T ( x_0 )+ 1/2l .\tag 2.1.3.1$$
By definition of $\underset T\to{sup}$ there exists $t_l$ such that
$$f_T (x_0) < f_{t_l} (x_0)+1/2l.$$
Thus
$$f_T (x_0 )< \sup\{f_{t_l} (x):x\in K_l\} +1/2l.\tag 2.1.3.2$$
The inequalities (2.1.3.1) and (2.1.3.2) imply that for any $l$ there exists
 $t_l$ such that
$$\sup\{ f_T (x): x\in K_l\}\leq \sup\{f_{t_l} (x): x\in K_l \}+1/l.\tag
2.1.3.3$$
Now set $T_0 =\{t_l\}$. Evidently, $f_{T_0} (x) \leq f_T (x)$ and thus
$$f_{T_0}^* (x)\leq f_T^* (x).\tag 2.1.3.4 $$

Let us prove the opposite inequality.

Let $x\in {\Bbb R}^m$. Choose a subsequence $\{K_{l_j}\}$ that tends to $x$.
From (2.1.3.3) we obtain
$$f^*_T (x)\leq \limsup\limits_{j\rightarrow \infty} \underset{x'\in K_{l_j}}
\to{\sup} f_T (x') \leq$$
$$\limsup\limits_{j\rightarrow \infty} \underset{x'\in K_{l_j}}
\to{\sup} f_{t_{l_j}} (x')\leq$$
$$ \limsup\limits_{j\rightarrow \infty} \underset {x'\in K_{l_j}}
\to{\sup} f_{T_0} (x')=f^*_{T_0}(x).\tag 2.1.3.5 $$
(2.1.3.4) and (2.1.3.5) imply the assertion of the theorem.  \qed
\enddemo
\newpage
\centerline{\bf 2.2.Measures and integrals}
\subheading{2.2.1}Let $G$ be an open set in ${\Bbb R}^m$ and  $\sigma
(G)$  a $\sigma$-algebra of Borel sets containing all the compact sets
$K\subset
G$.

Let $\mu$ be a countably additive nonnegative function on $\sigma (G)$,
which is finite on all  $K\subset G$. We will call it  a {\it measure} or a
{\it mass distribution}.

Let $G_0 (\mu)$ be the largest open set for which $\mu$ is zero. It is
the union of all the open sets $G'$ that $\mu (G') =0.$

 The set supp $ \mu :=G\backslash G_0 (\mu)$ is called
the {\it support} of $\mu$. It is closed in $G.$

We say that $\mu$ is {\it concentrated} on $E\in \sigma (G)$
if $\mu (G\backslash E) =0$.

\proclaim{ Theorem 2.2.1.1(Support)} The support of a measure $\mu$ is the
smallest closed set on which
the measure $\mu$ is concentrated.
\endproclaim

{\bf Exercise 2.2.1.1} Prove this.

 A measure $\mu$ can be concentrated on a non-closed set $E$ and then
$E\Subset\roman {supp}\ \mu.$
\subheading{Example 2.2.1.1}Let $E$ be a countable set dense in $G$. Then
$\roman{supp}\ \mu =G$ and , of course, $E\neq G$.

The set of all measures on $G$ will be denoted a $\Cal M (G).$

The measure $\mu _F (E):=\mu (E\cap F)$ is called the {\it restriction} of $\mu$ onto
$F\in \sigma (G)$.
It is easy to see that $\mu_F$ is concentrated on $F$ and supp$\mu \subset
\overline{F}.$

A countably additive function $\nu $ on $\sigma (G)$   that is finite for  all
$K \subset G$  is called a {\it charge}. We consider only real charges.

{\bf Example 2.2.1.2} $\nu :=\mu _1 -\mu _2,\ \ \mu _1 ,\mu _2\ \in
\Cal M (G).$

The set of all the charges will be denoted  ${\Cal M}^d.$
\proclaim {Theorem 2.2.1.2(Jordan decomposition)}Let $\nu \in \Cal M^d (G).$
Then there exist two sets $G^+,G^-$ such that

a) $G=G^+\cup G^-,\ G^+\cap G^-=\varnothing;$

b) $\nu(E)\geq 0\ for\ E \subset G^+;\ \ \nu (E) \leq 0\ for\ E\subset G^-.$
\endproclaim
  One can find the proof in \cite {Ha,Ch.VI Sec.29}

 The measures $\nu ^+:=\nu _{G^+}$ and $\nu ^-:=\nu _{G^-},$ where
$\nu _{G^+},\  \nu _{G^-}$ are restrictions of $\nu$ to $G^+,G^-,$ are called
the {\it positive} and {\it negative}, respectively, {\it variations} of $\nu$. The measure
$|\nu|:=\nu _+ + \nu _-$ is called the {\it full variation} of $\nu$ or just a
{\it variation}.
\proclaim {Theorem 2.2.1.3(Variations)}The following holds:
$$\nu ^+ (E)= \underset {E'\subset E}\to{\sup} \nu (E');
\ \nu ^- (E)= \underset {E'\subset E}\to{\inf} \nu (E');
\ \nu =\nu^+ + \nu^-.$$
\endproclaim
The proof is easy enough.

{\bf Exercise 2.2.1.2} Prove this.

{\bf Example 2.2.1.3} Let $\psi (x) $ be a locally summable function
 with respect to the Lebesgue's measure. Set $\nu (E):=\int_E \psi (x)dx$.

Then
$$\nu ^+ (E)=\int_E \psi ^+ (x)dx,\ \nu ^- (E)=\int_E \psi ^- (x)dx;\
|\nu| (E)=\int_E |\psi| (x)dx,$$
where
$$ \psi ^+ (x)=\max (0,\psi (x));\ \psi ^- (x)=-\min (0,\psi (x)).
\tag 2.2.1.1$$
\subheading{2.2.2}The function $f(x),\ x\in G $ is called a {\it Borel\
function} if the set \linebreak $E^A :=\{f(x)>A\}$ belongs to $\sigma (G)$
for any $A$ in $\Bbb R.$

Let $K \Subset G $ be a compact set and $f$ a Borel function. Then
 the Lebesgue-Stieltjes integrals of the form
 $\int_K f^+ d\mu,\ \int_K f^- d\mu$ with respect to a measure $\mu \in \Cal M (G)$ are
defined, and $ \int_K f d\mu:=\int_K f^+ d\mu-\int_K f^- d\mu$ is defined
 if at least one of the terms is finite.

We say that a property holds $\mu -almost\ everywhere\ on\ E$ if the set  $E_0$
 of  $x$ for which it does not hold satisfies the condition $\mu (E_0) =0.$

We will denote all the compact sets  in $G$ as $K$ (sometimes with indexes).
The following theorems hold:
\proclaim{Theorem 2.2.2.1(Lebesgue)}Let $\{f_n,\ n\in \Bbb N$\}  be  a
sequence of Borel functions on   $K$ and $g(x)\geq 0$
a function on $K$,that is summable with respect to  $\mu$ (i.e., its integral is finite),
 $|f_n (x)|\leq g(x)$ for $x$ in $K,$ and $f_n \rightarrow f$ when
$n\rightarrow\infty$.

Then $ \lim\limits_{n\rightarrow\infty} \int_K f_n d\mu =\int_K f d\mu.$
\endproclaim
\proclaim {Theorem 2.2.2.2(B.Levy)}Let $f_n\downarrow f$ when $n
\rightarrow\infty$, and $f$ be a summable function on $K$.

Then $ \lim\limits_{n\rightarrow\infty} \int_K f_n d\mu =\int_K f d\mu.$

\endproclaim
\proclaim{Theorem 2.2.2.3(Fatou's Lemma)}Let $f_n (x)\leq const <\infty$ for
$x$ in $K$.

Then $ \limsup\limits_{n\rightarrow\infty} \int_K f_n d\mu \leq \int_K
\limsup\limits_{n\rightarrow\infty} f_n d\mu.$
\endproclaim
The proofs can be found in \cite {Ha,Ch.V, Sec.27}.

Let $L(\mu)$ be the  space of functions that are summable with respect to  $\mu$.
We say that $f_n\rightarrow f$ in  $L(\mu)$ if  $f_n,f\in L(\mu)$ and
$$\|f_n-f\|:=\int |f_n-f|(x)d\mu\rightarrow 0$$
\proclaim {Theorem 2.2.2.4(Uniqueness in $L(\mu)$)} Let
$f_n\rightarrow f$ in  $L(\mu)$ and
$$\int f_n\psi d\mu \rightarrow \int g\psi d\mu$$
for any $\psi$ continuous on supp$\mu.$ Then $\|g-f\|=0.$
\endproclaim
For proof see, e.g, \cite {H\"o,Th.1.2.5 }.

\subheading{2.2.3} Let $\phi (x)$ be a Borel function on $G$. The set
$\roman{supp}\ \phi:=\overline {\{x:\phi(x)\neq 0\}}$ is called the
$support\ of\
\phi (x)$. A function $\phi$ is called
$finite\ in\ G$ if $\roman{supp}\ \phi\Subset G$

We say that  a sequence $\mu_n \in \Cal M$ converges $weakly$ to $\mu \in
\Cal M$ if the condition
 $\int \phi d\mu _n \rightarrow \int \phi d\mu $ holds for any continuous function $\phi$.

We will not show the integration domain, because it is always supp $\phi .$

The weak (it is called also $C^*$-) convergence will be denoted as
$\overset{*}\to\ri$.

\proclaim {Theorem 2.2.3.1($C^*$-limits)} If $\mu_n \overset*\to\ri \mu $, then for $E\in \sigma
(G)$ the following assertions hold
$$\underset n\rightarrow \infty \to{\limsup}\mu_n (\overline E )\leq
\mu (\overline E);$$
$$\underset n\rightarrow \infty \to{\liminf}\mu_n (\overset \circ \to E )
\geq \mu (\overset \circ \to E);$$
where $\overset \circ \to E$ is the interior  of $E$, $\overline E$ is the closure
of $E.$
\endproclaim
\demo {Proof}Let $\chi _{\overline E}$ be the characteristic function of
the
set $\overline E. $ It is upper semicontinuous. Thus there exists  a
decreasing sequence $\varphi _m $ of continuous functions finite in $G$
that converges to $\chi _{\overline E}$ as   $m\rightarrow \infty$.
Then we have
$$\mu_n (\overline E)=\int \chi _{\overline E}d \mu_n \leq \int \varphi _m
d \mu_n .$$
Passing  to the limit as $n\rightarrow\infty $ we obtain
$$\underset n\rightarrow \infty \to{\limsup} \mu_n (\overline E )\leq
\int \varphi _m d \mu .$$

Passing to the limit as $m\rightarrow\infty$ we obtain by Th.2.2.2.2
$$\underset n\rightarrow \infty \to{\limsup}\mu_n (\overline E )\leq
\int \chi _{\overline E}d \mu =\mu (\overline E).$$
The proof for $\overset \circ \to E$  is analogous.\qed
\enddemo
\proclaim {Theorem 2.2.3.2.(Helly)} Let $\{\mu_{\alpha} : \alpha \in A\}$
be a family of
measures uniformly bounded  on any compact  set $K\subset G$, i.e.,
$\exists C=C(K):\mu_\alpha (K)\leq C(K),\ for\ K\Subset G.$

Then this family is $weakly\ compact$, i.e., there exists a sequence
$\{\alpha_j :\alpha_j \in A \}$ and a measure $ \mu $ such that
$\mu_{\alpha_j} \overset*\to\ri \mu.$
\endproclaim
The proof can be found in \cite {Ha} .

A set $E$ is called $squarable$ with respect to measure $\mu$
($\mu$-squarable  ) if \linebreak $\mu (\partial E) =0.$

\proclaim {Theorem 2.2.3.3.(Squarable Ring)} The following holds

 sqr1) if $E_1,E_2$ are $\mu$ -squarable, the sets $E_1\cap E_2,E_1\cup E_2,
E_1\backslash E_2$ are $\mu$ -squarable;

 sqr2) for any couple: an open set $G$ and a compact set $K\subset G$ there
exists a $\mu$-squarable set $E$ such that $K\subset E \subset G$.
\endproclaim
\demo{Proof}The assertion sqr1) follows from
$$\partial (E_1\cup E_2) \bigcup \partial (E_1\cap E_2) \bigcup
\partial (E_1 \backslash E_2) \subset \partial E_1 \cup \partial E_2 .$$

Let us prove sqr2).
Let $K_t:=\{x:\exists y\in K: |x-y|<t\}$ be a $t$-neighborhood of the $K.$
It is clear that for all the small $t$ we have
$K\Subset K_t \Subset G.$
The function $a(t):=\mu(K_t)$ is monotonic on $t$ and thus has no more than
a countable set of jumps.

Let $t$ be a point of continuity of $a(t).$Then
$$\mu(\partial K _t)\leq \lim\limits_{\epsilon \rightarrow 0}
[\mu (K_{t+\epsilon}) -\mu(K_{t-\epsilon})] =0.$$
Thus it is possible to set $E:=K_t$ for this $t$.\qed
\enddemo

 A family $\Phi$ of sets is called a $dense\ ring$ if the following
conditions hold:

   dr1) $\forall F_1,F_2 \in \Phi \Longrightarrow F_1\cup F_2 ,F_1\cap F_2
\in \Phi;$

   dr2) $\forall K,G : K\Subset G \ \exists F\in \Phi :K\subset F
\subset G.$

The previous theorem shows that the class of $\mu$-squarable sets is a dense
ring. The following theorem shows how one can extend a measure from a dense
ring the Borel's algebra.

Let $\Phi$ be a dense ring and $\Delta (F),\ F\in \Phi$ a function of a set
which satisfies the conditions

$\Delta$1) $monotonicity$ on $\Phi$: $F_1\subset F_2 \Longrightarrow
\Delta (F_1) \leq   \Delta (F_2);$

$\Delta$2) $additivity$ on $\Phi$: $\Delta (F_1 \cup F_2) \leq \Delta (F_1) +
\Delta (F_2)$

and

$\Delta (F_1 \cup F_2) = \Delta (F_1) + \Delta (F_2)$ if
$F_1 \cap F_2 =\varnothing$

$\Delta$3) $continuity$ on $\Phi$: $\forall F\in \Phi$ and $ \epsilon >0$
there exists a
compact set $K$ and an open set $G\supset K$ such that $\forall F' \in \Phi :
K\subset F'\subset G$ the inequality $|\Delta (F) -\Delta (F')|<\epsilon$
holds.
\proclaim {Theorem 2.2.3.4.(N.Bourbaki)}There exists  a measure $\mu$ such that
$\mu (F)=\Delta (F),\ \forall F\in \Phi$ iff the conditions $\Delta$1) - $\Delta$3) hold.
\endproclaim
\proclaim {Theorem 2.2.3.5.(Uniqueness of Measure)} Under the conditions
$\Delta$1) - $\Delta$3) the measure is defined uniquely by  the formulae:
$$\mu (K)=\inf \{\Delta (F): F\in \Phi,\ F\supset K\};\tag 2.2.3.1$$
$$\mu (G)=\sup \{\Delta (F): F\in \Phi,\ F\subset G\};\tag 2.2.3.2$$
$$\mu (E)=\sup \{\mu (K):K\subset E\} =\inf \{\mu (G):G\supset E\},
\tag 2.2.3.3 $$
and every $F\in \Phi$ is $\mu$-squarable.
\endproclaim

For proof see \cite {Bo,Ch.4,Sec 3,it.10}. The squarability follows from  (2.2.3.3).

The following theorem connects the convergence of measures on any dense ring
and on the ring of sets squarable with respect to the limit measure.
\proclaim {Theorem 2.2.3.6.(Set-convergences)} If $\mu_n (F) \rightarrow
\mu (F)$ for all
$F$ in a dense ring $\Phi$, then $\mu_n (E) \rightarrow \mu (E) $ for any
$\mu$-squarable set $E.$
\endproclaim
\demo {Proof} Suppose $\overset \circ \to E \neq \varnothing.$

Let $\epsilon >0$ .By (2.2.3.3) one can find a compact set $K$  such that
$$\mu (K)+\epsilon \geq \mu (\overset \circ \to E )=\mu (E).\tag 2.2.3.4$$
One can also find an open set $G$ such that
$$\mu (G) -\epsilon \leq \mu (\overline E) =\mu (E).\tag 2.2.3.5$$
By property dr2) of a dense ring one can find $F,F'\in \Phi$ such that
$$K\subset F \subset \overset \circ \to E \subset E \subset \overline E
\subset F'\subset G .$$
Thus
$\mu_n (F) \leq \mu_n (E) \leq \mu_n (F')$ and hence
$$\mu (F) \leq \underset {n\rightarrow \infty }\to {\underline \lim}\mu_n (E)
\leq \underset {n\rightarrow \infty }\to {\overline \lim}\mu_n (E) \leq
 \mu (F').\tag 2.2.3.6 $$
From (2.2.3.4) and (2.2.3.5) we obtain
 $0\leq \mu (F') -\mu (F) \leq \mu (G)-\mu (K)\leq 2\epsilon $
for arbitrary small $\epsilon.$ Thus from (2.2.3.6) we obtain
$$\underset {n\rightarrow \infty }\to {\underline \lim}\mu_n (E) =
\underset {n\rightarrow \infty }\to {\overline \lim}\mu_n (E) =\mu (E).\tag
2.2.3.7 $$

That is to say that $\mu_n (E)\rightarrow \mu (E).$

If $\overset \circ \to E =\varnothing,$ then $\mu (\overline E) =0$ by
the definition of  a squarable set. One can show in the same way that
$\mu_n (E) \rightarrow 0.$ \qed
\enddemo

Now we connect the weak convergence to the convergence on squarable sets.
\proclaim {Theorem 2.2.3.7.(Set- and C*-convergences)} The conditions
$$\mu_n \overset*\to\ri \mu \tag 2.2.3.8 $$
and $\mu_n (E) \rightarrow \mu (E)$ on $\mu$-squarable sets $E$ are
equivalent.
\endproclaim
\demo{Proof} Sufficiency of (2.2.3.8) follows from Th.2.2.3.1.

{\bf Exercise 2.2.3.1} Prove this.

Let us prove necessity.

For any compact set one can find a $\mu$-squarable $E$ such that $K\subset E.$
Hence $\mu_n (K)\leq \mu (E)+1:=C(K) $ when $n$ is big enough.

By Helly's theorem (Th.2.2.3.2) there exists a measure $\mu' $ and a
subsequence $\mu_{n_j} \overset*\to\ri \mu'$.By the proved necessity $\mu'(E)=
\mu (E)$ on a dense ring of the squarable sets. Thus $\mu '=\mu$ by Uniqueness
theorem 2.2.3.5. And thus $\mu_n \overset*\to\ri \mu.$ \qed
\enddemo
Denote by
$$\mu_E(G):=\cases \mu(G\cap E)&\text { if }G\cap E\neq\varnothing\\0&\text { if }G\cap E=\varnothing \endcases$$
 the {\it restriction} of $\mu$ on the set $E.$
\proclaim {Corollary 2.2.3.8} Let $\mu_n \overset*\to\ri \mu$ and $E$ be a squarable set for $\mu.$
Then $(\mu_n)_E \overset*\to\ri (\mu)_E.$
\ep
 Indeed, if $E$ is a squarable set for $\mu$
it is a squarable set for $\mu_E$. So Th.2.2.3.7
implies the corollary.
\subheading {2.2.4} Let $\sigma ({\Bbb R}^{m_1} \times {\Bbb R}^{m_2})$ be
the $\sigma$-algebra of all the Borel sets,
$\Phi _i \subset \sigma ({\Bbb R}^{m_i}),\ i=1,2 ,$  be dense rings,
$\Phi :=\Phi_1 \otimes \Phi_2 \subset \sigma ({\Bbb R}^{m_1} \times
{\Bbb R}^{m_2})$ be a ring generated by all the sets of
form $F_1\times F_2,\ \ F_i\in \Phi_i .$
\proclaim {Theorem 2.2.4.1.(Product of Rings)}If $\Phi_i ,\ i=1,2 $ are
dense rings, then $\Phi_1 \otimes \Phi_2$ is a dense ring; if they consist of
squarable sets, then $\Phi$ consists of squarable sets.
\endproclaim
\demo{Proof} Let $K\subset G \subset {\Bbb R}^{m_1} \times {\Bbb R}^{m_2}$.
For every point $x\in K$ one can (evidently) find  $F_1\times F_2$ such that
$x\subset F_1\times F_2 \subset G.$ One can find a finite covering and obtain
a finite union F of sets of such form. Thus $F\in \Phi_1 \otimes \Phi_2 $
and
$F\subset G.$

The second assertion follows from the formula:
$$\partial (F_1\times F_2)= (\partial F_1 \times F_2) \cup( F_1 \times
\partial F_2).\ \qed $$
\enddemo

Let $\mu_i$ be a measure on $\sigma ({\Bbb R}^{m_i}),\ i=1,2,$ and
$\mu:=\mu_1\otimes \mu_2$  the {\it product of measures}, i.e., a measure on
$\sigma ({\Bbb R}^{m_1} \times {\Bbb R}^{m_2})$ such that
$\mu (E_1\times E_2)=\mu_1(E_1)\mu_2(E_2)$ for all $E_i\in
\sigma ({\Bbb R}^{m_i}) ,\ i=1,2.$

\proclaim {Theorem 2.2.4.2.(Product of Measures)} A measure
$\mu_1\otimes \mu_2$ is uniquely defined  by its values on
$\Phi_1 \otimes \Phi_2.$
\endproclaim
The assertion follows from Theorem 2.2.4.1 and Uniqueness Theorem 2.2.3.5.

\proclaim {Theorem 2.2.4.3.(Fubini) }Let $f(x_1,x_2)$ be a Borel function
on ${\Bbb R}^{m_1} \times {\Bbb R}^{m_2}.$ Then
$$\underset {{\Bbb R}^{m_1} \times {\Bbb R}^{m_2}} \to \int
f(x_1,x_2)d(\mu_1\otimes\mu_2) =
\underset {\Bbb R}^{m_1} \to \int d\mu_2\underset {\Bbb R}^{m_2}\to \int
f(x_1,x_2)d\mu_2 =
\underset {\Bbb R}^{m_2} \to \int d\mu_2\underset {\Bbb R}^{m_1}\to \int
f(x_1,x_2)d\mu_1 ,\tag 2.2.4.1$$
if at least one of parts of (2.2.4.1) is well defined.
\endproclaim
The proof can be found in \cite {Ha, Ch.VII, Sec.36}.
\newpage

\centerline {\bf 2.3.Distributions}
\subheading {2.3.1}Let us consider the set $\Cal D (G)$ of all infinitely
differentiable functions \linebreak$\varphi (x),\ x\in G\subset {\Bbb R}^m .$
\par
It is a {\it linear space} because for any constants $c_1 ,c_2 $
$$ \varphi_1 ,\varphi_2 \in \Cal D (G) \Longrightarrow c_1\varphi_1 +
c_2\varphi_2
\in \Cal D (G). \tag D1 $$
\par
It is a {\it topological} space with convergence  defined by
$$\varphi_n \overset {\Cal D} \to {\rightarrow} \varphi :\cases
a)\ \roman{supp}\ \varphi_n \subset K \Subset {\Bbb R}^m \\
for\ some\ compact\ K\\
and\\
b)\ \varphi_n \rightarrow \varphi \ uniformly\ on \ K\\
with\ all\ their\ derivatives.\endcases \tag D2 $$
\par
We consider some examples of functions $\varphi \in \Cal D.$
\par
Denote
$$ \alpha (t) =\cases
Ce^{-\frac {1}{1-t^2}},\ for\ t\in (-1;1)\\
0,\ \ \ \ \ \ for\ t \ \overline \in \ (-1;1).\endcases \tag 2.3.1.1 $$
Evidently $\alpha (|x|) \in \Cal D (\Bbb R^m)$ and  supp $\a\sbt \{x:|x|\leq 1 \}.$

 {\bf Exercise 2.3.1.1} Check this.

Let us find $C$ such that
$$\int \alpha (|x|) dx =\sigma_m\int_0^1 \alpha (t)t^{m-1}dt=1
\tag 2.3.1.2 $$
where $\sigma_m$ is  area of the unit sphere $\{|x|=1\}.$
Set
$$\alpha _\varepsilon (x):=\varepsilon ^{-m} \alpha
\fracwithdelims ( ) {|x|}{\varepsilon}.\tag 2.3.1.3$$
 For any $\varepsilon$ we have  $\alpha_\varepsilon \in \Cal D$ and
supp $\alpha_\varepsilon \subset \{x:|x|\leq \varepsilon \}.$
\par
Let $\psi (y),y\in K \subset G $ be a Lebesgue summable function. Consider
the function
$$\psi_\varepsilon (x):=\int \limits_K \psi (y) \alpha_\varepsilon
(x-y) dy .\tag 2.3.1.4 $$
The function belongs to $\Cal D (G)$ for $\varepsilon$ small
enough and its support is contained in the $\varepsilon$-neighborhood of $K$.
\par
\subheading {2.3.2}
Let $f(x),\ x\in G \subset \Bbb R ^m $ be a locally summable function in $G.$
The formula
 $$<f,\varphi>:=\int f(y) \varphi (y)dy,\ \varphi \in {\Cal D} (G)
\tag 2.3.2.1 $$
defines a {\it linear continuous functional} on $\Cal D$, i.e., one that
satisfies the conditions
$$ <f,c_1\varphi_1 +c_2\varphi_2>= c_1 <f,\varphi_1>+ c_2 <f,\varphi_2>;
     \tag D'1$$
$$(\varphi_n \overset {\Cal D} \to \rightarrow \varphi)\Longrightarrow
<f,\varphi_n>\rightarrow <f,\varphi>.\tag D'2$$
\par
However, (2.3.2.1) does not exhaust all the linear continuous functionals
as one can see further. An arbitrary linear continuous functional on $\Cal D
 (G)$ is called an L.Schwartz {\it distribution} and the linear topological
space of the functionals is denoted as $\Cal D'(G)$.
\par
Following are some examples of functionals that do not have the form of
(2.3.2.1):
$$<\delta _x ,\varphi>:=\varphi (x) \tag 2.3.2.2$$
This distribution is called the {\it Dirac delta-function}. Further,
$$<\delta^{(n)} _x ,\varphi>:=\varphi ^{(n)} (x). \tag 2.3.2.3$$
This distribution is called {\it the n-th derivative} of the Dirac
delta-function.

 {\bf Exercise 2.3.2.1}Check that the functionals (2.3.2.2) and (2.3.2.3) are
both distributions.
\par
A distribution of the form (2.3.2.1) is called {\it regular}.
\proclaim {Theorem 2.3.2.1.(Du Bois Reymond )} If two locally
summable functions $f_1$ and $f_2$ define the same distribution,  then they
coincide almost everywhere.
\endproclaim
For proof see, e.g., \cite {H/"o, Th.2.1.6}.
\par
Remark that the converse assertion is obvious.
\par
A distribution $\mu$ is called {\it positive} if  $<\mu,\varphi> \geq 0$
for any
$\varphi \in \Cal D (G)$ such that $\varphi (x) \geq 0$ for all $x \in
\Bbb R ^m  .$ We shall write this as $\mu >0$ in $\Cal D'.$
\par
\subheading {Example 2.3.2.1}Let $\mu (E)$ be a measure in $G.$
Then the distribution
$$<\mu,\varphi>:=\int \varphi (x)d\mu \tag 2.3.2.4$$
 is positive.
\par
This formula represents all the positive distributions as one can see from
\proclaim {Theorem 2.3.2.2.(Positive Distributions)} Let $\mu >0$ in
$\Cal D (G)$.Then there exists a unique measure $\mu (E)$ such that
the distribution $\mu$ is given by (2.3.2.3).
\endproclaim
For proof see, e.g., \cite {H/"o,Th.2.1.7}.

\subheading {2.3.3}Let us consider operations on distributions.
\par
A {\it product } of a distribution $f$ by an {\it infinitely differentiable}
function $\alpha (x)$ is defined by
$$<\alpha f,\varphi >:=<f,\alpha \varphi >.\tag 2.3.3.1 $$
It is well defined because $\alpha \varphi \in \Cal D$ too.
\par
A {\it sum} of distributions $f_1$ and $f_2$ is defined by
$$<f_1+f_2,\varphi>:=<f_1,\varphi>+<f_2,\varphi>, \tag 2.3.3.2$$
and the {\it partial derivative} $ \frac{\partial}{\partial x_k} $ is defined by
the equality
$$<\frac{\partial}{\partial x_k} f,\varphi >:=
<f,-\frac{\partial}{\partial x_k} \varphi >. \tag 2.3.3.3$$
\par
  These definitions  look reasonable because of the following
\proclaim {Theorem 2.3.3.1.(Operations on Distributions)}The sum
 of regular distributions corresponds to the sum of the functions;
the  product of a regular
distribution by an infinitely differentiable function corresponds to the
product of the functions;  the derivative of a regular distribution that
generated by a differentiable function corresponds to the derivative of that
function.
\endproclaim
\demo {Proof}We have, for example,
$$<\alpha\cdot (f),\varphi>:=\int f(x)[\alpha (x) \varphi (x)]dx =
\int [\alpha (x)f(x)] \varphi (x)dx: = <(\alpha f),\varphi > $$
For the sum we have
$$<(f_1) +( f_2),\varphi >:= <f_1 ,\varphi> + <f_2 ,\varphi>=
\int f_1 (x)\varphi (x)dx +\int f_2 (x) \varphi (x) dx =$$
$$=\int [f_1 (x) + f_2 (x)] \varphi (x)dx = <(f_1 + f_2 ),\varphi >.$$
Let $f(x)$ have the derivative $ \frac{\partial}{\partial x_k} f .$ Then
$$<\frac{\partial}{\partial x_k} f,\varphi >:=
<f,-\frac{\partial}{\partial x_k} \varphi > =$$
$$ =\int f(x_1,x_2,...x_m ) [-\frac{\partial}{\partial x_k} \varphi
(x_1,x_2,...x_m)]dx_1 dx_2,...dx_m =$$
$$=\int dx_1 ...dx_{k-1} dx_{k+1} ... dx_m \int f(x_1,x_2,...x_m )
[-\frac{\partial}{\partial x_k} \varphi
(x_1,x_2,...x_m)]dx_k.$$
Now we shall do integrating by parts and all the substitution will vanish,
because $\varphi$ is finite. So we obtain
$$<\frac{\partial}{\partial x_k} f,\varphi >=
\int \frac{\partial}{\partial x_k} f(x) \varphi (x)dx .$$
That is to say the derivative of the distribution corresponds to the
function derivative.\qed \enddemo
\par
\subheading {2.3.4}We say that a sequence of distributions $f_n$ {\it converges} to
a distribution $f$  if
$$<f_n,\varphi> \rightarrow <f,\varphi >\ \forall \varphi \in \Cal D (G).
\tag 2.3.4.1$$
\proclaim {Theorem 2.3.4.1.(Completeness of $\Cal D'$)} If the
sequence of  numbers  $<f_n,\varphi >$ has a limit for every
$\varphi \in \Cal D (G),$ then
this functional is a linear continuous functional on $\Cal D (G)$, i.e.,
a distribution.
\endproclaim
For proof see, e.g., \cite {H\"o,Th.2.1.8}.
\par
 Differentiation is continuous with respect to convergence of
distributions.
\proclaim {Theorem 2.3.4.2.(Continuity of Differential Operators)} If
$f_n \rightarrow f$ in $\Cal D (G),$ then
$\frac{\partial}{\partial x_k} f_n \rightarrow
\frac{\partial}{\partial x_k} f.$\endproclaim
\demo {Proof} Set in (2.3.4.1) $\varphi := -\frac{\partial}{\partial x_k}
\varphi .$ Then
$$<\frac{\partial}{\partial x_k} f_n ,\varphi >=
<f_n ,-\frac{\partial}{\partial x_k} \varphi > \rightarrow
<f ,-\frac{\partial}{\partial x_k} \varphi >=
<\frac{\partial}{\partial x_k} f ,\varphi >.$$
\qed
\enddemo
The following theorem shows that the $\Cal D'$- convergence is the weakest
of the convergences considered earlier.

\proclaim {Theorem 2.3.4.3.(Connection between Convergences)} Let $f_n$ be
a sequence of  Lebesgue summable functions on domain $G$ such that at least one of
the following conditions holds:

Cnvr1) $f_n \rightarrow f$ uniformly on any compact set $K\Subset G $
and $f$ is a locally summable function;

Cnvr2) $f_n \rightarrow f$  on any $K\Subset G ,$
satisfying the conditions of the Lebesgue theorem (Th.2.2.2.1);

Cnvr3) $f_n \downarrow f$ monotonically and $f$ is a locally summable function.

Then $f_n \rightarrow f$ in $\Cal D' (G).$
\endproclaim
\demo {Proof}All the assertions are  corollaries of the section 2.2.2 of
passing to the limit under an integral.
\par
Let us prove, for example, Cnvr3). Let $f_n\downarrow f$. Then
$$<f_n , \varphi >=\int f_n (x)\varphi (x)dx =
\int f_n (x)\varphi ^+ (x)dx -\int f_n (x)\varphi ^- (x)dx \tag 2.3.4.2$$
where $\varphi ^+$ and $\varphi ^-$ are defined in (2.2.1.1).
\par
Both last integrals in (2.3.4.2) have a limit by the B.Levy
theorem \linebreak( Th.2.2.2.2), and thus
$$\lim \limits_{n\rightarrow \infty} <f_n , \varphi >=
\int f (x)\varphi ^+ (x)dx -\int f (x)\varphi ^- (x)dx =
\int f (x)\varphi (x)dx= <f , \varphi >. \tag 2.3.4.3$$
(2.3.4.3) means that $f_n \rightarrow f$ in $\Cal D'. $ \qed
\enddemo

{\bf Exercise 2.3.4.1} Prove Cnvr 1) and 2).

\proclaim {Theorem 2.3.4.4.($\Cal D'$ and $C^*$ convergences)} Let $\mu _n ,
\mu $ be measures in $G$. The conditions $\mu_n \rightarrow \mu$ in
$\Cal D'(G)$ and $\mu_n \overset*\to\ri \mu$ are equivalent.
\endproclaim
It is clear that the first condition is necessary for the second one. The
sufficiency holds, because every continuous function can be approximated with
functions that belong to $\Cal D$. For more details see, e.g.,
\cite {H\"o,Th.2.1.9}.

Let $\alpha_\epsilon(x)$ be defined as in (2.3.1.3).For any $f\in \Cal D'(D)$ we can consider the
function $f_\epsilon (x):=<f,\alpha_\epsilon(x+\bullet)>.$ It is called
a {\it regularization} of the distribution $f.$
\proclaim {Theorem 2.3.4.5.(Properties of Regularizations)}The following
holds:

reg1) $f_\epsilon (x)$ is an infinitely differentiable function in any
$K\Subset D$ for sufficiently small $\epsilon;$

reg2) $f_{\epsilon}(x)\rightarrow f$ in $\Cal D'(D)$ as
$\epsilon \downarrow 0;$

reg3) if $f_n\rightarrow f$ in $\Cal D'(D),$ $(f_n)_\epsilon\rightarrow
f_\epsilon$ uniformly with all its derivatives on any compact set in $D.$
\endproclaim
 The property reg1) follows from the formula
$$\frac {\partial}{\partial x_j}f_\epsilon =
<f,\frac {\partial}{\partial x_j}\alpha_\epsilon (x+\bullet)>.$$
 The property reg2) follows from the assertion
$$ \phi_\epsilon (x):=\int \phi (y)\alpha_\epsilon (x+y)dy\rightarrow\phi (x)
\ in\ \Cal D (D)$$
as $\epsilon\downarrow 0.$

For the proof of reg3) see \cite {H\"o,Theorem's 2.1.8, 4.1.5}.

Let us note the following assertion;
\proclaim {Theorem 2.3.4.6.(Continuity $<\bullet,\bullet>$)} The function
$$<f,\phi>:\Cal D'(G)\times \Cal D (G)\mapsto \Bbb R$$
is continuous in the according topology,
\endproclaim
i.e., $f_n\rightarrow f$ in $\Cal D'(G)$ and $\phi_j\rightarrow \phi$ in
$ \Cal D (G)$ imply $<f_n,\phi_j>\rightarrow <f,\phi>.$

For proof see \cite {H\"o,Th.2.1.8}.

\subheading {2.3.5} Let $G_1\subset G.$ Then
$\Cal D' (G)\subset\Cal D' (G_1)$ , because every functional on $\Cal D (G)$
can be considered as a functional on $\Cal D (G_1).$
\par
A distribution $f\in \Cal D' (G)$ considered as a distribution in
$\Cal D' (G_1)$ is called the {\it restriction } of $f$ to $G_1$ and is denoted
 $f\mid _{G_1}.$
\proclaim {Theorem 2.3.5.1.(Sewing Theorem)} Let
$G_\alpha \subset {\Bbb R}^m $
be a family of  domains and in  every of them let there be a distribution
$f_\alpha \in \Cal D (G_\alpha),$ such that

 If $G_{\alpha_1} \cap G_{\alpha_2} \neq \varnothing$, the equality
     $$ f_{\alpha_1} \mid _{G_{\alpha_1} \cap G_{\alpha_2}} =
f_{\alpha_2} \mid _{G_{\alpha_1} \cap G_{\alpha_2}} \tag 2.3.5.1 $$
holds.

Then there exists one and only one distribution $f\in \Cal D' (G)$ where
$G=\underset {\alpha} \to {\bigcup} G_\alpha $ such that
$f\mid _{G_\alpha} =f_{\alpha}.$
\endproclaim
In particular, it means that every distribution is defined uniquely by its
restriction to  a neighborhood of  every point.
\par
Let $\Cal D (S_R)$ be a space of infinitely differentiable functions on the
 sphere \linebreak $S_R:=\{x:|x|=R\}.$ The corresponding distribution space
is denoted
as $\Cal D' (S_R).$ The sewing theorem holds for this space in the following
form:
\proclaim {Theorem 2.3.5.2.($\Cal D'$ on Sphere)} Let a family of domains
 $\Omega_\alpha $ cover $S_R$ and in every of them let there be a distribution
$f_\alpha \in \Cal D (\Omega_\alpha),$ such that

If  $\Omega_{\alpha_1} \cap \Omega_{\alpha_2} \neq \varnothing$ ,the equality
     $$ f_{\alpha_1} \mid _{\Omega_{\alpha_1} \cap \Omega_{\alpha_2}} =
f_{\alpha_2} \mid _{\Omega_{\alpha_1} \cap \Omega_{\alpha_2}} \tag 2.3.5.1 $$
holds.

Then there exists one and only one distribution $f\in \Cal D' (S_R)$
such that $f\mid _{\Omega_\alpha} =f_{\alpha}.$
\endproclaim
\subheading {2.3.6} Let
$$ L:=\dsize\sum \limits_{i,j}\frac{\partial}{\partial x_i}
 a_{i,j} (x) \frac{\partial}{\partial x_j } +q(x) \tag 2.3.6.1 $$
be a differential operator of  second order with infinitely differentiable
coefficients $a_{i,j} ,q .$

We will consider only three types of differential operators:
one dimensional operator with constant coefficients,the  Laplace operator and
the so called spherical operator (see Sec.2.4 ).
\par
For all these operators we have the following assertion which follows from
the general theory (see, e.g., \cite {H\"o,Th.11.1.1} ):
\proclaim {Theorem 2.3.6.1.(Regularity of Generalized Solution)} If the
equation \linebreak $Lu =0 $ has a solution $u \in \Cal D'(G) ,$  then $u$ is a regular
 distribution and can be realized as an infinitely differentiable function .
\endproclaim
\par
 A distribution that satisfies the equation
$$Lu = \delta_y \ \ in \ \Cal D' (G),\tag 2.3.6.2$$
 where $\delta_y $ is a Dirac delta function (see (2.3.2.2)),
is called a {\it fundamental solution } of $L$ at the point $y.$
\par
Every differential operator that we are going to consider has a fundamental
solution (see, e.g., \cite {H\"o,Th.10.2.1}).
\par
A restriction of the equation (2.3.6.2) to the domain $G_y :=G\backslash y$
is a homogeneous equation $Lu =0$ in $\Cal D' (G_y)$. Thus we have
\proclaim {Theorem 2.3.6.2.(Regularity of Fundamental Solution)} The
fundamental solution is an infinitely differentiable function outside the
point $y.$
\endproclaim
\subheading {2.3.7}We will need further also
the Fourier coefficients for the distribution
on the circle .

Let $\Cal D(S^1)$ be a set of all infinitely
differentiable function on the unit circle
$S^1.$ The set of all linear continuous functionals
over $\Cal D(S^1)$ with the corresponding topology (see 2.3.2) is the corresponding space of
distributions $\Cal D'(S^1)$ for which all the previous properties of distributions holds.

The functions $\{e^{ik\phi},\ k=0,\pm 1,\pm 2,...
\}$ belong to $\Cal D(S^1).$
The Fourier coefficients of $\nu\in  \Cal D'(S^1)$ are defined by
$$\hat \nu (k):=<\nu, e^{-ik\phi}> .\tag 2.3.7.1$$
The inverse operator is defined by
$$<\nu,g>=\frac {1}{2\pi}\sum\limits_
{k=-\iy}^{\iy}\hat\nu (k)<g,e^{ik\phi}> ,\tag 2.3.7.2$$
and the series converges, in any case, for those $\nu$ that are finite derivatives of summable functions, because Fourier coefficients of $g$
decrease faster then every power of $x.$

The convolution of distribution $\nu\in \Cal D'(S^1)$ and $g\in \Cal D(S^1)$ is defined by
$$\nu*g(\phi)=<\nu,g(\phi-\bullet)>,
\tag 2.3.7.3.$$
This is a function from $\Cal D(S^1).$

The convolution of distributions $\nu_1,\ \nu_2  \in \Cal D'(S^1)$ is defined by
$$<\nu_1*\nu_2, g>=\nu_1*(\nu_2*g).\tag 2.3.7.4
$$
In spite of the view it is commutative and
$$\widehat {\nu_1*\nu_2}(k)=\hat \nu_1 (k) \cdot
\hat\nu_2 (k).$$

 {\bf Exercise 2.3.7.1} Count the Fourier
coefficients of the functions
$$\ G(re^{i\phi})=\log |1-re^{i\phi}|\tag 2.3.7.5$$
for $r>1,\ r=1, r<1;$ the function defined by
$$\widetilde {\cos\r}(\phi):=\cos\r\phi,\  -\pi<\phi<\pi,\ \r\in (0,\iy)\tag 2.3.7.6$$
and $2\pi$-periodically extended; the function
$$\tilde \phi \sin p\phi,\ p\in \BN \tag 2.3.7.7$$
where $\tilde \phi$ is the $2\pi$-periodical
extension of the function $f(\phi)=\phi,\ \phi\in
[0,2\pi).$

 {\bf Exercise 2.3.7.2}Denote
$$
P_{p-1}(re^{i\phi}):=\Re\left\{\suml_{k=1}^{p-1}
\frac {r^ke^{ik\phi}}{k}\right\},\ \ p\in \BN.
\tag 2.3.7.8$$
Prove that for every distribution $\nu:$
$$ (P_{p-1}(\widehat{re^{i\bullet})*\nu})(p)=0\tag 2.3.7.9$$
The same for the function
$$G_{p}(re^{i\phi}):=G(re^{i\phi})+P_p(re^{i\phi})$$
for $r<1.$

\newpage

\centerline {\bf 2.4.Harmonic functions}
\subheading {2.4.1}We will denote as $\Delta$ the Laplace operator in
$\Bbb R ^m$:
$$\Delta :=\frac {\partial ^2}{\partial x_1^2} + ...+\frac {\partial ^2}{\partial x_m^2}.$$
We introduce in $\Bbb R ^m$ the spherical coordinate system by the formulae:
$$ \align
 x_1=& r\sin \phi _0\sin \phi _1 ...\sin \phi _{m-2};\ \ \ \ \ \ \ \ \\
 x_2 =& r\cos \phi _0\sin \phi _1 ...\sin \phi _{m-2};\\
 x_3 = &r\cos \phi _1\sin \phi _2 ...\sin \phi _{m-2}; \\
 ......&................................................... \\
 x_k = &r\cos \phi _{k-2}\sin \phi _{k-1} ...\sin \phi _{m-2};\\
 ......&.....................................................\\
 x_m = &r\cos \phi _{m-2},\endalign $$
where
$$ 0<\phi_0\leq 2\pi;\ \ 0\leq \phi_j <\pi,\ j=\overline {1,m-2};\ \
0 < r <\infty.$$

Passing  to the coordinates $(r,\phi_0,\phi_1,\  ...\  \phi_{m-2})$ in
the Laplace operator we obtain
$$\Delta = \frac {1}{r^{m-1}}\frac {\partial}{\partial r} r^{m-1}
\frac {\partial}{\partial r}+ \frac{1}{r^2}\Delta _{\bold x ^0}.$$
The operator $\Delta _{\bold x ^0}$ is called {\it spherical}, and has the form
$$\Delta _{\bold x ^0 }:=
\sum\limits_{i=0}^{m-2} \frac {1}{\Pi}\frac {\partial}{\partial \phi _i}
\frac{\Pi}{\Pi_i}\frac {\partial}{\partial \phi _i},$$
where
$$\Pi := \prod\limits_{j=1}^{m-2} \sin ^{j}\phi _j;\
\Pi_i := \prod\limits_{j=i+1}^{m-2} \sin ^{2}\phi _j;\
\Pi_{m-2}:=1.$$
In particular, for $m=2$, i.e., for the polar coordinates,
$$\Delta = \frac {1}{r}\frac {\partial}{\partial r} r
\frac {\partial}{\partial r}+ \frac{1}{r^2}
\frac {\partial ^2}{\partial \phi ^2}.$$
A distribution $H \in \Cal D '(G)$ is called {\it harmonic} if it satisfies
the equation $\Delta H =0.$

The next theorem follows from Theorem 2.3.6.1.
\proclaim {Theorem 2.4.1.1.(Smoothness of harmonic functions)} Any harmonic
distribution is equivalent to an infinitely differentiable function.
\endproclaim
 This function, of course, satisfies the same equation and is a {\it harmonic
function} in the ordinary sense.
 A direct proof can be found, e.g., in \cite {Ro,Ch.1,\S2 (1.2.5),p.60}.

Let $f(z),\ z=x+\imath y$ be a holomorphic function in a domain $G\subset \Bbb C$. Then the functions
$u(x,y):=\Re f(z)$ and $v(x,y):= \Im f(z)$ are harmonic in $G.$
In particular, the functions $r^n\cos n\varphi$ and $r^n\sin n\varphi$ where
$r=|z|,\ \varphi= \roman {arg}z$ are harmonic.

Set
$$\Cal E_m (x):=\cases -|x|^{2-m},\ &\text{for $ m\geq 3 $}\ \ \ \ \ \ \\
                        \log |z| ,  &\text{for $ m=2 $}
\endcases \tag 2.4.1.1 $$
(We will often denote points of the plane as $z$).

It is easy to check that  $\Cal E_m (x)$ is a harmonic function for
$|x|\neq 0.$

Set
$$\theta _m :=\cases (m-2)\sigma_m,\ \ &\text {for $ m\geq 3;$}\ \ \ \ \ \ \ \\    2\pi,            &\text {for $m=2 ;$}\endcases $$
where $\sigma_m$ is the surface area of the unit sphere in $\Bbb R^m.$

\proclaim {Theorem 2.4.1.2.(Fundamental Solution)}The function $\Cal E_m (x-y)$
satisfies in $\Cal D'(\Bbb R^m)$ the equation \footnote {$\Cal E_m $ is slightly a little different from the fundamental solution (see,(2.3.6.2)), but this is traditional in Potential Theory}
$$\Delta_x \Cal E_m (x-y)=\theta_m \delta (x-y),\tag 2.4.1.2$$
where $\delta (x)$ is the Dirac $\delta$ -function (see 2.3.2).
\endproclaim
\demo {Proof}Let us prove the equality (2.4.1.2) for $y=0.$ Suppose
$\phi \in\Cal D(\Bbb R^m)$ and \linebreak supp $\phi \subset K \Subset \Bbb R^m.$
We have
$$<\Delta \Cal E_m ,\phi >:=
\int \Cal E_m (x)\Delta\phi(x) dx =
\lim_{\epsilon \rightarrow 0}\int\limits_{|x|\geq \epsilon} \Cal E_m (x)\Delta
\phi(x) dx.$$
Transforming this integral by the Green formula and using the fact that
$\phi$ is finite we obtain
$$\int\limits_{|x|\geq \epsilon} \Cal E_m (x)\Delta\phi(x) dx =
\int\limits_{|x|\geq \epsilon} \Delta\Cal E_m (x)\phi(x) dx +
\int\limits_{|x|=\epsilon} \Cal E_m \frac {\partial \phi}{\partial n}ds -
\int\limits_{|x|=\epsilon} \phi \frac {\partial \Cal E_m}{\partial n}ds,$$
where $ds$ is an element of  surface area and $\frac {\partial }{\partial n}$
is the differentiation in the direction of the external normal.

Use the harmonicity of $\Cal E_m$. Then the first integral is equal to zero.
Further we have
$$\int\limits_{|x|=\epsilon} \Cal E_m \frac {\partial \phi}{\partial n}ds =
\epsilon\left .\left(\int\limits_{|x^0|=1} \frac {\partial}{\partial r}
\phi (r x^0)ds\right )\right|_{r=\epsilon}=O(\epsilon),\text {for $\epsilon\rightarrow 0$}$$
For the third term we have
$$\int\limits_{|x|=\epsilon} \phi \frac {\partial \Cal E_m}{\partial n}ds=
\frac{m-2}{\epsilon ^{m-1}} \epsilon ^{m-1}
\int\limits_{|x^0|=1}\phi (r x^0)ds = [\phi (0)+o(1)](m-2)\sigma_m .$$
Thus we obtain $<\Delta \Cal E_m ,\phi > = \phi (0) \theta_m $, and this proves
(2.4.1.2) for $y=0.$

It is clear that by changing $\phi (x)$ for $\phi (x+y)$ we obtain (2.4.1.2)
in the general case.\qed
\enddemo

We will consider now a domain $\Omega$ with a {\it Lipschitz} boundary
({\it Lipschitz domain}).
It means that every part of $\partial \Omega $ can be represented in some
local coordinates $(x,x'),\ x\in \Bbb R,\ x'\in \Bbb R^{m-1}$ in the form
$x=f(x'),$ where $f$ is a Lipschitz function, i.e.,
$$|f(x'_1)-f(x'_2)|\leq M_{\partial {\Omega}}|x'_1-x'_2|$$
where  $M$ depends only on the whole $\partial \Omega$ and does not depend  on this
local part.

Let $G(x,y,\Omega)$ be the  Green function of a Lipschitz domain $\Omega.$

It is known (see, e.g. \cite {Vl,Ch.V,\S28} ) that the Green function
has the following properties:
$$G(x,y,\Omega)<0,\ \text {for $(x,y)\in \Omega \times \Omega;\
G(x,y,\Omega)=0$ for $(x,y)\in \Omega \times \partial\Omega;$} \tag g1$$
$$G(x,y,\bullet)=G(y,x,\bullet);\tag g2$$
$$G(x,y,\bullet)-\Cal E_m (x-y)= H(x,y),\tag g3$$
where H is harmonic on $x$ and on $y$  within $\Omega.$
$$- G(x,y,\Omega_1)\leq - G(x,y,\Omega_2)\text { for } \Omega_1\subset
\Omega_2 \tag g4$$
From (g3) follows
\proclaim {Theorem 2.4.1.3.(Green Function)}The equality
$$\Delta_x G(x,y,\Omega)=\theta_m \delta (x-y),\tag 2.4.1.3$$
holds in $\Cal D'(\Omega).$
\endproclaim

Let $f(x)$ be a continuous function  on $\partial \Omega$. It is known
(see, e.g.\cite {Vl,Ch.V,\S 29}) that the function
$$H(x,f):=\int\limits_{\partial\Omega}f(y)\frac {\partial}{\partial n}
G(x,y,\Omega)ds_y \tag 2.4.1.4 $$
is the only harmonic function that coincides with $f$ on $\partial \Omega.$

The unique solution of the Poisson equation
$$\Delta u =p,\ u|_{\partial \Omega}=f $$
for a continuous function $p$ is given by the formula
$$ u(x,f,p):=\int\limits_{\partial\Omega}f(y)\frac {\partial}{\partial n_y}
G(x,y,\Omega)ds_y +\theta_m^{-1}\int\limits_{\Omega}G(x,y,\Omega)p(y)dy.\tag 2.4.1.5$$

Let $D$ be an {\it arbitrary} open domain.We can define a $G(x,y,D)$
in the following way. Consider a sequence $\Omega _n$ of a Lipschitz domains
such that $\Omega_n\uparrow D.$ The sequence of the corresponding Green
functions $G(x,y,\Omega _n)$ monotonically decreases. If it is bounded from
below in some point, it is bounded everywhere while $x\neq y$ (as it follows
from Theorem 2.4.1.7). It can be shown that the limit exists for any domain
the boundary of which have  positive {\it capacity} (see 2.5 and references there).
We will  mainly use the Green function for the Lipschitz domains.

Let $G(x,y,K_{a,R})$ be the Green function of the ball
$K_{a,R}:=\{|x-a|<R\}.$
\proclaim {Theorem 2.4.1.4.(Green Function for Ball)}
$$G_(x,y,K_{a,R})=\cases -|x-y|^{2-m}-(\frac {|y-a||x-y^*_{a,R}|}{R})^{2-m},\ &
\text {for $m\geq 3$}\\
 \log \frac {|\zeta -z|R}{|\zeta - a||z-\zeta ^*_{a.R}|}\ \ \ \ \ \ &\text {for
$m=2$}\endcases$$
where
$y^*_{a.R}:=a+(y-a)\frac {R^2}{|y-a|^2}$
is the inversion of $y$ relative to the sphere $\{|x-a|=R\}.$
\endproclaim
For the proof see, e.g.,  \cite {Br},Ch.6,\S3 .
\proclaim {Theorem 2.4.1.5.(Poisson Integral)} Let $H$ be a harmonic function
in $K_{a,R}$ and continuous in its closure. Then
$$ H(x)=\frac {1}{\sigma_m R}\int\limits_{|x-a|=R} H(y)\frac {R^2-|x-a|^2}
{|x-y|^m}ds_y,\ x\in K(a,R). \tag 2.4.1.6$$
In particular, for $m=2$
 $$H(a+re ^{\imath\phi})= \frac {1}{2\pi}\int\limits_{0}^{2\pi}
H(a+Re ^{\imath\psi})\frac {R^2 - r^2}{R^2 -2Rr\cos (\phi - \psi)+r^2}d\psi.$$
\endproclaim
This theorem  follows from (2.4.1.4).
\proclaim {Theorem 2.4.1.6.(Mean Value)} Let $H$ be harmonic in $G\subset
\Bbb R^m$. Then
$$H(x)=\frac {1}{\sigma_m R^{m-1}}\int\limits_{|x-a|=R} H(y)ds_y,\tag 2.4.1.7
$$
where $x\in G$ and $R$ is taken such that $K(x,R)\Subset G.$
\endproclaim
We must only set  $a:=x$ in (2.4.1.6).

We can rewrite (2.4.1.7) in the form
$$H(x)=\frac {1}{\sigma_m }\int\limits_{|y|=1} H(x+Ry)ds_y.$$
\proclaim {Theorem 2.4.1.7.(Harnack)}Suppose the family $\{ H_\alpha) ,\
\alpha \in A\}$ of harmonic functions in $G$ satisfies the conditions
$$H_\alpha (x)\leq C(K),\ \text {for $x\in K$;}\tag Har1$$
$$ H_\alpha (x_0)\geq B>-\infty,\ \text {for } x_0\in K \tag Har2 $$
for every compact  $K\Subset G$ and $C(K),B$ are constants not
depending on $\alpha$.

Then the family is precompact in the uniform topology, i.e., there exists
such a
sequence $H_{\alpha_n},$ and a function $H$ harmonic in the interior of $K$ and continuous in K
 such that $H_{\alpha_n}\rightarrow H$ uniformly in every $K.$
\endproclaim
 One can prove by using (2.4.1.6) that $|\grad H_\alpha|$  are bounded on every
compact set by a constant  not depending on $\alpha.$ Thus the family is
uniformly continuous and thus it is  precompact by the Ascoli theorem  .

For details see, e.g., \cite {Br,Supplement,\S7 }.

\proclaim {Theorem 2.4.1.8.(Uniform and $\Di'$-convergences)}Suppose the sequence
$H_n$ satisfies the conditions of the Harnack theorem and converges to a
function $H$ in $\Cal D'(G).$ Then $H_n$ converges to $H$ uniformly on every
compact $K\Subset G.$
\endproclaim
Of course, $H$ is harmonic in $G.$
\demo {Proof}By the Harnack theorem the family is precompact.Thus we must  only
prove the uniqueness of $H.$ Suppose there exist two subsequences such
that
$H^1_k \rightarrow H^1$ and $H^2_k \rightarrow H^2$ uniformly on every compact
$K\Subset G.$

By Connection between Convergences (Theorem 2.3.4.3) $H^1_k \rightarrow H^1$
and $H^2_k \rightarrow H^2$ in $\Cal D'.$ Hence, $H^1 =  H^2$ in $\Cal D'(G).$
By the De Bois Raimond theorem (Theorem 2.3.2.1) $H^1 =  H^2$ almost everywhere
and hence everywhere because these functions are continuous.\qed
\enddemo
Let $D$ be a domain with a smooth boundary $\partial D$ and let $F\subset
\partial D.$ Set
$$\omega (x,F,D):=\int \limits_F \frac {\partial G}{\partial n_y}(x,y)ds_y. $$
It is called a {\it harmonic measure} of $F$ with respect to $D.$
A harmonic measure can be defined for an arbitrary domain $D$ by a limit
process similar to the one we had for the Green function.In this case the formula (2.4.1.4) has the form
  $$ H(x,f):=\int\limits_{\partial D}f(y)d\omega (x,y,D). $$
However we can not assert that $H(x,f)$ coincides with $f$ in any point
$x\in D.$
We can only consider it as an operator that maps a function defined on
$\partial D$ to a harmonic function in $D.$

By (2.4.1.3) we obtain
\proclaim {Theorem 2.4.1.9.(Two Constants Theorem)} Let $H$ be harmonic in $D$
and satisfy the conditions
$$H(x)\leq A_1\ for\ x\in F;H(x)\leq A_2\ for\ x\in \partial D\backslash F$$
where $A_1$ and $A_2$ are constants.

Then
$$H(x)\leq  A_1\omega (x,F,D)+A_2\omega (x,\partial D\backslash F,D)\ for
\ x\in D.$$
\endproclaim

Let $y^*_{a,R}$ be the inversion from Green Function for Ball
(Theorem 2.4.1.4). Set $y^*:=y^*_{0,1}$, i.e., the inversion relative to
a unit sphere with the center in the origin. Let
$G^*:=\{y^*: y\in G\}$ be the inversion of a domain $G.$
\proclaim {Theorem 2.4.1.10.(Kelvin's Transformation)} If $H$ is harmonic in
$G,$ then
$$H^*(y):=|y|^{2-m}H(y^*)\tag 2.4.1.8$$
is harmonic in $G^*.$
\endproclaim
For the proof you must honestly compute  Laplacian of $H^*$.''The computation
is straightforward but tedious`` (\cite {He,Theorem 2.24}). It is not
so tedious if you  use the spherical coordinate system.

{\bf Exercise 2.4.1.1} Do this.

\subheading {2.4.2} Denote as $S_1:=\{x^0:|x^0|=1\}$ the unit sphere with
 center in the origin. A function $Y_\rho (x^0),\ x^0 \in \Omega \subset
S_1$ is
called a {\it spherical function } of {\it degree} $\rho$ if it satisfies
the equation
$$ \Delta _{\bold x ^0}Y + \rho (\rho+ m-2)Y = 0.\tag 2.4.2.1$$
For $m=2$ (2.4.2.1) gets the form
$$Y''(\theta)+ \rho ^2 Y(\theta)=0,$$
i.e.,
$$ Y(\theta)=a\cos \rho\theta +b\sin\rho \theta,$$
Spherical functions are obtained if we solve the equation $\Delta H =0$
 by the change
$H(x)=|x|^\rho Y(x^0).$
\proclaim {Theorem 2.4.2.1.(Sphericality and Harmonicity)} Let $Y_\rho (x^0)$
be  spherical in a domain $\Omega \subset S_1$  if and only if the functions
$H(x)=|x|^\rho Y_\rho (x^0)$ and \linebreak
$H^* (x)=|x|^{-\rho -m +2} Y_\rho (x^0)$ are
harmonic in the cone
$$Con (\Omega):=\{x=rx^0:x^0 \in \Omega,\ 0<r<\infty \}.\tag 2.4.2.2$$
If $\rho=k,\ k\geq 0,\ k\in \Bbb Z,$ and only in this case, $Y_k
(x^0)$ is spherical on the whole $S_1,$ $H(x)$ is a homogeneous
harmonic polynomial of degree $k$ and $H^*$ is harmonic in $\Bbb
R^m \backslash 0.$
\endproclaim

For the proof see, e.g.\cite {Ax}

The spherical functions of an integer degree $k$ form a finite-dimension space
 of dimension
$$dim (m,k)=\frac {(2k+m-2)(k+m-3)!}{(m-2)! k!}$$
In particular, $d(2,k)=2$ for any $k.$

For different $k$ the spherical functions $Y_k (x^0)$ are orthogonal on  $S_1.$
In particular, for $m=2$, it means the orthogonality of the trigonometric
functions system.
\proclaim {Theorem 2.4.2.2.(Expansion of Harmonic Function )}Let $H(x)$ be a
harmonic function in the ball $K_R :=\{|x|<R\}.$ There exists an orthonormal
system of spherical functions $Y_k(x^0),\ k=\overline {0,\infty},$ depending
on $H$ such that
$$H(x)=\sum\limits_{k=0}^\infty c_k Y_k (x^0)|x|^k,\ for\ |x|<R .
\tag 2.4.2.3$$
For any such system we have
$$c_k=\frac {1}{R^k}\int\limits_{S_1} H(Rx^0)Y_k (x^0)ds_{x^0}.\tag 2.4.2.4 $$
\endproclaim
For proof see, e.g., \cite {Ax,Ch.10},\cite {TT, Ch.4,\S 10} .
\proclaim {Theorem 2.4.2.3.(Liouville)}Let $H$ be harmonic in $\Bbb R^m$ and
 suppose
$$\liminf \limits _{R\rightarrow\infty}R^{-\rho}\max \limits_{|x|= R}H(x)
<\infty.\tag 2.4.2.5$$
holds.

Then $H$ is a polynomial of a degree $q\leq \rho.$
\endproclaim
\demo {Proof}We can suppose $H(0)=0$ because $H(x)-H(0)$ is harmonic and also
satisfies (2.4.2.5). Let $R_n\rightarrow\infty$ be a sequence for which
$$ R_n^{-\rho}\max \limits_{|x|= R_n}H(x) \leq const <\infty \tag 2.4.2.6$$
From (2.4.2.4) we obtain
$$|c_k|\leq A_k R^{-k}\int\limits_{S_1} |H(Rx^0)|ds_{x_0},\tag 2.4.2.7$$
where
 $A_k= \max_{S_1} |Y_k (x^0)|.$

From the mean value property (Theorem 2.4.1.6)
$$ \int\limits_{S_1} H(Rx^0)ds_{x_0}=H(0)\sigma_m=0.$$
Thus
$$\int\limits_{S_1} |H(Rx^0)|ds_{x_0}=2\int\limits_{S_1} H^+ (Rx^0)ds_{x_0}
\leq 2\sigma _m \max \limits_{|x|= R}H(x).\tag 2.4.2.8$$
 From (2.4.2.8) and (2.4.2.7) we have
$$|c_k|\leq 2A_k R^{-k}\sigma _m \max \limits_{|x|= R}H(x).\tag 2.4.2.9$$
Set $R:=R_n$ and $k>\rho.$ Passing  to the limit when
$n\rightarrow\infty,$ we obtain $c_k =0$ for $k>\rho.$ Then
(2.4.2.3.) implies that $H$ is a harmonic polynomial of degree
$q\leq \rho.$ \qed
\enddemo

\newpage
\centerline {\bf 2.5.Potentials and Capacities}
\subheading {2.5.1}Let $G(x,y.D)$ be the Green function of a Lipschitz
domain $D.$ We will suppose it is extended as zero outside of $D.$
$$\Pi (x,\mu,D):= -\int G(x,y,D)d\mu_y $$
is called the {\it Green potential } of $\mu$ relative to $D.$

The domain of integration will always be $\Bbb R^m.$
\proclaim {Theorem 2.5.1.1.(Green Potential Properties)}The following holds:

GPo1) $  \Pi (x,\mu,D)$ is lower semicontinuous;

GPo2) it is summable over any $(m-1)$- dimensional hyperplane or smooth
 hyper-surface;

GPo3) $ \Delta \Pi (\bullet,\mu,D) = -\theta_m \mu \text { in } \Cal D'(D);$

GPo4) the {\bf  reciprocity law } holds:
 $$ \int \Pi (x,\mu_1,D)d\mu_{2x}= \int \Pi (x,\mu_2,D)d\mu_{1x}.$$

GPo5) {\bf semicontinuity} in $\mu$:
if $\mu_n\rightarrow \mu$ in $\Cal D'(\Bbb R^m),$ then
$$\liminf_{n\rightarrow \infty}\Pi (x,\mu_n,D)\geq \Pi (x,\mu,D).$$

GPo6) {\bf continuity} in $\mu$ in $\Cal D'$:
if $\mu_n\rightarrow \mu$, then $\Pi(\bullet,\mu_n,D)\rightarrow
 \Pi(\bullet,\mu,D)$ in $\Cal D'(\Bbb R^m)$ and in $\Cal D'(S_R)$ ,where
$S_R$ is the sphere $\{|x|=R\}.$
\endproclaim
\demo {Proof}Let us prove GPo1). Let $N>0.$ Set
$\  G_{N}(x,y):=\max (G(x,y), -N ),$
 a truncation of the function $G(x,y).$

The functions $G_{N}$ are continuous in $\Bbb R^m\times\Bbb R^m $ and
$G_{N}(x,y)\downarrow G(x,y)$ for every $(x,y)$ when
$N\rightarrow\infty.$
Set
$$\Pi_N (x,\mu,D):= -\int G_N (x,y,D)d\mu_y .$$
The functions $\Pi_N$ are continuous and $\Pi_N (x,\bullet )\uparrow
\Pi(x,\bullet)$ by the B.Levy theorem (Theorem 2.2.2.2). Then
$\Pi_N (x,\bullet )$ is lower semicontinuous by the Second Criterion of
semicontinuity (Theorem 2.1.2.9).

Let us prove GPo5).
Since Theorem 2.3.4.4.($\Cal D'$ and $C^*$ convergences)\newline
 $\lim_{n\rightarrow \infty} \Pi_N (x,\mu_n,D)
=\Pi_N (x,\mu,D).$ Further $\Pi (x,\mu_n,D)\geq \Pi_N (x,\mu_n,D)$,
 hence $\liminf_{n\rightarrow\infty}\Pi(x,\mu_n,D)\geq\Pi_N (x,\mu,D)$
Passing to limit while $N\rightarrow \infty,$ we obtain GPo5)

The assertion GPo2) follows from the local summability of the function
$|x|^{2-m}$ that can be checked directly.

Let us prove GPo3). For $\phi \in \Cal D (D)$ we have
$$\align <\Delta \Pi,\phi>:&=<\Pi ,\Delta \phi> =-\int d\mu_y\int G(x,y,D)
\Delta \phi (x)dx \\
&=-\int <\Delta _x G( \bullet,y,D), \phi> d\mu _y =- \theta _m
\int \phi(y)d\mu _y \\ &=-\theta _m <\mu,\phi>.\endalign $$
since
$$ <\Delta _x G( \bullet,y,D), \phi>= \theta _m \phi (y).$$
 by Theorem 2.4.1.3.
The property GPo4) follows from the symmetry of \linebreak $G(x,y,\bullet)$
(property (g2)).

Let us prove GPo6). Note that integral $\int |x|^{m-1}dx$ converges locally in
$\Bbb R^m$ and in $\Bbb R^{m-1}.$ From this one can obtain by some simple
estimates that functions
$\Psi (y):=\int G(x,y,D)\psi (x)dx$ while $\psi\in \Cal D (\Bbb R^m)$ and
$\Theta (y):=\int_{S_R} G(x,y,D)\theta (x)ds_x$ while $\theta\in
\Cal D (S_R)$ are continuous on $y\in \Bbb R^m.$

Now we have
$$<\Pi (\bullet,\mu_n,D),\psi>=\int \Psi (y)d\mu_n (y)\rightarrow
\int \Psi(y)d\mu (y)=<\Pi(\bullet,\mu,D),\psi>.$$
Thus the first assertion in GPo6) is proved. The second one can be
 proved in
the same way.
\qed
\enddemo
Set $\nu :=\mu_1 -\mu_2$, and let
$\Pi (x,\nu,D):= \Pi (x,\mu_1,D)-\Pi (x,\mu_2,D)$
be a potential of this charge. Consider the boundary problem of the form
$$\align \Delta u &=\mu_1 -\mu_2 ,\text{ in } \Cal D'(D)\\
 u|_{\partial D}&=f, \tag 2.5.1.2 \endalign$$
where $f$ is a continuous function.
\proclaim {Theorem 2.5.1.2.(Solution of Poisson Equation)}The solution of
the boundary problem (2.5.1.2) is given by the formula
$$u(x)=H(x,f)-\theta_m^{-1}\Pi (x,\nu,D),$$
where $H(x,f)$ is the harmonic function from (2.4.1.4).
\endproclaim
\demo {Proof}Since $\Pi (x,\nu,D)|_{\partial D}=0$ the function $u(x)$
satisfies the boundary condition. Using GPo3) we obtain
$$\Delta u = \Delta H -[\theta _m]^{-1} \Delta \Pi =\mu_1 -\mu_2 .$$
\qed
\enddemo
A potential of the form
$$\Pi (x,\mu):= \int \frac {d\mu_y}{|x-y|^{m-2}}$$
is called a {\it Newton} potential. It is the Green potential for
$D=\Bbb R^m.$ The potential
$$\Pi(z,\mu)= -\int \log |z-\zeta|d\mu_\zeta$$
is called {\it logarithmic}.
\subheading {2.5.2}Let $K\Subset D$.The quantity
$$\text {{\bf cap}} _G (K,D):=\sup \mu (K)\tag 2.5.2.1$$
where the supremum is taken over all the mass distributions $\mu$ for which
the
following conditions are satisfied:
$$\Pi (x,\mu,D )\leq 1 \tag 2.5.2.2$$
$$\roman {supp}\mu\subset K, \tag 2.5.2.3$$
is called the {\it Green capacity} of the compact set $K$ relative to the domain
$D.$
\proclaim {Theorem 2.5.2.1.(Properties of cap$_G$)}For
$\text{\bf {cap}}_G$ the
following properties hold:

capG1) monotonicity with respect to $K$:\newline
$K_1\subset K_2$ implies
$ \text{\bf cap} _G (K_1,D)\leq \text{\bf cap} _G (K_2,D).$

capG2) monotonicity with respect to $D$:\newline
$K\Subset D_1\subset D_2$ implies $\text{\bf cap} _G (K,D_1)\geq
\text{\bf cap} _G (K,D_2)$

capG3) subadditivity with respect to $K:$
$$\text{\bf cap} _G (K_1\cup K_2,D)\leq
\text{\bf cap} _G (K_1,D)+\text{\bf cap} _G (K_2,D)$$

\endproclaim
\demo {Proof}The set of all the mass distributions that satisfy (2.5.2.2) for
$K=K_1$ is not less than the analogous set for $K=K_2.$ Thus capG1) holds.

By the Green function property (g3) (see \S 2.4.1)
$-G(x,y,D_1)\leq -G(x,y,D_2).$ Thus the set of all $\mu$ that
satisfy (2.5.2.2)
for $D=D_1$ is wider than for $D=D_2.$ Hence capG2) holds.

Let supp $\mu \subset K_1\cup K_2$ and let $\mu_1,\mu_2$ be the restrictions of
$\mu$ to $K_1,K_2$ respectively.

If $\mu$ satisfies the (2.5.2.2) for $K:=K_1\cup K_2$ then
$\mu_1,\mu_2$ satisfy (2.5.2.2) for $K:=K_1,K_2$ respectively.

 From the inequality
$$\mu (K_1\cup K_2)\leq\mu (K_1)+\mu (K_2)$$
we obtain that
$$ \mu (K_1\cup K_2)\leq
\text{\bf cap} _G (K_1,D)+\text{\bf cap} _G (K_2,D)$$
 for any $\mu$ with supp $\mu \subset K_1\cup K_2.$ Thus capG3) holds.
\enddemo

The equivalent definition of the Green capacity is done by
\proclaim {Theorem 2.5.2.2.(Dual Property)}The following holds
$$\text{\bf cap} _G (K,D)=[\inf_\mu \sup_ {x\in D}
\Pi (x,\mu, D)]^{-1} \tag 2.5.2.4$$
where the infimum is taken over all the mass distributions $\mu$ such that
$\mu (K)=1$
\endproclaim
For proof see, e.g., \cite {La,Ch.2,\S4 it.18}.
For $D= \Bbb R^m,\ m\geq 3,$ the Green capacity is called {\it
Wiener } capacity ({\bf cap}$_m (K)$). It has the following properties in
addition to those of the Green capacity:

capW1) invariance with respect to translations and rotations, i.e.
$$\text {\bf cap}_m (V(K+x_0))= \text {\bf cap}_m (K),$$
where $VK$ and $K+x_0$ are the rotation and the translation of $K$
respectively.

The presence of the properties brings  the notion of capacity closer to
the notion of measure. Thus it is natural to extend the capacity to the
Borel algebra of sets.

The Wiener capacity of an open set is defined as
$$\text {\bf cap}_m (D):=\sup_K \text {\bf cap}_m (K),$$
where the supremum is taken over all  compact $K\Subset D.$

The {\it outer} and {\it inner} capacity of any set $E$ can be defined
by the equalities
$$\overline {\text {\bf cap} }_m (E):=
\inf _{D\supset E}\text {\bf cap}_m (D);\
\underline {\text {\bf cap} }_m (E):=
\sup _{K\subset E}\text {\bf cap}_m (D).$$
A set $E$ is called {\it capacible} if $\overline {\text {\bf cap} }_m (E)=
\underline {\text {\bf cap} }_m (E)$

\proclaim {Theorem 2.5.2.3.(Choquet)}Every set $E$ belonging to the Borel
ring
is capacible.
\endproclaim
For proof see, e.g., \cite {La,Ch2, Th.2.8}.

Sets which have  ``small size'' are sets of zero capacity.
We emphasize the following properties of these sets:

capZ1)If ${\text {\bf cap} }_m (E^j)=0,\ j=1,2, ...$ then
$\text {\bf cap} _m (\cup_1^\infty E^j)=0;$

capZ2) The property to have the zero capacity does not depend of type of
the capacity: Green, Wiener or logarithmic capacity that we define below.

\subheading {Example 2.5.2.1} Using Theorem 2.5.2.2 we obtain that
any point has  zero capacity, because for every mass distribution
concentrated in the point the potential is equal to infinity. The
same holds for any set of zero $m-2$ Hausdorff  measure (see
2.5.4). \subheading {Example 2.5.2.2} Any (m-1)- hyperplane or
smooth hypersurface has positive capacity, because the potential
with masses uniformly distributed over the surface is bounded.

The Wiener 2-capacity can be defined naturally only for sets with
diameter less then one, because  the logarithmic potential
 is positive only when this condition holds.

Instead,  one can use the {\it logarithmic} capacity which is defined by
the formulae
$$\text {\bf cap}_l (K):= \exp [-\text {\bf cap}_2 (K)]\tag 2.5.2.5$$
for $K\subset \{|z|<1\}$ and
$$\text {\bf cap}_l (K):=t^{-1}\text {\bf cap}_l (tK)$$
for any other bounded $K,$ where t is chosen in such a way that
$tK\subset  \{|z|<1\}.$

One can check that this definition is correct, i.e. it does not depend on $t.$
 \subheading {2.5.3}
\proclaim {Theorem 2.5.3.1.(Balayage)} Let $D$ be a domain such that
$\partial D\Subset \Bbb R^m$ ,and supp $\mu \Subset D$.
Then there exists a mass distribution $\mu _b$ such that for $m\geq 3$,
or for $m=2$ and for $D$ which is a bounded domain,
the following holds:

bal1) $\Pi (x,\mu _b)< \Pi (x,\mu)$ for $x\in D;$

bal2) $\Pi (x,\mu _b)= \Pi (x, \mu )$ for $x\notin \overline D;$

bal3) supp $\mu _b \subset \partial D;$

bal4) $\mu _b (\partial D) = \mu (D).$

If $m=2$ and the domain is unbounded, a
potential of the form
$$\hat \Pi (z,\mu):=-\int \log|1-z/\zeta|d\mu_{\zeta} .$$
satisfies all the properties.
\endproclaim
\demo {Proof} We will prove this theorem when
$\partial D$ is smooth enough. For
$y\in D, x\in \Bbb R^m \backslash \overline D$ the function
$|x-y|^{2-m}$ is a harmonic function  of $y$ on $D.$

Since $|x-y|^{2-m}\rightarrow 0$ as $y\rightarrow \infty$ we can apply the
Poisson formula (2.4.1.4)  even if $D$ is unbounded. Thus
$$|x-y|^{2-m} =\int_{\partial D} |x-y'|^{2-m}
\frac {\partial G}{\partial n_{y'}} (y,y') ds_{y'}\tag 2.5.3.1.$$
where $G$ is the Green function of $D.$
From this we have
$$\int_D |x-y|^{2-m}d\mu_y =\int_{\partial D}|x-y'|^{2-m}ds_{y'}
 \left (\int_D \frac {\partial G}{\partial n_{y'}} (y,y')d\mu_y\right ).$$
The inner integral is nonnegative, because
$\frac {\partial G}{\partial n} >0$ for $y'\in \partial D$.
Let us denote
$$d\mu_{b,y'}:=\left (\int_D \frac {\partial G}{\partial n} (y,y')
d\mu_{y}\right )ds_{y'}.$$
 Then we  obtain the properties bal2) and bal3).

The potential $\Pi (x,\mu_b)$ is harmonic in $D.$ Thus
the function
$$u(x):=\Pi (x,\mu_b)-\Pi (x,\mu)$$
is a {\it subharmonic} function (see Theorem 2.6.4.1).
Every subharmonic
function satisfies the {\it maximum principle} (see Theorem 2.6.1.2),
 i.e.
$$u(x)< \sup_{y\in\partial D} u(y)=0.$$
Thus the property bal1) is fulfilled.
To prove bal4) we can write the identity
$$\int_{\partial G} d\mu_{b,y'}=\int_G d\mu_y
\int_{\partial G} \frac {\partial G}{\partial n_{y'}}(y,y')dy'.$$
 The inner integral is equal to one identically, because the function
$\equiv 1,\ y\in G$ is
 harmonic. Thus bal4) is true.

Consider now the  special case when $m=2$, and $D$ is an unbounded domain.
Since $\log |1-z/\zeta|\rightarrow 0$ when $\zeta\rightarrow
\infty,$ we obtain an equality like
(2.5.3.1). Repeating the previous reasoning we obtain the last assertion
for $D$ with a smooth boundary.

{\bf Exercise 2.5.3.1.} Check this in details.
\qed
\enddemo
For the general case see \cite {La Ch.4,\S 1};
\cite {Ca,Ch.3,Th.4}.

Pay attention that the swept potential $\Pi (x,\mu _b)$ is also
a solution
of the Dirichlet problem in the domain $D$ and the boundary function
$f(x)=\Pi (x,\mu)$ in the following sense:
\proclaim {Theorem 2.5.3.2.(Wiener)}The equality bal2) holds in the
points $x\in \partial D$ which can be reached by the top of a cone
placed
outside $D.$ For $m=2$ it can fail only for isolated points.
\endproclaim
For proof see \cite {He},\cite {La,Ch.4,\S1,Th.4.3.}

The points of $\partial D$ where the equality bal2) does not hold are
called {\it irregular}.
\proclaim {Theorem 2.5.3.3.(Kellogg's Lemma)}The set of all the irregular points of
$\partial D$ has zero capacity.
\endproclaim
For proof see, e.g., \cite {He}, \cite {La,Ch.4,\S2,it.10}.

One can often compute the capacity using the following
\proclaim {Theorem 2.5.3.4.(Equilibrium distribution)}For any
compact $K$ with  $\text {\bf cap} _m (K) >0$
there exists  a mass distribution $\lambda _K$ such that the
following holds:

eq1) $\Pi (x,\lambda)=1,\ x\in
\overline {D}\backslash E,\ \text {\bf cap} _m (E) =0$;

eq2) $\text {supp} \lambda _K \subset \partial K;$

eq3) $\lambda _K (\partial K) =\text {\bf cap}_m (K).$
\endproclaim
 For proof see \cite {He}, \cite {La,Ch.2,\S1,it.3,Th.2.3}.

 Let us note that the set $E$ in the previous theorem is a set of irregular
points.

The mass distribution $\lambda _K$ is called {\it equilibrium distribution} ,
and the
corresponding potential is called {\it equilibrium potential}.
\subheading {2.5.4} Let $h(x),\ x\geq 0$ be a positive continuous,
 monotonically increasing function which satisfies the condition
$h(0)=0.$ Let $\{K_j^\epsilon\}$ be a family of balls such that their
diameters $d_j:=d(K_j^\epsilon)$ are no bigger then $\epsilon.$

Let us denote
$$m_h (E,\epsilon) :=\inf \sum h(\frac {1}{2} d(K_j^\epsilon)),$$
where the infimum is taken over all the coverings of the set $E$ by the
families $\{K_j^\epsilon\}.$

The quantity
$$m_h (E):=\lim_{\epsilon\rightarrow 0} m_h (E,\epsilon)$$
is called {\it h -Hausdorff measure}  \cite {Ca,Ch.II} .
\proclaim {Theorem 2.5.4.1.(Properties of $m_h$)}
The following properties hold:

h1) monotonicity: $$E_1\subset E_2 \Longrightarrow
m_h (E_1)\leq m_h (E_2);$$

h2) countable additivity:
$$m_h (\cup E_j)=\sum m_h (E_j);\ E_j\cap E_i =\varnothing ,\text { for} i\neq j;\
E_j\in \sigma (\Bbb R ^m).$$
\endproclaim
We will quote two conditions (necessary and sufficient) that connect
the $h$-measure to the capacity (see,\cite {La,Ch.3,\S4,it. 9,10}.

\proclaim {Theorem 2.5.4.2} Let $\text {\bf cap} E= 0.$Then
$m_h (E)=0$ for all $h$ such that
$$\int \limits_0 \frac {h(r)}{r^{m-1} }dr<\infty.$$
 \endproclaim
\proclaim{Theorem 2.5.4.3} Let $h(r)=r^{m-2}$ for $m\geq 3$
and $h(r)=(\log 1/r)^{-1}$ for $m=2.$ If the $h$-measure of a set $E$
is finite, then $\text {\bf cap}_m (E) =0.$
\endproclaim
Side by side with the Hausdorff measure the {\it Carleson measure}(see,
\cite {Ca,Ch.II},
is often considered. It is defined by
$$m_h ^C (E):=\inf \sum h(0.5 d_j),$$
where infimum is taken over all the coverings  of the set $E$
with balls of radii $0.5d_j.$ The inequality
$m_h^C(E) \leq m_h(E)$ obviously holds.
 Let $\beta-mes_C E$ be the Carleson measure for $h=r^\beta.$
The following assertion connects the $\beta-mes_C$ to capacity.
\proclaim {Theorem 2.5.4.4} The following inequalities hold
$$\beta-mes_C E\leq N(m) (\text {\bf cap}_m (E))^{\beta/m-2},\
for\ m\geq 3,\ \beta >m-2;$$
$$ \beta-mes_C (E)\leq 18 \text {\bf cap}_l (E),\
for\ m=2,\ \beta >0,$$
where $N$ depends only on the dimension of the space.
\endproclaim
For proof see \cite {La Ch.III,\S 4,it.10,Corollary 2}.
\subheading {2.5.5}Now we will formulate an analog of the Luzin
theorem for potentials.
\proclaim {Theorem 2.5.5.1} Let supp $\mu =K$ and let the potential
$\Pi (x,\mu) $ be bounded  on $K.$ Then for any $\delta >0$ there exists
a compact $K'\subset K$ such that  $\mu (K\backslash K')<\delta$ and the
potential $\Pi (x,\mu')$ of the measure $\mu':=\mu\mid _K$ (the
restriction of $\mu$ to $K$) is continuous.
\endproclaim
For proof see, e.g., \cite {La,Ch.3,\S 2,it.3,Th.3.6}.
Let us prove the following assertion:
\proclaim {Theorem 2.5.5.2} Let $\text {\bf cap} K>0.$ Then for
arbitrary small $\epsilon >0 $ there exists a measure $\mu$ such that
supp $\mu \subset K$, the potential $\Pi (x,\mu)$ is continuous and
\linebreak $\mu (K) >\text {\bf cap} (K) -\epsilon.$
\endproclaim
 \demo {Proof}Consider the equilibrium distribution $\lambda _K$ on
$K.$ Its potential is bounded by Theorem 2.5.3.4. By Theorem 2.5.5.1
we can find a mass distribution $\mu$ such that $\Pi (x,\mu)$ is
continuous, supp $\mu \subset K$ and
$\mu (K)>\lambda _K (K)-\epsilon = \text {\bf cap} (K) -\epsilon.$\qed
\enddemo

\newpage
\centerline {\bf 2.6.Subharmonic functions}
\subheading {2.6.1}
Let $u(x),\ x\in D\subset \Bbb R^m $ be a measurable
function bounded from above which can be $-\infty$ on a set of no more than zero measure.

Let us denote as
$$\Cal M (x,r,u):=\frac {1}{\sigma _m r^{m-1}}\int_{ S_{x,r}} u(y)ds_y
\tag 2.6.1.1$$
the {\it mean value} of $u(x)$ on the sphere $S_{x,r}:=\{y:|y-x|=r\}.$

The function $\Cal M (x,r,u)$ is defined if $S_{x,r}\subset D$, but it
can be $-\infty$ a priori.

A function $u(x)$ is called {\it subharmonic} if it is upper
semicontinuous , $\not\equiv -\infty$, and for any $x\in D$ there
exists $\epsilon = \epsilon (x)$ such that the inequality
$$u(x)\leq \Cal M (x,r,u) \tag 2.6.1.2$$
holds for all the $r<\epsilon.$

The class of functions  subharmonic in $D$  will be denoted as $SH (D).$
\subheading {Example 2.6.1.1} The function
$$u(x):=-|x|^{2-m}, x\in \Bbb R^m$$
belongs to $SH (\Bbb R^m)$ for $m\geq 3$, and the function
$$u(z):=\log |z|,\ z\in \Bbb R^2$$
is subharmonic in $\Bbb R^2.$

\subheading {Example 2.6.1.2} Let $f(z)$ be a holomorphic function in
a plane domain $D.$ Then $\log|f(z)|\in SH (D).$

\subheading {Example 2.6.1.3} Let $f=f(z_1,z_2,...,z_n) $ be a
holomorphic function of $z=(z_1,...z_n).$ Then
$u(x_1,y_1,...x_n,y_n) :=\log |f(x_1+iy_1,...,x_n +iy_n)|$ is
subharmonic in every pair $(x_j,y_j)$, and, as one can see in future,
in all the variables.

\subheading {Example 2.6.1.4} Every harmonic function is subharmonic as
it follows from Theorem 2.4.1.6.(Mean Value).
\proclaim {Theorem 2.6.1.1.(Elementary Properties)} The following
holds:

sh1) if $u\in SH (D)$ then $Cu \in SH (D)$ for any constant $C\geq 0;$

sh2) if $u_1,u_2 \in SH (D),$ then $u_1+u_2,\ max [u_1,u_2] \in SH (D);$

sh3) suppose $u_n\in SH (D),\ n=1,2,..$ , and the sequence converges to $u$
monotonically decreasing or uniformly on every compact set in $D.$Then
$u\in SH (D);$

sh4) suppose $u(x,y)\in SH(D_1)$ for all $y\in D_2$, and be upper
semicontinuous in $D_1\times D_2.$
Let $\mu$ be a measure in $D_2$ such that
$\mu (D_2)<\infty.$ Then the function
$u(x):=\int u(x,y)d\mu _y$
is subharmonic in $D_1.$

sh5) let $V\in SO(m) $ be an orthogonal transformation of the space $\Bbb R^m$ and $u\in SH
(\Bbb R^m).$Then $u(V\bul)\in SH (\Bbb R^m).$
\endproclaim
All the assertions follow directly from the definition of subharmonic
functions, properties of semicontinuous functions and properties of
the Lebesgue integral. For detailed proof see, e.g.,  \cite {HK,Ch.2}.
\proclaim {Theorem 2.6.1.2.(Maximum Principle)} Let $u\in SH (D),
\ G\subset \Bbb R^m$ and $u(x)\not \equiv const.$ Then the inequality
$$u(x)< \sup_{x'\in \partial D} \limsup_{y\rightarrow x',y\in D} u(y),
\ x\in D$$
holds,
\endproclaim
i.e. the maximum is not attended inside the domain.

The assertion follows from (2.6.1.2) and the upper semicontinuity of
$u(x).$ For  details see \cite {HK,Ch.2}.

Let $K\Subset D$ be a compact set with nonempty interior $\overset
\circ \to K,$ and let $f_n $ be a decreasing
sequence of functions continuous in $K$  that tends to $u\in SH (D).$
Such a sequence exists by Theorem 2.1.2.9.(The second criterion of
semicontinuity).

Consider a sequence $\{H(x,u_n)\}$ of functions which are harmonic in
$\overset \circ\to K$ and $H\mid _{\partial K} = f_n.$
The sequence converges monotonically to a  function $H(x)$ harmonic in $\overset \circ
\to K$ by Theorem 2.3.4.3.(Connection between
convergences), Theorem 2.4.1.8.(Uniform and $\Cal D'$ -convergences)
and Theorem 2.6.1.2. The limit depends only on $u$ as one can see, i.e. it
does not depend on the sequence $f_n.$ This harmonic function
$H(x):=H(x,u,K)$ is called {\it the least harmonic majorant} of $u$ in
$K.$

This name is justified  because of the following
\proclaim {Theorem 2.6.1.3.(Least Harmonic Majorant)}
Let $u\in SH(D)$.Then for any $K\Subset D$
$u(x)\leq H(x,u,K),\ x\in K.$
If $h(x)$ is harmonic in $\overset \circ \to K$ and satisfies the
condition $h(x)\geq u(x),\ x\in \overset \circ\to K,$ then
$H(x,u,K)\leq h(x),\ x\in \overset \circ \to K.$
\endproclaim
For proof see \cite {HK,Ch.3}.
\subheading {2.6.2}Let us study properties of the  mean values of
 subharmonic functions. Let $\Cal M (x,r.u)$ be defined by (2.6.1.1)
and $\Cal N (x,r,u)$ by
$$\Cal N (x,r,u):=\frac {1}{\omega _m r^m }\int_{K_{x,r}}u(y)dy  ,$$
where $\omega _m$ is the volume of the ball $K_{0,1}.$
 \proclaim {Theorem 2.6.2.1.(Properties of Mean Values)}The
following holds:\newline
me1) $ \Cal M (x,r,u)$ and $\Cal N (x,r,u)$ non-decreases  in $r$
monotonically;
\newline
me2) $u(x)\leq \Cal N (x,\bullet)\leq \Cal M (x,\bullet)$;
\newline
me3) $\lim_{r\rightarrow 0} \Cal M (x,r,u)=
\lim_{r\rightarrow 0}\Cal N(x,r,u) =u(x).$
\endproclaim
\demo {Proof} For simplicity let us prove me1) for $m=2.$ We have
$$\Cal M (z_0,|z|,u)=\frac {1}{2\pi}\int_0^{2\pi} u(z_0 + ze^{i\phi})d\phi$$
Since $u(z,\phi):=u(z_0+ze^{i\phi})$ is a family of subharmonic functions
that satisfies  the condition sh4) of Theorem 2.6.1.1, $\Cal M (z_0,|z|,u)$
is subharmonic in $z$ on any $K_{0,r}$. By Maximum Principle (Theorem 2.6.1.2)
we have
$$\Cal M (z_0,r_1,u)=\max_{S_{0,r_1}}\Cal M (z_0,|z|,u)\leq
 \max_{S_{0,r_2}}\Cal M (z_0,|z|,u)= \Cal M (z_0,r_2,u)$$
for $r_1<r_2.$

Monotonicity of  $\Cal N (x,r,u)$ follows from the equality
$$ \Cal N (x,r,u)=m\int_0^1 s^{m-1}\Cal M (x,rs.u)ds\tag 2.6.2.1$$
and monotonicity of $\Cal M (x,r,u).$

The property me2) follows now from the definition of a subharmonic function
and (2.6.2.1).

Let us prove me3). Let $M (u,x,r)$ is defined by (2.1.1.1). We have
$\Cal M (x,r,u)\leq M (u,x,r)$ and $M(u,x,r)\rightarrow u(x)$ because of
upper semicontinuity of $ u(x).$ Thus me2) implies me3).\qed
\enddemo

It is clear from me2) that a subharmonic function is locally summable.
From me3) we have the corollary
\proclaim {Theorem 2.6.2.2.(Uniqueness of subharmonic function)}If $u,v\in
SH(D)$ and $u=v$ almost everywhere, then $u\equiv v.$
\endproclaim
 Let $\alpha (t)$ be defined by the equality (2.3.1.1),
$\alpha_\epsilon (x)$ by (2.3.1.3).

 For a Borel set $E$ let
$$E^\epsilon :=\{x:\exists y\in D:|x-y|<\epsilon\}.$$
This is the $\epsilon$-extension of $E;$ this is ,of course,  an open set.

For an open set $D$ we set
$$D^{-\epsilon}:=\bigcup\limits _{K^\epsilon \subset D}K^\epsilon$$
This is the maximal set  such that its $\epsilon$-extension is a subset
of $D.$

One can see that $D^{-\epsilon}$ is not empty for  small $\epsilon$ and
$D^{-\epsilon}\uparrow D$ when $\epsilon \downarrow 0.$ Therefore
for any $D_1\Subset D$ there exists $\eps$ such that
$D_1\Subset D^{-\epsilon}.$

For $u\in SH (D)$ set
$$u_\epsilon (x):=\int u(x+y)\alpha_\epsilon (y)dy \tag 2.6.2.2$$
which is defined in $D^{-\epsilon}.$

\proclaim {Theorem 2.6.2.3.(Smooth Approximation)}The following holds:

ap1) $u_\epsilon$ is an infinitely differentiable  subharmonic function in any
open set $D_1\subset D^{-\epsilon}.$

ap2) $u_\epsilon\downarrow u(x)$ while $\epsilon \downarrow 0$ for all
$x\in D.$
\endproclaim
\demo {Proof}The property ap1) follows from sh4) (Theorem 2.6.1.1) and
the following equality that one can obtain from (2.6.2.2):
$$u_\epsilon (x)=\int u(y)\alpha_\epsilon (x-y)dy.\tag 2.6.2.3$$

{\bf Exercise 2.6.2.1} Prove this.

Let us prove ap2). From (2.6.2.2) we obtain
 $$u_\epsilon (x)
=\int_0^1\alpha (s) s^{m-1}\Cal M (x,\epsilon s,u)ds\tag 2.6.2.4$$
It follows from the property me1) (Theorem 2.6.2.1) that
$u_{\epsilon_1}\leq u_{\epsilon_2}$ while $\epsilon_1<\epsilon_2.$
Now we pass to limit in (2.6.2.4). Using me3) we have
$\Cal M (x,\epsilon s,u)\downarrow u(x)$. We can pass to the limit under the
integral because of Theorem 2.2.2.2.Thus
$$\lim_{\epsilon\downarrow 0}u_\epsilon (x)=
\int_0^1\alpha (s) s^{m-1}u(x)ds= u(x)$$\qed
\enddemo
\proclaim {Theorem 2.6.2.4.(Symmetry of $u_\epsilon$)} If $u(x)$ depends only
on $|x|$ then
$u_\epsilon $ depends only on $|x|.$
\endproclaim
\demo {Proof} Let $V\in SO(m)$ be a rotation  of $\Bbb R^m$.Then
$$ u_\epsilon(Vx)=\int u(y)\alpha_\epsilon (Vx-y)dy.$$
Set $y=Vy'$ and change the variables. We obtain
$$ u_\epsilon(Vx)=\int u(Vy')\alpha_\epsilon (V(x-y'))dy.$$
Since $\alpha_\epsilon =\alpha_\epsilon(|x|)$ and $u=u(|x|)$,
$\alpha_\epsilon (Vy)=\alpha_\epsilon (y)$ and  $u(Vy)=u(y)$.
Thus $u_\epsilon(Vx)=u_\epsilon(x)$ for any $V$ and thus
$u_\epsilon(x)=u_\epsilon(|x|).$\qed
\enddemo
\subheading {2.6.3}Since a subharmonic function is locally summable and
defined uniquely by its values almost everywhere,every $u\in SH (D)$  corresponds to a (unique) distribution
$$<u,\phi>:=\int u(x)\phi (x)dx,\ \phi\in \Cal D'.$$
\proclaim {Theorem 2.6.3.1.(Necessary Differential Condition for
Subharmonicity)}

If $u\in SH (D)$, then $\Delta u$ is a positive distribution in $\Cal D'(D).$
\endproclaim
\demo {Proof}Suppose for beginning that $u(x)$ has  second continuous
derivatives. By using (2.4.1.5) and (2.4.1.6) we can represent $u(x)$ in the
form
$$u(x)=\Cal M (x,r,u) + \int_{K_{x,r}}G(x,y,K_{x,r})\Delta u (y)dy,
\tag 2.6.3.1$$
where $G$ is negative for all $r.$
Suppose $\Delta u (x)<0$. Then it is negative in $K_{x,r}$ for some $r.$ Thus
the integral in (2.6.3.1) is positive and we obtain that
$u(x)-\Cal M (x,r,u)>0.$ This contradicts the subharmonicity of $u(x).$

Now suppose $u(x)$ is an arbitrary subharmonic function. Then
$\Delta u_\epsilon (x)\geq 0$ for every $x\in D$ when $\epsilon$ is small
enough.For each  $x$ there is a neighborhood $D_x$ such that every
 $u_\epsilon$ defines a distribution from $\Cal D'(D_x).$
Hence $\Delta u_\epsilon (x)$ defines a positive distribution from
$\Cal D'(D_x)$ . Passing to the limit in $u_\epsilon$ when
$\epsilon \downarrow 0$ we obtain in $\Cal D'(D_x)$ a distribution that is
defined by function $u(x).$ Since the Laplace operator is continuous in
any $\Cal D'$ (Theorem 2.3.4.2), $\Delta u >0$ in  $\Cal D'(D_x).$
From Theorem 2.3.5.1 we
obtain that $\Delta u$ is a positive distribution in $\Cal D' (D).$
\qed
\enddemo

The distribution $\Delta u$ can be realized as a measure by
Theorem 2.3.2.2.The measure $(\theta_m)^{-1} \Delta u$ is called the {\it Riesz}
measure of   the subharmonic function $u.$
\proclaim {Theorem 2.6.3.2.(Subharmonicity and Convexity)}
Let $u(|x|)$ be subharmonic in $x$ on $K_{0,R}.$ Then $u(r)$ is convex with
respect to $-r^{2-m}$ for $m\geq 3$ and with respect to $\log r$ for $m=2.$
\endproclaim
\demo {Proof}By Theorem 2.6.2.4 $u_\epsilon (x)$ depends on $|x|$
only, i.e., $u_\epsilon (x)=u_\epsilon(|x|)$, and the function $u_\epsilon(r)$
is smooth. Passing to the spherical coordinates we obtain
$$\Delta u_\epsilon =\frac {1}{r^{m-1}}\frac {\partial}{\partial r}r^{m-1}
\frac {\partial}{\partial r}u_\epsilon (r)\geq 0.$$
By changing variables $r=e^v$ for $m=2$ or $r=(-v)^{\frac {1}{2-m}}$for
$m\geq 3$ we obtain $[u_\epsilon (r(v)]''\geq 0,$ i.e., $u_\epsilon (r(v))$
is convex in $v.$

Passing to the limit on $\epsilon \downarrow 0$ we obtain that $u(r(v))$ is
convex too, as a monotonic limit of convex functions.\qed
\enddemo
\subheading {2.6.4} Now we  will consider the connection between
subharmonicity and potentials.
\proclaim {Theorem 2.6.4.1.(Subharmonicity of -$\Pi$)}$-\Pi (x,\mu,D)\in
SH(D)$
\endproclaim
It is because of GPo1) and GPo3) (Th.2.5.1.1).

The following theorem is inverse to Theorem 2.6.3.1.
\proclaim {Theorem 2.6.4.2.(Sufficient Differential Condition of
Subharmonicity)}

 Let $\Delta u\in \Cal D'(D)$ be a positive distribution.Then there
exists $u_1\in SH (D)$  that realizes $u.$
\endproclaim
\demo {Proof}Set $\mu:=\theta_m^{-1}\Delta u.$ Let $\Omega_1\Subset
\Omega\Subset D$ and $\Pi (x,\mu_{\Omega})$ be the Newtonian (or
logarithmic) potential of $\mu\mid _\Omega.$
By GPo5) (Th.2.5.1.1) the difference
$H:=u+\Pi$ is a harmonic distribution in $\Cal D'(\Omega_1).$ Hence
there exists a ``natural'' harmonic function $H_1$ that realizes $H$
(Theorem 2.4.1.1). Thus the function $u_1:=H_1-\Pi \in SH (\Omega_1)$ and
realizes $u$ in $\Cal D'(\Omega ).$ Since $\Omega$ and $\Omega_1$ can be
chosen such that a neighborhood of any $x\in D$ belongs to $\Omega_1$,
the assertion holds for $D.$\qed
\enddemo
By the way, we showed in this theorem that every subharmonic function can be
represented inside its domain of subharmonicity as a difference of a harmonic
function and a Newton potential. Thus all the smooth properties of a
subharmonic function depend on the smooth properties of the potential only
because any harmonic function is infinitely differentiable.

The following representation determines the harmonic function completely.
\proclaim {Theorem 2.6.4.3.(F.Riesz representation)}Let $u\in SH (D)$ and let $K$ be a compact
Lipschitz subdomain of $D.$ Then
$$u(x)=H(x,u,K)-\Pi(x,\mu_u,K)$$
where $\mu_u$ is the Riesz measure of $u$ and $H(x,u,K)$ the least
subharmonic majorant.
\endproclaim
\demo {Proof}We can prove as above that the function
$H(x):=u(x)+\Pi (x,\mu_u,K)$ is harmonic in $\overset \circ \to K.$ Since
$H(x)\geq u(x)$ we have $H(x)\geq H(x,u,K).$ So we need the reverse
inequality.

Let us write the same equality for $u_\epsilon$ that is smooth.
    $$u_\epsilon :=H(x,u_\epsilon)-\Pi (x,\mu_{u_\epsilon},K).$$
Passing to the limit as $\epsilon\downarrow 0$ we obtain
$$u(x)=H(x,u,K)-\lim_{\epsilon\downarrow 0}\Pi(x,\mu_{u_\epsilon},K),$$
and the  potentials converge because other summands converge. By Gpo5) $\lim_{\epsilon\downarrow 0}\Pi(x,\mu_{u_\epsilon},K)\geq
\Pi (x,\mu_u,K).$ Hence $H(x)\leq H(x,u,K)$\qed
\enddemo
\subheading {2.6.5} In this item we will consider subharmonic functions
in the ball $K_R:=K_{0,R}$  which are  h a r m o n i c  in some neighborhood of the
origin and write $u\in SH(R).$

Set
$$\align M(r,u)&:=\max \{u(x): |x|=r\}\\
\mu (r,u)&:=\mu_u (K_r)\\
\Cal M (r,u)&:=\Cal M (0,r,u)\\
N(r,u)&:=A(m)\int_0^r \frac {\mu (t,u)}{t^{m-1}}dt,
\text {where} A(m)=\max (1,m-2).\tag 2.6.5.1\endalign$$
\proclaim {Theorem 2.6.5.1.(Jensen-Privalov)}For $u\in SH(R)$
$$\Cal M(r,u)-u(0)=N(r,u),\ for\  0<r<R.\tag 2.6.5.2$$
\endproclaim
\demo {Proof} By Theorem 2.6.4.3 we have
$$u(x)=\frac {1}{\sigma_m r}\int_{|y|=r}u(y)\frac {r^2-|x|^2}{|x-y|^m}ds_y +
\int_{K_r}G(x,y,K_r)d\mu_y.$$
For $x=0$ we obtain
$$u(0)=\cases -\int_0^r \left (\frac {1}{t^{m-2}}-\frac {1}{r^{m-2}}\right )
d\mu(t,u)+\Cal M(r,u),&\text {for }m\geq 3;\\
-\int_0^r \log\frac {r}{t}d\mu(t,u)+ \Cal M (r,u),&\text {for }m=2.
\endcases $$
Integrating by parts gives
$$u(0)-\Cal M(r,u)=\cases -\mu (t,u) \left (\frac {1}{t^{m-2}}-\frac {1}{r^{m-2}}\right )\mid _0^r +(m-2)\int_0^r \frac {\mu (t,u)}{t^{m-1}}dt,&\text {for }
m\geq 3;\\
-\mu (t,u) \log\frac {r}{t}\mid _0^r +\int_0^r \frac {\mu (t,u)}{t}dt,&\text
{for }m=2.\endcases\tag 2.6.5.3$$
We have $\mu(t,u)=0$ for  small $t$ because of harmonicity of $u(x).$
Thus (2.6.5.3) implies (2.6.5.2).\qed
\enddemo

\proclaim {Theorem 2.6.5.2.(Convexity of $M (r,u)$ and $\Cal M (r,u)$)} These
functions increase monotonically and are convex with respect to $\log r$ for
$m=2$ and $-r^{2-m}$ for $m\geq 3.$
\endproclaim
\demo {Proof}Consider the case $m=2.$ Set
$M(z):=\max_\phi u(ze^{i\phi}).$
One can see that $M(r)=M(r,u).$

Let $u$ be a continuous subharmonic function.Then $M(z)$ is
subharmonic (Theorem 2.6.1.1, sh5) and continuous because the family
$\{u_\phi (z):= u(ze^{i\phi})\}$ is uniformly continuous.
The function $M (z)$ depends only on $|z|.$ Thus it is convex with respect to
$\log r$ by Theorem 2.6.3.2.

Let $u(z)$ be an arbitrary subharmonic function and $u_\epsilon
\downarrow u$ while $\epsilon \downarrow 0.$ Then $M(r,u_\epsilon)
\downarrow M(r,u)$ by Prop. 2.1.2.7 and is convex with respect to
$\log r$ by sh3, Theorem 2.6.1.1.

If $m\geq 3$ you should consider the function
$M(x):=\max_{|y|=|x|}u(V_y x)$ where $V_y$ is a rotation of $\Bbb R^m$
transferring $x$ into $y.$

The convexity of $\Cal M (r,u)$ is proved analogously.

{\bf Exercise 2.6.5.1} Prove it.

The monotonicity of $M(r,u)$ follows from the Maximum Principle (Th. 2.6.1.2).
The monotonicity of $\Cal M (r,u)$ was proved in Theorem 2.6.2.1.\qed
\enddemo
The following classical assertion is a direct corollary of Theorem 2.6.5.2.
\proclaim {Theorem 2.6.5.3.(Three Circles Theorem of Hadamard)}Let $f(z)$ be
a holomorphic function in the disc $K_R$ and let $ M_f (r)$ be its maximum on the
circle $\{|z|=r\}.$  Then
$$M_f (r)\leq ([M_f(r_1)]^{\log \frac {r_2}{r}} [M_f(r_2)]^
{\log \frac {r}{r_1}})^{ \frac {1}{\log \frac {r_2}{r_1}}}.$$
for $0<r_1\leq r\leq r_2 < R.$
\endproclaim
For proof you should write down the condition of convexity  with respect to $\log r$ of the function
$\log M_f(r)$ which is the maximum of the subharmonic function $\log |f(z)|$
.

{\bf Exercise 2.6.5.2.} Do this.

 For details see \cite {PS, Part I, Sec.III,Ch.6,Problem 304}.
\newpage

\centerline {\bf 2.7. Sequences of subharmonic functions}
\subheading {2.7.1}We will formulate the following analogue for the Montel
theorem of  normal families of holomorphic functions.

The family $$\{u_\alpha,\ \alpha \in A\}\subset SH (D)\tag 2.7.1.1 $$ is called
{\it precompact} in $\Cal D'(D)$ if, for any sequence
$\{ \alpha_n,\ n=1,2,...\} \subset A,$ there exists a subsequence
${\alpha_{n_j},\ j=1,2,...}$ and a function $u\in SH (D)$ such that
$u_{\alpha_{n_j}}\rightarrow u$ in $\Cal D'(D).$
\subheading {Example 2.7.1.1}$u_\alpha :=\log |z-\alpha|,\ |\alpha|<1$ form
a precompact family.
\subheading {Example 2.7.1.2}$u_\alpha:=\log |f_\alpha|$ where
$\{f_\alpha \}$ is a family of holomorphic functions bounded in a domain $D$
form a precompact family.

A criterion of  precompactness is given by
\proclaim {Theorem 2.7.1.1.(Precompactness in $\Cal D'$)} A family (2.7.1.1) is
 precompact iff the conditions hold:\newline
comp1) for any compact $K\subset D$ a constant $C(K)$ exists such that
$$u_\alpha (x)\leq C(K)\tag 2.7.1.2$$
for all $\alpha\in A$ and $x\in K;$\newline
comp2) there exists a compact $K_1\Subset D$
  such that
$$\inf_{\alpha\in A} \max \{u_\alpha (x):x\in K_1\}> -\infty.\tag 2.7.1.3$$
\endproclaim

For proof see \cite {H\"o,Th.4.1.9}.

\proclaim {Theorem 2.7.1.2} Let $u_n\rightarrow u$ in $\Cal D'(K_R).$ Then
$u_n\rightarrow u$ in $\Cal D'(S_r)$ for any $r<R.$
\endproclaim
\demo {Proof}We have $\mu_n\rightarrow \mu.$ Let us choose $R_1$ such that
$r<R_1<R.$ Then
$$u_n(x)=H(x,u_n,K_{R_1})-\Pi (x,\mu_n,K_{R_1})$$
by F. Riesz theorem (Theorem 2.6.4.3).

Now, we have $ \Pi (x,\mu_n,K_{R_1})\rightarrow \Pi (x,\mu,K_{R_1})$ in
$\Cal D'(R_1)$ by GPo6), Theorem 2.5.1.1. Thus
$H(x,u_n,K_{R_1})\rightarrow H(x,u,K_{R_1})$ in $\Cal D'(R_1).$

By Theorem 2.4.1.8 $H(x,u_n,K_{R_1})\rightarrow H(x,u,K_{R_1})$ uniformly
on any compact set in $ K_{R_1}$ , in particular, on $S_r.$ Hence
 $H(x,u_n,K_{R_1})\rightarrow H(x,u,K_{R_1})$ in $\Cal D'(S_r).$ Also
$ \Pi (x,\mu_n,K_{R_1})\rightarrow \Pi (x,\mu,K_{R_1})$ in
$\Cal D'(S_r)$ by GPo6), Theorem 2.5.1.1. Hence,
$u_n\rightarrow u$ in $\Cal D'(S_r).$\qed
\enddemo
We say that a sequence $f_n$ of locally summable functions {\it
converges in }$L_{loc}$ to a locally summable function $f$ if for
any $x\in D$ there exists a neighborhood $V\ni x$ such that
$\int_V|f_n-f|dx\ri 0.$ \proclaim {Theorem 2.7.1.3 (Compactness in
$L_{loc}$)} Under conditions of Theorem 2.7.1.1 the family
(2.7.1.1) is precompact in $L_{loc}.$ \ep For the proof see
\cite{H\"o, Theorem 4.1.9}.
 \proclaim {Theorem 2.7.1.4} Let
$u_n\rightarrow u$ in $\Cal D'(K_R).$ Then $u^+_n\rightarrow u^+$
in $\Cal D'(K_R)$. \ep This is because $u^+_n(x)\leq M,\ x\in K,$
for all compacts $K\Subset K_R.$ \subheading {2.7.2} The following
theorem shows that a subharmonic function is much more
``flexible'' that harmonic or analytic functions. \proclaim
{Theorem 2.7.2.1} Let $D\Subset \Bbb R^m$ be a Lipschitz domain
and let $u \in SH(D)$ satisfy the condition $u(x)<C$ for $x\in D.$
Then for any closed domain $D_1\Subset D$ there exists a function
$\tilde u (x):=\tilde u(x,D_1)$ such that \newline ext1)
$u(x)=\tilde u (x)\ for\ x\in D_1 ;$ \newline ext2) $\tilde u
(x)=C$ for $x\in \partial D;$\newline ext3) $\tilde u \in SH (D)$
and is harmonic in $D\backslash \overline D _1;$
\newline
ext4) $u(x)\leq \tilde u (x)\ for\ x\in D.$\newline
The function $\tilde u$ is defined uniquely.
\endproclaim
\demo {Proof}We can suppose without loss  of generality that $C=0, $
because we can consider the function $u-C.$

Let $u(x)$ be continuous in $\overline D_1.$ Consider a harmonic
function $H(x)$ which is zero on $\partial D$ and  $u(x)$ on $\partial D_1.$
We have $H(x)\geq u(x)$ for $x\in D\backslash D_1$ because of Theorem 2.6.1.3.
Set
$$\tilde u(x)=\cases H(x),&\ x\in  D\backslash D_1;\\
u(x), &\ x\in D_1.\endcases$$
The function $\tilde u(x)$ is  subharmonic in $D$. For $x\notin \partial D_1$
it is obvious, and for $x\in \partial D_1$ it follows from
$$u(x)=\tilde u (x)\leq \Cal M (x,r,u)\leq \Cal M (x,r,\tilde u)$$
for $r$ small enough.

It is easy to check that all the assertions of the theorem are fulfilled
for the function $\tilde u.$

{\bf Exercise 2.7.1.1} Check this.

Let $u(x)$ be an arbitrary subharmonic function. Consider the family
$u_\epsilon$ of smooth subharmonic functions that converges to $u(x)$
decreasing monotonically in  a neighborhood of $\overline D_1.$
 The  sequence $\widetilde {(u_\epsilon)}$ converges monotonically to
a subharmonic function that has all the properties ext1) - ext4).\qed
\enddemo

\proclaim {Theorem 2.7.2.2.(Continuity of $\tilde \bullet$)}
Let $u_n\rightarrow u$ in $\Cal D'(D)$ and $u_n(x)<0$ in $D.$ Then for any
$K\Subset D$ with a smooth boundary $\partial K$
$\widetilde {u_n}(\bullet,K)\rightarrow \tilde u(\bullet,K)$ in
$\Cal D'(D).$
\endproclaim
For proving, we need the following auxiliary statement:
\proclaim {Theorem 2.7.2.3} Let $u_n\rightarrow u$ in $ \Cal D'(D).$ Then
for any smooth surface $S\Subset D$ and any function $g(x)$ continuous in a neighborhood  of $S$  the assertion
$$\int_S u_n(x)g(x)ds_x\rightarrow \int_S u(x)g(x)ds_x\tag 2.7.2.1$$
holds.
\endproclaim
\demo {Proof}Since $u_n\rightarrow u$ in $\Cal D'(D)$ also the Riesz measures
of the functions converge. Hence $\mu_n (K)\leq C(K)$ for some
$K\Supset S.$
 Thus, for the sequence of
 potentials $\Pi (x,\mu_n),$ we have
$$\int_S \Pi (x,\mu_n)g(x)ds_x =\int d(\mu_n)_y \int_S \frac {g(x)ds_x}
{|x-y|^{m-2}}.$$

The inner integral is a continuous function of $y$ as  can be seen by simple estimates.Thus the assertion (2.7.2.1) holds for potentials.
Now, one can represent $u_n$ in the form
$$u_n(x)=H_n(x)-\Pi(x,\mu_n)$$
in K.
The sequence $H_n$ convergences in $\Cal D'$ and, hence, uniformly on $S.$
Thus (2.7.2.1) holds for every $u_n.$\qed
\enddemo
\demo {Proof of Theorem 2.7.2.2}Let $\phi\in \Cal D(D)$ and supp $\phi\subset
\overset \circ \to K.$ Then
$$<\widetilde {u_n},\phi>=<u_n,\phi>\rightarrow <u,\phi>=<\tilde u,\phi>.$$         Let $x\in D\backslash K.$ Then
$$\widetilde u_n (x)=\int _{\partial K}\frac {\partial G}{\partial n_y}(x,y)
u_n(y)ds_y.$$
By Theorem  2.7.2.3. $\widetilde {u_n}(x)\rightarrow \tilde u(x)$ for
$x\in D\backslash K.$ The sequence $\widetilde u_n$  is precompact in
$\Cal D'(D).$ Thus
 every limit $u_0$ of the $\widetilde u_n$ coincides with
$\tilde u (x)$ in $\overset \circ \to K$ and in $D\backslash K.$ Hence,
$u_0\equiv\tilde u$ in $\Cal D'(D).$\qed
\enddemo

\subheading {2.7.3} The property sh2), Theorem 2.6.1.1, shows that maximum of
any finite number of subharmonic functions is a subharmonic function too.
However, it is not so if the number is  not finite.
\subheading {Example 2.7.3.1} Set $u_n(z)=\frac {1}{n}\log |z|,\ n=1,2....$
The functions $u_n\in SH (K_1).$ Taking supremum in $n$ we obtain
$$u(z)=:\sup_{n}u_n (z)=\cases 0,\ &for\ z\neq 0;\\
-\infty\ &for\ z=0.\endcases$$
The function is not semicontinuous, thus it is not subharmonic. However,
it differs from a subharmonic function on a set of zero capacity.
The following theorem shows that this holds in general.
\proclaim {Theorem 2.7.3.1.(H.Cartan)} Let a family $\{u_\alpha \in SH (D),
 \ \alpha\in A\}$ be bounded from above and $u(x):=\sup_{\alpha\in A}u_\alpha (x).$ Then  $u^*\in SH (D)$ and the set \linebreak
$E:=\{x: u^*(x)>u(x)\}$ is  a zero
capacity set.
\endproclaim
For proving this theorem we need an auxiliary assertion
\proclaim {Theorem 2.7.3.2}Let $\Pi (x,\mu_n,D)$ be a monotonically decreasing
sequence of  Green potentials and supp $\mu_n\subset K\Subset D.$
Then there exists a measure $\mu$ such that the inequality
$$\lim_{n\rightarrow\infty}\Pi (x,\mu_n,D)\geq\Pi (x,\mu,D)$$
holds for all $x\in D$ with  equality  outside
 some set
of zero capacity.
\endproclaim
\demo {Proof}The sequence $\Pi (x,\mu_n,D)$ converges monotonically and thus
in $\Cal D'$ (Theorem 2.3.4.3). Then $\mu_n\rightarrow \mu$ in $\Cal D'$
(Theorem 2.2.4.2.) and thus in $C^*$- topology (Theorem 2.3.4.4).
By GPo5)(Theorem 2.5.1.1) we have
$$\lim_{n\rightarrow\infty}\Pi (x,\mu_n,D)\geq\Pi (x,\mu,D).$$
Suppose that the strict inequality holds on some set $E$ of a positive
capacity. By  Theorem 2.5.2.3 one can find a compact set $K\subset E$ such
that {\bf cap}$(K)>0 .$
Then  there exists a measure $\nu$ concentrated
on $E$ such that its potential $\Pi (x,\nu,D)$ is continuous
(Theorem 2.5.5.2). Thus we have
$$\int \Pi (x,\mu,D)d\nu<
\int\lim_{n\rightarrow\infty}\Pi (x,\mu_n,D)d\nu=
\lim_{n\rightarrow\infty}\int\Pi (x,\mu_n,D)d\nu$$
$$\lim_{n\rightarrow\infty}\int\Pi (x,\nu,D)d\mu_n =
\int\Pi (x,\nu,D)d\mu =\int\Pi (x,\mu,D)d\nu.$$
The equalities use Theorem 2.2.2.2.(B. Levy), reciprocity law ( GPo4), Theorem
2.5.1.1, $C^*$- convergence of $\mu_n$ and once more reciprocity law
respectively. So we have a contradiction.\qed
\enddemo
\demo {Proof Theorem 2.7.3.1}Suppose for beginning $u_n(x)\uparrow u(x).$
We can suppose also that $u_n<0.$ For any domain $G\Subset D$ the
sequence $\tilde u_n (x)\rightarrow u(x)$ for $x\in G$ (see Theorem 2.7.2.1),
because $u_n(x)=\tilde u_n (x)$ for $x \in G.$
Since $\tilde u_n =\Pi (x,\tilde \mu_n,D)$ for   $x\in D$,
$\tilde u (x) =\Pi (x,\tilde \mu,D)= u(x)$ for $x\in G$ and coincides with
$\lim_{n\rightarrow\infty}u_n (x)$ outside some set $E_G$ of zero
capacity. Consider a sequence of domains $G_n$ that exhaust $D.$ Then
$u(x)= \lim_{n\rightarrow\infty}u_n (x)$ outside the set
$E:=\cup_{n=1}^{\infty} E_{G_n}$ which has zero capacity by capZ1) (see item
2.5.2).

Now let $\{u_n,\ n=1,2...\}$ be a countable set that satisfies the
conditions of the theorem. Then the sequence
$v_n:=\max\{u_k:k=1,2...n\}\in SH (D) $ and $v_n\uparrow u.$ Applying
the previous reasoning we obtain the assertion of the theorem also in
this case.

Let $\{ u_\alpha,\ \alpha\in A\}$ be an arbitrary set satisfying the condition
of the theorem. By Theorem 2.1.3.2.(Choquet's Lemma) one can find a countable
set $A_0\subset A$ such that
$$(\sup_{A_0}u_\alpha)^* =(\sup_A u_\alpha)^*.$$
 Since $\sup_{A_0}u_\alpha\leq\sup_A u_\alpha$ , we have
$$E:=\{x:(\sup_A u_\alpha)^*>\sup_A u_\alpha\}\subset E_0:=
\{x:(\sup_{A_0} u_\alpha)^*>\sup_{A_0} u_\alpha\}.$$
Thus {\bf cap }$(E)\leq${\bf cap }$(E_0)=0.$\qed
\enddemo
 Corollary of Theorem 2.7.3.1 is
\proclaim {Theorem 2.7.3.3.(H.Cartan +)}Let $\{u_t,\ t\in (0;\infty)\}
\subset SH (D)$ be
a bounded from above family, and $v:=\limsup_{t\rightarrow\infty }u_t.$
Then $v^* \in SH(D)$ and the set $E:=\{x:v^*(x)>v(x)\}$ has zero capacity.
\endproclaim
\demo {Proof}Set $u_n:=\sup_{t\geq n}u_t,\ E_n:=\{x:(u_n)^*>u_n\},\
E:=\cup E_n.$ Since {\bf cap}$(E_n)=0,\ $ {\bf cap}$E =0$ too.

Let $x\notin E.$ Then
$$v(x)=\lim_{n\rightarrow\infty}\sup_{t\geq n} u_t (x)=
\lim_{n\rightarrow\infty}(u_n)^* (x).$$
The function
$$v^*:=\lim_{n\rightarrow\infty}(u_n)^* (x)$$
is the upper semicontinuous regularization of $v(x)$ for all $x\in D.$\qed
\enddemo
In spite of Example 2.7.3.1 we have
\proclaim {Theorem 2.7.3.4.(Sigurdsson's Lemma)}\text {\cite {Si}} Let $S\subset SH(D)$ be
compact in $\Di'.$ Then
$$v(x):=\sup\{u(x):u\in S\}$$
is upper semicontinuous
\endproclaim
and, hence, subharmonic.
\demo {Proof}Note that
$$u_\epsilon (x)=<u,\alpha (x-\bullet)>$$
(see (2.6.2.3),(2.3.2.1)); and it is continuous in $(u,x)$ with respect to the product
topology on $(SH(D)\cap \Di')\times \Rm$ (Theorem 2.3.4.6).

Let $x_0\in D,a\in \Bbb R$ and assume that $v(x_0)<a.$  We have to prove that
there exists a neighborhood $X$ of $x_0$ such that
$$v(x)<a,\ x\in X.\tag 2.7.3.1$$
We choose $\delta >0$ such that $v(x_0)<a-\delta$. If $u^0\in SH(D)$ and
$\epsilon$ is chosen sufficiently small, then
$$u^0 (x_0)\leq u^0_\epsilon (x_0)<a-\delta$$
by Theorem 2.6.2.3.(Smooth Approximation).

Since $u_\epsilon (x)$ is continuous, there exists an open neighborhood
 $U_0 $ of $u^0$ in $SH(D)$ and an open neighborhood $X_0$ of $x_0$ such
that
$$u_\epsilon (x)<a-\delta,\ u\in U_0,\ x\in X_0.$$
The property ap2) (Theorem 2.6.2.3) implies
$$u(x)<a-\delta,\ u\in U_0,\ x\in X_0.\tag 2.7.3.2$$
Since $u^0$ is arbitrary and $S$ is compact, there exists a finite covering
 $U_1,U_2,...,U_n$ of $S$ and open neighborhoods $X_1,X_2,...,X_n$ of $x_0$
such that (2.7.3.2) holds for all $(u,x):u\in U_j, \ x\in X_j , \ j=1,...,n.$
Set $X:=\cap_j X_j.$ Then (2.7.3.1) holds.\qed
\enddemo

\subheading {2.7.4}Now we are going to connect $\Cal D'$-convergence to
 convergence outside  a  zero capacity set, the so called {\it
quasi-everywhere} convergence.
\proclaim {Theorem 2.7.4.1.($\Cal D'$ and Quasi-everywhere Convergence)}
Let $u_n,\ u \in SH (D)$ and $u_n\rightarrow u$ in $\Cal D'(D).$ Then
$u(x) =\limsup_{n\rightarrow\infty }u_n (x)$ quasi-everywhere and
$u(x) =(\limsup_{n\rightarrow\infty }u_n (x))^*$ everywhere in $D.$
 \endproclaim
For proof we need the following assertion in the spirit Theorem 2.7.3.2.
\proclaim {Theorem 2.7.4.2} Let $\mu_n\rightarrow\mu$ in $\Cal D'(D)$
and supp $\mu_n \subset K\Subset D.$ Then
 $$\liminf_{n\rightarrow\infty}\Pi (x,\mu_n,D)\geq\Pi (\mu,D)$$
with equality  quasi-everywhere.
\endproclaim
\demo {Proof}The inequality was in GPo5), Theorem 2.5.1.1.

Suppose the set
$$E:=\{x:\liminf_{n\rightarrow\infty}\Pi (x,\mu_n,D)>\Pi (x,\mu,D)$$
has a positive capacity. By Theorem 2.5.2.3 one can find a compact set
$K\subset E$ such that {\bf cap}$(K)>0.$ By Theorem 2.5.5.2 one can find
a measure $\nu$ concentrated on $K$ with  continuous potential. As in proof
of Theorem 2.7.3.2 we have
$$\int \Pi (x,\mu,D)d\nu<
\int\liminf_{n\rightarrow\infty}\Pi (x,\mu_n,D)d\nu\leq
\liminf_{n\rightarrow\infty}\int\Pi (x,\mu_n,D)d\nu=$$
$$\liminf_{n\rightarrow\infty}\int\Pi (x,\nu,D)d\mu_n =
\int\Pi (x,\nu,D)d\mu =\int\Pi (x,\mu,D)d\nu.$$
The second inequality uses Theorem 2.2.2.3.(Fatou's Lemma). The equalities use
the  reciprocity law ( GPo4), Theorem 2.5.1.1),
$C^*$- convergence of $\mu_n$ and
once more reciprocity law respectively. So we have a contradiction.\qed
\enddemo
\demo {Proof of Theorem 2.7.4.1}Let $D_1\Subset D.$ Then the sequence
$u_n$ is bounded in $D_1$ by Theorem 2.7.1.1.  We can assume that
$u_n (x)<0 $ for $x\in D_1.$

 For any domain $G\Subset D_1$ the
sequence $\tilde u_n (x,G)\rightarrow u(x)$ in $\Cal D'(D_1 )$ by Theorem
2.7.2.2. We  also have the equality
$\tilde u_n =-\Pi (x,\tilde \mu_n,D_1).$
Thus $ \tilde \mu_n\rightarrow \tilde \mu$ in $\Cal D'(D_1).$
By Theorem 2.7.4.2
 $\liminf_{n\rightarrow\infty}\Pi (x,\tilde \mu_n,D_1)=\Pi (x,\tilde \mu,D_1)$
quasi-everywhere in $D_1$. Hence
$$\limsup_{n\rightarrow\infty}u_n= u\tag 2.7.4.1$$
 quasi-everywhere in $G$  because $u_n(x)=\tilde u_n (x)$ for $x \in G.$

Consider a sequence of domains $G_n$ that exhaust $D.$
Then (2.7.4.1) holds outside a set $E_n$ of  zero capacity and (2.7.4.1)
holds in $D$ outside the set
$E:=\cup_{n=1}^{\infty} E_n$ which has zero capacity by capZ1)
(see item 2.5.2), i.e., quasi-everywhere.\qed
\enddemo
\subheading {2.7.5}Now we connect the convergence of subharmonic functions in
$\Cal D'$ to the convergence relative to
the Carleson measure (see 2.5.4).

We say that a sequence of  functions $u_n$ {\it converges} to a function $u$
{\it relative} to
the $\alpha$- Carleson measure if the sets
$E_n:=\{x:|u_n(x)-u(x)|>\epsilon\}$ possess the property
$$\alpha -mes_C E_n \rightarrow 0.\tag 2.7.5.1$$
\proclaim {Theorem 2.7.5.1.($\Cal D'$ and $\alpha-mes_C $ Convergences)}
Let $u_n,u \in SH (D)$ and $u_n\rightarrow u$ in $\Cal D'(D).$ Then for an
every
$\alpha>0$ and  every domain $G\Subset D$ $u_n\rightarrow u$
relative to the $(\alpha + m -2)$-Carleson measure.
\endproclaim
  For proving this theorem we need some auxiliary definitions and assertions.

Let $\mu$ be a measure in $\Bbb R^m.$ We will call a point $x\in \Bbb R^m$
$(\alpha,\alpha',\epsilon)$-{\it normal} with respect to the measure $\mu$ ,
$(\alpha<\alpha')$ if the inequality

$$\mu _x (t):=\mu (K_{x,t})<\epsilon ^{-\alpha '}t^{\alpha + m - 2 }
$$ holds for all $t<\epsilon.$
\proclaim {Theorem 2.7.5.2} In  any $(\alpha,\alpha',\epsilon)$-normal point
the following inequality holds
$$\align -\int_{K_{z,\epsilon}}[\log |z-\zeta|-\log \epsilon]d\mu_\zeta&\leq
C\epsilon^{\alpha-\alpha'},\ for\ m=2;\\
\int_{K_{x,t}}[|x-y|^{2-m}-\epsilon ^{2-m}]d\mu_y&\leq
C\epsilon^{\alpha-\alpha'},\ for\ m\geq 3;\endalign$$
while $C=C(\alpha,m)$ depends on $\alpha$ and $m$ only.
\endproclaim
\demo {Proof}Let us consider the case $m=2.$ We have

$$\int_{K_{z,\epsilon}}\log \frac {\epsilon}{|z-\zeta|}d\mu_\zeta =
\int_0^\epsilon \log \frac {\epsilon}{t}d\mu_z (t).$$
Integrating by parts we obtain
$$\int_{K_{z,\epsilon}}\log \frac {\epsilon}{|z-\zeta|}d\mu_\zeta =
 \log \frac {\epsilon}{t}\mu_z (t)\mid_0^\epsilon +
\int_0^\epsilon \frac {\mu_z (t)}{t}dt\leq$$
$$\leq\epsilon^{-\alpha'}\int_0^\epsilon t^{\alpha -1}dt=
\frac {1}{\alpha}\epsilon^{\alpha-\alpha'}.$$
 Let us consider the case $m\geq3.$ We have

$$\int_{K_{x,t}}[|x-y|^{2-m}-\epsilon ^{2-m}]d\mu_y=
\int_0^\epsilon(t^{2-m}-\epsilon ^{2-m})d\mu_x (t)=$$
$$=(t^{2-m}-\epsilon ^{2-m})\mu_x (t)\mid_0^\epsilon +
(m-2)\int_0^\epsilon \frac {\mu_x (t)}{t^{m-2}}dt\leq$$
$$\leq\frac {m-2}{\epsilon ^{\alpha'}}\int_0^\epsilon t^{\alpha -1}dt=
\frac {m-2}{\alpha}\epsilon^{\alpha-\alpha'}$$
\qed
\enddemo
\proclaim {Theorem 2.7.5.3.(Ahlfors-Landkof Lemma)}Let a set
$E\subset \Bbb R^m$ be covered by balls with bounded radii such that
every point is a center of a ball. Then there exists an at most countable
subcovering of the same set with  maximal multiplicity $cr=cr(m),$
\endproclaim
i.e.,  every point of $E$ is covered no more than $cr$ times.

For proof see \cite {La, Ch.III, \S 4, Lemma 3.2}.
\proclaim {Theorem 2.7.5.4}Let $K\Subset D.$ The set
$E:=E(\alpha,\alpha',\epsilon,\mu )$ of points that belong to $K$ and are not
$(\alpha,\alpha',\epsilon)$-normal with respect to $\mu$ satisfies the
condition
$$(\alpha +m-2)-mes_C E\leq cr(m)\epsilon^{\alpha '}\mu (K^\epsilon)
\tag 2.7.5.2$$
where $K^\epsilon$ is the $2\epsilon$ -extension of $K.$
\endproclaim
\demo {Proof} Let $x\in E.$ Then there exists $t_x$ such that
$$\mu_x (t_x)\geq t_x ^{\alpha +m-2}\epsilon ^{-\alpha '}.$$
Thus every point of $E$ is covered by a ball $K_{x,t_x}.$
By the Ahlfors-Landkof lemma (Theorem 2.7.5.3) one can find a no more than
$cr$-multiple subcovering  $\{K_{x_j,t_{x_j}}\}.$ Then we have
$$\sum_j t_{x_j}^{\alpha+m-2}\leq cr(m)\epsilon ^{\alpha '}\mu(K^\epsilon).$$
By definition of the Carleson measure we obtain (2.7.5.2).\qed
\enddemo
\proclaim {Theorem 2.7.5.5} Let $\mu_n\rightarrow \mu$ in
$\Cal D'(\Bbb R^m)$
and supp $\mu_n\subset K\Subset \Bbb R^m.$ Then for every $\alpha>0$
and $G\Subset \Bbb R^m$
$\Pi(x,\mu_n)\rightarrow \Pi(x,\mu)$ relative to the
$(\alpha+m-2)$-Carleson measure.
\endproclaim
 \demo {Proof} Let $m=2.$ Set
$$\log_\epsilon |z-\zeta|=\cases \log |z-\zeta|,\ &for\ |z-\zeta|>\epsilon\\
\log \epsilon,\ &for\ |z-\zeta|\leq \epsilon\endcases.$$
 This function is continuous for $(z,\zeta)\in K\times K.$

Set $\nu_n:=\mu_n -\mu.$ Then we have
$$-\int \log |z-\zeta|d(\mu_n)_\zeta +\int \log |z-\zeta|d\mu_\zeta =
-\int \log |z-\zeta|d(\nu_n)_\zeta=$$
$$=-\int\log_\epsilon |z-\zeta|d\nu_n -
\int_{K_{z,\epsilon}}[\log |z-\zeta|-\log \epsilon]d\nu_n.$$
The function $\log_\epsilon |z-\zeta|$ is continuous in $\zeta$ uniformly
over $z\in K.$ Thus the sequence
 $$\Pi_\epsilon (z):= \int\log_\epsilon |z-\zeta|d\nu_n$$
converges uniformly to zero on $K.$
Suppose now that $z\notin E(\alpha,\alpha',\epsilon,\mu)\cup
 E(\alpha,\alpha',\epsilon,\mu_n),$ i.e., it is an $(\alpha,\alpha',\epsilon)$-
normal point for $\mu$ and $\mu_n.$ By Theorem 2.7.5.2 we have
$$\int_{K_{z,\epsilon}}[\log |z-\zeta|-\log \epsilon]d\nu_n<2C\epsilon ^
{\alpha-\alpha '}.$$
Thus for sufficiently large $n>n_0 (\epsilon) $
$$|\Pi (z,\mu_n)-\Pi (z,\mu)|=$$
$$=|\int \log |z-\zeta|d(\mu_n)_\zeta -\int \log |z-\zeta|d\mu_\zeta|<\delta
=\delta (\epsilon)$$
while  $z\notin E(\alpha,\alpha',\epsilon,\mu)\cup
 E(\alpha,\alpha',\epsilon,\mu_n):=E_n(\epsilon).$

By  Theorem 2.7.5.3 the Carleson measure of $E_n(\epsilon)$ satisfies the
inequality
$$\alpha-mes_C E_n (\epsilon)\leq cr(m)\epsilon ^{\alpha '}
[\mu(K)+\mu_n(K)]\leq
C\epsilon ^{ \alpha '}:=
\gamma (\epsilon)$$
where $C=C(K)$ does not depend on $n$ because $\mu_n (K)$ are bounded
uniformly.

Hence, for any $\epsilon>0$ the set
$$ E_n' (\epsilon):=\{z:|\Pi (z,\mu_n)-\Pi (z,\mu)|>\delta (\epsilon)\}$$
satisfies the condition
$$\alpha-mes_C E_n' (\epsilon)\leq \gamma (\epsilon).\tag 2.7.5.2$$
while $n>n_0=n_0(\epsilon).$

Let us show that $\Pi (z,\mu_n)\rightarrow\Pi (z,\mu)$ relative to
$\alpha-mes_C$ on $K.$ Let $\gamma_0,\delta_0$ be arbitrary small. One can
find $\epsilon$ such that $\delta(\epsilon)<\delta_0,\gamma (\epsilon)
<\gamma_0.$ One can find $n_0=n_0(\epsilon)$ such that (2.7.5.2) is fulfilled.
Now the set
$$ E_{n,\delta_0} :=\{z:|\Pi (z,\mu_n)-\Pi (z,\mu)|>\delta _0\}$$
is contained in $E_n'(\epsilon).$
Thus $\alpha-mes_ C E_{n,\delta_0}<\gamma_0$ and this implies the convergence
relative to $ \alpha-mes_ C.$
An analogous reasoning works  for $m\geq3.$
\qed
\enddemo
\demo {Proof of Theorem 2.7.5.1}Let $u_n\rightarrow u$ in $\Cal D'.$
 One can assume that $u_n,\ u$ are potentials
on any compact set(Theorem
2.7.2.2). Hence, by Theorem 2.7.5.5 it converges relative
$(\alpha + m-2) -mes_C.$\qed
\enddemo
\newpage

\centerline {\bf 2.8.Scale of growth. Growth characteristics of subharmonic
functions}
\subheading {2.8.1} Let $A$ be a class of nondecreasing functions
$a(r),\ r\in (0,\infty)$ such that $a(r)\geq 0$ and $a(r)\rightarrow \infty$
when $r\rightarrow \infty.$ The quantity
$$\rho[a]:=\limsup_{r\rightarrow \infty}\frac {\log a(r)}{\log r}\tag 2.8.1.1$$
is called the {\it order} of $a(r).$

Suppose $\rho:=\rho[a]<\infty.$ The number
$$\sigma [a]:=\limsup_{r\rightarrow \infty}\frac { a(r)}{ r^\rho}\tag 2.8.1.2$$
is called the {\it type number}.

If $\sigma [a]=0$, we say $a(r)$ has {\it minimal type}. If
$0<\sigma [a]<\infty$, $a(r)$ has {\it normal type}. If $ \sigma [a]=\infty,$
it has {\it maximal type}.
\subheading {Example 2.8.1.1}Set $a(r):=\sigma_0 r^{\rho_0}.$ Then
$\rho[a]=\rho_0,\  \sigma [a]=\sigma_0 .$
\subheading {Example 2.8.1.2}Set $a(r):=(\log r)^{-1}r^{\rho_0}.$ Then
$\rho[a]=\rho_0,\  \sigma [a]=0 .$
\subheading {Example 2.8.1.3}Set $a(r):=(\log r) r^{\rho_0}.$ Then
$\rho[a]=\rho_0,\  \sigma [a]=\infty .$
\proclaim {Theorem 2.8.1.1(Convergence Exponent) }The following equality holds:
$$ \rho[a]=\inf\{\lambda:\int^\infty \frac {a(r)dr}{r^{\lambda + 1}}<\infty \}.
\tag 2.8.1.3$$
If the integral converges for $\lambda= \rho[a]$, $ a(r)$ has minimal type.
\endproclaim

{\bf Exercise 2.8.1.1.} Prove this.

For proof see, e.g., \cite {HK,\S 4.2}.

\subheading {Example 2.8.1.4}Let $r_j,\ j=1,2,...$ be a nondecreasing
sequence of positive numbers. Let us concentrate the unit mass in every point
$r_j$ and define a mass distribution
$$n(E):=\{ \text {the number points of the sequence $\{r_j\}$ in $E$}\},\
E\subset \Bbb R .$$
Then
$$\int_0^\infty \frac {dn}{r^\lambda}=\sum_1^\infty \frac {1}{r_j^\lambda}.
\tag 2.8.1.4$$
The infimum of $\lambda$ for which the series in (2.8.1.4) converges is
usually called the {\it  convergence exponent} for the sequence $\{r_j\}$
\cite { PS P.I,Sec.1,Ch.III,\S 2}.
Using integrating by parts one can transform the integral in (2.8.1.4) to
the integral of the form (2.8.1.3) where $a(r)=n((-\infty,r)).$
Theorem 2.8.1.1 shows that the convergence exponent coincides with the
order of this $a(r).$

A function $\rho(r)$ is called a {\it proximate order} with respect to order
$\rho$ if

po1) $\rho(r)\geq 0$

po2) $\lim_{r\rightarrow \infty}\rho(r)=\rho$

po3) $\rho (r)$ has a continuous derivative on $(0,\infty)$

po4) $\lim_{r\rightarrow \infty}r \log r \rho'(r)=0.$\newline
Two proximate orders $\rho_1 (r)$ and $\rho_2 (r)$ are called
{\it equivalent}, if
$$ \rho_1 (r)-\rho_2 (r)=o\left(\frac {1}{\log r}\right).\tag 2.8.1.5 $$
For $a\in A$ set
$$\sigma [a,\rho (r)]:=\limsup_{r\rightarrow \infty}\frac { a(r)}
{ r^{\rho(r)}}.\tag 2.8.1.6$$
It is called  a {\it type number  with respect to a proximate order}
$\rho(r).$ It is clear that this type number is the same for  equivalent
proximate orders.
\proclaim {Theorem 2.8.1.2.(Proper Proximate Order)}Let $a\in A$ and
$\rho [a]=\rho<\infty.$ Then there exists a proximate order $\rho(r)$ such
that
$$0<\sigma [a,\rho(r)]<\infty.\tag 2.8.1.7$$
\endproclaim
For proof see \cite {L(1980), Ch.1, Sec.12, Th.16}.

If a proximate order satisfies  the condition (2.8.1.7), we will call it
the {\it proper} proximate order of $a(r)$ (p.p.o.).
The function $r^{\rho (r)}$ inherits a lot of useful properties of the power
function $r^\rho.$
\proclaim {Theorem 2.8.1.3.(Properties of P.O)}The following holds:

ppo1) the function $V(r):=r^{\rho (r)}$
increases monotonically for  sufficiently large values of $r.$

ppo2) for $q<\rho + 1$,
$$\int_1^r t^{\rho (t)-q}dt \sim \frac {r^{\rho (r)+1-q}}{\rho +1-q}$$
 and for $q>\rho+1,$
$$\int_r^\infty t^{\rho (t)-q}dt \sim \frac {r^{\rho (r)+1-q}}{q-\rho -1}$$
as $r\rightarrow \infty.$

ppo3) the function $L(r):=r^{\rho (r)-\rho}$ satisfies the condition
$\forall \delta >0,\  L(kr)/L(r)\rightarrow 1$ when $r\rightarrow
\infty$  uniformly for $k\in [\frac {1}{\delta},\delta ].$
\endproclaim

{\bf Exercise 2.8.1.2.} Prove these properties.

For proof see, e.g., \cite {L(1980), Ch.2, Sec12}.
The following assertion allows to replace any p.o. for
a smooth one.
\proclaim {Theorem 2.8.1.4.(Smooth P.O)}Let $\rho(r)$ be an arbitrary p.o.
There exists an infinitely differentiable equivalent p.o. $\rho_1 (r)$ such
 that
$$r^k \log r \rho_1 ^{(k)}(r)\rightarrow 0,\ k=1,2,...\tag 2.8.1.8$$
when $r\rightarrow \infty.$
\endproclaim
\demo {Proof}Let $\alpha_\epsilon$ be defined by (2.3.1.3). Set
$\epsilon:=0.5,\ po(x):=\rho(e^x)$ and
$$po_1(x):=po(n)+[po(n+1)-po(n)]\int_n^x \alpha_{0.5}(t+0.5) dt$$
for $x\in [n,n+1).$
The function $po_1(x)$ is continuous and infinitely differentiable due to
properties of  $\alpha_\epsilon$ and $ po_1(n)=po(n)$ for $n=1,2,... $ .
By property po3) of p.o. we have
$$(n+1)|po(n+1)-po(n)|\leq \frac {n+1}{n}\max_{y\in [n,n+1]}|y\cdot po'(y)|
\rightarrow 0$$
as $n\rightarrow \infty.$ Thus
$$\max_{y\in [n,n+1]}|y\cdot po^{(k)}_1(y)|\leq const\cdot
(n+1)|po(n+1)-po(n)|\rightarrow 0$$
as $n\rightarrow \infty.$

So $\rho_1 (r):=po_1 (\log r)$ is a p.o. that satisfies (2.8.1.8). Let us show
that it is equivalent to $\rho(r).$ Indeed
$$|po(x)-po_1(x)|=|\int_n^x [po(y)-po_1(y)]'y \frac {dy}{y}|\leq$$
$$\leq\max_{y\in [n,n+1]}[|y\cdot po'(y)|+|y\cdot po'_1(y)|]\log
\frac {n+1}{n}=o\left(\frac {1}{x}\right ),$$
when $x\in [n,n+1]$ and $n\rightarrow \infty.$
\qed
\enddemo
We will further need (in 2.9.3) the following assertion
\proclaim {Theorem 2.8.1.5.(A.A.Goldberg)}Let $\rho (r)\rightarrow \rho$ be
a p.o., and let $f(t)$ be a function that is locally summable on $(0,\infty)$  and such that
$$ \lim_{t\rightarrow 0}t^{\rho +\delta}f(t)=
\lim_{t\rightarrow \infty}t^{\rho + 1 + \gamma}f(t)=0\tag 2.8.1.9$$
for some $0<\delta,\gamma<1.$

Then
$$ \align \lim_{r\rightarrow \infty}
r^{-\rho(r)}\int_{cr^{-1}}^x (rt)^{\rho (rt)}f(t)dt&=
\int_{0}^x t^\rho f(t)dt\\
\lim_{r\rightarrow \infty}
r^{-\rho(r)}\int_{x}^\infty (rt)^{\rho (rt)}f(t)dt&=
\int_{x}^\infty t^\rho f(t)dt\tag 2.8.1.10\endalign$$
for any $c>0$ and any $x\in (0,\infty).$
\endproclaim
\demo {Proof}Set
$$I(r):=\int_{cr^{-1}}^\infty \frac {(rt)^{\rho (rt)}}{r^{\rho(r)}}f(t)dt.$$
It will be enough to prove that
$$\lim_{r\rightarrow \infty}I(r)=\int_{0}^\infty t^\rho f(t)dt\tag 2.8.1.11$$
because  both functions
$$f_0 (t,x):=\cases f(t),\ &for\ t\in (0,x)\\ 0 &for\ t\in [x,\infty)
\endcases$$
and $f_\infty (t,x):=f(t)-f_0 (t,x)$ also satisfy the condition of the theorem.

Let us represent the integral as the following sum:

$$I(r):=\int_{cr^{-1}}^\infty \frac {(rt)^{\rho (rt)}}{r^{\rho(r)}}f(t)dt=
I_1 (r,\epsilon)+I_2 (r,\epsilon)+I_3 (r,\epsilon),\tag 2.8.1.12$$
where
$$I_1 (r,\epsilon):=\int_{cr^{-1}}^\epsilon
\frac {(rt)^{\rho (rt)}}{r^{\rho(r)}}f(t)dt$$
$$I_2 (r,\epsilon):=\int_{\epsilon}^{\epsilon^{-1}}
\frac {(rt)^{\rho (rt)}}{r^{\rho(r)}}f(t)dt$$
$$I_3 (r,\epsilon):=\int_{\epsilon^{-1}}^{\infty}
\frac {(rt)^{\rho (rt)}}{r^{\rho(r)}}f(t)dt.$$
We can represent $I_2 (r,\epsilon)$ in the form
$$I_2 (r,\epsilon)=\int_{\epsilon}^{\epsilon^{-1}}
\frac{L(rt)}{L(r)}t^\rho f(t)dt.$$
By ppo3) (Theorem 2.8.1.3),
$$\lim_{r\rightarrow \infty}I_2 (r,\epsilon)=
\int_{\epsilon}^{\epsilon^{-1}} t^\rho f(t)dt.\tag 2.8.1.13$$
Let us estimate the ``tails''. From (2.8.1.9) we have
$$|f(t)|\leq Ct^{-\rho -\delta}$$
for $0<t\leq \epsilon$ where $C$ does not depend on $\epsilon$ and
$$|f(t)|\leq Ct^{-\rho -1 - \gamma}$$
for $t\geq \epsilon ^{-1}.$
We have
$$|I_1 (r,\epsilon)|\leq C\int_{cr^{-1}}^\epsilon
\frac {(rt)^{\rho (rt)}}{r^{\rho(r)}}t^{-\rho-\delta}dt:= CJ_1(r,\epsilon)$$
and
$$\limsup_{r\rightarrow \infty}|I_1 (r,\epsilon)|\leq
C\lim_{r\rightarrow \infty}J_1(r,\epsilon)\tag 2.8.1.14$$
Let us calculate the last limit.We perform the change $x=tr$:
$$J_1(r,\epsilon)=r^{-\rho(r)+\rho + \delta -1}
\int_{c}^{\epsilon r}t^{-\rho(x)-(\rho+\delta)}dx.$$
Now we use ppo2) for $q=\rho+\delta$ and ppo3):
$$\lim_{r\rightarrow \infty}J_1(r,\epsilon)=\frac {1}{1-\delta}
\lim_{r\rightarrow \infty}
\frac {(\epsilon r)^{\rho(\epsilon r)-(\rho+\delta)+1}}
{r^{\rho(r)-(\rho + \delta) +1}}=$$
$$=\frac {\epsilon ^{1-\delta}}{1-\delta}
\lim_{r\rightarrow \infty}\frac {L(\epsilon r)}{L(r)}=
\frac {\epsilon ^{1-\delta}}{1-\delta}$$
Substituting in (2.8.1.14) we obtain
$$\limsup_{r\rightarrow \infty}|I_1 (r,\epsilon)|\leq
C\frac {\epsilon ^{1-\delta}}{1-\delta}.\tag 2.8.1.15$$
Analogously one can obtain
$$\limsup_{r\rightarrow \infty}|I_3 (r,\epsilon)|\leq
C\frac {\epsilon ^{\gamma}}{\gamma}.\tag 2.8.1.16$$
 Using (2.8.1.13), (2.8.1.15) and (2.8.1.16)
 one can pass to the limit in  (2.8.1.12) as $r\rightarrow \infty$
 then let $\epsilon \rightarrow 0$, and obtain (2.8.1.11).  \qed
\enddemo

\subheading {2.8.2}Let
$$u(x):=u_1(x)-u_2(x)\tag 2.8.2.1$$
where $u_1,u_2 \in SH(\Bbb R^m),\ u_1(0)>-\infty,\ u_2(0)=0$ and
$\mu_1:=\mu_{u_1},\ \mu_2:=\mu_{u_2}$ are concentrated on disjoined sets.

Let $m=2,\ u_j(z):=\log |f_j(z)|,\ j=1,2$ where $f_j(z),\ j=1,2$ are entire
functions. Then the function $u(z)=\log |f(z)|,$ where $f(z):=f_1(z)/f_2(z),$
is
 meromorphic. The condition for masses means that $f_1$ and $f_2$ have no
common zeros, $u_2 (0)=0$ corresponds to $f_2(0)=1$ and $u_1(0)>-\infty$ means
$f_1(0)\neq 0$.

The class of such functions is denoted as $\delta SH(\Bbb R^m).$
In spite of the standardization conditions the representation (2.8.2.1) is not
unique. However for any pair of representations $u_1-u_2$ and $u'_1-u'_2$
$$u_j(x)- u'_j(x)=H_j(x),\ j=1,2\tag 2.8.2.2$$
where $H_j$ are harmonic and $H_j(0)=0.$

Really, from the equality
$u_1-u_2=u'_1-u'_2$ we obtain $\mu_1-\mu_2=\mu'_1-\mu'_2.$ Using the
Theorem 2.2.1.2.(Jordan decomposition) we obtain
$\mu_1=\mu'_1,\ \mu_2=\mu'_2.$ Thus (2.8.2.2) holds. Obviously $H_2(0)=0.$

Set
$$T(r,u):=\frac {1}{\sigma_m}\int_{|y|=1}\max (u_1,u_2)(ry)dy\tag 2.8.2.3$$
where $ \sigma_m$ is the surface square of the unit sphere.
It is called the{\it Nevanlinna characteristic} of $u\in \delta SH(\Bbb R^m).$

The Nevanlinna characteristic does not depend on the representation (2.8.2.1).
Indeed,
$$\int_{|y|=1}\max (u_1,u_2)(ry)dy=\int_{|y|=1}[(u_1-u_2)^+(ry)-u_2(ry)]dy=$$
$$=\int_{|y|=1}[(u'_1-u'_2)^+(ry)-u'_2(ry)+H_2(rx)]dy=$$
$$\int_{|y|=1}[\max (u'_1,u'_2)(ry)+H_2(rx)]dy=$$ $$=\int_{|y|=1}
\max (u'_1,u'_2)(ry)dy+H_2(0)=\int_{|y|=1}
\max (u'_1,u'_2)(ry)dy.$$

Note also that the class $\delta SH(\Bbb R^m)$ is linear.

Actually, let $u\in \delta SH(\Bbb R^m).$ Then $\lambda u\in \delta SH(\Bbb
R^m)$ for $\lambda>0.$  $-u\in \delta SH(\Bbb R^m),$ since
$$-u(x)=[u_2(x)-u_1(0)]-[u_1(x)-u_1(0)].$$

Let us show that $u_1+u_2\in\delta SH(\Bbb R^m)$ if
$u,v\in \delta SH(\Bbb R^m).$

Set $\nu:=\nu_u +\nu_v,$  where $\nu_u ,\nu_v$ are the corresponding charges.
By Theorem 2.2.1.2.(Jordan decomposition) $\nu=\nu^+ - \nu^-,$  where
$\nu_u ,\nu_v$ are measures concentrated on disjoint sets.

Let $u_1$ be a subharmonic function in $\Bbb R^m $  the mass distribution
of which coincides with $\nu^+.$ \footnote {We will give a construction of such function
for the case of finite order (item 2.9.2), but it is possible actually always, see ,for
example, \cite {HK,Th.4.1}}
Then $u_2:=u_1-(u+v)$ is a subharmonic function with the mass distribution
 $\nu^-.$ Hence
$u(x)+v (x)=[u_1(x)-u_2(0)]-[u_2(x)-u_2(0)].$

\proclaim {Theorem 2.8.2.1.(Properties $T(r,u)$)}The following holds

t1) $T(r,u)$ increases monotonically and is convex with respect to
$-r^{m-2}$ for
$m=2$ and with respect to $\log r$ for $m=2$

t2) For $u\in SH(\Bbb R^m)$, (i.e. $u_2\equiv 0$)
$$T(r,u)= \frac {1}{\sigma_m}\int_{|y|=1} u^+(ry)dy$$

t3) $T(r,u)=T(r,-u)-u_1(0)$

t4) $T(r,u+u')\leq T(r,u)+T(r,u'),\ T(r,\lambda u)=\lambda T(r,u)$ for
$\lambda >0.$
\endproclaim
\demo {Proof}Since $v(x):=\max (u_1,u_2)(x)$ is  subharmonic, t1) follows from Theorem 2.6.5.2.(Convexity of $M(r,u)$ and $\Cal M(r,u)$).

The property t2) is obvious, t3) follows from the equality
$-u(x)=u_2(x)-[u_1(x)-u_1(0)] -u_1(0).$

The properties t4) follow from the properties of maximum and t3).\qed
\enddemo
Set $\rho_T [u]:=\rho[a]$ (see, (2.8.1.1)) where $a(r):=T(r,u).$
It is called the {\it order of $u(x)$ with respect to $T(r).$}
\proclaim {Theorem 2.8.2.2.($\rho_T$-property)}For
$ u_1,u_2\in \delta SH(\Bbb R^m)$ the following inequality holds:
$$\rho_T[u_1+u_2]\leq \max(\rho_T[u_1],\rho_T[u_2]),\tag 2.8.2.4$$
 Equality in (2.8.2.4) is attained if $\rho_T[u_1]\neq\rho_T[u_2].$
\endproclaim
\demo {Proof}Set $u:=u_1+u_2.$ From t3) and t4)
$$T(r,u)\leq T(r,u_1)+T(r,u_2)+O(1)\leq 2\max [ T(r,u_1),T(r,u_2)]+ O(1).$$
From the definition of $\rho_T$ we obtain (2.8.2.4).

Suppose, for example, $\rho_T[u_1]>\rho_T[u_2].$ Let us show that
$\rho_T[u ]=\rho_T[u_1].$ From the equality $u_1=u+(-u_2)$ we obtain
$\rho_T[u_1]\leq \max(\rho_T[u],\rho_T[u_2]$
If $\rho_T[u ]<\rho_T[u_1],$ then from the previous inequality we would have
the contradiction  $\rho_T[u_1 ]<\rho_T[u_1].$ \qed
\enddemo

Let us define $\sigma_T [u]$ by (2.8.1.2) while $\rho:=\rho_T[u].$
Set also $\sigma_T [u,\rho(r)]:= \sigma [a,\rho(r)]$ (see (2.8.1.6)), where
$a(r):=T(r,u).$

The characteristics $\rho_T[u],\ \sigma_T [u],\ \sigma_T [u,\rho(r)]$ are
defined for $u\in \delta SH(\Bbb R^m).$ For the class of subharmonic
function we have the inclusion $ SH(\Bbb R^m)\subset \delta SH(\Bbb R^m)$ and,
of course, all these characteristics can be applied to a subharmonic function.
However, for the class $ SH(\Bbb R^m)$ the standard characteristic of
growth is $M(r,u)$ that we can not apply  to a $\delta-subharmonic$ function
 $u\in \delta SH(\Bbb R^m).$ Thus for  $u\in  SH(\Bbb R^m)$ we define new
{\it characteristics} $\rho_M[u],\ \sigma_M [u],\ \sigma_M [u,\rho(r)]$ in
the same
way by replacing $T(r,u)$ for $M(r,u).$ The following theorem shows
that there is not a big difference between characteristics with respect $T$
and $M$ for $u\in SH(\Bbb R^m).$
\proclaim {Theorem 2.8.2.3.(T and M -characteristics)}Let $u\in SH(\Bbb R^m)$
and $\rho(r)(\rightarrow \rho)$ any p.o.Then

$\rho$MT1)  $\rho_T[u]$ and $\rho_M[u]$ are finite simultaneously and
 $\rho_T[u]=\rho_M[u]:=\rho[u]$

$\rho$MT2) there exists $A:=A(\rho,m)$ such                                    that
$$A\sigma_M [u,\rho(r)]\leq\sigma_T [u,\rho(r)]\leq\sigma_M [u,\rho(r)]$$

\endproclaim
In particular, the last property means that the types with respect to $T(r)$
and $M(r)$ for the same p.o. are minimal, normal or maximal at the same time.
\demo{Proof}From t2), Theorem 2.8.2.1 we have $T(r,u)\leq M(r,u)$ for
 $u\in  SH(\Bbb R^m).$ Thus $ \rho_T[u]\leq\rho_M[u]$, proving the second part
of
$\rho$MT2).

Let $H(x)$ be the least harmonic majorant of $u(x)$ in the ball $K_{2R}.$
By the Poisson formula (Theorem 2.4.1.5) and Theorem 2.6.1.3
$$M(R,u)\leq M(R,H)= \max_{|x|=R}\frac {1}{\sigma_m 2R}\int_{|y|=2R}u(y)
\frac {(4R^2 -|x|^2)}{|x-y|^{m}}ds_y\leq\tag 2.8.2.5$$
$$\leq \frac {2^{m-2}}{\sigma_m}\int_{|y|=1}|u(2Ry)|ds_y=
2^{m-2}[T(2R,u)+T(2R,-u)]=2^{m-2}[2T(2R,u)-u(0)].\tag 2.8.2.5$$
From here one can obtain  $ \rho_T[u]\geq\rho_M[u]$ . The left side of
$\rho$MT2) with $A(\rho,m):=2^{-\rho-m+2}$ follows from the properties
of p.o.\qed
\enddemo

{\bf Exercise 2.8.2.1} Prove the first inequality from $\rho$MT2).

\subheading {2.8.3}
Let $\mu$ be a mass distribution (measure) in $\Bbb R^m$ ($\mu\in
\Cal M(\Bbb R^m)$).The characteristic
$$\rho [\mu]:=\rho[a]-m+2$$
for $a(r):=\mu (K_r)$ (see (2.8.1.1)) is called the {\it convergence
exponent} of $\mu$, and
$$\bar \Delta [\mu]:=\sigma [a]$$
for the same $a$ (see (2.8.1.2)) is called the {\it upper density} of $\mu.$

The least integer number $p$ for which the integral
$$\int^\infty \frac {\mu(t)}{t^{p + m}}dt\tag 2.8.3.1$$
converges is called the {\it genus} of $\mu$ and is denoted  $p[\mu].$

\proclaim {Theorem 2.8.3.1.(Convergence Exponent and Genus)}The following
holds:

ceg1) $p[\mu]\leq \rho[\mu] \leq p[\mu]+1$

ceg2) for $\rho[\mu]=p[\mu]+1,\  \bar\Delta [\mu]=0 .$
\endproclaim
\demo {Proof}From Theorem 2.8.1.1 (Convergence Exponent) we have
$ \rho[\mu] +1+m-2\leq p[\mu]+m.$ Thus $\rho[\mu] \leq p[\mu]+1.$ The same
theorem implies $ \rho[\mu]+m-2+1\geq p[\mu]+m-1.$ Thus
$p[\mu]\leq \rho[\mu],$ and ceg1) is proved.

Let  $\rho(\mu)=p[\mu]+1.$ Then the integral (2.8.3.1) converges for
$p[\mu]=\rho[\mu]-1.$ We use the inequality
$$\int_r^\infty \frac {\mu(t)}{t^{\rho[\mu] + m-1}}dt\geq
\mu (r)\int_r^\infty \frac {dt}{t^{\rho[\mu] + m-1}}dt= \frac {\mu(r)}
{r^{\rho[\mu] + m-2}}(\rho[\mu] + m-2)^{-1}.$$
Since the left side of the inequality tends to zero we obtain
$$\bar\Delta [\mu]=\lim_{r\rightarrow\infty }\frac {\mu(r)}
{r^{\rho[\mu] + m-2}}=0.$$ \qed
\enddemo
Set $$ \bar\Delta [\mu,\rho(r)]:=\sigma [a,\rho(r)+m-2],\tag 2.8.3.2$$
where $a(r):= \mu (r)$ (see (2.8.1.6)).
It is clear that $\rho(r)+m-2$  is also a p.o.
Set as in (2.6.5.1)
$$N(r,\mu):=A(m)\int_0^r \frac {\mu (t)}{t^{m-1}}dt,$$
 where $A(m)=\max (1,m-2).$
Set also
$$\rho_N[\mu]:=\rho[a],\ \bar\Delta_N [\mu,\rho(r)]:=\sigma [a,\rho(r)],$$
where $a(r):=N(r,\mu)$ .
This is the {\it N-order of $\mu$} and the {\it N-type of  $\mu$} with respect
to p.o. $\rho(r).$
\proclaim {Theorem 2.8.3.2.(N-order and Converges Exponent)} The following
holds:

Nce1) $\rho_N[\mu]$  and $\rho[\mu]$ are finite simultaneously and
$\rho_N[\mu]=\rho[\mu]$

Nce2) for $\rho>0$ there exists such $A_j:=A_j(\rho,m),\ j=1,2,$ that
$$A_1\bar\Delta[\mu,\rho(r)]\leq\bar\Delta_N [\mu,\rho(r)] \leq
A_2\bar\Delta[\mu,\rho(r)].$$
\endproclaim
\demo {Proof}We have the inequality
$$N(2r,\mu)\geq A(m)\int_r^{2r} \frac {\mu (t)}{t^{m-1}}dt\geq
 A(m)\mu(r)\int_r^{2r} \frac {dt}{t^{m-1}}\geq$$
$$\geq A(m)B(m)\frac {\mu(r)}{(2r)^{ m-2}},$$
where $B(m):=1-2^{2-m}$ for $m\geq 3$ and $B(2):=\log 2.$

From here one can obtain the inequality
$\rho[\mu]\geq \rho_N[\mu]$ and the left side of Nce2) for
$A_1(\rho,m):=A(m)B(m)2^{-\rho}.$
For proving the opposite inequalities we use the l'H\^opital rule (slightly
improved):
$$\limsup_{r\rightarrow\infty}\frac {N(r,\mu)}{r^{\rho(r)}}\leq
\limsup_{r\rightarrow\infty}\frac {N'(r,\mu)}{(r^{\rho(r)})'}=$$
$$=\limsup_{r\rightarrow\infty}\frac {\mu(r)r^{2-m}}
{r^{\rho(r)}[\rho(r)+r\log r \rho '(r)]}=\frac {1}{\rho}\bar \Delta [\mu].$$
Thus $\rho_N[\mu]\leq \rho [\mu]$ and the right side of Nce2) holds.\qed
\enddemo
We shall denote as
$\delta\Cal M (\Bbb R^m)$ the set of charges (signed measures) of the
form $\nu:=\mu_1-\mu_2$ where
$\mu_1,\mu_2\in \Cal M (\Bbb R^m) .$ Let us remember that $|\nu|\in
\Cal M (\Bbb R^m)$ is
the full variation of $\nu$ (see 2.2.1).
\proclaim{Theorem 2.8.3.3.(Jensen)}Let $u:=u_1-u_2 \in \delta SH (\Bbb R^m)$
and $\nu:=\mu_1-\mu_2 $ be a corresponding charge.Then

J1) $\rho [|\nu|]\leq \max (\rho [\mu_1],\rho [\mu_2])\leq \rho [u]$

J2) $\bar \Delta [|\nu|,\rho(r)]\leq\bar \Delta [\mu_1,\rho(r)]+
\bar \Delta [\mu_2 ,\rho(r)]
\leq A  \sigma _T [u,\rho(r)]$ \newline
for some $A:=A(\rho,m).$
\endproclaim
\demo {Proof}We can suppose without loss of generality that $u(0)=0$ because
the function $u(x)-u(0)$ has the same order and the same number type if
$\rho>0.$ We apply the Jensen-Privalov formula (Theorem 2.6.5.1) to the
functions $u_1,u_2$ and obtain
$$N(r,\mu_j)\leq \Cal M (r,u_j)\leq T(r,u).$$
Thus $N(r,|\nu|)\leq  N(r,\mu_1)+N(r,\mu_2)\leq 2T(r,u).$ From here one can
obtain J1) and J2) for
$\rho_N [|\nu|]$ and $\bar \Delta _N [|\nu|,\rho(r)].$ However, we can delete
the subscript $N$ because of Theorem 2.8.3.2.\qed
\enddemo
\newpage

\centerline {\bf 2.9.The representation theorem of subharmonic functions in $\Bbb R^m.$
}
\subheading {2.9.1}
Set
$$H(z,\cos\gamma,m):=\cases
-\frac {1}{2}\log(z^2-2z\cos\gamma +1),\ &for\ m=2\\
(z^2-2z\cos\gamma +1)^{-\frac {m-2}{2}},\ &for\ m\geq 3
\endcases\tag 2.9.1.1$$
The function $H(z,\cos\gamma,m)$ is holomorphic on $z$ in the disk $\{|z|<1\}.$
It can be represented there in the form
$$H(z,\cos\gamma,m)=\sum_{k=0}^{\infty}C^{\frac {m-2}{2}}_k
(\cos\gamma) z^k \tag 2.9.1.2$$ where every coefficient
$C^{\beta}_k(\bullet),\ k=0,1,...\ $ is a polynomial
 of degree $k.$

Such polynomials are called the {\it Gegenbauer } polynomials
. Note that $C^{\frac
{1}{2}}_k(\bullet)$ are the Legendre
 polynomials and
$$C^{0}_k(\lambda)=\cases 0,\ &for\ k=0\\
\frac {1}{k}\cos(k\arccos \lambda),\ &for\ k\geq 1,\endcases $$
i.e. they are proportional to the { \it Chebyshev} polynomials.

Thus for $m=2$ we have the equality
$$-\frac {1}{2}\log(z^2-2z\cos\gamma +1)=
\sum_{k=1}^{\infty}{\frac {\cos k\gamma}{k}
z^k} $$
that can be checked directly.

Let $x\in \Bbb R^m.$ Set $x^0:=x/|x|.$ Then the scalar product $(x^0,y^0)$ is
equal to $\cos\gamma $ where $\gamma$ is the angle between $x$ and $y.$

Let $\Cal E_m (x)$ be defined by (2.4.1.1). For $m\geq 3$ the function
$\Cal E_m (x-y)$ is
the Green function for $\Bbb R^m.$ One  can see that it is represented in the
form
$$G(x,y,\Bbb R^m):=\Cal E_m (x-y)=-|y|^{2-m}H( |x|/|y|,\cos\gamma,m)$$
where $\cos\gamma=(x^0,y^0).$

For $m=2$ the function $-H(|x|/|y|,\cos\gamma,2)$ plays the same role. Thus we
will denote it as $G(x,y,\Bbb R^2).$
\proclaim {Theorem 2.9.1.1.(Expansion of $G(x,y,\Bbb R^m)$)}The following
holds:
$$G(x,y,\Bbb R^m)=-\sum_{k=0}^{\infty}C^{\frac {m-2}{2}}_k
(\cos\gamma) \frac {|x|^k}{|y|^{k+m-2}}, \tag 2.9.1.3$$
for $|x|<|y|$ , and the functions
$$D_k (x,y):=C^{\frac {m-2}{2}}_k (\cos\gamma) \frac {|x|^k}{|y|^{k+m-2}
}\tag 2.9.1.4$$
are  homogeneous harmonic functions in $x$ and harmonic in $y$ for $y\neq 0.$
\endproclaim
\demo {Proof}The expansion (2.9.1.3) follows from (2.9.1.2).
The function $G(zx,y,\Bbb R^m)$ is harmonic for $|x|<|y|$ and, hence, for
any real $0\leq z<1.$ Hence, for any $\psi\in \Cal D (K_r)$ while $r:=0.5|y|$
the function $g(z):=<G(z\bullet,y,\Bbb R^m),\Delta\psi>=0$ for $z\in (0,1).$
The function $g$ is holomorphic for all the complex $z\in \{|z|<1\}$ because
$G(zx,y,\Bbb R^m)$ is holomorphic.  Thus
$g(z)\equiv 0,$ i.e. all its coefficients are zero.

From the expansion (2.9.1.3) we can see that the coefficients of
$G(zx,y,\Bbb R^m)$ are $D_k (x,y).$ Hence, $<D_k (\bullet,y),\Delta\psi>=0$
for every $\psi\in \Cal D (K_r).$ Thus $D_k (\bullet,y)$ is a harmonic
distribution . By Theorem 2.4.1.1 it is an ordinary harmonic function
for $|x|<0.5|y|. $

$C^{\frac {m-2}{2}}_k (\cos\gamma)$ is a polynomial of degree $k$
with respect to $(x^0,y^0).$ Thus $D_k (x,y)$ is a homogeneous
polynomial of $x$ and is harmonic for all $x.$

Let us prove the harmonicity in $y.$  By Theorem  2.4.1.10
 the function\linebreak
$D_k (y^*,x^0)|y|^{2-m}$ ($*$ stands for inversion) is harmonic in $y.$ We have
$$D_k (y^*,x^0)|y|^{2-m}=|y|^{2-m}D_k (y/|y|^2,x^0)=D_k (x^0,y).$$ \qed
\enddemo
Set
$$H(z,\cos\gamma,m,p)=H(z,\cos\gamma,m)-\sum_{k=0}^{p}C^{\frac {m-2}{2}}_k
(\cos\gamma)z^k \tag 2.9.1.5$$
\proclaim {Theorem 2.9.1.2} The following holds:
$$|H(z,\cos\gamma,m,p)|\leq A_1 (m,p)|z|^{p+1} \tag 2.9.1.6 $$
for $|z|\leq 1/2$, and
 $$|H(z,\cos\gamma,m,p)|\leq A_2 (m,p)|z|^p\tag 2.9.1.7 $$
 for $|z|\geq 2,\ -\pi<\arg z\leq \pi.$
The factor $|z|^p$ should be replaced by $\log |z|$ if \linebreak $m=2,\ p=0.$
\endproclaim
\demo {Proof}Consider the function $\phi (z):=H(z,\cos\gamma,m,p)z^{-p-1}.$
It is holomorphic in the disk $\{|z|\leq1/2\}.$ We apply the maximum
principle and obtain (2.9.1.6) where
$$A_1 (m,p)=2^{p+1}\max_{|z|=1/2} |\phi (z)|.$$
For proving (2.9.1.7) we consider the function
 $\psi (z):=H(z,\cos\gamma,m,p)z^{-p}$ that is holomorphic in the domain
$D:=\{z:|z|\geq 2,\ -\pi<\arg z\leq \pi\}$ and continuous in its closure.
Applying the maximum principle we obtain (2.9.1.7) where
 $$A_2 (m,p)=2^{p}\max_{z\in \partial D} |\psi (z)|.$$
\qed
\enddemo
Set
$$G_p(x,y,m):=-|y|^{2-m}H(|x|/|y|,\cos\gamma,m,p)$$
where $\cos\gamma=(x^0,y^0).$

Note the equality
$$G_p(x,y,m)=G (x,y,\Bbb R^m)+\sum_{k=0}^{p}D_k (x,y).$$

{\bf Exercise 2.9.1.1.} Check this using (2.9.1.3),(2.9.1.4) and (2.9.1.5).

It looks like a Green function for $\Bbb R^m$ but it tends more quickly to zero
at infinity and generally speaking it is not negative.

For $m=2$ it can be represented in the form
$$G_p(z,\zeta,2)=\log |E(z/\zeta, p)|$$
where $E(z/\zeta, p)$ is the primary Weierstrass factor:
$$E(z/\zeta, p) :=\left(1-\frac {z}{\zeta}\right)\exp
\left [\left(\frac {z}{\zeta}\right )+
\frac {1}{2}\left(\frac {z}{\zeta}\right )^2 + \cdot\cdot\cdot +
\frac {1}{p}\left(\frac {z}{\zeta}\right )^p\right ].$$

We will call it the {\it primary kernel}  analogously to the primary factor.
\proclaim {Theorem 2.9.1.3.(Estimate of Primary Kernel)
}The following holds:
$$|G_p(x,y,m)|\leq A(m,p)\frac {|x|^{p+1}}{|y|^{p+m-1}}\tag 2.9.1.8$$
for $|x|<2|y|,$
$$|G_p(x,y,m)|\leq A(m,p)\frac {|x|^p}{|y|^{p+m-2}}\tag 2.9.1.9$$
for $|y|<2|x|,$ and
$$G_p(x,y,m)\leq A(m,p)\min \left
(\frac {|x|^{p+1}}{|y|^{p+m-1}},\frac {|x|^{p}}{|y|^{p+m-2}}\right )
\tag 2.9.1.10$$
for all $x,y\in\Bbb R^m,$ where $A(m,p)$ does not depend on $x,y.$

For $m=2,\ p=0$ we have $G_p(z,\zeta,2)\leq A(0,2)\log (1+\frac {|z|}{|\zeta|}).$
\endproclaim
\demo {Proof}The inequality (2.9.1.8) follows directly from (2.9.1.6) and
(2.9.1.9) follows from (2.9.1.7). By the
condition $2\leq |x|/|y|$ (2.9.1.10) follows from (2.9.1.9).

Suppose $1/2\leq |x|/|y|\leq 2.$ Since all the summands in (2.9.1.5) are
bounded from below, for $1/2\leq z\leq 2$ we have
$$G_p(x,y,m)\leq A_1(m,p)|y|^{2-m}\leq A(m,p)\min \left
(\frac {|x|^{p+1}}{|y|^{p+m-1}},\frac {|x|^{p}}{|y|^{p+m-2}}\right )$$
also under these conditions.

The case  $m=2,\ p=0$ is obvious.\qed
\enddemo

\subheading {2.9.2} Let $\mu\in \Cal M (\Bbb R^m).$ We suppose below that
 its support does not contain the origin.

We will say that the integral $\int_{\Bbb R^m}f (x,y)d\mu_y$ converges
uniformly on $x\in D$ if
$$\sup_{x\in D}\left|\int_{|y|>R}f (x,y)d\mu_y\right|\rightarrow 0$$
when $R\rightarrow \infty.$

Hence,  the integral is permitted to be equal to infinity for some finite
 $x.$

Let $\mu$ have  genus $p$ (see, 2.8.3). Set
$$\Pi (x,\mu,p):=\int_{\Bbb R^m}G_p (x,y,m)d\mu_y.\tag 2.9.2.1$$
It is called the {\it canonical potential}.

In particular, let {m=2} and $\mu:=n$ be a {\it zero distribution}, i.e. it
has  unit masses concentrated on a discrete point set $\{z_j:j=1,2,...\}.$
Then
$$\Pi (z,n,p)=\log \left|\prod_{j=1}^{\infty}E\left(\frac {z}{z_j},
p\right )\right|$$
where
$$\prod_{j=1}^{\infty}E\left(\frac {z}{z_j},p\right )$$
is the {\it canonical Weierstrass product}.
\proclaim {Theorem 2.9.2.1.(Brelot-Weierstrass)}The canonical potential
(2.9.2.1) converges uniformly on any bounded domain. It is a subharmonic
function with $\mu$ as its Riesz measure.
\endproclaim

\demo {Proof} Let    $ |x|<R_0$ and $|y|>R$ . From the estimate of the
primary kernel (Theorem 2.9.1.3) we have
$$|\int_{|y|>R}G_p (x,y,m)d\mu_y|\leq
A(m,p)|x|^{p+1}\int_{|y|>R}|y|^{-p-m+1}d\mu_y=$$
$$=A(m,p)|x|^{p+1}\int_R^\infty t^{-p-m+1}d\mu(t).$$
Integrating by part we obtain
$$\int_R^\infty t^{-p-m+1}d\mu(t)=\frac {\mu(R)}{R^{p+m-1}}+
(p+m-1)\int_R^\infty\frac {\mu(t)}{t^{p+m}}dt.$$
The last integral converges since the genus of $\mu$ is $p.$ Hence,
 both summands tend to zero when $R\rightarrow \infty$.
Thus
$$ \sup_{|x|<R_0}|\int_{|y|>R}G_p (x,y,m)d\mu_y|\rightarrow 0$$
while $R_0$ is fixed and $R\rightarrow \infty$, i.e. the canonical potential
converges uniformly  on any bounded domain.

Let us represent the canonical potential for $R>R_0$ in the form
$$\Pi (x,\mu,p)=\int_{|y|<R}G (x,y,\Bbb R^m)d\mu_y +
\int_{|y|<R}\sum_{k=0}^{p}D_k (x,y)d\mu_y + \int_{|y|>R}G_p (x,y,m)d\mu_y .$$
The first summand is a potential, hence a subharmonic function and its
Riesz measure coincide with $\mu.$ The other summands are harmonic for
$|x|<R_0.$ \qed
\enddemo
The following proposition estimates the growth of the canonical potential
in terms of its masses.
\proclaim {Theorem 2.9.2.2.(Estimation of Canonical Potential)}
The following inequality holds:
$$M(r,\Pi(\bullet,\mu,p))\leq A \left [\int _{0}^{\infty}
\frac {\mu(r\tau)}{r^{m-2}}\frac {\min (1,\tau ^{-1})}{\tau ^{p+m-1}}d\tau
+\frac {\mu (r)}{r^{m-1}}\right ]
\tag 2.9.2.2 $$
where $A:=A(m,p)$ does not depend on $r$ and $\mu.$
\endproclaim
\demo {Proof}From (2.9.1.10)
$$\Pi(x,\mu,p)\leq A(m,p)\int_{\Bbb R^m}\min \left
(\frac {|x|^{p+1}}{|y|^{p+m-1}},\frac {|x|^{p}}{|y|^{p+m-2}}\right )d\mu_y$$
Set $r:=|x|,\ t:=|y|$. Then we have
$$M(r,\Pi(\bullet,\mu,p))\leq A \int _{0}^{\infty}
\min \left
(\frac {r^{p+1}}{t^{p+m-1}},\frac {r^{p}}{t^{p+m-2}}\right )d\mu(t)
\tag 2.9.2.3$$
The integral on the right side of (2.9.2.3) can be represented in the form
$$\int _{0}^{r}\frac {r^{p}}{t^{p+m-2}}d\mu(t)+
\int _{r}^{\infty}\frac {r^{p+1}}{t^{p+m-1}}d\mu(t)$$
Integrating every integral by parts  we obtain
$$ (p+m-2)\int _{0}^{r}\frac {r^{p}}{t^{p+m-1}}\mu(t)dt +
(p+m-1)\int _{r}^{\infty}\frac {r^{p+1}}{t^{p+m}}\mu(t)dt +
\frac {\mu (r)}{r^{m-1}}\leq$$
$$\leq(p+m-1)\int _{0}^{\infty}\min \left (1,\frac {r}{t}\right )
\frac {r^{p}}{t^{p+m-1}}\mu(t)dt +\frac {\mu (r)}{r^{m-1}}.$$
After the change $t=r\tau$ we obtain (2.9.2.2) where the new $A(m,p)$ is
equal to $A(m,p)(p+m-1).$ \qed
\enddemo
\proclaim {Theorem 2.9.2.3.(Brelot - Borel)}The order of the canonical
potential is equal to the  convergence exponent of its mass distribution, i.e.
$$\rho [\Pi (\bullet,\mu,p)]=\rho [\mu],$$
 if the genus of $\mu$ is equal to $p.$
\endproclaim  \demo {Proof}First assume $\rho[\mu]<p+1.$ Let us choose $\lambda$
such that $\rho[\mu]<\lambda<p+1.$

For some constant $C$ that  does not depend
on $t$ we have $\mu (t)\leq Ct^{\lambda +m-2}.$

Actually,
$\mu (t)/t^{\lambda +m-2}\rightarrow 0,$ because $\lambda > \rho[\mu].$
 Since $\mu (t)=0$ for small $t,$
this function is bounded and we can take its lower bound as $C.$

Now we have
$$f(r,\tau):=\frac {\mu(r\tau)}{r^{\lambda+ m-2}}\frac
{\min (1,\tau ^{-1})}{\tau ^{p+m-1}}\leq C\tau ^{\lambda -p-1}\min (1,1/\tau).
\tag 2.9.2.4$$
for all $\tau \in (0,\infty).$

We also have
$$\lim_{r\rightarrow \infty}f(r,\tau)=0 \tag 2.9.2.5$$
because  of $\lambda >\rho[\mu].$

Let us divide (2.9.2.2) by $r^\lambda$ and pass to the upper limit. By
Fatou's lemma (Theorem 2.2.2.3)
$$\limsup_{r\rightarrow\infty}\frac {M(r,\Pi (\bullet,\mu,p))}{r^\lambda}\leq
A(m,p)\left [\int_0^\infty \limsup_{r\rightarrow\infty}f(r,\tau)d\tau
+\limsup_{r\rightarrow\infty}\frac {\mu (r)}{r^{\lambda +m-1}}\right ]
=0\tag 2.9.2.6$$
Hence,
$$\lambda \geq \rho [\Pi (\bullet,\mu,p)].\tag 2.9.2.7$$

Since this holds for any $\lambda >\rho[\mu],$ we have $\rho[\mu]\geq
\rho [\Pi (\bullet,\mu,p)]$ under the assumption $\lambda <p[\mu]+1.$

Let $\rho [\mu]=p[\mu]+1$. By Theorem 2.8.3.1 $\bar\Delta[\mu]=0.$ Hence,
$\mu (t)t^{-p-m+1}\leq C$ and
$$f(r,\tau):=\frac {\mu(r\tau)}{(r\tau)^{p+m-1}}\min (1,\tau ^{-1})\leq
C\min (1,1/\tau).$$
The function $\min (1,1/\tau)$ is not summable on $(0,\infty).$ Therefore we
will act in a slightly different way. From Theorem 2.9.2.2 we have
$$\limsup_{r\rightarrow\infty}\frac {M(r,\Pi (\bullet,\mu,p))}{r^{p+1}}\leq
A(m,p)\left [\int_0^1 \limsup_{r\rightarrow\infty}f(r,\tau)d\tau\right ] +$$
$$+A(m,p)\left [\limsup_{r\rightarrow\infty}\int_r^\infty \frac {\mu (t)}{t^{p +m}}dt +
+\limsup_{r\rightarrow\infty}\frac {\mu (r)}{r^{p +m}}\right ]$$
The first integral is equal to zero because  $\bar\Delta[\mu]=0.$ The second
addend vanishes since the integral converges.Thus we have $p+1=\rho [\mu]\geq \rho [\Pi (\bullet,\mu,p)].$

The reverse inequality holds for any subharmonic function in $\Bbb R^m$
by the Jensen theorem (Theorem 2.8.3.3).\qed
\enddemo
\subheading {2.9.3}Let us denote as $\delta SH (\rho)$  the class of
functions $u\in \delta SH (\Bbb R^m)$ for which $\rho _T [u]\leq \rho.$
\proclaim {Theorem 2.9.3.1.(Brelot - Hadamard)}Let
$u=u_1-u_2 \in\delta SH (\rho),$ and let $p_1,p_2$ be the genuses of the
mass distributions $\mu_j:=\mu_{u_j},\ j=1,2.$ Suppose
$\supp [\mu_1-\mu_2]\cap \{0\}=\varnothing.$

Then the following equality
holds:
$$u(x)=\Pi (x,\mu_1,p_1) - \Pi (x,\mu_2,p_2) + \Phi _q (x)$$
where $\Phi _q (x)$ is a harmonic polynomial of degree
$q\leq\rho.$
\endproclaim
\demo {Proof}The function $v(x):=u(x)-\Pi (x,\mu_1,p_1)+\Pi (x,\mu_2,p_2)$ is
harmonic by the Brelot - Weierstrass theorem (Theorem 2.9.2.1). We also have
the inequality
$$\rho_T [v]\leq \max (\rho_T [u],\rho_T [\Pi (\bullet,\mu_1,p_1)],
\rho_T [\Pi (\bullet,\mu_2,p_2)])\tag 2.9.3.1$$
by Theorem 2.8.2.2 ($\rho_T$ - properties).
The property $\rho$MT1) (Theorem 2.8.2.3 ) implies
$$\rho_T [\Pi (\bullet,\mu_j,p_j)]=\rho_M [\Pi (\bullet,\mu_j,p_j)]:=
\rho [\Pi (\bullet,\mu_j,p_j)],\ j=1,2.$$
The Brelot-Borel theorem (Theorem 2.9.2.3) implies
$$\rho [\Pi (\bullet,\mu_j,p_j)]=\rho [\mu_j],\ j=1,2.$$
The Jensen theorem (Theorem 2.8.3.3) implies
$$\max (\rho [\mu_1],  \rho [\mu_2])\leq \rho_T [u].$$
From (2.9.3.1) we have
$$\rho_T [v]\leq\rho_T [u]\leq \rho.$$
Since $v$ is subharmonic, $\rho_T [v]=\rho_M [v]:=\rho [v]$ by Theorem 2.8.2.3,
and $ \rho [v]\leq \rho.$ Therefore
$$\lim_{r\rightarrow \infty}\frac {M(r,v)}{r^{\rho +\epsilon}}=0$$
for arbitrary small $\epsilon >0.$

By the Liouville theorem (Theorem 2.4.2.3) $v(x)$ is a harmonic
polynomial of degree $q\leq \rho+\epsilon,$ and thus $v(x)=\Phi _q
(x)$ for $q\leq \rho.$ \qed
\enddemo
For a non-integer $\rho$ the Brelot-Hadamard theorem allows to
connect  the growth of functions and masses more tightly than in the Jensen
theorem.
\proclaim {Theorem 2.9.3.2.(Sharpening of Jensen) }Let $\rho >0$ and be
non-integer, $u=u_1-u_2\in \delta SH (\Bbb R^m)$ with $\rho_T [u]=\rho,$
and let $\nu_u=\mu_1-\mu_2$ the corresponding charge.Then

pJ1) $\rho[\nu_u]=\max (\rho[\mu_1],\rho[\mu_2])=\rho$

pJ2) $A_1\sigma_T [u,\rho(r)]\leq \bar\Delta [\nu_u,\rho(r)]\leq
\bar\Delta [\mu_1,\rho(r)]+\bar\Delta [\mu_2,\rho(r)]\leq
A_2\sigma_T [u,\rho(r)],$\newline
where $A_j=A_j(m,\rho)$ and $\rho(r)$ is an arbitrary proximate order such
that $\rho(r)\rightarrow \rho $ when $r\rightarrow \infty.$
\endproclaim
For proving this theorem we need \proclaim {Theorem 2.9.3.3}Let $
\Pi (x,\mu,p)$ be a canonical potential with   non-integer $\rho
[\mu]:=[\rho]$, and let$\rho(r)(\rightarrow \rho)$
 be a proximate order.Then
$$\sigma [\Pi (\bullet,\mu,p),\rho(r)]\leq A(m,\rho,p)
\bar\Delta [\mu,\rho(r)]\tag 2.9.3.2$$
\endproclaim
\demo {Proof}We can suppose without loss of generality that
$\bar\Delta [\mu,\rho(r)]<\infty.$ By this condition and since
$\mu(t)=0,$ $0<t<c$ for some $c>0,$ we have the inequality
$$\mu (t)t^{-\rho(t)-m+2}\leq C$$
for all $t\in (0,\infty)$ and some $C>0$ that does not depend on $t.$
Set
$$I(r):=\int_{c/r}^\infty \frac {\mu (rt)}{r^{\rho(r)+m-2}}
\frac {\min (1,1/t)}{t^{p+m-1}}dt.$$
By Theorem 2.9.2.2 we have
$$\sigma [\Pi (\bullet,\mu,p),\rho(r)]=
\limsup_{r\rightarrow\infty}\frac {M(r,\Pi (\bullet,\mu,p)}{r^{\rho(r)}}
\leq A(m,p)\limsup_{r\rightarrow\infty}I(r).\tag 2.9.3.3$$

Let us choose $r_\epsilon$ such that
$$\sup_{r>r_\epsilon}\frac {\mu (r\epsilon)}{(r\epsilon)^{\rho(r\epsilon)+m-2}}
\leq \bar\Delta [\mu,\rho(r)]+\epsilon.$$
For such $r$ we have
$$I(r)=\int_{c/r}^\infty\frac{\mu (rt)}{(rt)^{\rho(rt)+m-2}}
\frac{(rt)^{\rho(rt)}}{r^{\rho(r)}}\frac {\min (1,1/t)}{t^{p+1}}dt\leq$$
$$\leq \sup_{c/r\leq t\leq\epsilon}\frac{\mu (rt)}{(rt)^{\rho(rt)+m-2}}
\int_{c/r}^\epsilon \frac{(rt)^{\rho(rt)}}{r^{\rho(r)}}
\frac {\min (1,1/t)}{t^{p+1}}dt +$$
$$\sup_{\epsilon\leq t\leq 1/\epsilon}...
\int_{\epsilon}^{1/\epsilon}...dt +
\sup_{1/\epsilon\leq t\leq \infty}...
\int_{1/\epsilon}^{\infty}...dt\leq $$
$$\leq C\int_{c/r}^\epsilon \frac{(rt)^{\rho(rt)}}{r^{\rho(r)}}
\frac {\min (1,1/t)}{t^{p+1}}dt +(\bar\Delta [\mu,\rho(r)]+\epsilon)
\int_{\epsilon}^{1/\epsilon}...dt+C\int_{1/\epsilon}^{\infty}...dt.$$
The function
$$f(t):=\frac {\min (1,1/t)}{t^{p+1}}$$
satisfies the conditions of Goldberg's theorem (Theorem 2.8.1.5) with
$p+1-\rho<\delta<1$ and $0<\gamma<p+1-\rho.$ Passing to the limit we have
$$\limsup_{r\rightarrow\infty}I(r)\leq
C\int_{0}^\epsilon t^{\rho -p}dt +(\bar\Delta [\mu,\rho(r)]+\epsilon)
\int_{\epsilon}^{1/\epsilon}t^{\rho -p -1}\min (1,1/t)dt+$$
$$C\int_{1/\epsilon}^{\infty}t^{\rho -p - 2}dt.$$
Passing to the limit as $\epsilon \rightarrow 0$ we obtain with the help of
(2.9.3.3)
$$\sigma [\Pi (\bullet,\mu,p),\rho(r)]\leq A(m,p)
\bar\Delta [\mu,\rho(r)]
\int_{0}^{\infty}t^{\rho -p -1}\min (1,1/t)dt.$$\qed
\enddemo
\demo {Proof of Theorem 2.9.3.2} The inequality
$\rho [\nu_u]\leq \rho$ and the last inequality in pJ2) follow from
the Jensen theorem (Theorem 2.8.3.3). Let us prove the reverse inequality
and the left side.

Since $\rho$ is non-integer, $q<\rho$ in the Brelot-Hadamard theorem (Theorem 2.9.3.1).
 Hence $M(r,\Phi _q)=o(r^\rho)$ and
$$T(r,u)\leq T(r,\Pi (\bullet,\mu_1,p))+ T(r,\Pi (\bullet,\mu_2,p))+o(r^\rho)
.$$
Thus
$$\rho_T [u]\leq \max (\rho[\Pi (\bullet,\mu_1,p)],
\rho[\Pi (\bullet,\mu_2,p]),$$ $$  \sigma_T [u,\rho(r)]\leq
\max ( \sigma_T[\Pi (\bullet,\mu_1,p) ,\rho(r)],
\sigma_T[\Pi (\bullet,\mu_2,p] ,\rho(r)]).$$
From Theorem 2.9.3.3 we obtain
$$\rho_T [u]\leq\max (\rho[\mu_1],\rho[\mu_2]);$$
$$\sigma_T [u,\rho(r)]\leq  A(m,\rho,p)
\max(\bar\Delta [\mu_1,\rho(r)],\bar\Delta [\mu_2,\rho(r)]
=$$ $$ =A(m,\rho,p)\bar\Delta [|\nu|,\rho(r)].$$
We can set $A_1:=A^{-1}(m,\rho,p)$ and obtain the left side of pJ2).\qed
\enddemo
\subheading {2.9.4}Let $u\in \delta SH (\Bbb R^m)$ and $\rho:=\rho_T[u]$ be
an integer number. We can always represent the function $u$ in the form
$$u(x)=\Pi (x,\nu,\rho)+\Phi_\rho (x) \tag 2.9.4.1$$
where $\Phi_\rho (x)$  is a harmonic polynomial of degree at most
$\rho.$ Actually, such a representation can be obtained from
Theorem 2.9.3.1
 by addition  and subtraction of terms of the form
$$ \Phi _{k_j}(x):=\int _{\Bbb R^m} D_{k_j}(x,y)d(\mu_j)_y,\ j=1,2;$$
where $p_j<k_j\leq \rho.$

All $\Phi _{k_j}(x)$ of such a kind are harmonic polynomials of
degree at most $\rho.$

Set
$$\Pi_<^R (x,\nu,\rho - 1):=\int_{|y|<R}G_{\rho -1} (x,y,m)d\nu_y,
\tag 2.9.4.2$$
$$\Pi_>^R (x,\nu,\rho):=\int_{|y|\geq R}G_{\rho} (x,y,m)d\nu_y,\tag 2.9.4.3$$
$$\delta _R(x,\nu,\rho):=\int_{|y|<R}D_\rho (x,y)d\nu_y.\tag 2.9.4.4$$
In particular, for $m=2$
$$\delta _R (z,\nu,\rho):=\frac {1}{\rho}\int_{|\zeta|<R}
\Re \left (\frac {z}{\zeta}\right )^\rho d\nu_\zeta.\tag
2.9.4.4a$$ Let $Y_\rho (x)$ be the homogeneous polynomial of
degree $\rho$ from the polynomial $\Phi_\rho$ in (2.9.4.1). Set
also
$$\delta_R (x,u,\rho):=\delta (x,\nu,\rho)+Y_\rho (x),$$
$$M(r,\delta):=\max_{|y|=1}|\delta (ry,u,\rho)|,\tag 2.9.4.5$$
$$\bar \Delta _\delta [u,\rho]:=\limsup_{r\rightarrow\infty}M(r,\delta)
r^{-\rho(r)}.$$ The functions $\delta _R(x,\nu,\rho)$ are
homogeneous polynomials that are determined completely by their
values on the unit sphere. Thus, by  the Harnack theorem (Theorem
2.4.1.7) we have \proclaim {Theorem 2.9.4.1}  $\bar \Delta _\delta
[u,\rho(r)]<\infty$ if and only if the family $\delta_R
(x,u,\rho)R^{\rho -\rho(R)},R>0$ is precompact in $\Cal D'(\Bbb
R^m).$
\endproclaim
Let $\rho$ be an integer number and $\rho(r)\rightarrow \rho$ be a p.o. Set
$$\Omega [u,\rho(r)]:=\max (\bar \Delta _\delta [u,\rho(r)],
\bar \Delta  [|\nu_u|,\rho(r)].$$
\proclaim {Theorem 2.9.4.2.(Brelot-Lindel\"of)}The following holds
$$A_1 \Omega [u,\rho(r)]\leq\sigma _T[u,\rho(r)]\leq
A_2\Omega [u,\rho(r)],$$
where $A_j:=A_j (m,\rho).$
\endproclaim
For proving this theorem we will first study the function $\Pi_<^R$ and
$\Pi _>^R.$ Set
$$T(r,\lambda,>):=T(r,\Pi_>^{\lambda r} (\bullet,\nu ,\rho)),$$
$$T(r,\lambda,<):=T(r,\Pi_<^{\lambda r} (\bullet,\nu ,\rho-1)).$$
\proclaim {Theorem 2.9.4.3.(Estimate of $T(\bullet,>,T(\bullet,<)$)}
The following holds:
$$T(r,\lambda,>)\leq A\left(\int_\lambda ^\infty
\frac {|\nu|(rt)}{r^{m-2}}\frac {\min(1,t^{-1})}{t^{\rho+m-1}}dt+
\frac {|\nu|(r)}{r^{m-1}}\right ),\tag 2.9.4.6$$
$$T(r,\lambda,<)\leq A\left (\int_0 ^\lambda
\frac {|\nu|(rt)}{r^{m-2}}\frac {\min(1,t^{-1})}{t^{\rho+m-2}}dt \right )+
\tag 2.9.4.7$$
$$+A\left (\frac {|\nu|(r\lambda)}{r^{m-2}}\int_\lambda ^\infty\frac {\min(1,t^{-1})}
{t^{\rho+m-2}}dt+\frac {|\nu|(r)}{r^{m-1}}\right ),$$
where $A:=A(m,\rho).$
\endproclaim
\demo {Proof}Let $\nu=\mu_1-\mu_2.$ Then $|\nu|=\mu_1+\mu_2$. We have
$$\Pi_<^R(x,\nu,\rho -1)=\Pi_<^R(x,\mu_1,\rho -1)-\Pi_<^R(x,\mu_2,\rho -1)
\tag 2.9.4.8$$
Since $\Pi_<^R(0,\mu_2,\rho -1)=0 ,$ we have (see t3),t4), Theorem 2.8.2.1)
$$T(r,\Pi_<^R(\bullet,\nu,\rho -1))\leq
T(r,\Pi_<^R(\bullet,\mu_1,\rho -1))+T(r,\Pi_<^R(\bullet,\mu_2,\rho -1))
\tag 2.9.4.9$$
Set $$\Pi_1:=\Pi_<^R(\bullet,\mu_1,\rho -1),\
\Pi_2:=\Pi_<^R(\bullet,\mu_2,\rho -1).$$
Let us estimate, for example, $T(r,\Pi_1).$ The masses of the canonical
potential $\Pi_1$ are concentrated in $K_R$. Applying Theorem 2.9.2.2
(Estimation of Canonical Potential) for $p=\rho-1$ we obtain
$$T(r,\Pi_1)\leq M(r,\Pi_1)\leq A\int_0^{\frac {R}{r}}
\frac {\mu_1(rt)}{r^{m-2}}\frac {\min(1,t^{-1})}{t^{\rho+m-2}}dt+$$
$$+A\frac {\mu_1 (R)}{r^{m-2}}
\int_{\frac {R}{r}}^{\infty}\frac {\min(1,t^{-1})}{t^{\rho+m-2}}dt +
\frac {\mu_1 (r)}{r^{m-1}}.$$
Set $R:=r\lambda.$ Then we obtain the inequality (2.9.4.7) for $\nu:=\mu_1.$
Analogously one can do for $\nu:=\mu_2.$ The inequality (2.9.4.9) allows to pass to  limit in (2.9.4.7) in the  general case.

Set $\Pi_1:=\Pi_>^R(\bullet,\mu_1,\rho).$ Applying (2.9.2.2) for $p=\rho$
we obtain
$$T(r,\Pi_1)\leq M(r,\Pi_1)\leq\int_{\frac {R}{r}}^{\infty}
\frac {\mu_1(rt)}{r^{m-2}}\frac {\min(1,t^{-1})}{t^{\rho+m-1}}dt +
\frac {\mu_1 (r)}{r^{m-1}}.$$
In the same way we obtain (2.9.4.6).\qed
\enddemo
Set
$$\sigma [\Pi_>,\rho(r)]:=\limsup_{r\rightarrow \infty}
\frac {T(r,\Pi_>^r(\bullet,\nu,\rho))}{r^{\rho(r)}},$$
$$\sigma [\Pi_<,\rho(r)]:=\limsup_{r\rightarrow \infty}
\frac {T(r,\Pi_<^r(\bullet,\nu,\rho))}{r^{\rho(r)}}.$$
\proclaim {Theorem 2.9.4.4} Let $\nu:=\mu_1-\mu_2 \in \delta \Cal M (\rho)$
and $\rho$ integer number. Then for any p.o. $\rho(r)\rightarrow \rho$
$$\max (\sigma [\Pi_>,\rho(r)],\sigma [\Pi_<,\rho(r)])\leq
A\bar \Delta [|\nu|,\rho(r)]$$
where $A:=A(m,\rho).$
\endproclaim
\demo {Proof}From (2.9.4.6) we have
$$T(r,\Pi_>^r(\bullet,\nu,\rho))=T(r,1,>)\leq
A\int_1^\infty \frac {|\nu|(rt)}{r^{m-2}}\frac {1}{t^{\rho +m}}dt+
\frac {|\nu|(r)}{r^{m-1}}.$$
Now we  repeat the reasoning of  Theorem 2.9.3.3 for
$\mu:=|\nu|$ and $p:=\rho.$ We will obtain
$$\sigma [\Pi_>,\rho(r)]\leq A\bar \Delta [|\nu|,\rho(r)]
\int_1^\infty t^{-2}dt.$$

For the other case we have from (2.9.4.7)
$$T(r,\Pi_<^r (\bullet,\nu,\rho-1))=T(r,1,<)\leq
A\int_0^1\frac {|\nu|(rt)}{r^{m-2}}\frac {1}{t^{\rho +m-1}}dt+$$
$$+A\frac {|\nu|(r)}{r^{m-2}}\left (\int_1^\infty t^{-\rho -m+1}dt +
r^{-1}\right ).$$
We divide this inequality by $r^{\rho(r)}$ and pass to the  upper limit while
$r\rightarrow\infty.$

The first summand of the right side gives
$$A\bar \Delta [|\nu|,\rho(r)]\int_0^1dt $$
by the reasoning of Theorem 2.9.3.3.

The second one can be computed directly, yielding
$$A\bar \Delta [|\nu|,\rho(r)]\int_1^\infty t^{-\rho -m+1}dt.$$

Combining all these inequalities we obtain the assertion of the theorem.\qed
\enddemo
\demo {Proof of Theorem 2.9.4.2}Let us represent $u(x)$ in the form
$$u(ry)=\Pi _<^r(ry,\nu_u,\rho -1)+ \Pi_>^r (ry,\nu_u,\rho)+
\delta_r (ry,u,\rho) +o(r^{\rho -1})\tag 2.9.4.10$$
where $|y|=1.$

Then we have
$$T(r,u)\leq T(r,\Pi _<^r(\bullet,\nu_u,\rho -1))+
T(r,\Pi_>^r (\bullet,\nu_u,\rho))+M(r,\delta)+o(r^{\rho-1}).$$
Let us divide this by $r^{\rho(r)} $ and pass to the upper limit. By
Theorem 2.9.4.4 we obtain
$$\sigma_T [u,\rho(r)]\leq A \max (\bar \Delta [|\nu|,\rho(r)],
\bar \Delta_\delta [u,\rho(r)])=A_2\Omega [u,\rho(r)]$$
where $A_2=A(m,\rho).$
Let us write (2.9.4.11) in the form
$$\delta_r (ry,u,\rho)=u(ry)-Pi _<^r(ry,\nu_u,\rho -1)-
\Pi_>^r (ry,\nu_u,\rho)+o(r^{\rho -1}).$$
We obtain
$$T(r,\delta_r (\bullet,u,\rho))\leq T(r,u)+
T(r,\Pi _<^r(\bullet,\nu_u,\rho -1))+
T(r,\Pi_>^r (\bullet,\nu_u,\rho)+o(r^{\rho-1}).$$
Since $\delta_R (\bullet,u,\rho)$ is harmonic and homogeneous, we have
by (2.8.2.5)
$$M(r,\delta_R)\leq 2^{m-1}T(2r,\delta_R)=2^{m-1+\rho}T(r,\delta_R).$$
Therefore we obtain the inequality
$$\bar \Delta_\delta [u,\rho(r)]\leq\sigma_T [u,\rho(r)]+
2A\bar \Delta [|\nu|,\rho(r)].$$
By the Jensen theorem (Theorem 2.8.3.3) we have
$$\Omega [u,\rho(r)]\leq A_1^{-1}\sigma_T [u,\rho(r)]$$
for some $A_1=A_1(m,\rho)$.\qed
\enddemo

 \newpage

\centerline {\bf 3.Asymptotic behavior of subharmonic functions of
finite order.}

\centerline {\bf 3.1.Limit sets}

\subheading {3.1.1}Let  $\{V_t:t\in (0,\iy)\}$ be a family of rotations of $\Bbb R^m$
that form a one-parametric group, i.e.,
$$V_{t_1}V_{t_2}=V_{t_1 t_2},\ V_{1}=I,\tag 3.1.1.0$$
where $I$ is the identity map.

The family of linear transformations
$$P_t:=tV_t\tag 3.1.1.1$$
is also a one-parametric group.

In particular, for $m=2$ the general form of the rotations is
$$V_t z=z\exp (i\gamma \log t),$$
where $\gamma$ is real.

The orbit $\{P_t z:t\in (0,\infty)\}$ of every point $z\neq 0$ is a logarithmic
spiral if $\gamma\neq 0$ and a ray if $\gamma=0.$

For $m\geq 3$ and $V_t\equiv I,\ t\in (0,\iy)$ the orbit of every point $x\neq 0$ is a ray
from the origin. For other $V_\bullet$ it is a spiral connecting the
origin to infinity.

It is clear that only one orbit $ \{P_t x:t\in (0,\infty)\}$  passes
through every $x\neq 0.$ The behavior of every point $y(t):=P_t x$ is
completely determined by a system of  differential equations with
 constant coefficients:
$$\frac {d}{dt} y=(I+V')y,\ \ V':=\frac {d}{dt} V_t\mid _{t=1}$$
with the initial condition of $y(1)=x.$
\subheading {3.1.2} Let $u\in SH (\rho)$ and $\sigma_M [u,\rho(r)]<\infty$
for some p.o. $\rho (r)\rightarrow \rho.$ We will write
$u\in SH (\Bbb R^m,\rho,\rho(r))$ or shorter $u\in SH (\rho(r)).$

For $u\in SH (\rho(r))$ set
$$ u_t (x) :=u(P_t x)t^{-\rho (t)}.\tag 3.1.2.1$$
We will denote this transformation as $(\bullet)_t.$
\proclaim {Theorem 3.1.2.1.(Existence of Limit Set)} The following holds:

els1) $u_t \in SH (\rho(r))$ for any $t\in (0,\infty)$

els2) the family $\{u_t\}$ is precompact at infinity,
\endproclaim
i.e., for any sequence $t_k\rightarrow \infty$ there exists a subsequence
$t_{k_j}\rightarrow \infty$ and a function $v\in SH(\Bbb R^m)$ such that
$u_{t_{k_j}}\rightarrow v$ in $\Cal D'(\Bbb R^m)$ (see, 2.7.1).
\demo {Proof}The functions $u_t$ are subharmonic by sh1) and sh5), Theorem
2.6.1.1.(Elementary Properties), and
$$M(r,u_t)=M(rt,u)t^{-\rho(t)}.$$
 Now we have
$$\sigma_M [u_t,\rho (r)]=t^{-\rho(t)}
\limsup_{r\rightarrow\infty } \frac {M(rt,u)}{(rt)^{\rho(rt)}}\cdot
\lim_{r\rightarrow\infty}\frac{(rt)^{\rho(rt)}}{r^{\rho(r)}}=
\sigma_M [u,\rho (r)]t^{\rho-\rho(t)},$$
because
$$\lim_{r\rightarrow\infty}\frac{(rt)^{\rho(rt)}}{r^{\rho(r)}}=
t^{\rho}\lim_{r\rightarrow\infty}\frac {L(rt)}{L(r)}=t^{\rho}\tag 3.1.2.2$$
(see, ppo3), Theorem 2.8.1.3.(Properties of P.O)).
Therefore els1) is proved.

Let us check the conditions of Theorem 2.7.1.1.(Compactness in $\Cal D'$).
We have
$$\limsup_{t\rightarrow\infty }M(r,u_t)=
\limsup_{t\rightarrow\infty } \frac {M(rt,u)}{(rt)^{\rho(rt)}}\cdot
\lim_{t\rightarrow\infty}\frac{(rt)^{\rho(rt)}}{t^{\rho(t)}}=
\sigma_M [u,\rho (r)]r^{\rho}.\tag 3.1.2.3$$
Thus, the family is bounded from above on every compact set and
$$\lim_{t\rightarrow\infty}u_t (0)=\lim_{t\rightarrow\infty}u(0)t^{-\rho(t)}
=0.$$
Therefore $u_t(0)$ are bounded from below for large $t$.\qed
\enddemo

 We will call the set of all functions $v$ from Theorem 3.1.2.1
the {\it limit set } of the function $u(x)$ with respect to $V_\bullet$ and
denote it $\bold {Fr}[u,\rho(r),V_\bullet,\Bbb R^m]$ or shortly $\Fr [u].$

The limit set does not depend on values of the subharmonic function on a
bounded set, hence, it is a characteristic of asymptotic behavior.

Set
$$U[\rho,\sigma]:=\{v\in SH (\Bbb R^m): M(r,v)\leq \sigma r^\rho,\ r\in
[0,\infty);\ v(0)=0\},$$
$$U[\r]:=\bigcup\limits_{\s>0}U[\rho,\sigma]\tag 3.1.2.4$$
and
$$v_{[t]}(x):=t^{-\rho}v(P_t x),\ t\in (0,\infty).\tag 3.1.2.4a$$
Let us emphasize that the transformation $(\bullet)_{[t]}$
coincides with $(\bullet)_t$ from (3.1.2.1) for $\rho
(r)\equiv\rho$ and satisfies the condition
$$(\bul)_{[t\tau]}=((\bul)_{[t]})_{[\tau]}\tag 3.1.2.4b$$
\proclaim {Theorem 3.1.2.2.(Properties of $Fr$)}The following holds:

fr1) $\Fr [u]$ is a connected compact set;

fr2) $\Fr[u]\subset U[\rho,\sigma] ,$ for $\sigma\geq \sigma_M [u];$

fr3) $ (\Fr [u])_{[t]}=\Fr [u],\ t\in (0,\infty),$

i.e., $v\in \Fr [u]$ implies $v_{[t]}\in \Fr [u].$

f4) if $\r_1(r)$ and $\r(r)$ are equivalent (see (2.8.1.5)), then
$\Fr [u,\r_1(r),\bul]=\Fr[u,\r(r),\bul].$
\endproclaim

We need  the following assertion
\proclaim {Theorem 3.1.2.3.(Continuity $u_t$)}
The functions
$$u_t,v_{[t]}:(0,\infty)\times \Cal D'(\Bbb R^m)\mapsto \Cal D'(\Bbb R^m)$$
are continuous in the natural topology.
\endproclaim
\demo {Proof}For any $\psi \in\Cal D(\Bbb R^m)$ consider
$$<u_t,\psi>:=\int u_t (x)\psi (x)dx = \int u (y)\psi (y/t)t^{m-\rho(t)}dy:=
<u,\psi(\bullet,t)>,$$
where $ \psi(y,t):=\psi (y/t)t^{m-\rho(t)}.$

The function $\psi(\bullet,t)$ is continuous in  $t$ in $\Cal D(\Bbb R^m).$
By Theorem 2.3.4.6 (Continuity $<\bullet,\bullet>$) $<u,\psi(\bullet,t)>$
is continuous in $(u,t).$\qed
\enddemo

\demo {Proof of Theorem 3.1.2.2} Let us denote as $clos\{\bullet\}$ the closure
in $\Cal D'$-topology.

The set $F_N:=clos\{u_t: t\leq N\}\supset \Fr [u]$
is compact in $\Cal D'$-topology. Indeed, let  $t_j\rightarrow t$ and
$t<\infty;$ then $u_{t_j}\rightarrow u_t$ because of Theorem 3.1.2.3.
If  $t_j\rightarrow \infty$ and $u_{t_j}\rightarrow v,$ then $v\in \Fr [u]$ by
its definition, hence, $v\in F_N.$ Since $\Fr [u]=\cap_{N=1}^{\infty}F_N,$ it
is  compact.

Let us prove the connectedness. Suppose $\Fr [u]$ is not connected.
Then it can
be written as a union of two disjoint non-empty closed sets $F^1$ and $F^2.$
 Let $V^1,V^2$ be  disjoint open neighborhoods of $F^1,F^2$ respectively in
$\Cal D'(\Bbb R^m).$ Since $F^1,F^2$ are nonempty there exist sequences
$\{s_j\},\{t_j\}$ such that
$s_j<t_j,\ s_j\rightarrow \infty,\ \ u_{s_j}\in V^1,\ u_{t_j}\in V^2.$
 Since the mapping $u_t:(0,\infty)\mapsto \Cal D'(\Bbb R^m)$ is continuous
by Theorem 3.1.2.3 its image is connected. Thus there exists a sequence
$\{p_j\}$ with $s_j<p_j<t_j$ such that $u_{p_j}\notin V^1\cup V^2.$ This
sequence has a subsequence that converges to a function $v\in \Fr [u]$ and
 $v\notin F^1\cup F^2.$ This is a contradiction. Hence, $\Fr [u]$ is connected
and fr1) is proved.

Set
$$\psi (r):=\limsup_{r\rightarrow \infty}M(r,u_t).$$
This function is convex with respect to $-r^{2-m}$ for $m\geq 3$ and
with respect to $\log r$ for $m=2$ and hence continuous.

 Indeed, $M(|x|,u_t)$ are subharmonic
(see Theorem 2.6.5.2.(Convexity $M(\bullet,u)$ and $\Cal M(r,u)$). By
Theorem 2.7.3.3.(H.Cartan +) the function $\psi^* (|x|)$ is subharmonic and
$\psi(|x|)=\psi^* (|x|)$ quasi-everywhere. However,
if $\psi(|x|)<\psi^* (|x|)$
at some point, the same inequality holds on a sphere which has a positive
capacity (see Example 2.5.2.2). Hence, $\psi(|x|)=\psi^* (|x|)$ everywhere,
and
 $\psi(|x|)$ is subharmonic. Thus  $\psi (r)$ is convex with respect to
$-r^{2-m}$ for $m\geq 3$ and with respect to $\log r$ for $m=2$ by
Theorem 2.6.3.2.(Subharmonicity and Convexity).

One can also see that for $u\in SH (\Bbb R^m)$
$$M(r,u_\epsilon)\leq M(r+\epsilon,u),$$
where $(\bullet)_\epsilon$ is defined by (2.6.2.3).

Let $v\in\Fr [u]$ and $u_{t_j}\rightarrow v$ in $\Cal D (\Bbb R^m).$ By
property reg3), Theorem 2.3.4.5 $(u_{t_j})_\epsilon \rightarrow v_\epsilon$
uniformly on any compact set. Thus
$$v_\epsilon(x)=\lim_{j\rightarrow \infty}(u_{t_j})_\epsilon \leq
\limsup_{t\rightarrow\infty)}M(|x|,(u_{t})_\epsilon) \leq \tag 3.1.2.5 $$
$$\leq \limsup_{t\rightarrow\infty)}M(|x|+\epsilon,u_{t})=\psi (|x|+\epsilon).
$$
If $\epsilon \downarrow 0$, then $v_\epsilon \downarrow v$ by Theorem 2.6.2.3
and $\psi(r+\epsilon )\rightarrow \psi (r)$ because of continuity.
Passing to the limit in (3.1.2.5) and using (3.1.2.3) we obtain
$$v(x)\leq \sigma_M [u,\rho (r)]|x|^{\rho}.\tag 3.1.2.6$$

Since $u(0)\leq u_\epsilon (0)$ we have
$u(0)t^{-\rho(t)}\leq (u_t)_\epsilon (0).$ Let us pass to the limit as
$t:=t_j\rightarrow \infty.$ We obtain $v_\epsilon (0) \geq 0 .$ Passing to the
limit as $\epsilon \downarrow 0$ we have
$$v(0)\geq 0.\tag 3.1.2.7$$
The inequalities (3.1.2.6) and (3.1.2.7) imply fr2).

 One can check the equality
$$(u_{t})_{[\tau]}=u_{t\tau}\cdot
\frac {(t\tau)^{\rho (t\tau)}}{t^{\rho(t)}\tau ^\rho}.\tag 3.1.2.8$$
 By using properties of p.o. we have
$$\lim_{t\rightarrow \infty}
\frac {(t\tau)^{\rho (t\tau)}}{t^{\rho(t)}\tau ^\rho}=1$$
(compare (3.1.2.2)).

Let $v\in \Fr [u]$ and $u_{t_j}\rightarrow v.$ Set
$t:=t_j, \ \tau :=t$  in (3.1.2.8) and pass to the limit. Then
$$v_{[t]}=\Cal D'-\lim _{j\rightarrow\infty}u_{t_jt}.$$
Thus $v_{[t]}\in \Fr [u].$
The property f3) is proved.

Let us prove f4).We have
$$\frac {u(P_t x)}{t^{\r_1(t)}}=\frac {u(P_t x)}{t^{\r(t)}}\times
e^{(\r_1(t)-\r(t))\log t}=\frac {u(P_t x)}{t^{\r(t)}}\times (1+o(1))$$
as $t\ri\iy$ because of (2.8.1.5).

This implies f4).

{\bf Exercise 3.1.2.1} Check this in details.

\qed
\enddemo
We can consider the limit sets as a mapping $u\mapsto \Fr[u]$.
The following theorem describes some properties of this mapping.

Set
$$U[\rho]:=\bigcup \limits_{\sigma>0} U[\rho,\sigma]\tag 3.1.2.9$$
where $ U[\rho,\sigma]$ is defined by (3.1.2.4).

Let $X,Y$ be subsets of a cone (i.e. a subset of a linear  space
that is closed with respect to sum and multiplication by a
positive number). The set $U[\rho]$ is such a cone. Set
$$X+Y:=\{z=x+y:x\in X,\ y\in Y\};\ \lambda X:=\{z=\lambda x:x\in X\}.
\tag 3.1.2.10$$
\proclaim {Theorem 3.1.2.4(Properties of $u\mapsto \Fr [u]$)}
The following holds:

fru1) $\Fr [u_1+u_2]\subset\Fr [u_1]+\Fr [u_2]$

fru2) $\Fr [\lambda u]=\lambda \Fr[u]$.
\endproclaim
\demo {Proof}Let $v\in \Fr [u_1+u_2].$ Then there exists
$t_j\rightarrow \infty$ such that $(u_1+u_2)_{t_j}\rightarrow v$ in
$\Cal D'.$ We can find a subsequence $t_{_{j_k}}$ such that
$(u_1)_{t_{_{j_k}}}\rightarrow v_1$ and $(u_2)_{t_{_{j_k}}}\rightarrow v_2.$
Then $v=v_1+v_2.$ The property fru1) has been proved.

The property fru2) is proved analogically.\qed
\enddemo
\subheading {3.1.3} We will write $\mu\in \Cal M (\Bbb R^m,\rho (r))$ or
shortly  $\mu\in \Cal M (\rho (r))$ if $\mu\in \Cal M (\Bbb R^m)$
(see 2.8.3) and $\bar \Delta [\mu,\rho(r)]<\infty$ (see 2.8.3.2).

Let us define a distribution $\mu _t$ for $\mu\in \Cal M (\rho (r))$ by
$$<\mu _t,\phi>:=t^{-\rho(t)-m+2}\int \phi (P^{-1}_t x)d\mu \tag 3.1.3.1$$
for $\phi\in \Di (\Bbb R^m).$

It is positive. Hence, it defines uniquely a measure $\mu _t.$
\proclaim {Theorem 3.1.3.1.(Explicit form of $\mu_t$)}For any $E\in \sigma (\Rm)$
the following holds:
$$\mu_t(E)=t^{-\rho(t)-m+2}\mu (P_t E)\tag 3.1.3.2$$
\endproclaim
\demo {Proof}It is enough to proof the assertion for some dense ring
(see Theorem 2.2.3.5), for example, for all compact sets.

Let $\chi_K$ be a characteristic function of a compact set  $K$ and let
$\phi_\epsilon \downarrow \chi_K$ be a monotonically converging sequence of
functions that belong to $\Di (\Rm)$ (see Theorems 2.1.2.1, 2.1.2.9 and
2.3.4.4).
Then
$$\int \phi_\epsilon (x) d \mu_t=t^{-\rho(t)-m+2}\int \phi_\epsilon (P^{-1}_t x)  d \mu
.$$
Since $\phi_\epsilon (P^{-1}_t x)\downarrow \chi_{P_t K}(x),$
$$\mu_t (K)=\int \chi_K (x)d \mu_t=t^{-\rho(t)-m+2}\int\chi_{P_t K}(x)d \mu =
t^{-\rho(t)-m+2}\mu (P_t K).$$
\qed
\enddemo
\proclaim {Theorem 3.1.3.2.(Existence of $\mu$-Limit Set)} The following holds:

mels1) $\mu_t\in \Cal M(\rho(r))$ for any $t\in (0,\infty)$;

mels2) the family $\{\mu_t\}$ is precompact in infinity,
\endproclaim
i.e.,
 for any sequence $t_k\rightarrow \infty$ there exists a subsequence
$t_{k_j}\rightarrow \infty$ and a measure $\nu\in \Cal M(\Bbb R^m)$ such that
$\mu_{t_{k_j}}\rightarrow \nu$ in $\Cal D'(\Bbb R^m)$ (see, 2.7.1).

\demo {Proof}We have
$$\mu_t (r)=\mu (rt)t^{-\rho(t)-m+2}.$$
Thus
$$\limsup_{r\rightarrow \infty}\frac {\mu_t(r)}{r^{\rho(r)+m-2}}=
\limsup_{r\rightarrow \infty}\frac {\mu (rt)}{(rt)^{\rho(rt)+m-2}}
\frac{(rt)^{\rho(rt)}}{ t^{\rho(t)+m-2}r^{\rho(r)}}=
t^{\rho-\rho(t)-(m-2)}\bar \Delta [\mu,\rho(r)].\tag 3.1.3.3$$
Therefore mels1) holds.

We also have
$$\limsup_{t\rightarrow \infty}\mu_t(r)=
\bar \Delta [\mu,\rho(r)]r^{\rho+m-2}.$$
Thus $\mu_t$ satisfies the assumption of the
Helly theorem (Theorem 2.2.3.2). Using also Theorem 2.3.4.4 we obtain
mels2). \qed
\enddemo

We will call the set of all measures $\nu$ from Theorem 3.1.2.1
the {\it limit set } of the mass distribution $\mu$ with respect to $V_\bullet$ and
denote it $\bold {Fr}[\mu,\rho(r),V_\bullet,\Bbb R^m]$ or shortly $\Fr [\mu].$

Set
$$\Cal M[\rho,\Delta]:=\{\nu: \nu(r)\leq \Delta r^{\rho+m-2},\  \forall r>0\}.
\tag 3.1.3.4$$
$$\Cal M[\rho]:=\bigcup _{\Delta>0}\Cal M[\rho,\Delta],$$
and
$$\nu_{[t]}(E):=t^{-\rho-m+2}\nu (P_t E)\tag 3.1.3.5$$
for $E\in \sigma (\Rm).$
\proclaim {Theorem 3.1.3.3.(Properties of $Fr[\mu]$)}The following holds:

frm1) $\Fr [\mu]$ is  connected and compact;

frm2) $\Fr[\mu]\subset \Cal M[\Delta,\rho] ,$ for $\Delta\geq
\bar \Delta [\mu,\rho(r)];$

frm3) $ (\Fr [\mu])_{[t]}=\Fr [\mu],\ t\in (0,\infty),$
\endproclaim
\demo {Proof}We will only prove frm2) because frm1) and frm3) are proved word by
 word as in Theorem 3.1.2.2.

Suppose $t_n\rightarrow \infty$ and $\mu_{t_n}\rightarrow \nu\in \Fr[\mu].$
Let us choose $r'>r$ such that the open ball $K_{r'}$ is squarable with
respect to
$\nu.$ It is possible because of Theorem 2.2.3.3, sqr2). By Theorems 2.2.3.7
and 2.3.4.4 $\mu_{t_n}(r')\rightarrow \nu (r').$ Thus (compare with (3.1.2.3))
$$\nu (r')=\lim_{t_n\rightarrow\infty}\mu_{t_n}(r')\leq
\limsup_{t\rightarrow\infty}\mu_{t}(r')=\bar \Delta [\mu,\rho(r)]
(r')^{\rho+m-2}.$$
 Choosing $r'\downarrow r$ we obtain
$$\lim_{r'\rightarrow r}\nu (r')=\nu(r)$$
because (2.2.3.3).
Thus frm2) holds.\qed
\enddemo
The following assertion is a ``copy'' of Theorem 3.1.2.4.
\proclaim {Theorem 3.1.3.4.(Properties of $\mu\mapsto \Fr [\mu]$)}
The following holds:

frmu1) $\Fr [\mu_1+\mu_2]\subset\Fr [\mu_1]+\Fr [\mu_2]$

frmu2) $\Fr [\lambda \mu]=\lambda \Fr[\mu]$.
\endproclaim
The proof is also a ``copy'' and we omit it.

\subheading {3.1.4}We are going to study the class $U[\rho]$ and obtain for
it  ``non-asymptotic'' analogies of Theorem 2.8.3.3 (Jensen), 2.9.2.3
(Brelot-Borel), 2.9.3.1 (Brelot - Hadamard)
\proclaim {Theorem 3.1.4.1.(*Jensen)} Let $v\in U[\rho].$ Then its Riesz
measure $\nu_v \in \Cal M[\rho].$
\endproclaim
\demo {Proof} As in Theorem 2.8.3.2 we have an inequality
$$\frac {\nu_v(r)}{r^{m-2}}\leq A(m)N(2r,\nu_v)\tag 3.1.4.1$$
Since $v(0)=0$ we have (Theorem 2.6.5.1.(Jensen-Privalov))
$$ N(2r,\nu_v)=\Cal M (2r,v)\leq M (2r,v)\leq 2^\rho \sigma r^\rho.\tag 3.1.4.2$$
Substituting (3.1.4.2) in (3.1.4.1) we obtain $\nu_v \in \Cal M[\rho,\Delta]$
for some $\Delta$. Thus $\nu_v \in \Cal M[\rho].$ \qed
\enddemo

Let $\rho$ be non-integer and $\nu\in \Cal M[\rho].$ Consider the
canonical potential $\Pi (x,\nu,p)$ where $p:=[\rho]$ (see
(2.9.2.1)). Let us emphasize that the support of $\nu$ may contain
the origin but $\nu(0)=0,$ i.e., there  is no concentrated mass in
the origin. Thus we must also check its convergence in the origin.
\proclaim {Theorem 3.1.4.2.(*Brelot-Borel)} Let $\rho$ be
non-integer and let $\nu\in \Cal M[\rho].$ Then $\Pi (x,\nu,p)$
converges and belongs to $U[\rho].$
\endproclaim
\demo {Proof} Using (2.9.1.9) we have
$$\left |\int_{|y|<2|x|}G_p (x,y,m)d\nu_y\right |\leq
A(m,p)|x|^p\int^{2|x|}_0 \frac {d\nu(t)}{t^{p+m-2}}.\tag 3.1.4.2$$
Let us estimate the integral in (3.1.4.2). Integrating by parts we obtain
$$I_< (x):=\int_0^{2|x|}\frac {d\nu(t)}{t^{p+m-2}}=
\frac {\nu(t)}{t^{p+m-2}}|_0^{2|x|}+(p+m-2)
\int_0^{2|x|}\frac {\nu(t)dt}{t^{p+m-1}}.$$
Since $\nu\in \Cal M[\rho,\Delta]$ for some $\Delta,$
$$I_< (x)\leq A(m,\rho,p)\Delta|x|^{\rho-p}.$$
Substituting this in (3.1.4.2) we obtain
$$\left |\int_{|y|<2|x|}G_p (x,y,m)d\nu_y\right |\leq A(m,\rho,p)\Delta
|x|^\rho .\tag 3.1.4.3$$
Analogously, using (2.9.1.8) we obtain
$$ \left |\int_{|x|<2|y|}G_p (x,y,m)d\nu_y\right |\leq A(m,\rho,p)\Delta
|x|^\rho .\tag 3.1.4.4$$
In particular, these estimates show that $\Pi(x,\nu,p)$ exists. Now using
(2.9.1.10) we have also
$$\int_{\frac {|x|}{2}\leq |y|\leq 2|x|}G_p (x,y,m)d\nu_y\leq
A(m,p)\int^{2|x|}_{\frac {|x|}{2}}d\nu (t)\min \left
(\frac {|x|^{p+1}}{t^{p+m-1}},\frac {|x|^{p}}{t^{p+m-2}}\right ).$$
The latter integral can also be easily estimated by
$\Delta A(m,p,\rho)|x|^\rho.$
Thus we have
$$\int_{\frac {|x|}{2}\leq |y|\leq 2|x|}G_p (x,y,m)d\nu_y\leq A(m,p,\rho)|x|^\rho.$$
Therefore by (3.1.4.4) and (3.1.4.3) we obtain $M(r,\Pi)\leq \sigma r^\rho$
for some $\sigma.$

Since $G_p (0,y,m)=0$ for all $y\neq 0$ and the integral converges,\linebreak
$\Pi(0,\nu,p)=0.$ \qed
\enddemo

We will need an assertion that looks like the Liouville theorem
(Theorem 2.4.2.3). \proclaim {Theorem 3.1.4.3.(*Liouville)}Let $H$
be a harmonic function in $\Rm$ and $H\in U[\rho].$ Then $H\equiv
0$ if $\rho$ is non-integer and $H$ is a homogeneous polynomial of
degree $p$ if $\rho=p$ is integer.
\endproclaim
In particular, for $m=2$ we have $H(re^{i\phi})=r^p \Re (ce^{ip\phi}).$
\demo {Proof}Like in the proof of the Liouville theorem we obtain the
inequality (2.4.2.9) and
$$|c_k|\leq A R^{-k}\max_{|x|=R}H(x)\leq A\sigma R^{\rho-k}$$
for some $\sigma>0.$

If $k>\rho,$ we will pass to the limit when $R\rightarrow\infty$ and obtain
$c_k=0.$ If $k<\rho$, we will do that when $R\rightarrow 0$ and obtain
$c_k=0.$ \qed
\enddemo
The following theorem can be considered as an analogy of the  Brelot-Hadamard
theorem (Theorem 2.9.3.1):
\proclaim {Theorem 3.1.4.4.(*Hadamard)} Let $\rho $ be non-integer and
$v\in U[\rho].$ Then
$$v(x)=\Pi (x,\nu_v,p)\tag 3.1.4.5$$
for $p=[\rho].$
\endproclaim
\demo{Proof}. Consider the function $H(x):=v(x)-\Pi(x,\nu_v,p).$ It is
harmonic.We also have by(2.8.2.5)
$$M(r,H)\leq A(m)T(r,H)\leq A(m) [T(r,v)+T(r,\Pi)]\leq \sigma r^\rho$$
for some $\sigma.$

Hence, $H(x)\equiv 0$ by Theorem 3.1.4.3. \qed
\enddemo
Let us consider the case of  integer $\rho.$

Let $\nu\in \Cal M [\rho]$ for an integer $\rho=p.$ Set
$$\Pi_< (x,\nu,\rho):=\int_{|y|<1}G_{p-1} (x,y,m)d\nu \tag 3.1.4.6$$
$$\Pi_> (x,\nu,\rho):=\int_{|y|\geq1}G_{p} (x,y,m)d\nu. \tag 3.1.4.7$$
Both potentials converge and belong to $U[\rho]$.
\proclaim {Theorem 3.1.4.5 (**Hadamard)}Let $\rho$ be integer and let
$v\in U[\rho]$. Then
$$v=H_\r(x)+\Pi_< (x,\nu,\rho)+\Pi_> (x,\nu,\rho),\tag 3.1.4.8$$
where $H_\r$ is a homogeneous harmonic polynomial of degree $\r.$
\endproclaim
The proof is exactly the same as in the *Hadamard theorem, but we
use the second case of  Theorem 3.1.4.3. We also note that the
polynomial may be equal to zero identically.

{\bf Exercise 3.1.4.1} Check this in details.

Let as check that $\nu$ from (3.1.4.8) has the following property that is
analogous to Theorem 4.9.4.2.
\proclaim {Theorem 3.1.4.6 (*Lindel\"of)} Let $\rho$ be integer and let
$v\in U[\rho]$. Then
$$\liml_{\eps\ri 0}\intl_{\eps\leq|y|<1}D_\r (x,y)\mu(dy)=H_\r(x)\tag 3.1.4.9$$
\ep
\demo {Proof} Consider the function
$$v^*_\eps(x):=v(x)+\cases \intl_{|y|<\eps}\frac {\nu (dy)}{|x-y|^{m-2}}
\ for\ m>2;\\
-\intl_{|y|<\eps}\log |x-y|\nu (dy),\ for\ m=2.\endcases
\tag 3.1.4.10$$
It is subharmonic with $\supp\ \nu\cap \{0\}=\varnothing.$
We represent this function like in (2.9.4.10) in the form
$$v^*_\eps(x)= \Pi^1_< (x,\nu^*_\eps,\rho)+\Pi^1_> (x,\nu^*_\eps,\rho) +
P^*_{\r-1} (x,v^*_\eps)+
\dl_1 (x,v^*_\eps,\r)$$
In this representation we can pass to limit as $\eps\ri 0$ in the left side
and in all the summands except may the last two from the right side.

{\bf Exercise 3.1.4.2} Check this, using that all the integrals converge for
$\nu\in\Cal M(\r,\Dl)$ and showing that the integral in (3.1.4.10) tends to
zero.

The last two summands forms a harmonic polynomial, the limit of
which is also a harmonic polynomial. Comparing the limit with the
representation (3.1.4.8), we obtain that
$P^*_{\r-1}(\bul,v^*_\eps)$ tends to zero and $\dl_1
(x,v^*_\eps,\r)$ tends to $H_\r (x).$ \qed
\enddemo
\proclaim {Theorem 3.1.4.7. (**Liouville)} If $v\in U[\r]$
satisfies inequality $v(x)\leq 0$ for $z\in \BR^m $ then
$v(x)\equiv 0.$\ep Otherwise it contradicts to subharmonicity in
0.

\subheading {3.1.5} Let us study the connection between $\Fr[u]$ and
$\Fr [\mu_u]$.

Note the following properties of the transformations $(\bul)_t$ and
$(\bul)_{[t]}.$
\proclaim {Theorem 3.1.5.0 (Connection between $u_t$ and $\mu_t$)}One has
$$(\mu_u)_t=\mu_{u_t};\ (\mu_v)_{[t]}=\mu_{v_{[t]}}.\tag 3.1.5.0$$
\ep
\demo {Proof} By the F.Riesz theorem (Theorem 2.6.4.3) and Theorem 2.5.1.1 ,GPo3)
we have for any $\psi\in \Di'(\Rm)$
$$<\mu_u,\psi>=\theta_m<\Dl u,\psi>=\theta_m< u,\Dl\psi>.$$
Using  the definition (3.1.3.1), we obtain
$$<(\mu_u)_t,\psi>=<(\mu_u)_t, \psi ((P_t)^{-1}\bul)>t^{-\r(t)-m+2}.$$
Thus
$$<(\mu_u)_t,\psi>=\theta_m< u,\Dl[\psi ((P_t)^{-1}\bul)]>t^{-\r(t)-m+2}.$$
Since the Laplace operator is invariant with respect to $V_t$ for any $t$
we have
$$\Dl[\psi ((P_t)^{-1}\bul)]=t^{-2}[\Dl\psi] ((P_t)^{-1}\bul).$$
Thus we obtain
$$<(\mu_u)_t,\psi>=\theta_m t^{-\r(t)-m}< u,[\Dl\psi] ((P_t)^{-1}\bul)>=$$
$$=\theta_m <u(P_t \bul)t^{-\r(t)},\Dl\psi>=\theta_m <u_t,\Dl\psi>=
<\mu_{u_t},\psi>.$$
\qed
\edm

{\bf Exercise 3.1.5.1.}Do this for $(\bul)_{[t]}.$

 We begin from the case of a non-integer $\rho$.
\proclaim {Theorem 3.1.5.1 (Connection between $\Fr$'s for non-integer $\rho$)}
Let $u\in U(\rho (r))$ and $\mu_u$ be its Riesz measure. Then
$$\Fr[\mu_u]=\{\nu_v:v\in \Fr[u]\},\tag 3.1.5.1$$
$$\Fr[u]=\{\Pi (\bullet ,\nu,p):\nu\in \Fr[\mu_u]\}.\tag 3.1.5.2$$
\endproclaim
\demo {Proof}Let $\nu\in \Fr[\mu_u].$ There exists $t_n\rightarrow \infty$
such that $(\mu_u)_{t_n}\rightarrow \nu$ in $\Di'.$ We can find a subsequence
$t'_n$ such that $u_{t'_n}\rightarrow v\in Fr [u].$ Thus
$(\mu_u)_{t'_n}\rightarrow \nu_v$ and therefore $\nu =\nu_v.$ Hence,
$Fr[\mu_u]\subset\{\nu_v:v\in \Fr[u]\}.$ Analogously we can prove that every
$\nu_v \in Fr [\mu_u]$ and hence (3.1.5.1) holds.

Let $\nu\in \Fr[\mu_u].$ We find a  sequence $t_n\rightarrow \infty$ such
that $(\mu_u)_{t_n}\rightarrow \nu$ in $\Di'.$ We find a subsequence
$t'_n$ such that $u_{t'_n}\rightarrow v\in Fr [u]$ and $\nu_v =\nu.$ By
the *Hadamard theorem (Theorem 3.1.4.4) $v= \Pi (\bullet ,\nu,p).$ Hence,
$\{\Pi (\bullet ,\nu,p):\nu\in \Fr[\mu_u]\}\subset \Fr[u].$ And vice versa,
since $ \Fr[u]\subset U[\rho]$ (Theorem 3.1.2.2, fr2)), every $v\in \Fr[u]$
is represented as $\Pi (\bullet ,\nu_v,p)$ and $\nu_v\in \Fr[\mu]$ by
(3.1.5.1).\qed
\enddemo

Let $\rho$ be integer and $u\in U(\rho(r))$. Let us consider the
precompact family of homogeneous polynomials $\delta_t
(x,u,\rho)t^{\rho-\rho(t)}$ from Theorem 2.9.4.1. For every
$t_n\rightarrow \infty$ we can find a subsequence $t'_n$ such that
the pair $(\delta_{t'_n} (\bullet,u,\rho){t'_n}^{\rho-\rho(t'_n)},
(\mu_u)_{t'_n})$ tends to a pair $(H_{\nu},\nu)$ where $H_\nu$ is
a homogeneous harmonic polynomial of degree $p.$ We denote the set
of all such pairs  as $(\Cal H, Fr)[u].$ Every $v\in U[\rho]$ can
be represented in the form (3.1.4.7). Thus for every $v$ the
polynomial $H^v:=H_p$ is determined. \proclaim {Theorem 3.1.5.2
(Connection between $\Fr$'s for integer $\rho$)} Let $u\in U(\rho
(r)).$ Then
$$(\Cal H, Fr)[u]=\{(H^v,\nu_v):v\in \Fr[u]\},\tag 3.1.5.3$$
$$\Fr[u]=\{v:=H_{\nu}+\Pi_< (\bullet ,\nu,\rho )+\Pi_> (\bullet ,\nu,\rho )
:(H_\nu,\nu)\in (\Cal H, Fr)[u]\}.\tag 3.1.5.4$$
\endproclaim
The proof is clear. \subheading {3.1.6} Up to now we supposed that
the family of rotations $V_\bullet$ was fixed. Now we take in
consideration that it can be vary and use the notation $\Fr
[u,V_\bullet].$ \proclaim {Theorem 3.1.6.1.(Dependence of $Fr$ on
$V_\bullet$)} Let $\Fr [u,V_\bullet]$ and $\Fr [u,W_\bullet]$ be
limit sets of $u$ with respect to rotation families $V_\bullet$
and $W_\bullet$ accordingly. Then for any $v\in \Fr[u,V_\bullet] $
there exist a rotation $V^v$ and $w^v \in \Fr [u,W_\bullet]$ such
that
$$v(x)=w^v(V^v x)
$$
for all $x\in \Bbb R^m.$
\endproclaim
\demo {Proof} Let $v\in \Fr[u,V_\bullet]$ and let $t_n\rightarrow\infty$ be a sequence
such that $t_n^{-\rho(t_n)} u(t_n V_{t_n}\bullet)\rightarrow v.$ Since the
family
$V_t$ is obviously precompact there exists a subsequence (for which we keep
the same notation),and a rotation $V^v$ such that
$W_{t_n}^{-1}V_{t_n}\rightarrow V^v$ and $w\in \Fr [u,W_\bullet]$ such that
$t_n^{-\rho(t_n)} u(t_n W_{t_n}\bullet)\rightarrow w.$
Now we have
$$v(\bullet)=\Di'-\lim t_n^{-\rho(t_n)} u(t_n V_{t_n}\bullet)=
\Di'-\lim t_n^{-\rho(t_n)} u(t_n W_{t_n}W_{t_n}^{-1}V_{t_n}\bullet)=
w(V^v \bullet).$$ \qed
\enddemo
\newpage
\centerline {\bf 3.2.Indicators}
\subheading {3.2.1}Let $u\in SH(\rho(r))$ and let $\Fr [u]$ be the limit set.
Set
$$h(x,u):=\sup \{v(x):v\in \Fr[u]\}\tag 3.2.1.1$$
$$\lh (x,u):=\inf \{v(x):v\in \Fr [u]\}.\tag 3.2.1.2$$
These functionals reflect the asymptotic behavior of $u$ along rays of the
form
$$l_{\bold x ^0}:=\{x=t\bold x ^0:t\in (0,\infty)\}\tag 3.2.1.3$$
and are called {\it indicator} of growth of $u$ and {\it lower indicator}
respectively.

Of course, the indicators depend on $\rho(r)$ and $V_t,$ but we will only
note that if necessary.

\proclaim {Theorem 3.2.1.1(Properties of Indicators)}The following holds

h1) $\lh$ is upper semicontinuous, $h$ is subharmonic;

h2) they are semiadditive and positively homogeneous, i.e.,
$$h(x,u_1+u_2)\leq h(x,u_1)+h(x,u_2);\tag 3.2.1.4$$
$$\lh(x,u_1+u_2)\geq \lh (x,u_1)+\lh (x,u_2);\tag 3.2.1.5$$
$$h,\lh (x,Cu)=Ch,\lh (x,u),\ C\geq 0;\tag 3.2.1.6$$

h3) invariance:
$$h,\lh _{[t]}(x,u)= h,\lh (x,u).\tag  3.2.1.7$$
\endproclaim
\demo {Proof}Semicontinuity of $\lh$ follows from Theorem 2.1.2.8.(Commutativity of inf and M(.). Semicontinuity and subharmonicity of $h$ follow from Theorem  2.7.3.4.(Sigurdsson's Lemma). The properties h2) follow from properties of infimum and supremum. The invariance follows from
invariance of $\Fr [u]$ (Theorem 3.1.2.2, fr3)). \qed

\enddemo
Set
$$x^0 (x):=P^{-1}_{|x|}(x)\tag 3.2.1.8$$
where $P_t$ is defined by (3.1.1.1).

This is an intersection of the orbit of $P_t$ that passes through a point $x$
 with the unit sphere.

If $V_t\equiv I$
$$x^0(x)=x/|x|:=x^0\tag 3.2.1.9$$
\proclaim {Theorem 3.2.1.2(Homogeneity $h,\lh$)}One has
 $$h,\lh (x,\bullet)=|x|^\rho h,\lh (x^0(x),\bullet)\tag 3.2.1.10$$
\endproclaim
Thus the indicators are determined uniquely by their values on the unit sphere,i.e., they are ``functions of direction.In particular, they are homogeneous
for $V_t\equiv I$:
$$h,\lh (x,\bullet)=|x|^\rho h,\lh (x^0,\bullet)\tag 3.2.1.11$$
The proof of (3.2.1.10) follows from h4), Theorem 3.2.1.1 if we set
$t:=|x|;\ x:=P_t ^{-1}x.$
\subheading {3.2.2}In this item we will suppose that $V_t\equiv I$
 and study the indicator.

Let $\Delta _{\bold x^0}$ be defined in 2.4.1. Its coefficients depend on a choice
of  the spherical coordinate system. However, one has
\proclaim {Theorem 3.2.2.1}Let $\psi(y)$ have continuous second derivatives on
the unit sphere $S_1.$ Then the differential form $\Deltax0 \psi (y)dy$ is
invariant with respect to the choice of  spherical coordinate system.
\endproclaim
\demo{Proof}Let $\phi(x)$ be a smooth function in $\Rm.$ Then
$\Delta\phi (x)dx$ is invariant with respect to the choice of an orthogonal
system because $\Delta$ (the Laplace operator) and an element of volume are
invariant. Set $\phi(x)=\psi (y),$ where $y:=x^0/|x|.$ Then
$$\Delta \phi dx =\Deltax0 \psi (y)dy\  m r^{m-3}dr.$$
Since $r$ is invariant with respect to rotations of the coordinate system,
$\Deltax0 \psi (y)dy$ is invariant with respect to the choice of a  spherical
coordinate system.\qed
\enddemo
Note that for $m=2$ this theorem is obvious because
$$\Deltax0 =\frac {d^2}{d\theta^2}\tag 3.2.2.1$$
and it does not depend on  translations with respect to $\theta.$

We define the operator $\Deltax0 $ on $f\in \Di' (S_1)$ by
$$<\Deltax0 f,\psi>:=<f,\Deltax0 \psi>,\ \psi\in \Di(S_1)$$
in a fixed spherical coordinate system.

The definition is correct. Indeed, suppose in a fixed system
$$ \roman {supp}\psi \subset S_1\backslash
\{\theta _j=0;\pi:j=1,2,...,m-2\}.\tag 3.2.2.2$$
Then all the coefficients of $\Deltax0$ are infinitely differentiable  and
$\Deltax0 \psi (y)\in \Di(S_1).$ By Theorem 3.2.2.1 we obtain that the condition of Theorem 2.3.5.2.($\Di'$ on Sphere) are
fulfilled.

Note that for $m=2$ the operator $\Deltax0$ is realized by the formula
(3.2.2.1) on functions of the form $f=f(e^{i\theta}),$ i.e., on $2\pi$-periodic
functions.
\proclaim {Theorem 3.2.2.2(Subsphericality of Indicator)}One has
$$[\Deltax0 +\rho(\rho+m-2)]h(y,u):=s >0 \tag 3.2.2.3$$
 in $ \Di'(S_1),$
\endproclaim
i.e., $s$ is a measure on $S_1.$
\demo {Proof} It is sufficient to prove this locally, in any spherical system.
Let $R(r)$ be finite, infinitely differentiable and non-negative in
$(0;\infty)$ and let $\psi\in \Di(S_1)$ be non-negative and satisfy
(3.2.2.2). Set $\phi (x):=R(|x|)\psi(x^0).$ Using the subharmonicity of
$h (x,u)$
(h1), Theorem 3.2.1.1 and (3.2.2.2), we have
$$0\leq \int h(x,u)\Delta \phi (x)dx =$$
$$=\int _{(y,r)\in S_1\times (0;\infty)}r^\rho h(y,u)
\left [\frac {1}{r^{m-1}}\frac {\partial}{\partial r}r^{m-1}
\frac {\partial}{\partial r}+\frac {1}{r^2}\Deltax0\right ]
\psi(y) r^{m-1}dydr.$$
Transforming the last integral we obtain
$$\int h(x,u)\Delta \phi (x)dx =$$
$$=\int_0^\infty r^\rho\left [\frac {1}{ r^{m-1}}
\frac {\partial}{\partial r}r^{m-1}
\frac {\partial}{\partial r}R(r)\right ]r^{m-1}dr
\int_{S_1}h(y,u)\psi (y)dy+\tag 3.2.2.4 $$
$$+\int_0^\infty r^{\rho-2}r^{m-1}R(r)dr
\int_{S_1}h(y,u)\Deltax0\psi(y)dy. $$
Integrating by parts in the first summand we obtain
$$ \int_0^\infty r^\rho\left [\frac {1}{ r^{m-1}}
\frac {\partial}{\partial r}r^{m-1}
\frac {\partial}{\partial r}R(r)\right ]r^{m-1}dr=
\int_0^\infty R(r)\rho(\rho+m-2)r^{\rho+m-3}dr.\tag 3.2.2.5$$
Substituting (3.2.2.5) into (3.2.2.4), we have
$$0\leq \int_0^\infty r^{\rho+m-3}R(r)dr
\int_{S_1}h(y,u)[\Deltax0 +\rho(\rho+m-2)]\psi(y)dy.$$
Since $R(r)$ is an arbitrary non-negative function,
$$\int_{S_1}h(y,u)[\Deltax0 +\rho(\rho+m-2)]\psi(y)dy\geq0$$
for arbitrary $\psi.$ \qed
\enddemo
We will call an upper semicontinuous function which satisfies (3.2.2.3)
a $\rho$-{\it subspherical} one.
Now we are going to study properties of these functions.
\subheading {3.2.3}We consider the case $m=2.$ A $\rho$-subspherical function
for $m=2$ is called  $\rho$-{\it trigonometrically convex} ($\rho$-t.c.). We
will obtain for such a function a representation  like in Theorems 3.1.4.4,
3.1.4.5.(*,** Hadamard). Set
$$T_\rho:=\frac {d^2}{d\phi^2}+\rho^2. $$

Let us find a fundamental solution of this operator.

Let $\rho$ be non-integer. Let us denote as $\widetilde {\cos\rho}(\phi)$
the periodic continuation of $\cos\rho\phi$ from the interval $(-\pi,\pi).$
\proclaim {Theorem 3.2.3.1.(Fundamental Solution of $T_\rho$)}One has
$$\frac {1}{2\rho \sin\pi\rho}T_\rho \widetilde {\cos\rho}(\phi-\pi)
=\delta (\phi)\ \text {in $\Di'(S_1)$ }$$
\endproclaim
\demo {Proof}Let $f\in \Di (S_1).$ We have
$$\int_0^{2\pi}\widetilde {\cos\rho}(\phi-\pi)[f''+\rho^2 f]d\phi =
\lim_{\epsilon\rightarrow 0}\int_\epsilon ^{2\pi-\epsilon}
 \cos\rho(\phi-\pi)[f''+\rho^2 f]d\phi\tag 3.2.3.1 $$
Integrating by parts  we obtain
$$\int_\epsilon ^{2\pi-\epsilon}
 \cos\rho(\phi-\pi)[f''+\rho^2 f]d\phi=
\cos\rho(\phi-\pi)\left .f'(\phi)\right |_\epsilon ^{2\pi-\epsilon}+
\rho\sin\rho(\phi-\pi)\left .f(\phi)\right |_\epsilon ^{2\pi-\epsilon}+$$
$$+\int_\epsilon ^{2\pi-\epsilon} f(\phi)T_\rho \cos\rho(\phi-\pi)d\phi.$$
However, $T_\rho \cos\rho(\phi-\pi)=0$ for $\phi\in (\epsilon,2\pi-\epsilon).$
Thus the limit in (3.2.3.1) is equal to $f(0)2\rho\sin\pi\rho.$ \qed
\enddemo
Let $s$ be a measure on the circle $S_1.$ Set
$$\Pi (\phi,s):=\int_0^{2\pi}\widetilde {\cos\rho}(\phi-\psi-\pi)s(d\psi).$$
\proclaim {Theorem 3.2.3.2.}One has
$$T_\rho \Pi (\phi,s)=(2\rho\sin\pi\rho) ds\ \roman {in} \Di'(S_1) .$$
\endproclaim
The proof is the same as GPo3) in Theorem 2.5.1.1.
\proclaim {Theorem 3.2.3.3(Representation of $\rho$-t.c.f for a non-integer
$\rho$)}Let $h$ be $\rho$-t.c. on $S_1$ for non-integer $\rho$ and let
$s:=T_\rho h.$ Then
$$h(\phi)=\frac {1}{2\rho\sin\pi\rho}\Pi (\phi,s).$$
\endproclaim
The proof is like in Theorem 3.1.4.4(*Hadamard).
\subheading {3.2.4}We will suppose in this item that $V_t=I,m=2,\rho$ is
integer.

\proclaim {Theorem 3.2.4.1(Condition on $s$)}Let $\rho$ be integer,
$h$ be $\rho$-t.c. and $T_\rho h=s$. Then
$$\int_0^{2\pi}e^{i\rho \phi}ds =0\tag 3.2.4.1$$
\endproclaim
\demo{Proof}We have for $f\in\Di(S_1):$
$$<s,f>=<T_\rho h,f>= <h,T_\rho f>.$$
Since $e^{i\rho \phi}\in\Di(S_1) $ for integer $\rho$ and
$T_\rho e^{i\rho \phi}=0,$ we have for $f:= e^{i\rho \phi}$
$$<s,e^{i\rho \bullet}>= <h,T_\rho e^{i\rho \bullet}>=0.$$
\qed
\enddemo
Let us denote the periodic continuation of the function $f(\phi):=\phi$ from
the interval $[0,2\pi)$ to $(-\infty,\infty)$ as $\tilde \phi.$
\proclaim {Theorem 3.2.4.2.(Generalized Fundamental Solution for $T_\rho$)}
One has
$$T_\rho [-\frac {1}{2\pi\rho} \tilde \phi\sin\rho\phi]=
\delta (\phi)-\frac{1}{\pi}\cos\rho\phi$$
in $\Di'(S_1).$
\endproclaim
\demo {Proof}Let $\phi\in (\epsilon, 2\pi -\epsilon).$ Then
$$T_\rho \tilde \phi\sin\rho\phi=2\rho \cos\rho\phi$$
because $\tilde \phi=\phi$ when $\phi\in (\epsilon, 2\pi -\epsilon).$
We have also
$$(\phi\sin\rho\phi)'=\sin\rho\phi+\phi\rho\cos\rho\phi.$$
Thus
$$<T_\rho\tilde \bullet\sin\rho\bullet,f>=\int_0^{2\pi} \tilde \phi\sin\rho\phi
\ T_\rho fd\phi=\lim_{\epsilon\rightarrow 0}\int_{\epsilon}^{2\pi-\epsilon}
\phi\sin\rho\phi\ T_\rho f d\phi.$$
Integrating by parts we obtain
$$\int_{\epsilon}^{2\pi-\epsilon}
\phi\sin\rho\phi\ T_\rho f d\phi=\phi\sin\rho\phi f'(\phi)
|_\epsilon^{2\pi-\epsilon}-$$
$$-f(\phi)[\sin\rho\phi+\phi\rho\cos\rho\phi|_\epsilon^{2\pi-\epsilon}+
\int_{\epsilon}^{2\pi-\epsilon}T_\rho[\phi\sin\rho\phi]f(\phi)d\phi.$$
Passing to the limit as $\epsilon\rightarrow 0$ and taking in account
that $f $ is periodic and continuous we obtain
$$<T_\rho[\tilde \bullet\sin\rho\bullet],f>=-2\pi\rho f(0)+2\rho\int_0^{2\pi}
\cos\rho\phi f(\phi)d\phi=-2\pi\rho f(0)+ <\cos\rho\bullet, f>.$$
\qed
\enddemo
Set
$$\hat\Pi (\phi,ds):=\int_0^{2\pi}\widetilde {(\phi -\psi)}
\sin\rho (\phi-\psi)s(d\psi).$$
\proclaim {Theorem 3.2.4.3.}One has
$$T_\rho\hat\Pi (\bullet,ds)=-2\pi\r ds\ \roman {in} \ \Di'(S_1)$$
for $s$ that satisfies (3.2.4.1).
\endproclaim
\demo {Proof}Using Theorem 3.2.4.2 we obtain
$$<T_\rho \hat\Pi (\bullet,ds),f>=<s,f>-
\frac {1}{\pi}<\int _0^{2\pi}
\cos\rho (\bullet-\psi)ds_\psi,f>.$$
The last integral is zero because of Theorem 3.2.4.1. \qed
\enddemo
\proclaim {Theorem 3.2.4.4.(Representation of $\rho$-t.c.f. for an integer
$\rho$)} Let $h$ be a $\rho$-t.c.f.for an integer $\rho$ and
$T_\rho h:=s .$ Then
$$h(\phi)=\Re ce^{i\phi}+\hat\Pi (\phi,ds).$$
for some complex constant $c.$
\endproclaim
\demo {Proof}The function
$H(\phi):=h(\phi)-\hat\Pi (\phi,ds)$
satisfies the equation $T_\rho H=0$ in $\Di'(S_1)$ because of Theorem 3.2.4.3
 and it is real.Thus $H(\phi)=\Re ce^{i\phi}.$ \qed
\enddemo
\subheading {3.2.5}The class $TC_\rho$ of $\rho$-t.c.functions has a number
of properties of subharmonic functions.

The function $\widetilde {\cos\rho}\phi$ is continuous and
$\widetilde \phi \sin\rho\phi$ is continuous for  integer $ \rho.$ Therefore
any $\rho$-t.c.f is continuous as  follows from Theorem 3.2.3.3 and
3.2.4.4.

Set
$$\Cal E (\phi):=\frac {1}{2\rho}\sin\rho|\phi|.$$
For any interval $I:=(\alpha,\beta)\Subset (-\pi,\pi)$ this function
satisfies the equality
$$ T_\rho \Cal E =\delta $$
in $\Di'(\alpha,\beta)$, where $\delta$ is the Dirac function in zero.

Let $G_{I}(\psi,phi)$ be the Green function of $T_\rho$ for the interval
$I.$ By definition it must be symmetric with respect to $\phi,\psi$ and
have the form
$$G_I (\phi,\psi):= \frac {1}{2\rho}\sin\rho|\phi-\psi|
 +A_I\cos\rho\phi\cos\rho\psi+ B_I \sin\rho\phi\sin\rho\psi,\tag 3.2.5.1$$
where $A_I,B_I$ are chosen such that $G_I(\phi,\psi)$ be equal to zero on
$\partial \{I\times I\}.$ An explicit form of $G_I$ is given by
$$G_I(\phi,\psi)=\cases \frac {\sin\rho(\beta - \phi)\sin\rho (\psi -\alpha)}
{\rho\sin\rho(\beta -\alpha)},&\text{ for } \psi<\phi;\\
\frac {\sin\rho(\beta - \psi)\sin\rho (\phi -\alpha)}
{\rho\sin\rho(\beta -\alpha)},&\text{ for } \phi<\psi.\endcases$$
The following assertion analogous to the Riesz theorem (Theorem 2.6.4.3):
\proclaim {Theorem 3.2.5.1.(Representation on $I$)}Let $h\in TC_\rho $ and
the let $I$ be an interval of length $mes I<\pi/\rho.$ Then
$$h(\phi)=Y_\rho(\phi,h)-\int\limits_\alpha^\beta G_I(\phi,\psi)s(d\psi),$$
where $Y_\rho(\phi,h)$ is the only solution of the boundary problem:
$$T_\rho Y=0,\ Y(\alpha)=h(\alpha),\ Y(\beta)=h(\beta)\tag 3.2.5.2$$
and $ds:=T_\rho h.$
\endproclaim
\demo{Proof}Set
$$\Pi_I(\phi,ds):= \int\limits_\alpha^\beta G_I(\phi,\psi)s(d\psi)$$
One can check like in Theorem 3.2.3.2 that $T_\rho\Pi_I=-ds$ in $\Di'(I).$
Then the function
$$Y_\rho(\phi):=h(\phi)+\Pi_I(\phi,ds)$$
satisfies the conditions (3.2.5.2).\ \qed
\enddemo
The explicit form of $Y_\rho (\phi)$ is
$$Y_\rho (\phi)=\frac {h(\alpha)\sin\rho (\beta -\phi)+
h(\beta)\sin\rho (\phi-\alpha)}{\sin\rho (\beta-\phi)}.\tag 3.2.5.3$$
Since $\Pi_I(\phi)\geq 0$ we have
\proclaim {Theorem 3.2.5.2. ($\rho$-Trigonometric Majorant)}Suppose $h\in TC_\rho $
and $Y_\rho (\phi)$ is the solution of (3.2.5.2). Then                         $$h(\phi)\leq Y_\rho (\phi),\ \phi\in I $$
if $\be -\a<\pi/\r.$
\endproclaim
This inequality can be written in the symmetric form
$$h(\alpha)\sin\rho(\beta-\phi)+h(\phi)\sin\rho(\alpha -\beta)+
h(\beta)\sin\rho(\phi-\alpha)\geq 0\tag 3.2.5.4$$
for $\max (\alpha,\phi,\beta)-\min(\alpha,\phi,\beta)<\pi/\rho.$
It is called the {\it fundamental relation of indicator.}
\proclaim {Theorem 3.2.5.3.(Subharmonicity and $\rho$-t.c.)}A function
$h(\phi)\in TC_\rho$  iff the function $u(re^{i\phi}):=h(\phi)r^\rho$ is
subharmonic in $\Bbb R^2.$
\endproclaim
\demo {Proof} Sufficiency follows from Theorem 3.2.2.2. Let us prove necessity.
The function $u_1(z):=r^\rho \sin \rho|\phi|$ is subharmonic. Actually, it is
harmonic for $\phi\neq 0 ,\ r\neq 0$ and can be represented in the form
$$ u_1(z)=\max (r^\rho\sin\phi,-r^\rho\sin\phi)$$
in a neighborhood of the line $\phi=0.$
Hence, it is subharmonic because of sh2), Theorem 2.6.1.1
(Elementary Properties).

The function $u_2(z):=r^\rho\Pi_I(\phi)$ is subharmonic because of sh5)
 and sh4),
Theorem 2.6.1.1. The function $r^\rho Y_\rho(\phi)$  is harmonic for
$r>0$. This can be checked directly. Hence, $u(z)$ is subharmonic for $r>0$
because of  Theorem 3.2.5.1. By Theorem 2.6.2.2 $u(z)$ is also subharmonic
for $r=0,$ because it is, obviously, continuous at $z=0.$ \qed
\enddemo
\proclaim {Theorem 3.2.5.4.(Elementary Properties of $\rho$-t.c.Functions)}One
 has

tc1) If $h\in TC_\rho,$ then $Ah\in TC_\rho$ for $A>0;$

tc2) If $h_1,h_2\in TC_\rho,$ then $h_1+h_2,\ \max (h_1,h_2)\in TC_\rho.$
\endproclaim
These properties follow from Theorem 3.2.5.3 and properties of subharmonic
functions.

{\bf Exercise 3.2.5.1} Prove Th.3.2.5.4.

Similarly to (usual) convexity, $\rho$-t.convexity of functions implies several analytic properties.
\proclaim {Theorem 3.2.5.5}Let $h\in TC_\rho.$ then there exist right
$(h'_+ )$ and left $(h'_- )$ derivatives and they coincide everywhere except,
maybe,for countable set of points.\ep

\demo {Proof} It is enough to proof these properties for the potential
$$\Pi(\phi):=\int\limits_\alpha ^\beta \sin\rho |\phi -\psi |ds_\psi,$$
because of (3.2.5.1) and Theorem 3.2.5.1.

We will prove the following
$$\Pi'_+(\phi)=\rho\int\limits_\alpha^{\phi -0}\cos\rho(\phi -\psi )ds_\psi
+\rho\mu (\phi) -\rho \int\limits^\beta_{\phi +0}\cos\rho(\phi -\psi )ds_\psi;
\tag 3.2.5.5$$
$$\Pi'_-(\phi)=\rho\int\limits_\alpha^{\phi -0}\cos\rho(\phi -\psi )ds_\psi
-\rho\mu (\phi) -\rho \int\limits^\beta_{\phi +0}\cos\rho(\phi -\psi )ds_\psi,
\tag 3.2.5.6$$
  where $\mu (\phi)$ is the measure, concentrated in the point $\phi.$

We have for $\Delta>0$:
$$\frac {\Pi(\phi +\Delta}{\Delta}=\int\limits_\alpha^{\phi -0}
\frac {\sin\rho |\phi +\Delta -\psi |-\sin\rho |\phi -\psi |}{\Delta}ds_\psi+$$
$$+\frac {\sin\rho\Delta}{\Delta}\mu(\phi)+\int\limits_{\phi +0}^{\phi +\Delta}
...+\int\limits^\beta_{\phi +\Delta}... .$$
Let us estimate the second integral.We have
$$\int\limits_{\phi +0}^{\phi +\Delta}
\left|\frac {\sin\rho |\phi +\Delta -\psi |-\sin\rho |\phi -\psi |}{\Delta}
\right|ds_\psi\leq $$
$$\leq\frac {2\sin\rho\Delta}{\Delta}[s(\phi +\Delta)-s(\phi +0)]=o(1)$$
when $\Delta\rightarrow +0.$

Passing to the limit, we obtain (3.2.5.5). The equality (3.2.5.6) is obtained
by the same way when $\Delta<0.$

Since $\mu(\phi)\neq 0$ at most in a countable set for all the other points
$\Pi'_+(\phi)=\Pi'_-(\phi).$ \ \qed
\enddemo
\subheading {3.2.6} Now we consider the case $m\geq 3.$ We will obtain for
the $\rho$-subspherical function a representation like for the
$\rho$ -trigonometrically convex functions.
\proclaim {Theorem 3.2.6.1 (Subharmonicity and Subsphericality)} Let $h$ be
 subspherical in a neighborhood of $y\in S_1.$ Then the function
$u(x):=h(y)r^\rho,\ x=ry$
is subharmonic in the corresponding neighborhood of the ray
$x=ry:0<r<\infty.$
\endproclaim
\demo {Proof} Let $f\in \Di'(\Rm\setminus 0).$ We can represent it in the
form $f:=f(rx),\ x\in S_1$ where $f(\bullet,x)\in \Di'(0,\infty)$ for any $x.$

Then
$$<u,\Delta f>=\int\limits_0^\infty \int\limits_{S_1} u(rx)
\Delta f(rx)r^{m-1}drds_x =$$
$$=\int\limits_0^\infty\int\limits_{S_1} u(rx)\frac {1}{r^{m-1
}}
\frac {\partial}{\partial r}r^{m-1}\frac {\partial}{\partial r}
f(rx)r^{m-1}drds_x+
\int\limits_0^\infty\int\limits_{S_1} u(rx)\frac {1}{r^2}
\Delta_x0 f(rx)r^{m-1}drds_x.$$
 Integrating by parts in the first integral, we obtain
$$ \int\limits_0^\infty\rho(\rho+m-2)r^{\rho+m-3}
\int\limits_{S_1}h(x)f(rx)drds_x.$$
Set
$$\Cal S_\rho :=\Deltax0 +\rho(\rho+m-2)\tag 3.2.6.1$$
Together with the second summand we obtain
$$<u,f>=\int\limits_0^\infty\left
[\int\limits_{S_1}h(x)\Cal S_\rho f(rx)ds_x\right ]
r^{\rho+m-3}dr >0$$
if $f(rx)\geq 0.$ \qed
\enddemo

Note that the Riesz measure for such $u$ has the form
$$\mu(r^{m-1}drds_x)=r^{\rho+m-3}dr\nu_h(ds_x),$$
where $\nu_h$ is a positive measure on $S_1,$ that is equal to
$\Cal S_\rho h$ in $\Di'(S_1).$

For a non-integer $\rho$ set
$$\Cal E_\rho (x,y):=\int\limits_0^\infty G_p(x,ry,m)r^{\rho+m-3}dr,$$
where $G_p$ is the primary kernel and $x,y\in S_1.$
\proclaim {Theorem 3.2.6.2} For non-integer $\rho$ and any $\rho$-subspherical
 function $h$ one has
$$h(x)=\int\limits_{S_1}\Cal E_\rho (x,y)d\nu_h.$$
\endproclaim
\demo {Proof}Set in (3.1.4.5) $v:=r^\rho h(x).$ It is clear that
$v\in U[\rho].$ We have
$$r^\rho h(x)=\int\limits_{S_1}\int\limits_0^\infty G_p(rx,ty,m)
t^{\rho+m-3}dt\nu_h(ds_x).$$
Now we make the change $t':=t/r$ and use the homogeneity of $G_p(rx,ty,m).$ \qed
\enddemo
\subheading {Exercise 3.2.6.1} Show that $\Cal E_\rho (x,y)$ is a fundamental
solution of the operator $\Cal S_\rho.$

For an integer $\rho=p$ set
$$\Cal E'_\rho (x,y):=\int\limits_0^1 G_{p-1}(x,ry)r^{\rho+m-3}dr+
\int\limits_1^\infty G_p(x,ry)r^{\rho+m-3}dr.$$
\subheading {Exercise 3.2.6.2} Prove.
\proclaim {Theorem 3.2.6.3} For any integer $\rho=p$ and any
$\rho$-subspherical function $h$ one has
$$h(x)=Y_p(x)+\int\limits_{S_1}\Cal E'_\rho (x,y)d\nu_h.$$
where $Y_p$ is some $p$-spherical function.

For any $p$-spherical function $Y$
$$\int\limits_{S_1}Y(x)d\nu_h=0.$$
\endproclaim
\subheading {3.2.7}We return to the general case when $x\in \Rm,\ V_t$ is
a one parametric group, $\rho(r)$ is a proximate order and
$u\in SH(\rho(r)).$  The following theorem represents indicators in a form
of limits in usual topology.
\proclaim {Theorem 3.2.7.1(Classic Indicators)}One has
$$h(x,u)=\underset T \to \sup
[\limsup\limits_{t_j\rightarrow\infty} u_{t_j}]^*(x)=
[\limsup_{t\rightarrow\infty}u_t(x)]^* \tag 3.2.7.1 $$
where * can be deleted outside a set of  zero capacity, and
$$\underline h (x,u)=\underset T \to \inf
[\limsup\limits_{t_j\rightarrow\infty} u_{t_j}]^*(x) ,\tag 3.2.7.2$$
where $T$ is the set of all the sequences that tend to infinity.
\endproclaim
\demo {Proof} Let us prove (3.2.7.1). Set
$$ h_1 (x,u,\{t_j\}):= \limsup_{t_j\rightarrow\infty}u_{t_j}(x).\tag 3.2.7.3$$
Let $v\in \Fr[u]$ and $u_{t_j}\rightarrow v$ in $\Di'$. Then
$$h_1^* (x,u,\{t_j\})=v(x)\tag 3.2.7.4$$
by Theorem 2.7.3.3.(H.Cartan+). Thus
$$\sup\limits_T h_1^* (x,u,\{t_j\})\geq h(x,u).\tag 3.2.7.5$$
Let $\epsilon >0$ be arbitrary small, and $t_j:=t_j(x)$ be a sequence such that
$$h_1^* (x,u,\{t_j\})\geq \sup\limits_T h_1^* (x,u,\{t_j\})-\epsilon$$
We can find a subsequence $\{t_j\}$ (we keep the same notation for it) and
$v\in \Fr[u]$ such that $u_{t_j}\rightarrow v$ in $\Di'.$ From (3.2.7.4)
we obtain
$$h(x,u)\geq v(x)\geq \sup\limits_T h_1^* (x,u,\{t_j\})-\epsilon.$$
Thus the reverse inequality to (3.2.7.5) holds.
Therefore
$$h(x,u)=\sup\limits_T h_1^* (x,u,\{t_j\})$$
Let us prove the second equality in (3.2.7.1).
Since
$$\sup\limits_T h_1 (x,u,\{t_j\})=\limsup_{t\rightarrow\infty}u_{t}(x)$$
we have
$$h(x,u)\geq[\limsup_{t\rightarrow\infty}u_{t}]^*(x)\tag 3.2.7.6$$
Let us prove the opposite inequality. Let $v\in \Fr [u].$ There exists
a sequence $t_j\rightarrow\infty$ such that
$u_{t_j}\rightarrow v $ in $\Di'(\Rm).$ By (3.2.7.4)
$$[\limsup_{t\rightarrow\infty}u_{t}]^*\geq h_1^* (x,u,\{t_j\})=v(x)$$
Since it holds for every $v\in\Fr [u]$ we have the reverse inequality  to
(3.2.7.6).Hence, (3.2.7.1) is proved completely.

Let us prove (3.2.7.2). From (3.2.7.4) we have
$$\inf\limits_{T}h^* (x,u,\{t_j\})\leq v(x)$$
for all $v\in\Fr[u].$ Therefore
$$\inf\limits_{T}h^* (x,u,\{t_j\})\leq\underline h(x,u).\tag 3.2.7.7$$

Let us prove the opposite inequality.
Let  $\{t_j\}$
be any sequence that tends to $\infty$. Let us find a subsequence
$\{t_{j'}\}$
such that $u_{t_{j'}}\rightarrow v $ in $\Di'(\Rm).$ Then
$$h(x,u,\{t_j\})\geq \limsup
\limits_{j'\rightarrow\infty}u_{t_{j'}}(x).$$
Taking * from the two sides of this inequality and using Theorem 2.7.3.3,
we obtain
$$h^*(x,u,\{t_j\})\geq [\limsup
\limits_{j'\rightarrow\infty}u_{t_{j'}}]^*(x)=v(x)\geq \underline h (x,u).$$
 This implies the reverse inequality  to (3.2.7.7). Hence (3.2.7.2)
holds. \qed
\proclaim {Corollary 3.2.7.2} If all the functions (3.2.7.3) are upper
semicontinuous, then
$$h (x,u)=\limsup\limits_{t\rightarrow\infty}u_{t}(x)$$
$$\underline h (x,u)=\liminf\limits_{t\rightarrow\infty}u_{t}(x)$$
\endproclaim
\demo {Proof}
We have
$h^*(x,u,\{t_j\})=h(x,u,\{t_j\})$
and thus
$$h (x,u)=\underset T \to \sup
[\limsup\limits_{t_j\rightarrow\infty} u_{t_j}](x)=
\liminf\limits_{t\rightarrow\infty}u_{t}(x).$$
$$\underline h (x,u)=\underset T \to \inf
[\limsup\limits_{t_j\rightarrow\infty} u_{t_j}](x)=
\liminf\limits_{t\rightarrow\infty}u_{t}(x).$$
\qed
\enddemo
\proclaim {Theorem 3.2.7.3(Indicators of Harmonic Function)}Let $u\in
SH(\rho (r))$ be  harmonic for all the
large $|y|$ in a ``cone'' of the form
$$Co_\Om:=\{y=P_t x:x\in \Om,\ t\in (0;\infty)\}$$
where $\Om\subset S_1.$ Then
$$h (x,u)=\limsup\limits_{t\rightarrow\infty}u_{t}(x)\tag 3.2.7.7$$
and
$$\underline h (x,u)=\liminf\limits_{t\rightarrow\infty}u_{t}(x)\tag 3.2.7.8$$
for $x\in Co_\Om.$
\endproclaim
\demo {Proof} The harmonicity of $u$ in $Co_\Om$ implies
$[u_t]_\epsilon (x)=u_t(x)$ for large $t$ and sufficiently small $\epsilon$
when $x\in Co_\Om.$

The family $[u_t]_\epsilon$ is uniformly continuous by reg3), Theorem 2.3.4.5
(Properties of Regularizations).Thus the function (3.2.7.5) is continuous.
Therefore we can use Corollary 3.2.7.2.\qed
\enddemo
\proclaim {Theorem 3.2.7.4.(Indicator for $m=2$)}Let $u\in SH(\Bbb R^2)$.Then
$$h(z,u)=\limsup\limits_{t\rightarrow\infty}u_{t}(x),\tag 3.2.7.9$$
\endproclaim
i.e., the star in (3.2.7.1) can be deleted.
\demo {Proof} Let as denote as $h_1(z,u)$ the right part of (3.2.7.9).The
``homogeneity'' of the indicator (3.2.1.10) and also of $h_1(z,u)$ implies the
following property: if the inequality
$h_1(z,u)<h(z.u)$ holds for some $z_0,$ it holds on the whole orbit
$$z=\{P_t z_0:0<t<\infty\}$$
that has a positive capacity in $\Bbb R^2.$ This
contradicts  Theorem 3.2.7.1.\qed
\enddemo
\newpage
\centerline {\bf 3.3. Densities}
3.3.1.In the sequel $G$ is an open set,  $K$ is a compact and  $E$ -
a bounded Borel set.

Let $\mu\in \Cal M (\rho(r))$ and $\Fr [\mu]:=\Fr [\mu,\rho(r),V_t,\Bbb R^m]$
 be the limit set of $\mu$. Set
$$\overline \Delta (G,\mu):=\sup \{\nu (G):\nu\in \Fr[\mu]\};$$
$$\overline \Delta (E,\mu):=\inf \{\overline \Delta (G,\mu):G\supset E\};$$
$$\underline \Delta (K,\mu):=\inf \{\nu (K):\nu\in \Fr[\mu]\};$$
$$\underline \Delta (E,\mu):=\sup \{\underline \Delta (K):K\subset E\}.$$
The quality $\overline \Delta (E,\mu),\ (\underline \Delta (E,\mu))$ is
called the {\it upper} ( {\it lower) density} of $\mu$ relative to the proximate
order $\rho (r)$ and the family $V_t.$
\proclaim {Theorem 3.3.1.1.(Properties of Densities)}The following properties hold

dens1) if $E=\varnothing,$ then $ \overline \Delta (E,\bullet)=
\underline \Delta (E,\bullet)=0$

dens2) $\forall E,\ \underline \Delta (E,\bullet)\leq
\overline \Delta (E,\bullet);$

dens3)  monotonicity:
$\underline \Delta ,\overline \Delta (E_1,\bullet)\leq
\underline \Delta ,\overline \Delta (E_2,\bullet)$
for $E_1\subset E_2$;

dens4)  generalized semi-additivity
\footnote {see Exercise 3.3.1.1} with respect to a set:
$$\overline \Delta (E_1\cup E_2,\bullet)+
\underline \Delta (E_1\cap E_2,\bullet)
\leq
\overline \Delta(E_1,\bullet)+\overline \Delta(E_2,\bullet)\tag 3.3.1.1$$
$$\underline \Delta(E_1\cup E_2,\bullet)+\overline\Delta (E_1\cap E_2,\bullet)
\geq
\underline \Delta (E_1,\bullet)+\underline \Delta (E_2,\bullet).\tag 3.3.1.2$$

dens5)   continuity from the right and from the left.
$$E_n \uparrow E\ \Longrightarrow
\overline \Delta (E_n,\bullet)\uparrow \overline \Delta (E,\bullet);\
K_n \downarrow K\ \Longrightarrow
\overline \Delta (K_n,\bullet)\downarrow \overline \Delta (K,\bullet);
\tag 3.3.1.3$$
$$ E_n \downarrow E\ \Longrightarrow
\underline \Delta (E_n,\bullet)\downarrow\underline \Delta (E,\bullet);\
G_n \uparrow G\ \Longrightarrow
\underline \Delta (G_n,\bullet)\uparrow\underline \Delta (G,\bullet).
\tag 3.3.1.4$$

dens6)  semi-additivity and  positive homogeneity with respect to $\mu,$
i.e.,
 $$ \overline \Delta (E,\mu_1 +\mu_2)\leq
\overline \Delta (E,\mu_1 )+\overline \Delta (E,\mu_2);\tag 3.3.1.5$$
$$\underline \Delta (E,\mu_1 +\mu_2)\geq
\underline \Delta (E,\mu_1 )+\underline \Delta (E,\mu_2);\tag 3.3.1.6$$
$$\overline \Delta,\underline \Delta (E,\lambda\mu)=
\lambda\overline \Delta,\lambda\underline \Delta (E,\mu)\tag 3.3.1.7$$
for $\lambda \geq 0;$

dens7) invariance with respect to the map $(\bullet)_{[t]}$ (see, 3.1.2.4a),
i.e.,
$$t^{-\rho -m +2}\overline \Delta,\underline \Delta (P_t E,\bullet)=
\overline \Delta,\underline \Delta (E,\bullet).$$
\endproclaim

\demo {Proof of Theorem 3.3.1.1}The property dens1) holds because the empty set
is open by definition. The properties dens2) and dens3) hold because of the
monotonicity of $\nu.$

Let us prove dens4). Since $\nu $ is a measure we have
$$\nu (G_1\cup G_2,\mu)+\nu (G_1\cap G_2,\mu)=\nu (G_1,\mu)+\nu (G_2,\mu)$$
for any $G_1\supset E_1$ and  $G_2\supset E_2.$

From this we obtain
$$\nu (G_1\cup G_2,\mu)+\nu (K_1\cap K_2,\mu)\leq\nu (G_1,\mu)+\nu (G_2,\mu)
\tag 3.3.1.8$$
for $K_1\subset E_1$ and  $K_2\subset E_2.$

The right side of (3.3.1.8) is no larger than
$\overline \Delta (G_1,\bullet)+\overline \Delta (G_1,\bullet).$
Now we can take supremum over $\nu \in \Fr \mu$ in the first summand of the
left side and infimum in the second summand. Thus we obtain
$$ \overline \Delta (G_1\cup G_2,\bullet)+
\underline \Delta (K_1\cap K_2,\bullet)\leq
\overline \Delta (G_1,\bullet)+\overline \Delta (G_1,\bullet)
\tag 3.3.1.9$$
Since $\overline \Delta (E,\bullet)$ and $\underline \Delta (E,\bullet)$ are
monotonic with respect to $E,$
$$\inf \{\overline \Delta (G_1\cup G_2,\bullet):
G_1\supset E_1,\ G_2\supset E_2\}=\overline \Delta (E_1\cup E_2,\bullet)$$
and
$$\sup\{\underline \Delta (K_1\cap K_2,\bullet):
K_1\subset E_1,\ K_2\subset E_2\}.$$
Thus we obtain the first inequality in dens4) from (3.3.1.9).
The second one can be obtained analogously\footnote {Exercise 3.3.1.2}.

Let us prove dens5). For arbitrary $G\supset K$ there exists $n_0$ such that
$K_n\subset G$ for $n>n_0$. According to dens3)
$$\vdelt(K,\bullet)\leq\vdelt(K_n,\bullet)\leq\vdelt(G,\bullet).$$
Hence,
$$\vdelt(K,\bullet)\leq\lim\limits_{n\rightarrow \infty}\vdelt(K_n,\bullet)
\leq \vdelt(G,\bullet).$$
Taking infimum over all $G\supset K,$ we obtain the second assertion in
(3.3.1.3).

For $G_n\uparrow G$ we have the equality
$$\lim\limits_{n\rightarrow \infty}\vdelt(G_n,\bullet)=
\sup\limits_{n}\vdelt(G_n,\bullet)=\vdelt(G,\bullet)\tag 3.3.1.10$$
because one can change the order of taking the supremum on $n$ and on
$\nu\in\Fr [\mu].$

Let $E_n\uparrow E$ and let $\epsilon$ be arbitrarily small. One can find
$G_n\supset E_n$ such that
$$\vdelt(G_n,\bullet)\leq\vdelt(E_n,\bullet)+\epsilon.$$
Since
$G:=\bigcup_1^\infty G_n \supset E $ we have
$$\vdelt(G_n,\bullet)-\epsilon\leq\vdelt(E_n,\bullet)\leq\vdelt(E,\bullet)
\leq\vdelt(G,\bullet).$$
Using (3.3.1.10), we obtain
$$\vdelt(E,\bullet)-\lim\limits_{n\rightarrow \infty}\vdelt(E_n,\bullet)\leq
\epsilon.$$
Since $\epsilon$ is arbitrary small
$$\vdelt(E,\bullet)\leq\lim\limits_{n\rightarrow \infty}\vdelt(E_n,\bullet)$$
and hence the first assertion in (3.3.1.8) holds.

The assertion (3.3.1.9) can be proved analogously.
\footnote {see Exercise 3.3.1.3}

Let us prove dens6).One has
$$\vdelt(G,\mu_1+\mu_2)=\sup\{\nu (G):\nu\in\Fr[\mu_1+\mu_2]\}.$$
Since
$$\Fr[\mu_1+\mu_2]\subset\Fr[\mu_1]+\Fr[\mu_2]$$
(see frmu1), Theorem 3.1.3.4 (Properties of $\mu\mapsto \Fr[\mu])$) one can
 continue the previous inequality as
$$\leq \sup\{\nu (G):\nu\in\Fr[\mu_1]+\Fr[\mu_2]\}=$$
$$=\sup\{\nu (G):\nu\in\Fr[\mu_1]\}+\sup\{\nu (G):\nu\in\Fr[\mu_2]\}=
\vdelt(G,\mu_1)+\vdelt(G,\mu_2).$$
Passing to the infimum over $G\supset E,$ we obtain (3.3.1.5). The
assertions (3.3.1.6) and (3.3.1.7) can be proved analogously.\footnote {see Exercise 3.3.1.4}

The properties dens7) follow from the invariance of $\Fr[\mu]$ (see frm3),
Theorem 3.1.3.3. (Properties of $Fr[\mu]$)). \qed
 \enddemo
{\bf Exercise 3.3.1.1.} Prove  the  {\it subadditivity } of $\overline \Delta (E,\bullet):$
$$\overline \Delta (E_1\cup E_2,\bullet)
\leq
\overline \Delta(E_1,\bullet)+\overline \Delta(E_2,\bullet)$$
and the {\it superadditivity} of $\underline \Delta(E,\bullet):$
$$\underline \Delta(E_1\cup E_2,\bullet)
\geq
\underline \Delta (E_1,\bullet)+\underline \Delta (E_2,\bullet).$$
from the Theorem 3.3.1.1.

{\bf Exercise 3.3.1.2.} Prove (3.3.1.2).

{\bf Exercise 3.3.1.3.} Prove (3.3.1.4).

{\bf Exercise 3.3.1.4.} Prove (3.3.1.6) and (3.3.1.7).

Set for $I\subset (0,\infty)$ and $\Omega \subset S_1$
$$Co_{\Om}(I):= \{x=P_t y:y\in \Omega,\ t\in I\}.$$
Also set $I_t:=(0,t).$
\proclaim {Theorem 3.3.1.2.(Cone's Densities)}One has
$$\vdelt,\ndelt (Co_{\Om}(I_t))=t^{\rho +m-2}\vdelt,\ndelt (Co_{\Om}(I_1)).$$
\endproclaim
We obtain this from dens7), Theorem 3.3.1.1, taking $E:=Co_{\Om}(I_1).$
\enddemo
\subheading {Exercise 3.3.1.5} Show that for $m=2, S_1=\{|z|=1\},
\Om =\{z=e^{i\phi}:\phi \in (\a,\be)\}$ $Co_\Om (I_t)$ is a sector of radius $t$ corresponding to the arc $(\a,\be)$ on the unit circle.

\subheading {3.3.2} Let $\delta (E)$ be a {\it monotonic} function of
$E\in\Rm.$ A set $E$ is called $\delta$-squarable if
$$\sup\limits_{K\subset E}\delta (K)=\inf \limits_{G\supset E}\delta (G).
\tag 3.3.2.1$$
\subheading {Example 3.3.2.1}Let $\delta (E)$ be a measure. Then (3.3.2.1)
implies  $\delta (\p E)=0,$ i.e., $E$ is $\delta$-squarable in sense of
item 2.2.3.
\subheading {Exercise 3.3.2.2} Prove.
\proclaim {Theorem 3.3.2.1} If  $\vdelt (\p E)=0$ then $E$ is
$\vdelt$-squarable. If $E$ is $\ndelt$-squarable, then $\ndelt (\p E)=0.$
\endproclaim
Set
 $$E_t:=\{x:\exists y\in E: |x-y|<t\}.$$
This is a $t$ -extension of $E.$

A family of sets $\Cal A_1$ is said to be {\it dense in} a family $\Cal A_2$
if for each set $E_2\in \Cal A_2$ and an arbitrary small $\epsilon >0$
there exists a set  $E_1\in \Cal A_1$ such that
$$\overline {E_1\Delta E_2}:=
\overline {(E_1\setminus E_2)\cup(E_2\setminus E_1)}\subset
(\p E_2)_\eps .\tag 3.3.2.2$$

\subheading {Exercise 3.3.2.3} Prove
\proclaim {Theorem 3.3.2.2}The relation ``to be dense in'' is reflexive
and transitive,
\endproclaim
i.e., $\Cal A_1$ is dense in $\Cal A_1,$ and
$$\{\Cal A_1\text { is dense in }\Cal A_2\}\wedge
\{\Cal A_2\text { is dense in }\Cal A_3\}\Longrightarrow
\{\Cal A_1\text { is dense in }\Cal A_3\}.\tag 3.3.2.3$$

There are  lots of squarable sets.
\proclaim {Theorem 3.3.2.3} For any monotonic $\delta (E)$ the class of
$\delta$-squarable sets is dense in the class of all the subsets of $\Rm.$
\endproclaim
\demo {Proof}For any $E\subset\Rm$ set
$$E(t):=E\cup(\p E)_t. \tag 3.3.2.4$$
One can check that
$$\overline {E\Delta E(t_1)}\subset (\p E)_{t_2}\tag 3.3.2.5$$
and
$$\overline {E(t_1)}\subset \overset {\circ} \to {E} (t_2)\tag 3.3.2.6$$
for $t_1<t_2.$

The function $f(t):=\delta (\overset \circ \to E(t))$ is monotonic. Hence,
its set of continuity points has a concentration point at $t=0$.

Suppose $\eps >0$ is arbitrarily small, and $t_0<\eps $ is a continuity point
 for $f(t).$ From (3.3.2.6) we have
$$\lim\limits_{t\rightarrow t_0 -\eps}\delta (\overline E (t))\leq
\sup\limits_{K\subset E_{t_0}}\delta (K)\leq
\inf\limits_{G\supset E_{t_0}}\delta (G)\leq
\lim\limits_{t\rightarrow t_0 +\eps}\delta ( \overset \circ \to E (t))$$
Hence, $E_{t_0}$ is $\delta$-squarable. From (3.3.2.5) we have
$$\overline {E\Delta E(t_0)}\subset (\p E)_{\eps}.$$

\qed
\enddemo
Set
$$\vdelt ^{cl}(E)=\limsup\limits_{t\rightarrow\infty}\mu_t(E);\
\ndelt ^{cl}(E)=\liminf\limits_{t\rightarrow\infty}\mu_t(E).$$
These are classic densities determined without $\Di'$-topology.
They are monotonic. The following assertion connects these densities
to $\vdelt$ and $\ndelt.$
\proclaim {Theorem 3.3.2.4.(Classic Densities)} For any
$\vdelt ^{cl}$-squarable set $E$
$$\vdelt ^{cl}(E)=\sup\{\nu (E):\nu\in \Fr [\mu]\}=\vdelt (E,\mu).\tag 3.3.2.7$$
 For any
$\ndelt ^{cl}$-squarable set $E$
$$\ndelt ^{cl}(E)=\inf\{\nu (E):\nu\in \Fr [\mu]\}=\ndelt (E,\mu).
\tag 3.3.2.7'$$
\endproclaim
The theorem follows obviously from the following assertion
\proclaim {Theorem 3.3.2.5} One has
$$\sup\limits_{K\subset E}\vdelt ^{cl}(K)\leq
\sup\limits_{\nu\in \Fr}\nu(E)\leq
\vdelt (E)\leq \inf\limits_{G\supset E}\vdelt ^{cl}(G);\tag 3.3.2.8$$
$$\sup\limits_{K\subset E}\ndelt ^{cl}(K)\leq
\inf\limits_{\nu\in \Fr}\nu(E)\leq
\ndelt (E)\leq \inf\limits_{G\supset E}\ndelt ^{cl}(G);\tag 3.3.2.9$$
\endproclaim
\demo {Proof}Let us prove, for example,  (3.3.2.9). Let us choose any $G$ and $K$ such that
$K\subset E\subset G.$ We can find a sequence $t_j\rightarrow\infty$ such that
$$\lim\limits_{j\rightarrow\infty}\mu_{t_j}(G)=\ndelt ^{cl}(G).$$
Choose a subsequence $t_{j_n}$ such that $\mu_{t_{j_n}}\rightarrow \nu$
in $\Di'$ for some $\nu\in \Fr.$

Using Theorems 2.3.4.4.(D'and C*) and 2.2.3.1.(C*-limits), we obtain
$$\nu(G)\leq \liminf\limits_{n\rightarrow\infty}\mu_{t_{j_n}}(G)=
\ndelt ^{cl}(G).\tag 3.3.2.10$$
By the same theorems
$$\ndelt ^{cl}(K)\leq \limsup\limits_{n\rightarrow\infty}\mu_{t_{j_n}}(K)\leq
\nu(K).\tag 3.3.2.11$$
From (3.3.2.10) and (3.3.2.11) we obtain
$$\ndelt ^{cl}(K)\leq\nu(E)\leq\nu(G)\leq \ndelt ^{cl}(G)\tag 3.3.2.12$$
because of monotonicity of $\nu (E).$
Taking supremum over all  $K\subset E$ and infimum over all $G\supset E,$
we obtain (3.3.2.9).\qed
\enddemo
 {\bf Exercise 3.3.2.4.} Prove (3.3.2.8).

\proclaim {Corollary 3.3.2.6}The following holds
$$\vdelt ^{cl}(K_t)=\vdelt (K_t,\mu)=t^{\r+m-2}\vdelt (K_1,\mu), t\geq 0.\tag 3.3.2.13$$

 $$\ndelt ^{cl}(K_t)=\ndelt (K_t,\mu)=t^{\r+m-2}\ndelt (K_1,\mu), t\geq 0.
\tag 3.3.2.13'$$
where $K_t= \{x:|x|<t\}$ is the ball.\ep
\demo {Proof}The right equalities follow from Theorem 3.3.1.2 with
$\Omega :=S_1.$ The left equalities hold at least for one $t$ because of
Theorem 3.3.2.4 and hence for all $t.$
\qed\edm
\subheading {3.3.3}Let us note generally speaking that  values of $\vdelt$ and $\ndelt$ on the
sets $Co_\Om (I_t)$  do not determine their values even
 on  the sets  $Co_\Om (I)$ for $I=(t_1,t_2).$ However the following assertion
holds.
\proclaim {Theorem 3.3.3.1.(Existence of Density)}Let $\Phi$ be a dense ring
(see, 2.2.3) on $S_1$. Then the conditions
$$\vdelt (Co_\Om (I_t))=\ndelt (Co_\Om (I_t))\tag 3.3.3.1$$
for $\Om\in \Phi$ and some $t$ determine uniquely a measure $\Delta (\Om)$
on $S_1.$  $\Fr[\mu]$ consists of one single measure $\nu$ and
$$\nu(Co_\Om (I_t))=t^{\rho+m-2}\Delta (\Om)\tag 3.3.3.2$$
for all the $t\in (0,\infty).$
\endproclaim
To prove this we need an assertion that is valuable by itself.
Set
$$\vdelt (\Om):=\vdelt (Co_\Om (I_1));
\ndelt (\Om):=\ndelt (Co_\Om (I_1))\text {for} \Om\in S_1.\tag 3.3.3.3$$
We will call them {\it angular } densities because for $m=2$ and $V_t\equiv I$
$\Om$ determines an angle in the plane.

Let $\Om^G$ denote an open set in $S_1$ and $\Om^K$ a closed one.
\proclaim {Theorem 3.3.3.2.(Angular Densities)}One has
$$\vdelt (\Om)=\sup\limits_{\Om^G\supset\Om}\vdelt (\Om^G);\
\ndelt (\Om)=\inf\limits_{\Om^K\subset\Om}\ndelt (\Om^K).\tag 3.3.3.4$$
\endproclaim
\demo {Proof}We need to prove two assertions
$$\forall \eps>0\  \exists \Om^G:\vdelt (\Om^G)<\vdelt (\Om)+\eps;\tag 3.3.3.5$$
$$\forall \eps>0 \ \exists \Om^K:\ndelt (\Om^K)>\vdelt (\Om)-\eps;\tag 3.3.3.6$$
Let us  prove (3.3.3.5). Set
$$\Om^G(\eps):=Co_\Om (I_{1+\eps})\cup \{|x|<\eps\}.$$

This is an open set that contains $Co_\Om (I_1).$ One can show the following:

{\bf Exercise 3.3.3.1.} For every
open set $G\supset Co_\Om (I_1)$ there exists $\eps>0$ and $\Om^G\subset S_1$
such that $\Om^G(\eps)\subset G.$

We will show
$$\vdelt (\Om^G(\eps))<\vdelt (\Om^G)+o(1) \tag 3.3.3.7$$
uniformly with respect to $\Om^G\subset S_1$ while $\eps\rightarrow 0.$

We have from Exercise 3.3.1.1
 $$\vdelt (\Om^G(\eps))\leq \vdelt(Co_\Om (I_{1+\eps}))+\vdelt (\{|x|<\eps\}).
\tag 3.3.3.8$$
The property dens7), Theorem 3.3.1.1, gives
$$\vdelt(Co_{\Om^G} (I_{1+\eps}))=\vdelt(Co_{\Om^G} (I_{1}))(1+\eps)^{\rho+m-2}.$$
Since
$ \vdelt(Co_{\Om^G} (I_{1}))\leq \vdelt (\{|x|<1\})$
we have
$$\vdelt(Co_{\Om^G} (I_{1+\eps}))=\vdelt(Co_{\Om^G} (I_{1}))+o(1) \tag 3.3.3.9$$
uniformly with respect to $\Om^G\subset S_1$ as $\eps\rightarrow 0.$

By dens7) we  also have
$$\vdelt (\{|x|<\eps\})=\vdelt (\{|x|<1\})\eps^{\rho+m-2}=o(1).\tag 3.3.3.10$$
From (3.3.3.10),(3.3.3.9) and (3.3.3.8) we obtain (3.3.3.7).
Hence (3.3.3.5) is proved.

Let us prove (3.3.3.6). Set
$$\Om^K(\eps):=Co_{\Om^K}(\bar I_{1-\eps})\setminus \{|x|<\eps\}$$
where $\bar I$ is the closure of $I.$

One can show the following:

{\bf Exercise 3.3.3.2.} For any compact $K\subset Co_\Om (I_1)$ there
exist $\Om^K\subset \Om$ and $\eps>0$ such that
$K\subset\Om^K (\eps)\subset Co_\Om (I_1).$

From the definition of $\vdelt (\Om)$ and the monotonicity we obtain (3.3.3.6).
\qed
\enddemo
\demo {Proof of Theorem 3.3.3.1}
Suppose (3.3.3.1) holds .The property dens7), Theorem 3.3.1.1,
implies (3.3.3.1) for all the $t\in (0,\infty).$

Set
$\Delta (\Om):=\vdelt (\Om)=\ndelt (\Om )$
for $\Om\in \Phi.$

Let us prove that $\Delta$ satisfies the conditions $\Delta1)-\Delta3)$
from 2.2.3. The conditions $\Delta1)$ and $\Delta2)$ follow from
dens3) and dens4), Theorem 3.3.1.1, Exercise 3.3.1.1.

Let us prove  $\Delta3).$
By Th.3.3.3.2 for arbitrary $\Om\in \Phi$ and $\eps>0$ we can choose  $\ \Om^G\supset \Om\ $ such that
$\ \vdelt (\Om)>\vdelt (\Om^G)-\eps\  $ and
$\ \Om^K\subset \Om\ $ such that $\ndelt (\Om)<\ndelt (\Om^K)+\eps .$

Suppose $\Om'\in \Phi$ satisfies the condition
$\Om^K\subset \Om'\subset\Om^G.$ Then
$$\Delta (\Om') =\vdelt (\Om')\leq \vdelt (\Om^G)\leq \vdelt (\Om)+\eps
=\Delta (\Om)+\eps$$
and
$$\Delta (\Om)-\eps=\ndelt(\Om)-\eps\leq \ndelt (\Om^K)\leq\ndelt(\Om')=
\Delta (\Om'),$$
 implying $\Delta3).$
\qed
\enddemo
\newpage
\centerline {\bf 4.Structure of Limit Sets }
\centerline {\bf 4.1. Dynamical Systems }
\subheading {4.1.1} The most complete and effective description of an arbitrary
limit set can be done in terms of dynamical systems (see, \cite {An}).

A family of the form
$$T^t:M\mapsto M,\ t\in \Bbb R,$$
on a compact metric space $(M,d)$ with a metric $d(\bullet,\bullet)$ is
{\it a dynamical system} $(T^\bul,M)$ if it satisfies
the condition
$$T^{t+\tau}=T^t\circ T^\tau,\ t,\tau \in \Bbb R$$
and the map $(t,m)\mapsto T^t m$
 is continuous with respect to $(t,m),$ for all $ t\in \Bbb R,m\in M.$****

 Let $m,m'\in M,$ and $\eps,s>0.$ An {\it
$(\eps,s)$-chain from $m$ to $m'$} is a finite sequence
$m_0=m,m_1,...,m_n=m',$
 satisfying the conditions $d(T^{t_j} m_j,m_{j+1})<\eps,\ j=0,1,...,n-1,$ for
some $t_j\geq s.$

A dynamical system $(T^t,M)$ is called {\it chain recurrent} (see, \cite {HS}),  an arbitrarily small $\eps >0$ and an arbitrarily large $s>0$ there
exists  an $(\eps,s)$-chain in $M$ from $m$ to $m .$

\proclaim {Theorem 4.1.1.1 (Properties of Chain Recurrence)}Let $(T^t,M)$ be
a dynamical system on a compact set.Then the following conditions are equivalent:

cr1) $M$ is connected and $(T^t, M)$ is chain recurrent;

cr2) for every open proper $U\subset M$ satisfying
$$T^t U\subset U,\ -\infty<t<0,\tag 4.1.1.1$$
the boundary $\p U$ contains a non-empty $T$-invariant subset of $M;$

cr3) for every closed proper  $K\subset M$ satisfying
$$T^t K\subset K,\ t\geq 0,\tag 4.1.1.2$$
the boundary $\p K$ contains a non-empty $T$-invariant subset of $M;$

cr4) there does not exist any open proper $V\subset M$ satisfying
$T^\tau \operatorname{clos}V\subset V$ for some $\tau>0;$

cr5) for any small $\eps>0$, large $s>0$, and every  pair of points $m,m'$
there exists an $(\eps,s)$-chain from $m$ to $m'.$
\endproclaim
\demo {Proof}The conditions cr2) and cr3) are equivalent. Let us prove, for
example, cr2)$\Longrightarrow $ cr3). Set $U:=M\setminus K.$ It is open.
Applying
to (4.1.1.2)
$T^{-t}$ and, using the invariance of $M,$ we obtain
(4.1.1.1) for $U$. Hence $\p U$ contains a non-empty invariant subset of $M.$ Since
$\p K=\p U$ we obtain cr2).

Let us prove the implication cr1)$\Longrightarrow $cr3). Let $K\subset M$ be
 closed, proper and satisfy (4.1.1.2). Since $M$ is proper $\p K$ is
non-empty.

Let $W$ denote the interior of $K$ in $M$.The continuity of $T$ and (4.1.1.2)
imply
$$T^tW\subset W \tag 4.1.1.3$$
 for $t\geq 0.$ Indeed, $T^t W\subset K.$ It must be
open.Thus it can not contain any point of $\p K$, since else it would contain some  neighborhood of this point, contradicting the definition of
$\p K.$

Suppose that $\p K$ does not contain any non-empty $T$-invariant set.

Let us show that there exists  $s>0$ such that
$$T^s K\subset W.\tag 4.1.1.4$$
For any $m\in \p K$ there exists $t=t(m)$ such that $T^t m\in W.$ There
exists a neighborhood $V_m$ of $m$ in $\p K$ that passes to $W$ under
$T^{t(m)}$-action
because of continuity of $T^tm$ on $m$.

We also have $T^t V_m \subset W$ for $t>t(m)$ because of (4.1.1.3).
   Since $\p K$ is compact we can
cover it by a finite number of neighborhoods and obtain $s$ such that
$$T^s \p K\subset W. \tag 4.1.1.5$$
(4.1.1.5) and (4.1.1.3) give (4.1.1.4).

Set $\eps:=0.5 d(\p K, T^s K).$ From (4.1.1.2) we see that
$T^t K\subset T^s K$ for $t>s.$
 Therefore there does not exist any $(\eps,s)$-chain
from a small neighborhood of a point $m\in \p K$ to itself. This contradicts
the chain recurrence of $M.$

Let us prove cr3)$\Longrightarrow$cr4).

Assume that there exists an open
proper
$V\subset M$ satisfying
$T^\tau \text {clos}V\subset V$ for some
$\tau>0.$

We will construct  $K$  that does not satisfy cr3).
Set
$W:=\bigcup\limits_{0\leq t\leq \tau} T^t V$ and
$K:=\text {clos} W.$

Then $$T^s W\subset W, \forall s\geq 0.\tag 4.1.1.6$$
Indeed, let $s=k\tau +s',\ s'\in [0,\tau),\ k\in \Bbb Z.$ Then
$$T^s W =\bigcup\limits_{t\in [0,\tau]}T^{t+s}V .\tag 4.1.1.7$$
Since $T^\tau V\subset V$ we have
$T^{t+k\tau}V\subset T^{t} V  $ for $t>0.$
From (4.1.1.7) we obtain
$$T^s W=\bigcup\limits_{t\in [0,\tau]}T^{t+s'+k\tau}V\subset
\bigcup\limits_{t\in [0,\tau]}T^{t+s'}V =
\bigcup\limits_{t'\in [s',\tau+s']}T^{t'}V=$$
$$ \bigcup\limits_{t\in [s',\tau]}T^{t}V\  \cup
\bigcup\limits_{t\in [\tau,\tau+s']}T^{t}V:=W_1\cup W_2.$$
 Further we have
$W_1\subset W$ by definition.$W_2$ can be represented in the form
  $$W_2= \bigcup\limits_{t\in [0,s']}T^{t+\tau}V.$$
Since
$$T^{t+\tau}V=T^{t}T^\tau V$$
and
$$T^\tau V\subset V$$
by the assumption we get:
$$W_2\subset W_1\subset W.$$

 This implies (4.1.1.6).

 The same holds for $K$ because
of continuity of $T^t,$ i.e. $K$ satisfies (4.1.1.2).

Let us prove the equality
$$K= \bigcup\limits_{0\leq t\leq \tau} T^t \text {clos}V.\tag 4.1.1.8$$
 Denote as $K'$ the right side of (4.1.1.8).

The set $K'$ is closed
because of compactness of $[0,\tau]$ . Indeed,
let the sequence $\{T^{t_j}v_j:j=1,2,...\} \in T^{t_j} (\clos V)$ converge to $w.$ Choose a subsequence
$t_{j_k}\rightarrow s\in [0,\tau].$ Then
$$v:=\lim \limits_{k\rightarrow \infty}v_{j_k}=
\lim \limits_{k\rightarrow \infty}T^{-t_{j_k}}w= T^{-s}w.$$
Since $\clos V$ is closed, $v\in \clos V.$ Thus $w=T^s v$ for some
$s\in [0,\tau]$ and some $v\in \clos V,$ i.e. $w\in K'.$

Now, $W\subset K'$ because
$$T^t V\subset \clos T^t V =T^t\clos V.$$
Hence,
$$K:=\clos W\subset \clos K'=K'.$$
We also have
$$(T^t V\subset W
 \ \forall t\in [0,\tau])\Longrightarrow  (\clos T^t V=T^t \clos V\subset \clos W =K,\ \forall t\in [0,\tau])  .$$
Hence, $K'\subset K.$ Therefore $K=K'$ , i.e. (4.1.1.8) holds.

From (4.1.1.8) and $T^{\tau}\clos V\subset V$ we obtain
$T^\tau \clos W\subset W.$ Hence $T^\tau\p K \subset W.$  This and
$\p K\cap W= \varnothing $ imply
$$T^\tau\p K\cap \p K=\varnothing.\tag 4.1.1.9$$

To obtain a contradiction and complete the proof of
cr3)$\Longrightarrow$ cr4)
 we have to show that $K$ is a proper subset, because  both  cases:
$\p K=\varnothing$ and $ \p K\neq\varnothing$ will contradict
cr3).

Since $V$ is proper $T^t V$ is proper for any $t\in (-\infty,\infty).$
Otherwise $T^t V=M$ implies $V=T^{-t}M=M$ that is a contradiction.

Since $V$ is a neighborhood of the compact set $T^\tau \clos V$ we can find
$\alpha>0$ such that
$T^t\circ T^\tau \clos V\subset V$ for $t\in [0,\alpha].$ Then
$T^t\clos V\subset T^{-\tau}V$ for $t\in [0,\a].$

By iteration of this inclusion we obtain $T^{jt}\clos V\subset T^{-j\tau}V$
for any integer $j.$
When $j\alpha>\tau$ it follows that $K\subset T^{-j\tau}V.$
The last set is proper because we mentioned already that
 $T^t V$ is proper for any $t\in (-\infty,\infty).$
Hence $K$ is proper.

So $K$ satisfies the conditions of cr3) but $\p K$ does not contain a non-empty
$T$-invariant set. This contradiction proves the implication
cr3)$\Longrightarrow$ cr4).

Let us prove cr4)$\Longrightarrow$cr5). Let $\eps>0$ be small
 and $s>0$ be
large. Let $V$ denote the set of all
$m'\in M$ such that there exists an
$(\eps,s)$-chain from $m$ to $m'.$ This
set is open and closed. Indeed,
 let $m'\in V.$ There
exists an $(\eps,s)$-chain $m=m_0,...,m_{n-1},m_n=m'$ from $m$ to $m'.$
Choose $\eps_1<\eps -d(m_n,m_{n-1})$ and consider the closed neighborhood
$W:=\{m'':d(m',m'')\leq \eps_1\}.$ Then for any $m''\in W$
 the chain $m=m_0,...,m_{n-1},m_n=m''$ is an $(\eps,s)$-chain from
$m$ to $m''.$ Hence, with every point $V$ contains its closed neighborhood.
Therefore it is open and closed.Therefore it is a connected component of $M.$

We also have $T^s m\in V$ because for that case $n=1, m_0=m,m_1=T^s m.$
Hence $T^s \clos V\subset V.$
If $V$ does not coincide with the whole $M$ the latter contradicts to cr4).
Hence $V=M.$

Finally, let us prove cr5)$\Longrightarrow$ cr1). If $M$ is a
 union of two non-empty disjoint sets $A$ and $B,$ then both of them are
open end
 closed. Since $M$ is compact, the distance $\eps$ between $A$ and $B$ is
positive . Hence every $(\eps/2,s)$-chain starting at a point of $A$ remains
in $A,$ contradicting cr5).

Since for every point $m\in M$ the set $V$ from the proof of
cr4)$\Longrightarrow$ cr5) coincides  with $M$, cr1) holds.
\qed
\enddemo
\proclaim {Theorem 4.1.1.2} Let $T^t$ be chain recurrent on $M_\a,\a\in A.$ Then $T^t$ is chain recurrent on $M=\bigcup_{\a\in A}M_\a.
$\ep
This is because every ($\eps, s$)- chain from $m$ to $m'$ in $M_\a$ is also
($\eps,s$)-chain in $M.$

\subheading {4.1.2}Here we prove two auxiliary assertions that will be used
further.
\proclaim {Theorem 4.1.2.1}Let $T$ be  chain recurrent on a connected compact
$M$ and let $\{q_j\}$ be a sequence in $M.$
Then there exist sequences $\{\alpha_{\nu}\}$ and $\{\om_{\nu}\}$ of real
numbers and a sequence $\{p_{\nu}\}$ in $M$ having $\{q_j\}$ as a subsequence,
such that
$$ \alpha_{\nu}\rightarrow -\infty;\ \om_{\nu}\rightarrow \infty\tag 4.1.2.1$$
and
$$d(T^{\om_{\nu}}p_\nu,T^{\alpha_{\nu+1}}p_{\nu+1})\rightarrow 0\tag 4.1.2.2$$
as $\nu\rightarrow\infty.$
\endproclaim
\demo {Proof}In addition to $\{\alpha_{\nu}\}$, $\{\om_{\nu}\}$ and
$\{p_{\nu}\}$ we define, by induction, a sequence $\{\eps_\nu\}$ of positive
real numbers, tending to zero, and an increasing sequence $\{\nu_j\}$ of
positive integers, such that $\{p_{\nu_j}\}=\{q_j\}$ and
$$d(T^{\om_{\nu}}p_\nu,T^{\alpha_{\nu+1}}p_{\nu+1})<\eps_\nu,
\nu=1,2,... .\tag 4.1.2.3$$

We start by setting $\alpha_1=-1,\eps_1=1, \nu_1=1, \om_1=5$ and $p_1=q_1.$
Assume now that $\alpha_\nu,\eps_\nu, \om_\nu$ and $p_\nu$ have been chosen
for $\nu=1,2,...,\nu_j.$ Set
$$\alpha =\alpha_{\nu_j}-1,\  \eps =\eps_{\nu_j}/2,\
 \om=\om_{\nu_j}.\tag 4.1.2.4$$
 By Theorem 4.1.1.1, cr5) there exists a sequence
$r_0:=T^\om q_j,r_1,...,r_m:=T^\alpha q_{j+1}$ such that
$d(T^{t_k}r_k,r_{k+1})<\eps$ for $k=0,1,...,m-1$, where $t_k\geq \om$. Now we
set $\nu_{j+1}=\nu_j+m+1.$ For $\nu=\nu_j+k+1,\ k=0,1,...,m-1,$ we set
$\alpha_\nu=-t_k/2,\om_\nu=t_k/2,p_\nu=T^{t_k/2}r_k,$ and finally, for
$\nu=\nu_{j+1}$ we set $\alpha_\nu =\alpha, \eps_\nu=\eps, \om_\nu= \om+1,
p_\nu=q_{j+1}.$

 Let us check that
with this setting the properties (4.1.2.1) hold .
Since $\om_{\nu_{j+1}}=\om_{\nu_{j}}+1$ we have
$\om_{\nu_{j}}\rightarrow\infty$ as $j\rightarrow\infty.$
From $t_k\geq \om=\om_{\nu_j}$ we obtain $\alpha_\nu\rightarrow -\infty$ and
$\om_\nu \rightarrow \infty.$ Hence (4.1.2.1) holds.

One can see from (4.1.2.4) that $\eps_\nu=\eps_{\nu_j}/2\rightarrow 0.$
To prove (4.1.2.2) it is enough to check (4.1.2.3).
For $k=0$ we have
$p_\nu=T^{t_0/2}r_0=T^{t_0/2+\om}q_j=T^{t_0/2+\om}p_{\nu_j}.$
Hence, $T^{\alpha_\nu} p_\nu=T^\om p_{\nu_j}=T^{\om_{\nu_j}} p_{\nu_j}.$
Thus
$$d(T^{\om_{\nu_j}} p_{\nu_j},T^{\alpha_\nu} p_\nu)=0\tag 4.1.2.5$$
for this case.

For $k=1,...,m-2$ and the corresponding $\nu$ we have
$$T^{\om_\nu} p_\nu=T^{t_k/2}\circ T^{t_k/2}r_k=T^{t_k}r_k$$
and
$$ T^{\alpha_{\nu+1}} p_{\nu+1}=T^{-t_{k+1}/2}\circ T^{t_{k+1}/2}r_{k+1}=
r_{k+1}.$$
Hence,
$$d(T^{\om_\nu} p_\nu,T^{\alpha_{\nu+1}} p_{\nu+1})=
d(T^{t_k}r_k,r_{k+1})<\eps=\eps_\nu\tag 4.1.2.6$$
Finally, for the last link of the chain we obtain
$$k=m -1,\ \nu=\nu_j+m,\ \nu+1=\nu_{j+1},\  \alpha_{\nu_{j+1}}=\alpha ,$$
$$T^{\alpha_{\nu+1}}p_{\nu+1}=T^{\alpha_{\nu_{j+1}}}p_{\nu_{j+1}}=
T^\alpha q_{j+1}=r_m$$
Thus (4.1.2.6) holds  for $k=m-1.$ Hence, (4.1.2.3) also holds.
Therefore (4.1.2.2) holds.
\qed
\enddemo

\proclaim {Lemma 4.1.2.2} Let $p_k,q_k\in M$ and $d(p_k,q_k)\ri 0$ as
$k\ri\iy.$ Then there exists a sequence $\{\g_k\uparrow \iy\}$ such that
$$d
(T^\tau p_k,T^\tau q_k)\ri 0\tag 4.1.2.7$$
uniformly with respect to $\tau\in [-\g_{k+1},\g_k].$
\ep
\demo {Proof} Let $[-\g,\g]$ be a fixed segment.Then
$d(T^\tau p_k,T^\tau q_k)\ri 0$ uniformly in this segment.

Indeed, suppose there exist sequences $\tau_j,k_j$ such that
$d(T^{\tau_j} p_{k_j},T^{\tau_j} q_{k_j})\geq \eps>0.$
Choosing a subsequence we can assume that $\tau_j\ri \tau\in [-\g,\g],\
p_{k_j}\ri p\in M$ and $q_{k_j}\ri q=p.$ Using continuity of
$T^\tau m$ on $(\tau,m)$ and continuity of $d(\bul,\bul)$ in both
arguments we obtain $0=d(p,p)\geq \eps >0.$

Denote
$$\eps(\g,k):=\max\limits_{\tau  \in [-\g,\g]}d(T^\tau p_k,T^\tau q_k).$$
This function increases monotonically in $\g$ and tends to zero for any
$\g$ as $k\ri\iy.$

Choose $l_n$ such that $\eps (n,k)\leq 1/n$ for $k\geq l_n.$ Set
$\g_{k+1}:=n $ for $l_n<k\leq l_{n+1}.$ One can see that
$\eps(\g_{k+1},k)\ri 0$ as $k\ri\iy.$ Since
$$\max\limits_{\tau  \in [-\g_{k+1},\g_k]}d(T^\tau p_k,T^\tau q_k)\leq
\eps(\g_{k+1},k),$$
 $\{\g_k\}$ satisfies (4.1.2.7).
\qed
\edm

\subheading {4.1.3} Let us consider some corollaries of the previous results.
\proclaim {Theorem 4.1.3.1} $(T^t,M)$ is chain recurrent iff $M$ is connected
and for any $m\in M,$ small $\eps>0$ and large $s>0$ there exist an
$(\eps,s)$-chain from $m$ to $m.$
\endproclaim
i.e. we can omit $V$ from the definition of the chain recurrence.
The assertion follows from cr5).
\subheading {Exercise 4.1.3.1}Prove Theorem 4.1.3.1.

 We connect the  property of being chain recurrent with other well known
 characteristics of  dynamical system .

A point $m_0\in M$ is called {\it non-wandering} (see \cite {An})
if for any neighborhood $\Cal O$ of $m_0$ and  arbitrarily large number
$s\in \Bbb R$ there exists $m\in \Cal O $ and $t\geq s$ such that
$T^t m \in  \Cal O. $

This means that the ``returns'' take place to an arbitrary small neighborhood
of the point $m_0.$ We shall denote as
 $\Om (T^\bul)$ the set of non-wandering
points. It is a closed invariant subset of $M.$

The set $A\subset M$ is called an {\it attractor} if it satisfies
the following conditions:

attr1) for any neighborhood $\Cal O\supset A$ there exists a neighborhood
$\Cal O',\ A\subset \Cal O'\subset \Cal O$ such that
$T^t \Cal O'\subset \Cal O \ t\in \Bbb R,$ where $T^t \Cal O'$ is the image
of $\Cal O';$

attr2) there exists a neighborhood $\Cal O\supset A$ such that
$T^t m\rightarrow A$ when $t\rightarrow\infty$ for $m\in \Cal O.$

\proclaim {Theorem 4.1.3.2 }If  $\Om (T^t)=M, $ then $(T^t,M)$ is
 chain recurrent; if $(T^t,M)$ has an attractor $A\neq M,$ it is not
chain recurrent.
\endproclaim

\demo {Proof}The property $\Om(T^t)=M$ obviously implies the chain
recurrence for $m=1$.Suppose there exists an attractor $A\neq M.$
Take a point $m_0$
that does not belong to $A$ and choose a neighborhood $\Cal O\supset A$ such
that $d(m_0,\clos \Cal O)=2\eps>0 .$ This is possible because an attractor is
closed.
Let $\Cal O'$ be chosen by attr1) and $s$ be such that $T^s m\in \Cal O'.$ Then
there does not exist any $(\eps,s)$-chain from a small neighborhood of $m_0$
itself. By Theorem
4.1.3.1 $(T^t,M)$ is not chain recurrent.\qed
\enddemo

Let us give  examples of dynamical systems on  connected compacts that are
chain recurrent.
\proclaim {Theorem 4.1.3.3} Let $M$ be a
connected compact and let $T^t$ be the
identity map. Then $(T^t,M)$ is chain recurrent.
\endproclaim
This theorem, of course, is trivial. However, if $M$ consists of a single
point this dynamical system determines
an important class of subharmonic and entire functions of {\it
completely regular growth} (see \cite {L(1980),Ch.III}).

Let $m\in M.$ Set
$$\Bbb C (m):=\text {clos}\{T^t m:-\infty<t<\infty\}\tag 4.1.3.1$$
It is closed, connected and invariant.
\subheading {Exercise 4.1.3.2} Prove this.

Let us denote as $ \Om (m)$ the set of all limits of the form
$$\Om (m):=\{m'\in M:(\exists t_k\rightarrow\iy)
(m'=\lim\limits_{k\rightarrow\infty} T^{t_k}m\}\tag 4.1.3.2.$$
This is a limit set as $t\rightarrow \infty.$ It is the ``tangle'' at the end
of the curve. Denote by $A(m)$ the
analogous set for $t\rightarrow -\infty.$
\subheading {Exercise 4.1.3.3} Prove that $A(m)$ and $ \Om (m)$ are invariant.
\proclaim {Theorem 4.1.3.4}$(T^t,\Bbb C (m)$) is  chain recurrent iff
 $$A (m)\cap \Om (m)\neq\varnothing.\tag 4.1.3.3$$
\endproclaim
\demo {Proof}Suppose $B:=A (m)\cap \Om (m)=\varnothing.$ Then $ \Om (m)$ is an
attractor and $(T^t,\Bbb C (m)$) is not chain recurrent by
Theorem 4.1.3.2.

Suppose $B\neq\varnothing.$ We will use cr2) from Theorem 4.1.1.1.

Let $U$ be an open proper subset of $\Bbb C (m)$ satisfying (4.1.1.1).
 Consider two  cases:

i) $B$ contains a point of $U.$ Thus $U$ contains a sequence of form
$T^{t_k}m,\ t_k\ri\iy.$ From (4.1.1.1) we obtain that $U$ contains $T^t m$ for
 all $t\in (-\iy,\iy).$ Thus $U\supset \Bbb C (m)$ and
$\clos U=  \Bbb C (m).$ Set $K= \Bbb C (m)\setminus U.$ One can
show that $K$ satisfies (4.1.1.2)(see the beginning of proof of Theorem 4.1.1.1).
Hence $K$ contains the set
$$K^*:=\bigcap\limits_{t\geq 0}T^t K\tag 4.1.3.4$$
that is invariant ( Exercise 4.1.3.4).

Therefore $K^*\subset K\subset \clos U\setminus U=\p U.$ By cr2)
$(T^t,\Bbb C (m))$ is chain recurrent.

ii) $B$ contains no point of $U.$ Then $B\subset A(m)\subset \p U.$ By cr2)
$(T^t,\Bbb C (m))$ is chain recurrent.
\qed
\enddemo
\subheading {Exercise 4.1.3.4}Let $U$ satisfy (4.1.1.1) and
$K:=M\setminus U.$ Prove that $K^*$ from (4.1.3.4) is invariant.
\subheading {4.1.4}The connectedness of $M$ is a necessary condition for
a dynamical system to be chain recurrent.

Let $M$ be a subset of a linear space. The set $M$ is called {\it
polygonally connected} if  every pair of points $m_1,m_2$ can be
connected by polygonal path.

Of course,  polygonal connectedness implies connectedness and even
arcwise connectedness. \proclaim {Theorem 4.1.4.1} Let $(T^t,M)$
be a dynamical system such that $M$ is a polygonally connected
set. Then $(T^t,M)$ is chain recurrent. \ep \demo {Proof}Let $U$
be an open proper subset of $M,$ satisfying (4.1.1.1). We choose
$m_1\in U$ and $m_2$ in an invariant subset $K^*$ of
$K:=M\setminus U.$ Then there exists a polygonal path from $m_1$
to $m_2:$
$$m_\th:= (j+1-\th)m'_{j}+(\th -j)m'_{j+1},\text { for } \th\in [j,j+1]$$
$$ j=0,1,...,l-1;\ m'_0:=m_1,\ m'_l:=m_2.$$
Now $M$ is invariant, so for each $t$ the continuous path
$\th\mapsto T^t m_\th$ lies in $M.$

If $t\in (-\iy,0)$ its initial point $T^t m_1$ belongs to $U$ and its endpoint $T^t m_2$ belongs to $K^*\subset K.$

For each $t\in (-\iy,0)$ we set
$$\th(t):=\min [\th\in [0;l]: T^t m_\th \in K].$$
Then $\th(t)>0,\ T^t m_{\th(t)}\in \p U$ and (4.1.1.1) implies that $t\mapsto
\th(t)$ is a decreasing function. Hence the limit
$$\th(-\iy):=\liml_{t\ri -\iy}\th(t)$$
exists and is positive.

Set $m_3:=m_{\th (-\iy)}.$ We claim that $A(m_3)\subset \p U$ ($A(\cdot)$ is
a set defined before Theorem 4.1.3.4).
If $\th(-\iy)\in (j,j+1]$ for some $j\in [0,l]$ then $\th (t)\in (j,j+1]$ for
 $t$ that is near to $-\iy$, and
$$T^t m_3= T^t m_{\th(t)}+(\th(t)-\th(-\iy))T^tm'_j +
(\th (-\iy)-\th(t))T^t m'_{j+1}.$$

The first term in the right hand side lies in $\p U.$ The set $M$ is compact
and invariant so the other terms tend to zero as $t\ri-\iy.$ Hence
$A(m_3)\subset \p U.$

Thus $\p U$ contains this invariant subset and $(T^t,M)$ is chain
recurrent by cr2), Theorem 4.1.1.1. \qed \edm We have the obvious
\proclaim {Corollary 4.1.4.2}Let $(T^t,M)$ be a dynamical system
such that $M$ is a convex set. Then $(T^t,M)$ is chain recurrent.
\ep This is because the polygonal path can be taken as a line
segment connecting every pair of points. \subheading {4.1.5} Let
$U(\r,\sigma)$ be a set of subharmonic functions defined in
(3.1.2.4). It is invariant with respect to the transformation
$(\bul)_{[t]}$ defined in (3.1.2.4a).

Set (subindex!)
$$T_t v:=v_{[e^t]}.\tag 4.1.5.1$$
Since $(\bul)_{[t]}$ has the property (3.1.2.4b*)
$$T_{t+\tau} v=(T_t\circ T_\tau) v,\ \forall t,\tau\in \BR.\tag 4.1.5.2$$
By Theorem 3.1.2.3  $\ T_\bullet$ is continuous in the appropriate topology
and hence
$(T_\bullet,U[\r,\sigma])$ is a dynamical system.

\proclaim {Theorem 4.1.5.1 (Universality of $U[\r,\s]$)}Let $(T^\bul,M)$ be
a chain recurrent dynamical system on a compact set $M.$ Then for any $\r,\s$ there exists
$U\sbt U[\r,\s]$ and a homeomorphism $imb:M\mapsto U$ such that
$imb\circ T^t=T_t\circ imb,\ t\in (-\iy,\iy).$
\ep
i.e., any dynamical system can be imbedded in $(T_\bul,U[\r,\s]).$

It is sufficient to prove the theorem by supposition
 $P_t x=tx$ because $(T^P_t, U[\r,\s])$ is a dynamical system for any $P_t$
and $imb :(T_t,U[\r,\s])\mapsto(T^P_t,U[\r,\s])$ where
$imb: u(x)\mapsto T_{-t}T^P_tu(x)$ is also a homeomorphism of dynamical systems.

{\bf Exercise 4.1.5.1.}Consider Theorem 3.1.6.1 from this point of view.

We need some auxiliary definitions and results. Let us denote as $\Cal M (S^{m-1})$
the set of measures $\nu$ with bounded full variation on the unit sphere
$S^{m-1}.$ Introduce the metric $d(\nu,0):=\var \nu$ and consider the set
$$K:=\{\nu:\nu>0,\ d(\nu,0)\leq 1\},$$
i.e., the intersection of the cone of positive measures with the unit ball.

The following assertion is a corollary of the Keller's theorem  (see,e.g.
\cite {BP, Th.3.1, p.100}).
\proclaim {Theorem 4.1.5.2 (Imbedding)} Every metric compact set can be
homeomorphically imbedded to $K.$
\ep
Thus we can assume below that for any $m\in M$ there exists a positive
measure
$$ Y(\bul,m)=Y(dx^0,m)\in K $$
such that
$$(Y(\bul,m_1)=Y(\bul,m_2))\Longrightarrow (m_1=m_2)\tag 4.1.5.3$$
and $Y(\bul,m)$ is continuous with respect to the metrics.

We also introduce a new coordinate system. For $x:=e^yx^0\in \Rm\setminus 0$ set
$Pol(x)=(y,x^0).$ This formula gives a one-to-one map from $\Rm\setminus 0$ onto
the cylinder $Cyl:=(-\iy,\iy)\times S^{m-1}.$ Thus, for any $(y,x^0)\in Cyl,$
$\ Pol^{-1}(y,x^0)=e^yx^0.$

For $m=2$ this is a common cylinder.

\subheading {4.1.6}\demo {Proof of Theorem 4.1.5.1} We consider separately the cases of integer
 and non-integer $\r.$

Let $\r$ be non-integer and $\s>0.$ For any $v\in U[\r,\s],$ one has
the representation of Theorem 3.1.4.4.
$$v(x)=\Pi (x,\mu,\r)\tag 4.1.6.1$$
where $\mu\in \Cal M [\r,\Dl]$ and $\Dl$ depends only on $\s$ (Theorem
2.8.3.3).

Vice versa, every $\mu\in\Cal M [\r,\Dl]$ generates $v$ by (4.1.6.1) and
$$v_{[t]}(x)=\Pi (x,\mu_{[t]},\r).$$

Let us ``transplant'' $\mu$ in Cyl. For $\mu$ that has a dense $f_\mu (rx^0),$
we set
$$\nu (dy\otimes dx^0):=f_\mu (e^yx^0)e^{(-\r-2)y}(dy\otimes dx^0).$$
i.e., the density $f_\nu$ of $\nu$ is defined by
$$f_\nu(x^0,y):=f_\mu (e^yx^0)e^{(-\r-2)y}.$$
Respectively
$$f_\mu (x^0,r)=f_\nu(x^0,\log r)r^{\r+2}$$
We can extend this equality for all $\mu\in \Cal M [\r,\Dl]$ by using
limit process in $\Di '$ topology.

{\bf Exercise 4.1.6.3.}Do that using, for example, Theorem 2.3.4.5.

We can also define $\nu$ as a distribution in $\Di'(Cyl).$ Namely, for
$\psi\in \Di(Cyl)$ we set
$$\psi^*(x^0,r):=\psi (Pol^{-1}(x^0,\log r))r^{-\r-m+2}$$
and
$$<\nu,\psi>:=\int \psi^* (x^0, r)\mu (dx^0\otimes r^{m-1}dr)$$

 {\bf Exercise 4.1.6.4.}Check that this definition gives the same $\nu.$

The transformation $P_t x=(x^0,tr),\ rx^0\in \Rm\setminus 0$ passes to
$$Pol\circ P_t\circ Pol^{-1}(x^0,y)=( x^0,y+\log t)$$
Thus $T_{e^\tau}\mu$ gives a transformation $S_\tau \nu$ defined by
$$S_t f_\nu (x^0,y):=f_\nu (x^0,y+t)$$
for densities or by
$$<S_t\nu,\psi>:=\int \psi (x^0,y-t)\nu (dx^0\otimes dy)\tag 4.1.6.2$$
for distributions ( $\psi \in \Di (Cyl).$)

{\bf Exercise 4.1.6.5.}Check the equivalence.

From $\mu\in \Cal M[\r,\Dl]$  we obtain
$$\intl_{y\leq 0}e^{(\r +m-2)y}S_t\nu (dy\otimes dx^0)\leq \Dl,\  t\in\BR,
\tag 4.1.6.3$$

{\bf Exercise 4.1.6.6.} Check this.

Let $X(t)$ be a positive function  satisfying the condition
$$\intl_{-\iy}^{\iy}X(t)dt=1$$
and such that the linear hull of its
translations are dense in $L^1(-\iy,\iy)$.
 We can chose, for example, the
function
$$X(t):=\frac {1}{\sqrt {2\pi}}e^{-\frac{t^2}{2}}$$
because its Fourier transformation does not vanish in $\BR$
(it is $e^{-\frac{ s^2}{2}}$).

{\bf Exercise 4.1.6.7.}Check these properties.

Let us define $\nu (\bul,m)$ by
$$<\nu(\bul,m) ,\psi):=\intl_{(x^0,y)\in Cyl} \psi (x^0,y)\left (\r\intl_{-\iy}^{\iy}Y(dx^0,T^{y-t}m)X(t)dt\right ) dy.\tag 4.1.6.4$$
 Now we check the property
$$S_\tau \nu(\bul,m)=\nu (\bul, T^\tau m)$$
Using (4.1.6.2), we obtain
$$<S_\tau\nu (\bul, m),\psi>=
\int \psi (x^0,y)
\left (\r\intl_{-\iy}^{\iy}Y(dx^0,T^{y+\tau-t}m)X(t)dt\right ) dy=$$
$$\int \psi (x^0,y)
\left (\r\intl_{-\iy}^{\iy}Y(dx^0,T^{y-t}( T^\tau m))X(t)dt\right ) dy=
<\nu (\bul, m),T^\tau m>.$$
We also check  the condition (4.1.6.3).
$$\intl_{y\leq 0}e^{\r y}S_t\nu (dy\otimes dx^0)=
\intl_{-\iy}^{\iy}X(t)dt\intl_{y\leq 0}e^{\r y})\r dy
\intl_{S^{m-1}}Y(dx^0,T^{y+\tau-t}m\leq$$
$$\leq\sup\limits_{\tau\in\BR}Y(S^{m-1},T^{\tau}m)\intl_{-\iy}^{\iy}X(t)dt\leq 1,$$
since $Y(\bul,\bul)\in K.$

Now we should ``transplant'' $\nu$  back to $\Rm\setminus 0$ such that $S_\tau$
passes to $(\bul)_{[e^\tau]}.$ Define $\mu(\bul,m)$ by
$$< \mu(\bul,m),\psi^*>:=<\nu (\bul,m),\psi>,\tag 4.1.6.5$$
where $\psi^*(rx^0)\in \Di(\Rm\setminus 0)$ and
$$\psi (x^0,y):=\psi^* (e^yx^0)e^{-(\r-m+2)y}\ \in \Di(Cyl).$$

Then
$$<(\mu)_{[e^\tau]},\psi^*>=<(\mu),T_{-\tau}\psi^*>=<\nu,S{-\tau}\psi>=
<S_\tau \nu,\psi>.$$
The condition $\mu (\bul,m) \in \Cal M[\r,\s]$ is also satisfied.

{\bf Exercise 4.1.6.8.}Check these properties.

Now we use (4.1.6.1) to transplant the dynamical system to $U[\r,\s].$
This completes a construction of an homomorphism $(T^t,M)\mapsto (T_t,U[\r,\s]).$

Let us check that it is an imbedding, i.e., we must check the
one-to-one correspondence. One-to-one correspondence of
$v(\bul,m)$ and $\mu (\bul,m)$ is known (Theorem 3.1.4.4).
One-to-one correspondence of $\mu (\bul,m)$ and $\nu(\bul,m)$ can
be also checked easily.

{\bf Exercise 4.1.6.9.} Check this in details.

So we should check the one-to-one correspondence of $\nu(\bul,m)$ and $Y(\bul,m).$

Suppose
$$\nu(\bul,m_1)=\nu(\bul,m_2).$$
Then
$$<\nu(\bul,m_1),\psi>=<\nu(\bul,m_2),\psi>\ \forall \psi\in \Di(Cyl).$$
In particular, set
$$\psi(x^0,y)=\phi (x^0)R(y),\ \phi\in \Di (S^{m-1}),\ R\in \Di (-\iy,\iy).$$
Then
$$<\nu(\bul,m_1),\psi>=\int R(y)dy
\intl_{-\iy}^\iy <Y(\bul,T^{y-t}m_1),\phi>_{S^{m-1}}
X(t)dt=\tag 4.1.6.6$$
$$=<\nu(\bul,m_2),\psi>=\int R(y)dy
\intl_{-\iy}^\iy <Y(\bul,T^{y-t}m_2),\phi>_{S^{m-1}}.$$
where
 $$ <Y(\bul),\phi>_{S^{m-1}}:=\intl_{S^{m-1}}\phi(x^0)Y(dx^0).$$
Set
$$F_j(y):=<Y(\bul,T^y m_j),\phi>_{S^{m-1}},\ j=1,2.$$
From  (4.1.6.6) we obtain for the convolutions
$$(F_1*X)(y)\equiv(F_1*X)(y),\ y\in (-\iy,\iy).$$
Thus
$$F_1(y)\equiv F_2(y),\ y\in (-\iy,\iy)$$
because of the property of $X.$

Hence
$$Y(\bul,T^y m_1)\equiv Y(\bul,T^y m_2),\ y\in (-\iy,\iy).$$
In particular, for $y=0$ we have
$$Y(\bul,m_1)=Y(\bul,m_2).$$
Hence $m_1=m_2$ because of (4.1.5.3), and this completes the proof of one-to-one
correspondence.

Consider the case of an integral $\r.$ For this case we can
use $v\in U[\r,\s]$ of the form
$$v(x)=\Pi_< (x,\mu,\rho)+\Pi_> (x,\mu,\rho)$$
instead of (4.1.6.1).

\qed
\edm
{\bf Exercise 4.1.6.10.}Check this.
\subheading {4.1.7}The most simple set satisfying the conditions of Theorem
4.1.3.4 is the set that is generated by a function $v\in U[\r]$ that has the
property
$$v_{[te^P]}=v_{[t]},\ t\in (0,\iy)$$
for some $P.$

Then
$$T_{t+P}v=T_t v,\ t\in (0,\iy)$$
i.e., the dynamical system $T_\bul$ is {\it periodic} with the period $P$
on the set
$$\Bbb C(v):=\{T_tv:0\leq t\leq P\}.$$
\proclaim {Theorem 4.1.7.1(Periodic Limit Set)} For all $P>0,\r>0,\s>0,$
there exists $v\in U[\r,\s]$ such that the dynamical system $(T_\bul,\BC (v))$
is periodic.
\endproclaim
\demo {Proof}Suppose $\r$ is non-integer. Let us take $\mu\in \Cal M[\r,\Dl]$
such that the canonical potential $\Pi(x,\mu,[\r])$ belongs to $U[\r,\s].$
This is possible because of Theorem 3.1.4.2.

Denote as $\mu^*_P$ the restriction of $\mu$ on the spherical ring
$\{x:1<|x|<e^P\}$ and set
$$\mu_P:=\sum\limits_{k=-\iy}^{\iy}T_{kP}\mu^*_P.$$
We have $\mu_P\in\Cal M[\r,\Dl]$ and
$$T_{t+P}\mu_P=T_t( \sum\limits_{k=-\iy}^{\iy}T_{(k+1)P}\mu^*_P)=T_t\mu_P.$$
Then $v:=\Pi(x,\mu_P,[\r])\in U[\r,\s]$ and $T_{t+P}v=T_t v$ because of
(3.1.5.0).

For an integer $\r$ we use the function
$$v(x):=\Pi_< (x,\mu_P,\r)+\Pi_> (x,\mu_P,\r).$$
\qed
\edm
\newpage

\centerline {\bf 4.2. Subharmonic function with prescribed limit set }
\subheading {4.2.1} The following two theorems  describe structure of limit sets
in terms of dynamical systems.
\proclaim {Theorem 4.2.1.1 (Necessity)} Let
$u\in SH(\BR^m,\r,\r (r)).$ Then the dynamical system
$(T_\bullet,\Fr [u,\bullet])$ is chain recurrent.
\ep
The chain recurrence is also sufficient.
\proclaim {Theorem 4.2.1.2 (Sufficiency)} Let $U$ be a compact connected and $T_\bullet$--
invariant subset of $U[\r,\sigma]$ for some $\sigma >0,$ such that the
dynamical system $(T_\bullet,U)$ is chain recurrent.

Then for any proximate order $\r (r)\rightarrow \r$ there exists
$u\in SH(\BR^m,\r,\r (r))$ such that
$$\Fr [u,\r (r),V_t,\BR^m]=U.$$
\ep
\demo {Proof of Theorem 4.2.1.1} We need the curve $u_t,\ t\geq 1,$ and
$\Fr[u,\bullet]$ to be contained in a common metric space $X.$ Thus we
set
$$X:=\{v\in SH(\BR^m):\sup\limits_{r\geq1}M(r,v)r^{-\r -1}\leq
\sup\limits_{r\geq1}M(r,u)r^{-\r -1}\}.$$
We want to use Theorem 4.1.1.1 cr 2). Let $U$ be an open proper subset of
$\Fr[u,\bullet]$ satisfying (4.1.1.1) and let $F$ be a $T_\bullet$--invariant subset
of $K:=\Fr[u,\bullet]\setminus U.$

Such $F$ exists. Indeed, $K$ is closed and $T_tK\subset K$ for $t>0$ (see
proof of Theorem 4.1.1.1, cr2)$\Longleftrightarrow$cr3)).Thus
$\Om (K)\subset K$ where $\Om (\bullet)$ was defined in (4.1.3.2). The set
$\Om (K)$ is invariant with respect to $T_t$ (see Exercise 4.1.3.2). So the
set of such sets $F$ is not empty.

If $F$ intersects $\p U$ at a point $v$, then $A(v)\subset F\cap\p U.$ Since
$A(v)$ is invariant (Exercise 4.1.3.3) $\p U$ contains a nonempty $T_\bullet$
--invariant set. So we obtain the assertion of the theorem  using
Theorem 4.1.1.1, cr2).

Suppose $F$ does not intersect $\p U.$ Let $U_0$ be an open set in $X$ such
that
$$U_0\cap\Fr [u,\bullet]=U,\ \clos U_0\cap\Fr [u,\bullet]=\clos U.\tag 4.2.1.1 $$
(see Exercise 4.2.1.1).
Since $\clos U_0 \cap F=\varnothing$ we can take a sequence of open neighborhoods $U_1,U_2,...$ of $F$ in $X$ such that the all sets $\clos U_j,\ j=1,2,...$ do not intersect
 $\clos U_0$ and $U_j\downarrow F.$

By definition of $\Fr [u,\bullet]$ we can find intervals
$a_j\leq t\leq b_j$ with $a_j\rightarrow\iy$ such that
$u_{e^{a_j}}\in \p U_j,\ u_{e^{b_j}}\in \p U_0,$ and
$\ u_{e^t}\not\in \clos U_0\cup\clos U_j$ for $a_j<t<b_j$.
We can pass to a subsequence and assume that
$$u_{e^{a_j}}\rightarrow w\in F.\tag 4.2.1.2$$
Let us use the following identity :
$$u_{e^{t+a_j}}=(u_{e^{a_j}})_{e^t}\frac {\r (e^t)\r (e^{a_j})}{\r (e^{t+a_j})}.
\tag 4.2.1.3$$
By (4.2.1.2),(4.2.1.3) and the property (3.1.2.2) of a proximate order we
obtain
$$u_{e^{t+a_j}}\rightarrow T_t w \in F$$
uniformly for any bounded set of $t.$

Thus $b_j-a_j\rightarrow \iy.$

 Passing to a subsequence we may assume that
$u_{e^{b_j}}\rightarrow v\in \Fr [u,\bullet]\cap\p U_0=\p U.$
Since $ u_{e^{t+b_j}}\rightarrow T_t v$ and $u_{e^{t+b_j}}\notin U_0$ when
$a_j-b_j<t<0$ we obtain that $T_tv\notin U$ when $t<0.$

Hence the whole backward orbit $\{T_tv:t<0\}$ lies in $\p U,$ which must
therefore contain the $T_\bullet$--invariant set $A(v).$
\qed
\edm

{\bf Exercise 4.2.1.1.} Prove that the set
$$U_0:=\bigcup\limits_{v\in U}\{w\in X:dist(v,w)<dist(v,K)/2\}$$
satisfies the conditions (4.2.1.1).
\demo {Proof}We have
$$U_0\supset U \Longrightarrow
U_0\cap\Fr [u,\bullet]\supset U\cap\Fr [u,\bullet]=U$$
Thus
$$U_0\cap\Fr [u,\bullet]\supset U\tag 4.2.1.4$$
From (4.2.1.4) we have
$$ \clos U_0\cap \Fr [u,\bullet]=\clos U_0\cap \clos\Fr [u,\bullet]=
\clos (U_0\cap\Fr [u,\bullet])\supset \clos U
\tag 4.2.1.5$$
 Finally $(4.2.1.4)\wedge  (4.2.1.5)\Longrightarrow (4.2.1.1).$
 \qed
\edm

\subheading {4.2.2}To prove Theorem 4.2.1.2 we need some preparation. Theorems of the next items form the basis of the construction that we will use in the proof.

Let $\beta$ be an infinitely differentiable function on $\Bbb R$ such that
$0\leq\beta(x)\leq 1,\ \beta(x)=0 $ for $x\leq 0$ and $\beta (x)=1$ for
$x\geq 1.$ We can set, for example,
$$\be (x):=A\intl_{-\iy}^x \a(y+1)dy$$
where $\a$ is taken from (2.3.1.1) and
$$A= \intl_{-\iy}^\iy \a(y+1)dy.$$
Suppose that the sequences $\{r_k,\sigma _k,\ k=0,1,..\}$ satisfy
the following conditions:
$$r_0=1;\ r_k<r_k\sigma_k<r_{k+1}/\sigma _{k+1}<r_{k+1},\ k=1,2,...
\tag 4.2.2.1$$
$$
\sigma
_k\uparrow\infty;\
  \frac {r_{k+1}}{\sigma_{k+1}r_k \sigma_k}\uparrow\infty .\tag 4.2.2.2
$$
Set
$$\psi _k (r):=\beta \left (\frac {\log r -\log (r_k/\sigma_k)}
{\log (\sigma_k r_k) -\log (r_k/\sigma_k)}\right )
-$$
$$\beta \left (\frac {\log r -\log (r_{k+1}/\sigma_{k+1})}
{\log (\sigma_{k+1} r_{k+1}) -\log (r_{k+1}/\sigma_{k+1})}\right )$$
$$\psi_0 (r):=1- \beta \left (\frac {\log r -\log (r_1/\sigma_1)}
{\log (\sigma_1 r_1) -\log (r_1/\sigma_1)}\right )$$
The sequence $\{\psi_k\},\ k=0,1,...$ forms a {\it partition of unity}
 with the following properties
\proclaim {Theorem 4.2.2.1(Partition of Unity)} One has
$$\sum_{k=0}^\infty \psi _k =1; \tag prtu1 $$
$$\roman {supp}\psi_k \subset (r_k/\sigma_k,r_{k+1}\sigma_{k+1});\tag prtu2$$
$$\psi_k (r)=1,\text { for }r\in (r_k\sigma_k,r_{k+1}/\sigma_{k+1});\tag prtu3$$
$$\roman {supp}\psi_k\cap\roman {supp}\psi_l=\varnothing  \text { for }
|k-l|>1;\tag prtu4$$
$$\lim\limits_{k\rightarrow \infty}\max\limits_r \psi'_k(r)r=
 \lim\limits_{k\rightarrow \infty}\max\limits_r \psi''_k(r)r^2=0.\tag prtu5$$
Moreover
$$\max\limits_r |\psi'_k(r)r|,\max\limits_r |\psi''_k(r)r^2|\leq \g_k\tag prtu6$$
where $\g_k$ can be made to tend to zero arbitrarily fast by choosing
the sequences $\{\sigma_k\}$ and $\{r_k\}.$
\ep
\demo {Proof} Set
$$\be_k(r):= \beta \left (\frac {\log r -\log (r_k/\sigma_k)}
{\log (\sigma_k r_k) -\log (r_k/\sigma_k)}\right ).$$
The functions $\be_k(r)$ and $\be_{k+1}(r)$ vanish for $r<r_k/\sigma_k$
because $\be (x)=0$ for $x\leq 0$, and both of them are equal to one for
$r\geq \sigma_kr_k$ because $\be (x)=1$ for $x\geq 1.$ Hence, (prtu2) holds.

One has for any $r\in (0,\iy)$
$$\sum_{k=0}^n \psi _k =1-\be_{n+1}(r).$$
As mentioned $\be_{n+1}(r)=0 $ for $n$ such that $r_{n+1}/\sigma_{n+1}>r.$
Thus (prtu1) holds.

Counting derivatives of $\psi_k,$ we have:
$$\max\limits_r |r\psi'_k (r)|\leq$$
$$\left [(\log (\sigma_k r_k) -\log (r_k/\sigma_k))^{-1}+
(\log (\sigma_{k+1} r_{k+1}) -\log (r_{k+1}/\sigma_{k+1}))^{-1}\right ]
\max\limits_x|\be'| (x).$$
Thus we can take the right side of the inequality as $\g_k$ and regulate its
vanishing by choice of the ratio in (4.2.2.2). The same holds for
$r^2\psi''(r).$ Hence (prtu5) and (prtu6) are proved.

{\bf Exercise 4.2.2.1.} Check (prtu4).
\qed
\edm

\subheading{4.2.3}Now we construct a function which is of zero type but has a
`` maximal possible'' mass density.
\proclaim {Theorem 4.2.3.1 (Maximal Mass Density Function)}Let $\r (r)\rightarrow \r,\ \r>0$ be a smooth
proximate
order (i.e., having properties (2.8.1.8)), and let $\g(r),\ r\in [0,\iy),$ satisfy
the conditions: $\g (r)>0$ and
$\g(r)\rightarrow 0$, as $r\ri \iy.$

Then there exists an infinitely differentiable  subharmonic function
$\Phi (x)$ such that
$$\Delta\Phi(x)\geq \g(x)|x|^{\r(r)-2}\tag 4.2.3.1$$
and
$$(\Phi)_t \ri 0\tag 4.2.3.2$$
in $\Di '$ as $t\ri\iy.$
\ep
To prove Theorem 4.2.3.1 we need an elementary lemma.
\proclaim {Theorem 4.2.3.2 (Convex Majorization)} Let $a(s),\ s\in [s_0,\iy)$ be a
function such that $a(s)\rightarrow -\iy$ as $ s\rightarrow \iy.$

Then there exists an infinitely differentiable, convex function $k(s)$ such that

k1) $k(s)\geq a(s);$

k2) $k(s)\downarrow -\iy$ as $s\rightarrow \iy$;

k3) $k^{(n)}(s)\ri 0$ for all $n=1,2,...\ .$
\ep

\demo {Proof}Set
$$a^*(s):=\sup \{a(t):t\geq s\}.$$
Then $a^*(s)\downarrow -\iy$ as $s\rightarrow \iy.$

Set $b_0:=-a^*(s_0)$ and denote as $s(b),\ b\in [b_0,+\iy) $ the
function inverse to the function $-a^*(s).$ Let us construct a convex function
 that majorates $s(b)$ and tends to infinity monotonically with all
its derivatives. It can be done in the following way.

First we construct a piece-wise linear convex function. Set
$$s_1(b):=s_0+1+\a _0(b-b_0),\ b\in [b_0,b_0+1],$$
and chose $\a _0$ such that the inequality $s_1(b)>s(b)$ holds for
$b\in [b_0,b_0+1].$

For this we choose
$$\a_0\geq \sup\limits_{b\in [b_0,b_0+1]}\frac{ s(b)-s_0-1}{b-b_0}$$
Since $s(b)-s_0-1<0$ the right side is finite.

For all the following intervals we set
$$ s_1(b):=s_1(b_0+j) +\a_j(b-b_0-j),\ b\in[b_0+j,b_0+j+1],$$
where $\a_j\geq\a_{j-1}$ and satisfies the condition
$$\a_j\geq\sup\limits_{b\in [b_0+j,b_0+j+1]}
\frac{ s_1(b)-s_1(b_0+j)}{b-b_0-j}$$                        To obtain a smooth function set
$$s_2(b):=\int \a(b-x)s_1(x)dx ,\tag 4.2.3.3$$
where $\a(x) $ is defined by (2.3.1.1).
Then $s_2(b)$ is infinitely differentiable, monotonic and convex .

{\bf Exercise 4.2.3.1.} Check this.

Set
$$k(s):=-s_2^{-1}(s),\tag 4.2.3.4$$
where $s_2^{-1}(s)$ is the inverse function to $s_2.$
One can check that $k(s)$ satisfies the properties k1), k2), k3).

\qed
\edm
{\bf Exercise 4.2.3.2.} Check that $k(s)$ satisfies k1),k2),k3).

\demo {Proof of Theorem 4.2.3.1}We are going to show that $\Phi$ can be taken
in the form
$$\Phi (x):=ce^{k(\log |x|^2)}|x|^{\r (|x|)}.\tag 4.2.3.5$$
where $c$ and $k(s)$ will be chosen later.

Note that $\Phi(x)=\Phi (|x|)$ depends only on $r=|x|$ and pass to the
variable $s:=\log r^2.$ Then for $\phi (s):=\Phi(e^{s/2})$ we have
$$\Delta \Phi (x)=r^{1-m}\frac {\p}{\p r}r^{m-1}\frac {\p}{\p r}
ce^{k(\log r^2)}r^{\r (r)}=$$
$$ce^{-s}\left (\frac {\p^2}{\p s^2 }+\frac {m-2}{2}\frac {\p}{\p s }\right )
\phi (s)\geq cme^{-s} \min [\phi''(s),\phi '(s)].\tag 4.2.3.6$$
Let us chose $k$ as in Theorem 4.2.3.2 with
$a(s):=\log\g (r)=\log\g (e^{\frac {s}{2}}).$
 Now we estimate the derivatives from below.
$$\phi '(s)=\phi (s) [k'(s)+\frac {1}{4}s\r '(e^{\frac {s}{2}}) +
\frac {1}{2}\r (e^{\frac {s}{2}})].$$
By k3) and k1)  $k'(s)\ri 0$ and $k(s)\geq a(s).$
Also $s\r '(e^{\frac {s}{2}})\ri 0$ and
$\r (e^{\frac {s}{2}})\ri \r$ by properties of proximate order
(Theorem 2.8.1.4).
Thus we can chose $c$ such that
$$\phi '(s)>\frac {1}{m}e^{\log \g (e^{\frac {s}{2}}) +
\frac {s}{2}\r (e^{\frac {s}{2}})}.\tag 4.2.3.7$$
Differentiating once again, we obtain
$$\phi ''(s)=\phi (s) [k'(s)+\frac {1}{4}s\r '(e^{\frac {s}{2}})+
\frac {1}{2}\r (e^{\frac {s}{2}})]^2 +k''(s)+
\frac {1}{2}\r '(e^{\frac {s}{2}})+\frac {1}{8}s\r'' (e^{\frac {s}{2}})].$$
From here we obtain by choosing $c$ :
$$\phi ''(s)>\frac {1}{m}e^{\log\g (e^{\frac {s}{2}}) +
\frac {s}{2}\r (e^{\frac {s}{2}})}.\tag 4.2.3.8$$
Using (4.2.3.6) ,(4.2.3.7) and (4.2.3.8) we obtain:
$$\Delta \Phi (s)>e^{\log\g (e^{\frac {s}{2}}) +
\frac {s}{2}\r (e^{\frac {s}{2}})}.$$
Returning to the variable $r$ we obtain (4.2.3.1).
Correctness of (4.2.3.2) can be checked directly using k2) and properties
of the proximate order (Theorem 2.8.1.3).

{\bf Exercise 4.2.3.3.} Check this.
\qed
\edm
\subheading {4.2.4}
We have already approximated distributions and subharmonic functions by infinitely differentiable functions (Theorems 2.3.4.5 and 2.6.2.3). Now we need to
make more precise this approximation.Namely, we are going  to make it uniform with respect to
$v\in U[\r,\s].$
 We will denote
$$\p^l:=\frac {\p^{|l|}}{(\p x_1)^{l_1}(\p x_2)^{l_2}...(\p x_m)^{l_m}}
\tag 4.2.4.1$$
where $l=(l_1,l_2,...l_m ),\ |l|=l_1+l_2+...+l_m.$

Set for $v\in U[\r,\s]$
$$R_\eps v (x):=\int \a _\eps (x-y)v(y)dy \tag 4.2.4.2$$
where $\a _\eps$ is taken from (2.3.1.3).

We have changed the notation from 2.3.1 and 2.6.2 because a subindex of
$v$ was already engaged for $t.$

For a fixed $0<\dl\leq 0.5$ set
$$Str(\dl):=\{x:\dl\leq|x|\leq\dl^{-1}\}\tag 4.2.4.3$$
\proclaim {Theorem 4.2.4.1 (Estimation of $R_\eps$)} Let $v\in U[\r,\s].$ Then

R1. for a fixed $g\in \Di(\Rm\sm 0)$ with  supp $g\sbt Str(\dl)$
$$|<R_\eps v-v,g>|\leq o(1,g)2\s \dl^{-\r} \tag 4.2.4.4$$
where $o(1,g)\ri 0$ as $\eps\ri 0;$

R2. the inequality
$$|\p^lR_\eps v(x)|\leq A(m)\s\eps^{-|l|-m+1}|x|^{-|l|+\r},\tag 4.2.4.5$$
with $A(m)$  depending only on $m,$ holds for $\eps <|x|/2$.
\ep
\demo {Proof} One has
$$<R_\eps v,g>=<v, R_\eps g>.\tag 4.2.4.6$$
Thus
$$<R_\eps v-v,g>=<v, R_\eps g -g>.\tag 4.2.4.7$$

{\bf Exercise 4.2.4.1.} Check (4.2.4.6) and (4.2.4.7).

Now
$$|<v, R_\eps g -g>|\leq \max\limits_{Str(\dl)}|R_\eps g -g|(x)
\intl_{Str(\dl)}|v|(x)dx \tag 4.2.4.8$$
The first factor is $o(1)$ because $g$ is smooth. For the second one we have
$$\intl_{Str(\dl)}|v|(x)dx\leq 2 \intl_{Str(\dl)}v^+(x)dx\leq 2\s \dl^{-\r}.
\tag 4.2.4.9$$
This and (4.2.4.8) imply R1).

Differentiating the equality
$$R_\eps v(x):=C_m\int \eps ^{-m}  \a (|x-y|/\eps)v(y)dy,$$
we have
$$|\p ^l R_\eps v(x)|\leq
 C_m\eps ^{-|l|-m}\max\limits_{\{|y|<\eps\}}|\p ^l\a (|y|)|
\intl_{\{|y|<\eps\}}|v(x-y)|dy.$$
Suppose $|x|=1.$ Then for $0<\eps\leq 0.5,$ we have
$$\intl_{|y|<\eps}|v|(x-y)dy\leq\intl_{1-\eps<|x|<1+\eps}|v|(x)\leq$$
$$2\intl_{1-\eps<|x|<1+\eps}v^+(x)dx\leq
\s_m2\cdot 2\eps\s(1+\eps)^\r\leq \s_m 6\s\eps$$
where $\s_m$ is the square of the unit sphere.
Hence for $|x|=1$
$$|\p^lR_\eps v(x)|\leq A(m)\s\eps^{-|l|-m+1}.\tag 4.2.4.10$$
with
$$A(m)=6\s_m\max\limits_{y\in \Rm }|\p ^l\a (|y|)|.$$
Set $t=|x|$. Apply the inequality (4.2.4.10) to $v:=v_{[t]}(y)$ with
$y:=x/|x|$.
 Then
$$|\p^lR_\eps v_{[t]}(y)|\leq A(m)\s \eps^{-|l|-m+1}.$$
Computing the derivatives, we obtain
 $$\p^lR_\eps v_{[t]}(x)=t^{-\r}t^{|l|}\p^lR_\eps v(x)|_{x=ty}$$
Thus one has R2.\qed
\edm

\subheading {4.2.5}In this item we describe the main part of a construction
that will be used in the proof of the Theorem 4.2.1.2.

Let $\{v_j\in U[\r,\s],\ j=1,2,...\}$ and  $\{\psi_j,\ j=1,2...\}$ be the partition of unity from Theorem 4.2.2.1. Let us chose
$\eps_j\downarrow 0$ such that the condition
$$\g_j\eps_j^{-m}\rightarrow\iy \tag 4.2.5.1$$
holds for $\g_j$  taken from Theorem 4.2.2.1, (prtu 6).

Set
$$v(x|t):=\sum_{j=0}^\iy \psi_j(t)(v_j)_{[t]}(x), \tag 4.2.5.2$$
where $(\cdot)_{[t]}$ defined by (3.1.2.4a).

One can see that $v(x|t)\in U[\r,3\s]$ for all $t.$

{\bf Exercise 4.2.5.1.} Show this, using properties of $\{\psi_j\}$ and
invariance of $U[\r,\s]$ with respect to $(\cdot)_{[t]}$.

 We can consider $v(x|t)$ as
a curve (a {\it pseudo-trajectory}) in $ U[\r,3\s].$

Set
$$ u(x):=\sum_{j=0}^\iy \psi_j(|x|)R_{\eps_j}(v_j)(x)|x|^{\r(|x|)-\r}.
\tag 4.2.5.3 $$
where $R_\eps$ is defined by (4.2.4.2).

It is infinitely differentiable function in $\Rm .$
\proclaim {Theorem 4.2.5.1(Construction )}One has
$$u_t -v(\bullet|t)\rightarrow 0\tag 4.2.5.4$$
in $\Di'(\Rm),$ and
$$\Delta u (x)=f(x)+\g (x)|x|^{\r(|x|)-2}\tag 4.2.5.5$$
with $f(x)\geq 0$ and $\g(x)=o(1)$ as $|x|\rightarrow\iy.$
\ep
Let us note that the function $u(x)$ is ``almost-subharmonic'' and can be
made subharmonic by summing with the function $\Phi$ from Theorem 4.2.3.1.

{\bf Exercise 4.2.5.2.} Prove this.

So we have \proclaim {Theorem 4.2.5.2(Pseudo-Trajectory
Asymptotics)} For any $v(x|t)$ of the form (4.2.5.2) there exists
an infinitely differentiable function $u\in SH(\r(r))$ that
satisfies (4.2.5.4). \ep \demo {Proof of Theorem 4.2.5.1} One has
$$ u_t(x):=\sum_{j=0}^\iy \psi_j(t|x|)(R_{\eps_j}(v_j))_{[t]}(x)a(x,t),$$
where
$$a(x,t):=\frac{|tx|^{\r(t|x|)-\r}}{t^{\r(t)-\r}}.$$
For any $0<\dl<0.5$ and $x\in Str(\dl)$
$a(x,t)\rightarrow 1$ uniformly in $|x|$as $t\ri\iy.$ This follows from
Theorem 2.8.1.3 , ppo3).

{\bf Exercise 4.2.5.3.} Check this in details.

We have
$$u_t(x)-v(x|t)=\sum_{j=0}^\iy
[\psi_j(t|x|)(R_{\eps_j}(v_j))_{[t]}(x)a(x,t)-\psi_j(t)(v_j)_{[t]}(x)],
\tag 4.2.5.6$$
and there are  no more than three summands in the sum for sufficiently large
$t=t(\dl)$ because of Theorem 4.2.2.1, prtu4. Let us estimate every summand.
One has
$$b_j(x,t):=
[\psi_j(t|x|)(R_{\eps_j}(v_j))_{[t]}(x)a(x,t)-\psi_j(t)(v_j)_{[t]}(x)]=$$
$$[\psi_j(t|x|)-\psi_j(t)](R_{\eps_j}(v_j)))_{[t]}a(x,t)+ $$
$$\psi_j(t)(R_{\eps_j}(v_j))_{[t]}(x)[a(x,t)-1]+
\psi_j(t)[(R_{\eps_j}(v_j))_{[t]}(x)- (v_j)_{[t]}(x)]:= $$
 $$:=(a_1+a_2+a_3)(x,t).$$
Let us estimate $<b_j(t,\bullet),g>$ for every $g\in \Di(\Rm\sm 0).$

We can assume that $\roman{supp}\ g\sbt Str (\dl).$
Set
$$M(g):=\max\limits_{x\in Str (\dl)}|g| (x).$$

 We have
$$|<a_1 (\bullet,t),g>|\leq M(g)\max\limits_{r\in (0,\iy)} |r\psi'_j(r)|\dl ^{-1}
\intl_{Str(\dl)}|(R_{\eps_j}(v_j))_{[t]}|(x)dx.$$
One can check that
$$\intl_{Str(\dl)}|(R_{\eps_j}(v_j))_{[t]}|(x)dx\leq
3\s \dl^{-\r}.$$

{\bf Exercise 4.2.5.4.} Check this using (4.2.4.9) and the invariance of
$U[\r,\s]$ with respect to $(\bullet)_{[t]}$ (see (3.1.2.4)).

Hence
$$|<a_1 (\bullet,t),g>|\leq C_1(g)\g_j.\tag 4.2.5.7$$
Let us estimate $a_2(x,t)$. We have
$$<a_2(\bul,t),g>\leq\max\limits_{Str(\dl)}
|a(x,t)-1|\psi_j(t)M(g)3\s\dl^{-\r}=C_2(g)o(1)\tag 4.2.5.8$$
where $o(1)\ri 0$ as $t\ri\iy$.

For estimating $a_3(x,t),$ we use Theorem 4.2.4.1(Estimation of $R_\eps$),
 (4.2.4.4):
$$|<a_3 (\bullet,t),g>|\leq o(\eps_j,g)2\s\dl^{-\r}.\tag
 4.2.5.9$$
where $o(\eps_j,g)\ri 0$ as $j\ri \iy.$

Hence (4.2.5.7), (4.2.5.8) and (4.5.5.9) imply
$$<b_j(\bul,t),g>\ri 0\tag 4.2.5.10$$
as $t\ri\iy$ and $j\ri\iy.$

Suppose, for a large fixed $t,$ the sum (4.2.5.6) contains $b_j(x,t)$ for
$j=j(t),j=j(t)+1$ and $j=j(t)+2.$ This implies that $j(t)\ri\iy$ as $t\ri\iy.$

Since
$$<u_t(\bul)-v(\bul|t),g>=<b_{j(t)}(\bul,t),g>+<b_{j(t)+1}(\bul,t),g>+
<b_{j(t)+2}(\bul,t),g>$$
 we obtain from (4.2.5.10) that $<u_t(\bul)-v(\bul|t),g>\ri 0$ as $t\ri\iy$
for any $g\in \Di(\Rm\sm 0).$

This is (4.2.5.4).

Let us prove (4.2.5.5). We have
$$\Dl u=\sum_{j=0}^\iy [\Dl(R_{\eps_j}v_j)(x)\psi_j (x) |x|^{\r(|x|)-\r}+
\sum_{l,m.k}\p^l(R_{\eps_j}v_j)(x)\p^n\psi_j(x)\p^k |x|^{(\r|x|)-\r}],
\tag 4.2.5.11$$
where $l,m,k$ are multi-indexes that satisfy the condition:
in any summand there are derivatives in the same variable, the  derivatives of
$\psi_j$ and $|x|^{\r (|x|)-\r}$ have no more than  second order and
the derivatives of $R_{\eps_j}v_j (x)$ have no more than first order.

{\bf Exercise 4.2.5.5.}Check this.

 As usual, the derivative of zero order is the function itself.

For any $x\in Str (\dl),$ the outside sum contains no more
then three summands. First we consider only the terms in the square brackets.

The first term is non-negative because of subharmonicity of $R_{\eps_j}v_j$
 and non-negativity of all the other factors. Set
$$f(x):=\sum_{j=0^\iy}[\Dl(R_{\eps_j}v_j)(x)\psi_j (x) |x|^{\r(|x|)-\r}\geq 0.\tag 4.5.2.12$$
 Using Theorem 4.2.4.1, R2) we obtain
$$|\p^l(R_{\eps_j}v_j)(x)|\leq$$
$$A(m)\s \eps_j ^{-|l|-m+1}|x|^{-|l|+\r}\tag 4.2.5.13$$
for $|l|=0$ or $|l|=1.$

From Theorem 4.2.2.1, prtu6), and inequality $|\p_{x_i}|x||\leq 1$
we obtain
$$|\p^n\psi_j (|x|)\leq |\psi^{(n)}(r)||_{r=|x|}\leq \g_j|x|^{-|n|}\tag 4.2.5.14$$
for $|n|=1,2.$

Using properties of the smooth proximate order (Theorem 2.8.1.4), one can
obtain
$$|\p ^{|k|} |x|^{\r(|x|)-\r}|=(|x|^{\r(|x|)-\r-|k|})|_{r=|x|}(1+o(1)),
\tag 4.2.5.15$$
as $|x|\ri\iy.$

{\bf Exercise 4.2.5.6.} Check in details (4.2.5.13), (4.2.5.14) and (4.2.5.15).

Thus, for every term of the inner sum, we have
$$|\p^l(R_{\eps_j}v_j)(x)\p^n\psi_j(x)\p^k |x|^{\r(|x|)-\r}|\leq$$
$$A(m)\s\g _j\eps_j ^{-|l|-m+1}|x|^{-2+\r}|x|^{\r(|x|)-\r}\leq \be_j |x|^{\r(|x|)-2},\tag 4.2.5.16$$
 where $\be_j\ri 0$ because of the condition (4.2.5.1).

Recall that for every large $x$ the outside sum contains no more then three
summands, say, $j=j(x),j=j(x)+1 $ and $j=j(x)+2.$ Thus $j(x)\ri\iy$ as
$|x|\ri\iy.$ Hence (4.2.5.12) and (4.2.5.16) imply (4.2.5.5).\qed
\edm

 \subheading {4.2.6}\demo {Proof of Theorem 4.2.1.2}Let $v(\bul|t)$ have the form (4.2.5.2). We denote as $\Om (v)$ a set of
the $\Di'$--limits of the form
$$ w:=\liml_{t_k\ri\iy}v(\bul|t_k).$$
We are going to construct some $v(\bul|t)$ for which
$$ \Om (v)=U,\tag 4.2.6.1$$
and at the next step to use Theorem 4.2.5.2 to obtain a subharmonic
function with the same limit set.

First we describe the construction of the function  $v(\bul|t).$
Let $\{r_k,t_k,\ k=1,2,...\}$ be an alternating sequences
$r_0=1 ,\ r_k<t_k<r_{k+1}$ such that
$$\liml_{k\ri\iy}\frac {t_k}{ r_k}= \liml_{k\ri\iy}\frac{r_{k+1}}
{t_k}=\iy.$$
Let us chose in $U$ a countable, dense set $\{g_j\}$ and form from it
a sequence $\{w_k\}$ such that every element $g_j$ is repeated infinitely often. For example,
$$w_1:=g_1,w_2:=g_1,w_3:=g_2,w_4:=g_1,w_5:=g_2,w_6:=g_3,...\ .$$
Set
$$q_k:=(w_k)_{[1/t_k]}=T_{-\log t_k}w_k$$
in the notation (4.2.1.1).

Now we use that $(T_{\bul},U)$ is chain recurrent. Set
$$\a_k:=\log \frac {r_k}{t_k};\om_k:=\log \frac{r_{k+1}}{t_k}$$
and find, by Theorem 4.1.2.1, a sequence $\{v_j\}\spt \{q_k\}$ such that
the condition (4.1.2.1) holds, i.e.,
$$T_{\om_k}v_k -T_{\a_{k+1}}v_{k+1}\ri 0 \tag 4.2.6.2$$
as $k\ri\iy.$

Set in Theorem 4.1.2.3
$$p_k:=T_{\a_{k+1}}v_{k+1},\ q_k:=T_{\om_k}v_k$$
and find $\g_k$ such that the condition
$$T_\tau\circ T_{\om_k}v_k -T_\tau\circ T_{\a_{k+1}}v_{k+1}\ri 0 \tag 4.2.5.3$$
holds uniformly for $\tau\in [-\g_{k+1},\g_k].$

Set
$$\s_k:=\min\left [e^{\g_k},\sqrt {\frac {t_k}{r_k}}\right].$$
These $\s_k$ satisfy the conditions (4.2.2.1) and (4.2.2.2).

{\bf Exercise 4.2.6.1.} Check this.

We define $v(\bul|t)$ by (4.2.5.2) with described $v_j$ and with $\psi_j$ from
Theorem 4.2.2.1, corresponding to the chosen $r_j$ and $\s_j.$
Let us  prove (4.2.6.1).

Consider for  fixed $k$ the following three cases.

1. $t\in [r_k\s_k,r_{k+1}/\s_{k+1})$;

2. $t\in [r_{k+1}/\s_{k+1}, r_{k+1});$

3. $t\in [r_k,r_k\s_k).$

For the first case we have
$$v(\bul|t)=(v_k)_{[t/t_k]}=T_{\log t/t_k}v_k.$$
For the second one
$$v(\bul|t)=\psi_k(t)(v_k)_{[t/t_k]}+\psi_{k+1}(t)(v_{k+1})_{[t/t_{k+1}]}=$$
$$(v_k)_{[t/t_k]}+\psi_{k+1}(t)[(v_{k+1})_{[t/t_{k+1}]}-(v_k)_{[t/t_k]}$$
We transform the expression in the square brackets
$$(v_{k+1})_{[t/t_{k+1}]}=T_{\log (t/t_{k+1})}v_{k+1}=
T_{\log (t/r_{k+1})}\circ T_{\log (r_{k+1}/t_{k+1})}v_{k+1}=
T_{\log (t/r_{k+1})}\circ T_{\a_{k+1}}v_{k+1}.$$
For the second term, we obtain
$$(v_k)_{[t/t_k]}=T_{\log (t/r_{k+1})}\circ T_{\om_{k}}v_{k}.$$

{\bf Exercise 4.2.6.2.} Check this.

Setting $\tau:=\log (t/r_{k+1}),$ we have
$$v(\bul|t)= (v_k)_{[t/t_k]}+\psi_{k+1}(t)[T_\tau\circ T_{\a_{k+1}}v_{k+1}-
T_\tau\circ T_{\om_{k}}v_{k}],\ \tag 4.2.6.4$$
where $\tau\in [-\log\s_{k+1},0)\sbt [-\g_{k+1},\g_k]$.
For the third case, set $\tau:=\log (t/r_k).$ Then
$$v(\bul|t)= (v_k)_{[t/t_k]}+\psi_{k}(t)[T_\tau\circ T_{\om_{k-1}}v_{k-1}-
T_\tau\circ T_{\a_{k}}v_{k}],\ \tag 4.2.6.5$$
where $\tau\in [0,\log\s_{k})\sbt [-\g_k,\g_{k-1}]$.

Let $t_N\ri\iy$ be an arbitrary sequence. Choosing a subsequence,  we may
suppose that there exist the limits
$(v_{k(t_N)})_{[t_N/t_{k(t_N)}]}\ri v^*\in U$
and
$v(\bul|t_N)\ri v_\iy.$

Choosing a subsequence, we may suppose that $t_N$ satisfies either 1. or 2.
or 3. For the case 1., we obtain at once $ v_\iy= v^*\in U.$

For the case 2., from  (4.2.6.4), (4.2.6.2) and Theorem 4.1.1.3 we
obtain that the superfluous addends tend to zero, and hence $
v_\iy\in U.$

The same holds for the case 3. Hence $\Om(v)\sbt U.$

Further, for $t=t_k$, we have  $v(\bul|t)=w_k.$ The sequence
$\{w_k\}$ contains the set $\{g_j\}$ that is dense in $U. $ Thus
$\Om(v)\spt U.$ Thus equality (4.2.6.1) has been proved.

As  already said, the application of Theorem 4.2.5.2 concludes the proof.
\qed
\edm
\newpage

\centerline {\bf 4.3. Further properties of Limit Sets }
\subheading {4.3.1} Let as mark the following property of the
pseudo-trajectory $v(\bul|t)$ defined in (4.2.5.2):
\proclaim {Theorem 4.3.1.1} One has
$$T_\tau v(\bul|e^t)-v(\bul|e^{t+\tau})\ri 0\tag 4.3.1.1$$
as $t\ri\iy$ uniformly with respect to $\tau\in [a,b]$ for any
$[a,b]\sbt (-\iy,\iy).$
\ep
\demo {Proof} Using definition of $(\bul)_t$ (see (3.1.2.1)) and (4.2.5.4)
the remainder in (4.3.1.1) can be represented in the form
$$b(t,\tau,\bul):=T_\tau v(\bul|e^t)-v(\bul|e^{t+\tau})=
T_\tau (u_{e^t})-u_{e^{t+\tau}}+o(1)$$
where $o(1)\ri 0$ uniformly with respect to $\tau\in [a,b]$ for any
$[a,b]\sbt (-\iy,\iy).$

{\bf Exercise 4.3.1.1.} Check this in details.

Then we obtain
$$b(t,\tau,\bul)=u_{e^{t+\tau}}[e^{\r (e^t)-\r (e^{t+\tau})}-1]+o(1)\ri 0$$
uniformly in the same sense due to precompactness of the family
$\{u_{e^t}\}$ and properties of the proximate order.

{\bf Exercise 4.3.1.2.} Check this in details.\qed
\edm

The property (4.3.1.1) shows that the pseudo-trajectory $v(\bul|t)$ behaves
 asymptotically like the dynamical system $T_\bul.$ Thus a
pseudo-trajectory with this property is called
an {\it asymptotically dynamical } pseudo-trajectory with
 {\it  dynamical asymptotics } $T_\bul$ (a.d.p.t.).

Theorem 4.2.5.1 shows that for any a.d.p.t. of the form (4.2.5.2) there exists
$u\in SH(\r (r))$ that satisfies the condition
$$u_{e^t}-v(\bul|e^t)\ri 0.\tag 4.3.1.2$$
as $t\ri\iy.$

The following assertion shows that we can suppose $v(\bul|\bul)$
to be an arbitrary, in some sense, a.d.p.t.

We call a pseudo-trajectory $w(\bul|\bul)$ {\it piecewise continuous} if the
property
$$w(\bul|t+h)-w(\bul|t)\ri 0 \tag 4.3.1.3$$
as $h\ri 0$ holds for all $t$ except may be a countable set  without  points
of condensation.

Let $U\sbt U[\r,\s]$ for some $\s>0.$ A pseudo-trajectory $w(\bul|\bul)$
is called $\om${\it --dense} in $U$ if $\Om (w)=U$ (see (4.1.3.2)),i.e.,
$$\{v\in U[\r]:(\exists t_j\ri\iy)\  v=\Di'-\lim w(\bul|e^{t_j})\}=U\tag 4.3.1.4$$

We have proved already that $v(\bul|\bul)$ defined by (4.2.5.2) has this
property (see (4.2.6.1)).

Now we consider again the dynamical system $(T_\bul,U)$ where $U\sbt U[\r,\s]$
for some $\s>0$ and $T_t$ is defined by (4.2.1.1).
\proclaim {Theorem 4.3.1.2 (A.D.P.T. and Chain Recurrence)}$(T_\bul,U)$ is
chain recurrent iff there exists an a.d.p.t. that is piecewise continuous and $\om$-dense in $U.$
\ep
Necessity has been proved already, because the pseudo-trajectory (4.2.6.2)
possesses this property. Sufficiency will be proved later.

The claim of piecewise continuity can be justified by
\proclaim {Theorem 4.3.1.3} For any $u\in SH(\r (r))$ there exists a
piecewise continuous pseudo-trajectory $w(\bul|\bul)$ such that
$$u_{t}-w(\bul|t)\ri 0.\tag 4.3.1.5$$
as $t\ri\iy.$
\ep
Of course, $w(\bul|\bul)$ is a.d.p.t.

{\bf Exercise 4.3.1.2.}Check this.
\subheading {4.3.2} \demo {Proof of Theorem 4.3.1.3}
Let $\{t_n\}$ be any sequence such that
$$t_n\ri\iy,\ t_{n+1}/t_n \ri 1,\tag 4.3.2.1$$
for example, $t_n=n.$

There exists a sequence $\{v_n\}\sbt \Fr [u]$ such that
$$u_{t_n}-v_n\ri 0.\tag 4.3.2.2$$
Set
$$w(\bul|t):=v_n,\ for\ t_n<t\leq t_{n+1}.\tag 4.3.2.3$$
This is a piecewise continuous function.

Let us prove that
$$u_t -w(\bul|t)\ri 0.\tag 4.3.2.4 $$
Assume the opposite; i.e., there exists a sequence $\{t'_k\}$ such that it is not
true. We can suppose that
$$u_{t'_k}\ri w_1\in \Fr [u], \ w(\bul|t'_k)\ri w_2\in \Fr [u],\ w_1\neq w_2.
\tag 4.3.2.5$$
Let us find a sequence $\{n_k\}$ such that $t_{n_k}<t'_k<t_{n_k +1}$ .
Then
$$t_{n_k}/t'_k\ri 1.\tag 4.3.2.6$$
From (4.3.2.5), (4.3.2.3) and (4.3.2.2) we have
$$u_{t_{n_k}}\ri w_2.\tag 4.3.2.7$$
Then we have, using properties of $(\bul)_t$ and the proximate order,
$$u_{t'_k}=(u_{t_{n_k}})_{[t'_k/t_{n_k}]}(1+o(\log(t'_k/t_{n_k}))\ri w_2
\tag 4.3.2.8$$
because of (4.3.2.6) and the continuity of $u_{[t]}$ on $(u,t).$

However (4.3.2.8) contradicts (4.3.2.5). Thus (4.3.2.4) holds.
\qed
\edm
\subheading {4.3.3}Now we will prepare the proof of Theorem 4.3.1.2.

Let $\{v_k,\ k=1,2,...\}\subset U[\rho,\s]$ for some $\s$ be a sequence of
functions and
$\{r_k,\ k=1,2,...\},\{t_k\  k=1,2,...\}$
be two sequences such that
$$0<r_k<t_k<r_{k+1},\ k=1,2,...\tag 4.3.3.1$$
and
$$\lim\limits_{k\rightarrow \infty}t_k/ r_k=
\lim\limits_{k\rightarrow \infty}r_{k+1}/t_k=\infty.\tag 4.3.3.2$$
Set
$$w^*(\bul|t):=(v_k)_{t/t_k},\text { for } t\in [r_k,r_{k+1})\tag 4.3.3.3$$
where $k=1,2,...$

\proclaim {Theorem 4.3.3.1}Let $w (\bul|\bul)\subset U$ be an arbitrary
 $\omega$-dense a.d.p.t. and
$\{p_j,\ j=1,2,...\}\subset U$ an arbitrary sequence. Then there exists
a sequence $\{v_k,\ k=1,2,...\}$$\supset\{p_j,\ j=1,2,...\}$ and sequences
$\{r_k,\ k=1,2,...\}$ and $\{t_k,\  k=1,2,...\}$ satisfying (4.3.3.1) and
(4.3.3.2)
such that for $w^*(\bul|\bul)$ determined by (4.3.3.3) the condition
$$w^*(\bul|t)-w (\bul|t)\rightarrow 0\tag 4.3.3.4$$
as $t\ri\iy$ is fulfilled.
\endproclaim
This proposition shows that any $\omega$-dense a.d.p.t. is equivalent to
 one constructed of long pieces of
trajectories of the dynamical system $T_\tau.$
\demo {Proof of Theorem 4.3.3.1}We can take sequences
$\{\eps_j\downarrow 0,\ j=1,2,...\}$ and \linebreak
$\{b_j\uparrow\infty,\ j=1,2,...\}$
and choose a sequence $\{\tau_j,\ j=1,2,..\}$ such that the inequalities
$$d (T_\tau p_j -T_\tau  w(\bul|\tau_j))<\eps_j/2 \tag 4.3.3.5$$
and
$$d (T_\tau w(\bul|t)-w(\bul|e^\tau t))< \eps_j/2,\ t>\tau_j
\tag 4.3.3.6 $$
are fulfilled uniformly with respect to $\tau\in[b_{j+1}^{-2},b_j^{2}].$

Indeed, $w(\bul|\bul)$ is $\omega$-dense in $U,$ and hence we can find
$\tau_n\rightarrow \infty$ such that
$$
p_n - w(\bul|\tau_n)\rightarrow 0.
$$
Set in Lemma 4.1.2.2
$$
p_n:=p_n,\ q_n:=w(\bul|\tau_n),\ \g_n:=2\log b_j
$$
Then for any $\eps_j$ we can find $\tau_j:=\tau_{n_j}$ such that (1.2.5)
holds uniformly with respect to  $\tau\in [b_{j+1}^{-2},b_j^{2}].$

The inequality (4.3.3.6) holds, because $w (\bul |\bul)$ is asymptotically
dynamical (see (4.3.1.1)).

We can also suppose without loss of generality that
$\tau_j>\tau_{j-1}b_{j-1}^2,$ i.e., that the sequence $\{\tau_j\}$ is rather
thin.

The inequality (4.3.3.5) shows that for intervals of $t$ that are determined by
the inequality $b_{j+1}^{-2}\leq t/\tau_j\leq b_j^2$ our pseudo-trajectory is
already close to some trajectories.

Now we  divide the spaces between such intervals into equal parts in
the logarithmic scale such that their logarithmic lengths would be between
$\log b_j$ and $\log b_{j+1}$,so that they tend to infinity.

To this and set
$$n_j:=\left [\frac {\log \tau_{j+1}- \log \tau_{j}}{b_j}\right]$$
where $[\cdot]$ means the entire part, and
$$\gamma _j:=(\tau_{j+1}/\tau_{j})^{\frac {1}{2n_j}}.$$
It is clear that
$b_j<\gamma _j<b_j^2.$
As centers of new intervals we take the points
$$\tau_{j,l}:=\tau_j\gamma_j^{2l},\ l=0,1,...,n_j.$$
Thus
$\tau_{j,0}=\tau_j$
and
$\tau_{j,n_j}=\tau_{j+1}.$
The ends of the intervals are $\tau_{j,l}/\gamma_j$ and $\tau_{j,l}\gamma_j.$
Now we complete the sequence $\{p_j\}$ by the values of the
pseudo-trajectory $w(\bul|t)$ in the centers of the intervals, i.e., we set
$$p_{j,l}:=w(\bul|\tau_{j,l}),\ l=1,...,n_j-1$$
For $t\in (\tau_{j,l}/\gamma_j,\tau_{j,l}\gamma_j),\ l=1,...,n_j-1$ we have
$$d ((p_{j,l})_{t/\tau_{j,l}}-w(\bul |t))<\eps_j/2\tag 4.3.3.7$$
because of (4.3.3.6).

For $l=0$ and $l=n_j$ we set accordingly
$$p_{j,0}:=p_{j};\ p_{j,n_j}:=p_{j+1}.$$
Using (4.3.3.5) and (4.3.3.6) we have the inequality like (4.3.3.7) for
$l=0,l=n_j$ but with $\eps_j$ instead of $\eps_j/2.$

We complete the proof re-denoting all the centers $\tau_{j,l}$ as $t_k,$
 all the ends as $r_k$ and all the $p_{j,l}$ as $v_k.$
\qed
\enddemo
\subheading {4.3.4}\demo {Proof of sufficiency in Theorem 4.3.1.2} A direct
corollary of the previous Theorem 4.3.3.1 is
$$w^*(\bul|r_k-0)-w^*(\bul|r_k)\ri 0\tag 4.3.4.1$$
as $k\ri\iy.$

Really, $ w^*(\bul|\bul)$ is an a.d.p.t.

{\bf Exercise 4.3.4.1.}Check this as in Theorem 4.3.1.1 using that
$w(\bul|\bul)$ is asymptotically dynamical.

For $\tau \in [-\eps,0]$ and $t=r_k$ we have uniformly on $\tau$
$$T_\tau w^*(\bul|t)-w^*(\bul|t)= T_\tau (v_k)_{r_k/t_k} - (v_{k+1})_{r_k/t_k}
\ri 0.$$
Setting $\tau=0$ we obtain 4.3.4.1.

Let $V\sbt U$ be an arbitrary open set, $\eps >0$ arbitrary small and $s>0$
arbitrary large. We show that there exists an $(\eps,s)$-chain from $V$ to
$V.$

Choose $s_1$ such that

i. for $r_k>s_1,\  d( w^*(\bul|r_k-0),w^*(\bul|r_k))<\eps.$

It is possible by virtue of (4.3.4.1).

ii.$w(\bul|s_1)\in V$.
This is possible because of $w(\bul|\bul)$ is $\om$ -dense.

iii. $d(w^*(\bul|t),w(\bul|t))< d(w(\bul|s_1),\p V)$ for $t>s_1$.
This is possible because of Theorem 4.3.3.1.

Choose $s_2>s_1$ such that $w(\bul|s_2)\in V.$ This is possible because
$w(\bul|\bul)$ is  $\om$-dense. Then the pseudo-trajectory
$w^*(\bul|e^t)$ for $s_1\leq e^t\leq s_2$ is an $(\eps,s)$-chain connecting
$w^*(\bul|s_1)$ and $  w^*(\bul|s_2$ that belong to $V.$

{\bf Exercise 4.3.4.1.} Check this in details.

Hence $(T_\bul,U)$ is chain recurrent.
\edm
\subheading {4.3.5}We will prove one more existence theorem that is a corollary
of Theorem 4.2.1.2.
\proclaim {Theorem 4.3.5.1} Let $\Lambda\sbt U[\r]$ be a compact connected and $T_\bullet$--
invariant subset of $U[\r,]$ Then for any proximate order $\r (r)\rightarrow \r$ there exists
$u\in SH(\BR^m,\r,\r (r))$ such that
$$h(x,u)=\sup \{v(x):v\in \Lambda\}\tag 4.3.5.1$$
$$\underline h (x,u)=\inf \{v(x):v\in \Lambda\}.\tag 4.3.5.2$$
\ep \demo {Proof}Let $U:=Conv \Lambda$ be the convex hull of
$\Lambda.$ It is linearly connected and hence polygonally
connected (see 4.1.4). By Theorem 4.1.4.1 it is chain recurrent
and by Theorem 4.2.1.2 for any proximate order $\r (r)\rightarrow
\r$ there exists $u\in SH(\BR^m,\r,\r (r))$ such that
$$\Fr [u,\r (r),V_t,\BR^m]=U.$$
Since every $v\in U$ can be represented in the form $v=av_1+(1-a)v_2$ for
$0\leq a\leq 1,\ v_1,v_2\in \Lambda$ we obtain (4.3.5.1) and (4.3.5.2) from
Theorem 3.2.1.1 (Properties of Indicators),h2).
\edm
\subheading {Exercise 4.3.5.1} Check this.
\subheading {4.3.6} In applications we need the following
\proclaim {Theorem 4.3.6.1} Let $p\in P\sbt \BR^m$ and let $P$ be a connected closed set. Let $U_P:=\{v(z,p):p\in P\sbt \BR^m\}$ be a family
of functions with parameter $p$  such that for every $p\in P$ $v(\bullet,p)\in U[\r]$ and satisfy the condition (4.1.3.3). Then there exists $u\in SH(\r(r)$ such that
 $\Fr [u]=U_P.$
\ep
This is a direct corollary of Theorems 4.1.1.2, 4.1.3.4 and 4.2.1.2.

{\bf Exercise 4.3.6.1} Explain this in details.
\subheading {4.3.7} In the next three \S\S\  we return to the periodic limit
sets (see Th.4.1.7.1). We show that the limit set $\bold {Fr}[u,\rho(r),V_\bul,\Bbb R^m]$ of every subharmonic function
$u\in SH(\r(r),\BR^m) ,\ \r(r)\ri\r$ for  noninteger $\r$ can be approximated in some sense by
periodic limit sets (\cite {Gi(1987)},\cite {GLO,Ch.3,\S2,Th.10}).

Here we give some definitions. Let $X_n \sbt U[\r],n=1,2,...$ be a sequence of compact sets. We say that $X_n$ {\it converges to a compact set} $Y\sbt U[\r],$ i.e.,
$$\Di'-\lim\limits_{n\ri\iy}X_n=Y\tag 4.3.7.1$$
if the following two conditions hold:

converg1) $\forall x_n\in X_n,\ n=1,2,... \, \exists x_{n_j}\in X_{n_j},j=1,2,...$ and
$y\in Y$ such that  $\Di'-\lim\limits_{j\ri\iy} x_{n_j}=y;$

converg2) $\forall  y\in Y \, \exists x_n\in X_n,\ n=1,2,..,$ such that $x_n\ri y.$

On every  compact set  $K$ in $\Di'$- topology one can introduce a metric $d(\bul,\bul)$
such that the topology generated by this metric is equivalent to $\Di'$- topology
(see,e.g.\cite {AG(1982)}).

Denote by
$$X_\eps:=\{ y\in K:\exists x\in X \text { such that } d(y,x)<\eps\}$$
the $\eps-neighborhood$ of $X.$

Let $X,Y$ be two compact sets.  Set
$$d(X,Y):=\inf\{\eps:X\sbt Y_\eps, Y\sbt X_\eps\}.$$

{\bf Exercise 4.3.7.1} Prove the assertion
$$(4.3.7.1)\Longleftrightarrow \{d(X_n,Y)\ri 0\}.\tag 4.3.7.2$$

We prove the following
\proclaim {Theorem 4.3.7.1.(Approximation by Periodic Limit Sets)} Let \newline
$u\in SH(\r(r),\BR^m) ,\ \r(r)\ri\r$ for  noninteger $\r.$ Then for every
$V_\bul$ there exists a sequence $u_n\in SH(\r(r),\BR^m)$ with periodic
limit sets $\bold {Fr}[u_n\bold ,\rho(r),V_\bul,\Bbb R^m]$ such that
$\bold {Fr}[u_n,\bul]\ri\bold {Fr}[u,\bul].$
\ep
This theorem is a corollary of the following
\proclaim {Theorem 4.3.7.2}Let
$\mu\in \Cal M(\r(r),\BR^m) ,\ \r(r)\ri\r$ for  noninteger $\r.$ Then for every
$V_\bul$ there exists a sequence $\mu_n\in SH(\r(r),\BR^m)$ with periodic
limit sets $\bold {Fr}[\mu_n\bold ,\rho(r),V_\bul,\Bbb R^m]$ such that
$\bold {Fr}[\mu_n,\bul]\ri\bold {Fr}[\mu,\bul].$
\ep
\demo {Proof of Theorem 4.3.7.1}The canonical potential $u(x):=\Pi (x,\mu,p)$ (see (2.9.2.1))
of a measure $\mu\in \Cal M(\r(r),\BR^m)$ belongs to $SH(\r(r),\BR^m)$ by Th.2.9.3.3 and has a limit set
$$\bold {Fr}[u ,\bul]=\{\Pi (\bullet ,\nu,p):\nu\in \Fr[\mu_u]\}$$
by Th.3.1.5.2. The potentials $u_n(x):=\Pi (x,\mu_n,p)$ have periodic limit sets
$$\bold {Fr}[u_n ,\bul]=\{\Pi (\bullet ,\nu,p):\nu\in \Fr[\mu_{u_n}]\}$$
by Th.3.1.5.0.
Let us prove that
$$\bold {Fr}[u_n ,\bul]=:X_n\ri Y:=\bold {Fr}[u ,\bul].$$
If $v_n\in \bold {Fr}[u_n ,\bul]$ then from the corresponding sequence of
$\nu_n:=\nu_{v_n}\in \bold {Fr}[\mu_n,\bul]$ we can find a subsequence
$\nu_{n_j}$ and $\nu\in \bold {Fr}[\mu,\bul]$ such that $\nu_n\ri \nu$
(by Th.2.2.3.2 (Helly)). It is easy to check, using Th.3.1.4.3 (*Liouville),
 that $v^\ast=\Di' -\liml_{j\ri\iy}\Pi (\bul,\nu_{n_j})$
exists and coincides with $v=\Pi (\bullet ,\nu,p)\in \Fr [u,\bul].$

So the condition converg1) is verified. In the same way one can check converg2).
\qed
\edm
{\bf Exercise 4.3.7.2} Prove this in details.
\subheading {4.3.8} Now we are going to prove Theorem 4.3.7.2. We begin from
\proclaim {Proposition 4.3.8.1} For any $\mu\in \Cal M(\r(r),\bul)$ there exists
$\hat\mu\in \Cal M(\r,\bul)$ such that
$$\bold {Fr}[\hat\mu ,\r,\bul]=\bold {Fr}[\mu ,\r(r),\bul].\tag 4.3.8.0$$
\ep
 In other words we can suppose further that $\r(r)\equiv \r.$
\demo {Proof} Set $L(r)=r^{\r(r)-\r}$ and
$$\hat\mu (dx):=L^{-1}(|x|)\mu (dx). \tag 4.3.8.1$$
Using properties of proximate order (2.8.1., po1)-po4)), it is easy to check  that

$$[L^{-1}(r)]'=L^{-1}(r)o(1)  \text{ and } [L(r)]'=L(r)o(1),\text{ as }\ri 0.\tag 4.3.8.2 $$

{\bf Exercise 4.3.8.1} Prove this.

Let us show that $\hat\mu\in\Cal M[\r,\Dl]$ for some $\Dl.$ Indeed
$$\frac{\hat\mu (R)}{R^{\r+m-2}}=R^{-\r-m+2}\intl_0^R\frac{\mu (dr)}{L(r)}=$$
$$\frac{\mu (r)}{R^{\r+m-2}L(r)}|^R_0+R^{-\r-m+2}\intl_0^R\mu (r)(L^{-1})'dr=$$
We suppose that $\mu (r)=0$ in some neigborhood of zero.
Using (4.3.8.2) we obtain further
$$=\mu (R)R^{-\r (r)}+R^{-\r-m+2}\intl_0^R\mu (r)(L^{-1})o(1)dr.$$
Using the l'Hospital rule, we obtain
$$\liml_{R\ri\iy}R^{-\r-m+2}\intl_0^R\mu (r)(L^{-1}(r))o(1/r)dr=
(-\r-m+2)\liml_{R\ri\iy}\mu (R)R^{-\r(R)}o(1/R).$$
Thus
$$\limsup\limits_{R\ri\iy}\frac{\hat\mu (R)}{R^{\r+m-2}}\leq
\limsup\limits_{R\ri\iy}\mu (R)(L^{-1}(R))[1+o(1/R)]=\overline \Delta [\mu,\r(r)]<\iy$$
Let us note that $\mu_t=L(t)\hat\mu_{[t]}$.  This implies equality (4.3.8.0) because
 $L(t)\ri 1$ as $t\ri\iy.$

{\bf Exercise 4.3.8.2.} Prove this in details.
\qed
\edm

\demo {Proof of Th.4.3.7.2} As we already said we can suppose that
$\mu\in \Cal M[\r].$ Let $\nu\in\Fr[\mu].$ We can suppose that
$$\nu(\{|x|=1\})=0.\tag 4.3.8.2a$$
 Otherwise  we can find $\tau$ such that
$\nu_{[\tau]}( \{|x|=1\})=0$ and if $\nu_n\ri\nu_\tau$ and are periodic then
$(\nu_n)_{[1/\tau]}$ are also periodic and $(\nu_n)_{[1/\tau]}\ri\nu.$

Let $r_n\ri\iy$ be such that $\mu_{[r_n]}\ri\nu.$ By passing to subsequences
we can make $r_{n+1}/r_n >r_n.$

Denote  $K_n:=\{x:r_n\leq|x|<r_{n+1}\}.$ Set for every $E\sbt K_n$
 $$\mu_n|_E:=\mu (E)$$
and continue it periodically with the period
$P_n=r_{n+1}/r_n$ by the equality
$$\mu_n(P_n^kE)=P_n^{k\r}\mu (E),\ k=\pm 1,\pm 2,...\tag 4.3.8.3 $$
Since every $X\in \Bbb R^m$ can be represented in the form
$$X=\bigcup\limits_{k=-\iy}^{\iy}\{X\cap P^k_nK_n\},$$
  we can define
$$\mu_n (X):=\sum\limits_{k=-\iy}^{\iy}\mu_n(\{X\cap P_k^n K_n\}).$$
It is easy to check that $\mu_n$ is periodic with the period $P_n$ and
$\mu_n\in \Cal M[\r,\Dl]$ with $\Dl$ independent of $n.$

{\bf Exercise 4.3.8.3} Check this.

Let us prove that
$$\Fr[\mu]=\liml_{n\ri\iy}\Fr[\mu_n].\tag 4.3.8.4 $$
Check the condition converg1). Let $\nu_{n_j}\in \Fr [\mu_{n_j}]$ and suppose
 $\Di'-\liml_{j\ri\iy}\nu_{n_j}:=\nu.$ Let us prove that $\nu \in \Fr[\mu ].$

Since  $\Fr [\mu_{n_j}]$ is a periodic limit set
$$\nu_{n_j}=(\mu_{n_j})_{[\tau_j]}.$$
Take $k_j$ such that
$$ \tau'_j:=\tau_jP^{k_j}_{n_j}\in [r_{n_j},r_{n_j+1}).$$
From periodicity $\mu_{n_j}$ we obtain
$$ \nu_{n_j}=(\mu_{n_j})_{[\tau'_j]}.$$
Passing to a subsequence if necessary, we can consider three cases:

i) $\liml_{j\ri\iy}\tau'_j/r_{n_j}=\iy,\liml_{j\ri\iy}\tau'_j/r_{n_j+1}=0;$

ii)$\liml_{j\ri\iy}\tau'_j/r_{n_j}=\tau;\ 1\leq\tau <\iy;$

In this case we have also $\liml_{j\ri\iy}\tau'_j/r_{n_j+1}=0.$

iii)$\liml_{j\ri\iy}\tau'_j/r_{n_j+1}=\tau;\ 0<\tau\leq 1;$

In this case we have also $\liml_{j\ri\iy}\tau'_j/r_{n_j}=\iy.$

Consider the case i). Let $\phi\in \Di (\Bbb R^m\setminus O).$ Then
$\supp \phi(x/\tau'_j)\sbt (r_{n_j},r_{n_j+1})$ for $j\geq j_0.$
It is easy to see that for $j\geq j_0$
$$<(\mu_{n_j})_{[\tau'_j]},\phi>=<\mu_{[\tau'_j]},\phi>.$$

{\bf Exercise 4.3.8.4} Check this.

Since $ \mu_{[\tau'_j]}\ri \nu \in \Fr [\mu]$ by definition
the condition converg1) holds
for the case i).

Consider the case ii).

Recall that $O\notin \supp \phi.$ Then there exists $1\leq c <\iy$ such that
$$\supp \phi\sbt \{x:|x|\in (1/c,c)\}.$$

Define
$$\phi_t(x):=\phi(x/t)(1/t)^\r.$$
Represent $\tau'_j$ in the form
$$\tau'_j:=e_j\tau r_{n_j}$$
where
$$e_j:=\frac {\tau'_j}{r_{n_j}\tau}.$$
The condition ii) means that
$$e_j\ri 1.\tag 4.3.8.4a$$
Compute
$$<\nu_j,\phi>:=<(\mu_{n_j})_{[\tau'_j]},\phi>=<\mu_{n_j},((\phi_\tau)_{e_j})_{r_{n_j}}>.$$

Note that
$$\supp \phi_\tau  \sbt \{x:|x|\in (\tau/c,\tau c )\}.$$
We can increase $c$ so that $1\in (\tau/c,\tau c ).$

Consider the following partition of unity. Choose the functions
$\eta_k\in \Di (\Bbb R^m),\ k=1,2,3$ so that
$$\eta_1(t)+\eta_2(t)+\eta_3(t)=1$$
for $t\geq 1$ and
$$\supp \eta_1\sbt \{x:|x|<1-\eps\},$$
$$\supp \eta_2\sbt \{x:|x|\in(1-2\eps,1+2\eps)\},$$
$$\supp \eta_3\sbt \{x:|x|>1+\eps\},$$
where $\eps$ is an arbitrary number, satisfying
$$\tau/c<1-2\eps<1+2\eps <\tau c.$$
Represent $\phi_\tau$ in the form
$$\phi_\tau =\psi_1 +\psi_2 +\psi_3,$$
where
$$\psi_k=\phi_\tau \eta_k,\ k=1,2,3.$$
In this notation
$$<\nu_j,\phi>=\sum\limits_{k=1}^{3}<\mu_{n_j}, (\psi_k)_{e_jr_{n_j}}>.\tag 4.3.8.5$$
Choose $j_\eps$ such that for $j\geq j_\eps$ the following inclusions hold:
$$\supp (\psi_1)_{e_jr_{n_j}}\sbt \{x:|x|\in (r_{n_j}\tau/c,r_{n_j}(1-\eps))\};$$
$$\supp (\psi_2)_{e_jr_{n_j}}\sbt \{x:|x|\in ((1-\eps)r_{n_j},r_{n_j}(1+\eps))\};$$
$$\supp (\psi_3)_{e_jr_{n_j}}\sbt \{x:|x|\in ((1+\eps)r_{n_j},\tau r_{n_j}\}.$$
Thus for $\psi_3$ we have
$$<\mu_{n_j},(\psi_3)_{e_jr_{n_j}} >=
\int (e_jr_{n_j})^{-\r}\psi_3(|x|/(e_jr_{n_j}))\mu_{n_j}(dx)=$$
$$
\int (e_jr_{n_j})^{-\r}\psi_3(|x|/(e_jr_{n_j}))\mu(dx)=<\mu_{[r_{n_j}]},
(\psi_3)_{e_j}>.$$
Since $\mu_{[r_{n_j}]}{  \overset{\Di'}\to {\ri}} \nu$ and
$(\psi_3)_{e_j}\overset{\Di}\to{\ri} \psi_3$ we have (see Th.2.3.4.6)
$$\liml_{j\ri\iy} <\mu_{n_j},(\psi_3)_{e_jr_{n_j}} >=<\nu,\psi_3>.\tag 4.3.8.6$$

Consider the addend with $\psi_1.$ Because of periodicity $\mu_{n_j}$ we have
$$<\mu_{n_j},(\psi_1)_{{e_j}r_{n_j}}>=<(\mu_{n_j})_{[P_{n_j}]},(\psi_1)_{{e_j}r_{n_j}}>.$$

Transforming the RHS we obtain
$$<(\mu_{n_j})_{[P_{n_j}]},(\psi_1)_{{e_j}r_{n_j}}> =<\mu_{n_j},
((\psi_1)_{P_{n_j}})_{{e_j}r_{n_j}}>.$$
Since $P_{n_j}=r_{n_j+1}/r_{n_j}$ the following inclusion holds for $j\geq j_\eps:$
$$\supp (\psi_1)_{P_{n_j}e_jr_{n_j}}\sbt
\{x:|x|\in (P_{n_j}r_{n_j}\frac {\tau}{c},P_{n_j}r_{n_j}(1-\eps))\}=$$
$$\{x:|x|\in (r_{n_j+1}\frac {\tau}{c},r_{n_j+1}(1-\eps))\}\sbt
\{x:|x|\in (r_{n_j},r_{n_j+1})\}.$$
Thus
$$<\mu_{n_j},((\psi_1)_{P_{n_j}})_{{e_j}r_{n_j}}>=<\mu,(\psi_1)_{P_{n_j}e_jr_{n_j}}>=$$
$$<\mu,(\psi_1)_{{e_j}r_{n_j}}>=<\mu_{[r_{n_j}]},(\psi_1)_{e_j}>.$$
Hence
$$\liml_{j\ri\iy}<\mu_{n_j},(\psi_1)_{{e_j}r_{n_j}}>=<\nu,\psi_1>\tag 4.3.8.7 $$
because $e_j\ri 1$ and $\mu_{[r_{n_j}]}\overset {\Di'}\to{\ri}\nu$ (see Th.2.3.4.6).

From (4.3.8.5), (4.3.8.6) and (4.3.8.7) we obtain
$$\liml_{j\ri\iy}<\nu_j,\phi>=<\nu,\psi_1+\psi_3>+
\liml_{j\ri\iy}<\mu_{n_j},((\psi_2)_{e_jr_{n_j}}>.$$

Let us estimate the last limit. We have
$$\liml_{j\ri\iy}<\mu_{n_j},((\psi_2)_{e_jr_{n_j}}>=
<(\mu_{n_j})_{[e_jr_{n_j}]},\psi_2>$$
Define
$$E_1(\eps):=\{x:|x|\in (1-2\eps,1)\};
E_2(\eps):=\{x:|x|\in [1,1+2\eps)\}$$
Suppose $\eps$ is chosen so that
$$\nu (\p E_k)=0,\ k=1,2.\tag 4.3.8.7a$$
Recall that $\nu$ satisfies the condition (4.3.8.2a), hence $E_1,E_2$ are
$\nu$-squarable and hence (see Th.2.2.3.7)
$$\liml_{n\ri\iy}\mu_{[r_n]}(E_k(\eps))=\nu (E_k(\eps)),\ k=1,2.$$
Define
$$C_\phi:=\max \{\phi(x):x\in \Bbb R^m\}.$$
Then for $j\geq j_\eps$
$$|<(\mu_{n_j})_{[e_jr_{n_j}]},\psi_2>|\leq
C_\phi(\mu_{n_j})_{[e_jr_{n_j}]}(E_1(\eps)\cup E_2(\eps))=$$
$$C_\phi( (\mu_{n_j})_{[e_jr_{n_j}]}(E_1(\eps))+ (\mu_{n_j})_{[e_jr_{n_j}]}(E_2(\eps)).$$
By definition
$$(\mu_{n_j})_{[e_jr_{n_j}]}(E_2(\eps))=\mu_{[e_jr_{n_j}]}(E_2(\eps)).$$
Because of (4.3.8.4a) we obtain
$$\liml_{j\ri\iy}\mu_{[e_jr_{n_j}]}(E_2(\eps))=\nu(E_2)$$
{\bf Exercise 4.3.8.5} Check in details.

To compute the limit of the first addend we use periodicity of $\mu_{n_j}:$
$$\mu_{n_j}(E_1(\eps))=P_{n_j}^{-\r}\mu_{n_j}(P_{n_j}E_1(\eps))=
(\mu_{n_j})_{[P_{n_j}]}(E_1(\eps)),$$
 where
$$r_{n_j}P_{n_j}E_1(\eps)=\{x:|x|\in (r_{n_j+1}(1-2\eps),r_{n_j+1})\}.$$
Thus
$$(\mu_{n_j})_{[e_jr_{n_j}]}(E_1(\eps))=(\mu_{n_j})_{[e_jr_{n_j}P_{n_j}]}(E_1(\eps))=$$
$$(\mu_{n_j})_{[e_jr_{n_j+1}]}(E_1(\eps))=\mu_{[e_jr_{n_j+1}]}(E_1(\eps)).$$
From this we obtain
$$\liml_{j\ri\iy}\mu_{n_j})_{[e_jr_{n_j}]}(E_1(\eps))=\nu (E_1(\eps)).$$
Therefore
$$\liml_{j\ri\iy}|<(\mu_{n_j})_{[e_jr_{n_j}]},\psi_2>|\leq C_\phi\nu(E_1(\eps)\cup E_2(\eps)).\tag 4.3.8.8$$
Because of (4.3.8.2a) we have
$$\nu (\{x:|x|\in (1-2\eps,1+2\eps)\})\ri 0\tag 4.3.8.9$$
as $\eps\ri 0$ over the set of $\eps$ satisfying (4.3.8.7a).
From (4.3.8.8) and (4.3.8.9) we obtain
$$<\liml_{j\ri\iy}<(\mu_{n_j})_{[e_jr_{n_j}]},\psi_2>\ri 0\tag 4.3.8.10$$
and
$$<\nu,\psi_2>\ri 0\tag 4.3.8.11$$
when $\eps\ri 0.$ Hence if $\eps\ri 0 $ satisfying (4.3.8.7a) we have
$$\liml_{j\ri\iy}< \nu_j,\phi>-<\nu_{[\tau]},\phi>=$$
$$\liml_{\eps\ri 0}[<\nu,\psi_1 +\psi_3>+ \liml_{j\ri\iy}<(\mu_{n_j})_{[e_jr_{n_j}]},\psi_2>-<\nu, \psi_1+\psi_2+\psi_3>]=$$
$$\liml_{\eps\ri 0}[\liml_{j\ri\iy}<(\mu_{n_j})_{[e_jr_{n_j}]},\psi_2>-<\nu, \psi_2>]=0$$
The last equality holds because every addend tends to zero.

The case iii) can be considered in an analogous way.

{\bf Exercise 4.3.8.6.} Consider it.

Thus the condition converg1) was checked.
\qed
\edm
\subheading {4.3.9} Now we should check the condition converg2).
We need
\proclaim {Lemma 4.3.9.1}Let $\mu\in \Cal M[\r],\ \nu\in \Fr[\mu]$ and
$r_n\ri\iy,\ n=1,2,..$ is  a sequence such that
$$\Di'-\liml_{n\ri\iy}\mu{[r_n]}=\nu_0.\tag 4.3.9.1$$
Then passing if necessary to a subsequence, we can find $\{r_n\}$ such that
for arbitrarily  $\nu\in\Fr[\mu]$ a sequence $t_j\ri\iy$ exists such that
$$\Di'-\liml_{j\ri\iy}\mu{[t_j]}=\nu.\tag 4.3.9.2$$
and for every  $n$ we can find $t_j\in[r_n,r_{n+1}].$
\ep
 \demo {Proof} Note that if the assertion of the Lemma is satisfied for the
sequence $\{r_n,\ n=1,2,...\}$ it is satisfied for every subsequence of
 $\{r_n,\ n=1,2,...\}.$

Let $M$ be a countable set that is dense in $\Fr [\mu].$ Since reduction $\Di'$ -
topology on $\Cal M[\r]$ is metrizable, it is sufficient to prove that we can
choose a subsequence $r_n$ for which assertion of the lemma is satisfied for
 all $\nu\in M. $ We can do it using a diagonal process.

Let $r^0_n\ri\iy$ be an arbitrary sequence such that
$$\Di'-\liml_{n\ri\iy}\mu_{[r^0_n]}=\nu_0$$
and let a sequence $t^1_j$ satisfies the condition
$$\Di'-\liml_{j\ri\iy}\mu_{[t^1_j]}=\nu_1.$$
Omitting in the sequence $\{r^0_n\}$ the ends of the segments $[r^0_n,r^0_{n+1}] $
that do not contain elements $\{t_j\}$ we obtain a subsequence
$\{r^1_n\}.$ Continuing in such way we obtain a subsequence $\{r^m_n\}$
satisfying  (4.3.9.1) and the sequences
$\{t^1_j\},\{t^2_j\}...\{t^m_j\},j= 1,2...$ satisfying
 $$\Di'-\liml_{j\ri\iy}\mu_{[t^l_j]}=\nu_l,\ l=1,2,...m.\tag 4.3.9.3$$
Taking a diagonal sequence $\{r^n_n\},\ n=1,2,...$ we
observe that it is a subsequence of every subsequence $\{r^m_n\}$ and hence satisfies
the assertion of the lemma.
\qed
\edm
\demo {Proof of converg2)} We can suppose that $\{ r_n\}$ from
the construction of $\mu_n$ with periodic limit sets satisfies the assertion of Lemma 4.3.9.1.
Let $\nu\in \Fr[\mu]$  and $\mu_{t_j}\ri\nu$ under condition
$t_j\in [r_n,r_{n+1}].$ We should consider as in the proof
of converg1) three cases i),ii) and iii). But all these cases were
already considered and hence it was proved that
$$\nu_n:=(\mu_n)_{r_n}\ri \nu.$$
{\bf Exercise 4.3.9.1.}Check this.

\qed
\edm

\newpage

\centerline {\bf 4.4. Subharmonic curves. Curves with prescribed  limit sets.}

\subheading {4.4.1}In this  paragraph we consider  subharmonic functions
$u\in SH(\r(r))$ in the plane of finite type with respect to some proximity order
$\r(r)\ri \r.$

The pair $\bold u:=(u_1,u_2),\ u_1,u_2 \in SH(\r(r))$ is called a {\it subharmonic
curve } or short by {\it curve}.

The family
$$(\bold u)_t:=((u_1)_t,(u_2)_t)$$
is precompact in the topology of convergence in $\Di '$ -topology on every component. The set of all limits
$$\Fr [\bold u]:=\{\bold v=(v_1,v_2):\exists t_j\ri\iy ,\bold v=\Di '-\liml_{j\ri\iy}
\bold u _{t_j} \}$$
is called the limit set of the curve $\bold u.$

Actually this set describes coordinated asymptotic behavior of pairs of subharmonic functions.

\proclaim {Theorem 4.4.1.1} $\Fr [\bold u]$ is closed, connected, invariant with respect to $(\bullet)_{[t]}$ (see 3.1.2.4a) and is contained in the set
$$\bold U [\r,\bold\s]:=\{\bold v =(v_1,v_2): v_n(z)\leq \s_n|z|^\r, v_n(0)=0, n=1,2.\}$$
where $\bold\sigma :=(\s_1,\s_2)$
\ep
\subheading {Exercise 4.4.1.1} Prove this by using Th.3.1.2.2.

Let us define $\bold \s >0$ as $\s_n >0,\ n=1,2.$ Set
$$\bold U[\r]:=\bigcup_{\bold \s > 0}\bold U [\r,\bold \s]$$
We will write $ \bold U\sbt \bold U[\r]$ if $\bold U \sbt \bold U [\r,\bold\s]$ for some $\bold \s.$

Since $(T_\bullet, \bold U[\r])$ is a dynamical system we have two theorems
analogous to Theorems 4.2.1.1 and 4.2.1.2.
\subheading {Exercise 4.4.1.2} Formulate and prove this theorems.

All the other assertions and definitions of \S\S 4.2,4.3 can be repeated for
subharmonic curves.

Let $\bold U\sbt \bold U[\r].$ Set
$$ U':=\{v':\exists v'':(v',v'')\in \bold U\}.$$
This is a projection of $\bold U.$ Set for $v'\in U'$
$$U''(v'):=\{v'':(v',v'')\in \bold U\}.$$
This is the fibre over $v'.$
\proclaim {Theorem 4.4.1.2} Let $\bold U\Subset \bold U[\r]$ be closed and invariant and assume that every fiber $U''(v')$ is convex.Let $U'=\Fr [u']$ for some
$u'\in U(\r(r)).$Then there exists $u''\in U(\r(r))$ such that
$\Fr (u',u'')=\bold U.$
\ep

We construct a pseudo -trajectory asymptotics in the form 4.2.5.2
replacing $u$ with $\bold u$ and $v$ with $\bold v.$ We can
directly check that this curve satisfies the assertion of the
Theorem.

{\bf Exercise 4.4.1.3} Check this.

\proclaim { Theorem 4.4.1.3.(Concordance Theorem)} Let $u\in U(\r(r))$ and $v^0\in \Fr [u],$  and suppose
$v\in U[\r]$ has the property
$$\liml_{\tau\ri -\iy}T_\tau v=\liml_{\tau\ri +\iy}T_\tau v=\tilde v$$
Then there exists a function $w\in U(\r(r))$ such that the limit set of the curve $\bold u=(u,w)$  $\Fr [\bold u]= (\Fr [u],\BC(v))$ and for every sequence $t_n\ri\iy$ such that $\liml_{n\ri\iy} w_{t_n}=v$
$$\liml_{n\ri\iy}\bold u_{t_n} =(v^0,v)\tag 4.4.1.3.$$
\ep
For proving this theorem we should use a.d.p.t. (4.2.5.2). If $v_j=v^0$ we replace $v_j$ by $\bold v_j:=(v^0,v).$ If $v_j\neq v^0$ we replace $v_j$ by $\bold v_j:=(v_j,\tilde v).$

\subheading {Exercise 4.4.1.3} Do that and exploit Th.4.3.1.2 and Th.4.2.1.2.

\proclaim {Corollary 4.4.1.4} Under conditions of Th.4.4.1.3, if
$\liml_{n\ri\iy}w_{t_n}=T_\tau v,$ then $\liml_{n\ri\iy}u_{t_n}=T_\tau v^0.$
\ep
We should apply $T_\tau$ to (4.4.1.3) and use its continuity  in $\Di'$-topology.

\qed

\newpage
\centerline {\bf  5. Applications to Entire functions}
\vskip .15in
\NoBlackBoxes
\centerline {\bf 5.1.Growth characteristics of entire  functions }
\baselineskip 20pt
\vskip .10in
\subheading {5.1.1}Let $f(z)$ be an entire function. The function $u(z):=\log |f(z)|$ is subharmonic in $\BR ^2 (=\BC).$
Hence the scale of growth subharmonic functions considered in
\S 2.8 is transferred  completely to entire functions.
We will mark passing to entire function by changing index
$u$ for index $f.$ For example,
$$M(r,f):=M(r,\log |f|), T(r,f):=T(r,\log |f|).$$
If $u(z):=\log |f(z)|$ has
 order $\r[u]=\r,$ then  $f(z)$ has order
$\r[f]:=\r$  and so on.

We will write $f\in A(\r,\r(r))$ and say ``$f$
{\it is an entire function of order $\r$ and normal type with respect to
proximate order} $\r(r)$'' if $\log |f|$ is a subharmonic
function of order $\r$ and normal type with respect to
the same proximate order. Shortly, if
$\log |f|\in SH(\r,\r(r),\BR ^2)$ then $f\in A(\r,\r(r)).$
\subheading {Exercise 5.1.1.1} Give definitions of
$$T(r,f),\ M(r,f),\ \r_T [f],\  \r_M [f],\ \s_T[f,\r(r)],\
\s_M[f,\r(r)]$$ and reformulate all the assertions of \S 2.8
in terms of entire and meromorphic functions.
\subheading {5.1.2}A divisor of zeros of an entire function
can be represented as an integer mass distribution $n$ on
a discrete set $\{z_j\}\sbt \BC.$ The multiplicity of a zero
$z_j$ is the mass concentrated at the point $z_j.$

The notation for characteristics of the behavior of zeros will mimic that of he behavior of the behavior of masses, replacing
$\mu$ for $n.$ For example, $n(K_r),n(r)$ is the number
of zeros (with multiplicities) in the disk $K_r,$
$\r[n]$ is the convergence exponent, $\vdelt [n]$ is the upper
density and so on.
\subheading {Exercise 5.1.1.2} Give definitions of
$N(r,n),\ \r_N [n],\ \vdelt_N [n],\ p[n].$
\subheading {5.1.3}The limit set $\Fr [f]$ of an entire
function
$f\in A(\r,\r(r))$ is defined as the limit set of the
subharmonic function
$u(z):=\log |f(z)|\in SH(\r,\r(r),\BR ^2)$ (see \S 3.1),
i.e.,
$$\Fr [f]:=\Fr [\log |f|].\tag 5.1.2.1$$
It possesses, of course, all the  properties described in
Ch.'s 3,4 but it is not clear now if there
exists an
{\it entire} function with prescribed limit set, i.e.,
whether the subharmonic function in Theorem 4.2.1.2 can be chosen to be $\log |f(z)|$ where
$f\in A(\r,\r(r)).$ It turns out that this is possible and we prove this    in \S 5.3.

As it was mention in 3.1.1 the general form of
$V_\bullet$ for the case of the plane
 is
$$V_t z=ze^ {i\gamma \log t},$$
where $\gamma$ is real.

The limit set $\Fr [n]$ of a divisor $n$ is the limit set of
the corresponding mass distribution $n$ (see 3.1.3).

Of course generally speaking $n_t$ (see (3.1.3.2)) is not an integer mass
distribution.
\subheading {Exercise 5.1.3.1} Give a complete definition
of $\Fr [f]$ and $\Fr [n],$
and reformulate all the theorems of \S \S 3.1.2, 3.1.3 in
terms of entire functions and their zeros.

The connection between  $\Fr [f]$ and  $\Fr [n]$ is preserved
completely (see \S 3.1.5).
\subheading {Exercise 5.1.3.2}Reformulate the theorems of
\S 3.1.5 for entire functions.

\subheading {5.1.4}Let $f=f_1/f_2$ be a meromorphic function, where $f_1,f_2$ have no common zeros.
If $f_2(0)=1,\ f_1(0)\neq 0$ and
$f_1,f_2\in A(\r,\r(r)),$ then
$u:=\log|f_1|-\log|f_2|\in \dl SH(\r,\r(r)),$
and we write $f\in Mer(\r,\r(r))$ and say
``$f$ is a {\it  meromorphic function of order
 $\r$ and normal type with respect to the proximate
order} $\r(r).$ For $f\in Mer(\r,\r(r))$ we use the following characteristics: $T(r,f), \r_T[f],
\s_T[f,\r(r)].$ The charge of $\log|f|$ consists of integer positive and negative masses.

\newpage
\subheading {5.2. $\Di'$ -topology and Topology of
exceptional sets }
\subheading {5.2.1}Let $\a-\mes$ be the Carleson measure
defined in 2.5.4. Set for $C\sbt \BR ^2$
$$\a-\overline {\mes}C:=\limsup\limits_{R\ri\iy}
[\a-\mes (C\cap K_R)]R^{-\a}.\tag 5.2.1.1$$
It is called the {\it relative} Carleson $\a$-measure.
\proclaim {Theorem 5.2.1.1 (Properties of the Relative
Carleson Measure)}
One has

rCm1) If $C$ is bounded $\a-\overline {\mes}C=0;$

rCm2)
$$\a-\overline {\mes}(C_1\cap C_2)\leq
\a-\overline {\mes}C_1 + \a-\overline {\mes}C_2,$$
 i.e., the relative Carleson measure is sub-additive;

rCm3)
$$C_1\sbt C_2\Rightarrow \a-\overline {\mes}(C_1)\leq
 \a-\overline {\mes}C_2,$$
 i.e., the relative Carleson measure is monotonic
with respect to sets;

rCm4)
$$\a_1>\a_2\Rightarrow \a_1-\overline {\mes}C\leq
 \a_2-\overline {\mes}C,$$
 i.e., the relative Carleson measure is monotonic with respect to $\a.$
\endproclaim

{\bf Exercise 5.2.1.1.} Prove this.

A set $C\sbt\BR^2$ for which $\a-\overline {\mes}C=0$ is
called a $C_0^\a -set.$ If $\a-\overline {\mes}C=0$ for all
$\a>0,$ $C$ is called a $C_0^0 -set.$

Let us recall that if $u_1,u_2\in SH(\r,\r(r),\BR ^2),$ then
$u=u_1-u_2 \in \dl SH(\r,\r(r),\BR ^2)$ (see 2.8.2).
\proclaim {Theorem 5.2.1.2 ($\Di'$ topology and Exceptional
sets)}
Let $u\in \dl SH(\r,\r(r),\BR ^2)$ In order that
$$u_t\ri 0\tag 5.2.1.2 $$
 in $\Di'$ as $t\ri \iy$ it is sufficient that
$$u(z)|z|^{-\r(|z|)}\ri  0\tag 5.2.1.3$$
as $z\ri\iy$ outside some $C_0^2 -set.$

If (5.2.1.2) holds, then (5.2.1.3) holds outside some
$C_0^0 -set.$
\ep
\subheading {5.2.2} To prove Theorem 5.2.1.2 we need
some auxiliary assertions. Recall that $dz$ is an element
of area following the notation of the previous
chapters.
\proclaim {Proposition 5.2.2.1}Let $u\in SH(\r,\r(r),\BR ^2),$
and  $C_{0,R}^2:= C_0^2\cap K_R.$ Then
$$\intl_{C_{0,R}^2}|u|(z)dz=o(R^{\r(R)+2})\tag 5.2.2.1$$
as $R\ri\iy.$
\ep
\demo {Proof} Suppose (5.2.2.1) does not hold. Then there
exists a sequence $R_j\ri\iy$ such that
$$\liml_{R_j\ri\iy}R_j^{-\r(R_j)-2}\intl_{C_{0,R_j}^2}
|u|(z)dz
=A>0\tag 5.2.2.2$$
Consider  the following family of $\dl$-subharmonic
functions:
$$u_j(\z):= R_j^{-\r(R_j)}u(\z R_j)\tag 5.2.2.2a$$
It can be represented as a difference $u_j=u_{1,j}-u_{2,j}$
of subharmonic functions of the same form.

Thus it is precompact in $L_{loc}$ (Theorem 2.7.1.3). Let us
choose a convergent subsequence for which we keep the same
notation. Its limit $v$ is a locally summable function.

Now let $\chi_j$ be the characteristic functions of the
sets
$$E_j:=R_j^{-1}C^2_{0,R_j}.$$
Since $\mes E_j\ri 0$ it is possible to choose a sequence
(for which we keep the same notation) such that
$\chi_j\ri 0$ almost everywhere. We  will also suppose that
$R_j$ are the same for $\chi_j$ and $u_j.$ Thus
$$\intl_{|\z|\leq 1}|\chi_j (\z)u_j(\z)-0\cdot v(\z)|d\z=
\intl_{|\z|\leq 1}|\chi_j (\z)u_j(\z)|d\z\ri 0.$$
By change of variables $z=R_j\z$ we obtain that
$$R_j^{-\r(R_j)-2}\intl_{C_{0,R_j}^2}
|u|(z)dz =\intl_{|\z|\leq 1}
|\chi_j (\z)u_j(\z)|d\z\ri 0.$$
Hence the limit in (5.2.2.2) is equal to zero. Contradiction.\qed
\edm
\proclaim {Proposition 5.2.2.2} Under condition (5.2.1.2)
the set
$$C:=\{z:|u(z)||z|^{-\r(|z|)}>\eps\}$$
is a $C^0_0-set$ for arbitrary $\eps.$
\ep
\demo{Proof}Assume the contrary; that is ,$\exists \a>0$
such that
$$\a-\overline{\mes}C=2\dl>0.\tag 5.2.2.3$$
One can see that for some $\eta >0$
$$\limsupl_{R\ri\iy}(\a-\mes K_{\eta R})R^{-\a}\leq \dl/2
.\tag 5.2.2.4$$
{\bf Exercise 5.2.2.1.} Check this.

(5.2.2.3) and (5.2.2.4) imply that there exists a sequence
$R_j\ri\iy$ such that
$$\liml_{R_j\ri\iy}\a-\mes[C\cap(K_{R_j}\setminus
K_{\eta R_j})]R_j^{-\a}\geq \frac {3}{2}\dl.$$
Set
$$E_j:=R_j^{-1}C\cap (K_{R_j}\setminus
K_{\eta R_j}).$$
It is clear that $E_j\sbt K_1\setminus K_\eta$ and for
sufficiently large $j$
$$\a-\mes E_j\geq\dl.\tag 5.2.2.5$$
Set $u_j$ as in (5.2.2.2a). We claim that for large $j$ and
$\z\in E_j$
$$|u_j|(\z)\geq \frac {\eps}{2}|\z|^\r.\tag 5.2.2.6$$
Indeed,
$$|u_j|(\z)=\frac {|u|(R_j\z)}{R_j^{\r(R_j)}}=
\frac {|u|(z)}{|z|^{\r(|z|)}}(1+o(1))|\z|^\r\geq \frac
{\eps}{2}|\z|^\r.$$ We used here properties of the proximate order
and the equivalence
$$z=R_j\z\in C\cap (K_{R_j}\setminus
K_{\eta R_j})\Leftrightarrow \z\in E_j.$$
{\bf Exercise 5.2.2.2.} Check this in details.

Now we will show that the condition (5.2.1.2) contradicts
 (5.2.2.6). Since $u\in \dl SH(\r,\r(r),\BR ^2)$ it is
a difference of $u_1,u_2\in SH(\r,\r(r),\BR ^2).$ The
corresponding sequences $u_{1,j}$ and $u_{2,j}$ are
precompact in $\Di'$ and there exist subsequences (with
the same notation) that converge to $v_1$ and $v_2,$
respectively.

By Theorem 2.7.5.1 these sequences converge to
$v_1$ and $v_2$ with respect to $\a-\mes$ on
$K_1\setminus K_\eta.$ Since $u_t\ri 0$ in $\Di',$ it
follows that $v_1=v_2.$ Thus $u_j\ri 0$  with respect to $\a-\mes$ on $K_1\setminus K_\eta.$ However, this contradicts
 (5.2.2.5) and (5.2.2.6).\qed
\enddemo
\proclaim {Proposition 5.2.2.3}Let $\{C_j\}_1^\iy$ be a sequence
of $C_0^0$-sets. There exists a sequence $R_j\ri\iy$ such
that the set
$$C=\bigcup\limits_{j=1}^{\iy}\{C_j\cap(K_{R_{j+1}}\setminus
K_{R_j})\}\tag 5.2.2.7$$
is a $C^0_0$-set.
\ep
\demo {Proof} Choose $\eps_j\downarrow 0$ and
$\a_j\downarrow 0.$ Set $R_0:=1.$ Suppose $R_{j-1}$
was already chosen. Take $R_j$ such that
$$\a_j-\mes [C\cap K_{R_{j-1}}]
<\eps_j R^{\a_j} \tag 5.2.2.8 $$
for $R>R_j.$

It is possible because of property rC1) Theorem 5.2.1.1.
We can also increase $R_j$ so that
$$\a_j-\mes [C_j\cap K_R]
<\eps_jR^{\a_j} \tag 5.2.2.9 $$
and
$$\a_j-\mes [C_{j+1}\cap K_R]
<\eps_jR^{\a_j} \tag 5.2.2.10 $$
for $R>R_j.$

It is possible because $C_j$ and $C_{j+1}$ are $C^0_0$-
sets.

Let us estimate $\a_j-\mes [C\cap K_R]$ for
$R_j\leq R<R_{j+1}.$ From (5.2.2.8),(5.2.2.9) and (5.2.2.10)
 we obtain
$$\a_j-\mes [C\cap K_R]\leq
3\eps_jR^{\a_j}. \tag 5.2.2.11 $$
Let $\a>0$ be arbitrary small. Find $\a_j<\a.$
For $R_{j+1}\geq R>R_j$ we have
$$\a-\mes [C\cap K_R]R^{-\a}\leq \a_j-\mes [C\cap K_R]
R^{-\a_j}\leq 3\eps_j.$$
Hence $\a -\overline\mes C=0.$ \qed
\edm
 \subheading {5.2.3}
\demo {Proof of Theorem 5.2.1.2}Let $\phi\in \Di(\BC)$ and
supp $\phi\sbt K_R.$ Then for any $\eps >0$
$$J(t):=\int \phi(z) u_t(z) dz=\left (\intl_{K_R\setminus K_\eps}+
\intl_{ K_\eps}\right )\phi(z) u_t(z)dz:=J_1 (t)+J_2(t)
\tag 5.2.3.1$$
We have for $J_2$ (see 2.8.2.3):
$$|J_2|(t)\leq \max\limits_{|z|\leq\eps}|\phi(z)|\times
const.\intl_0^\eps T(r,|u_t|)rdr\leq const.T(\eps,|u_t|)\eps ^2.
\tag 5.2.3.2$$
Further (see Theorem 2.8.2.1)
$$T(r,|u_t|)\leq 2T(r,u_t)+O(t^{-\r(t)})\leq 2[T(r,u_{1,t})+T(r,u_{2,t})]+O(t^{-\r(t)})\leq$$
$$\leq
2[M(r,u_{1,t})+M(r,u_{2,t})]+O(t^{-\r(t)}\tag 5.2.3.3$$
Using (5.2.3.2), (5.2.3.3) and (3.1.2.3) we obtain
$$\limsup\limits_{t\rightarrow\iy}|J_2(t)|\leq const.
\eps^{\r+2}\tag 5.2.3.4 $$
To estimate $J_1(t)$  write
$$|J_1(t)|\leq const.\left (\intl_{\tilde K_t\setminus
C^2_{0,Rt}}|u(z)|dz +\intl_{C^2_{0,Rt}}|u(z)|dz\right )
t^{-\r(t)-2}:=J_{1,1}(t)+J_{1,2}(t),\tag 5.2.3.5$$
where $\tilde K_t:=\{z:\eps t\leq |z|\leq Rt\}.$

{\bf Exercise 5.2.3.1.} Check this using the change
of variable $z=t\z.$

The summand $J_{1,1}$ is $o(1)$ as $t\rightarrow\iy$ by
(5.2.1.3).

{\bf Exercise 5.2.3.2.} Check this  using the properties of
the proximate order (Theorem 2.8.1.3, ppo3).

The summand $J_{1,2,}=o(1)$ by Theorem 5.2.2.1. Thus
$$\limsup\limits_{t\ri\iy}|J(t)|\leq const.\eps ^{\r+2}$$
for any $\eps.$ Hence it is equal to zero and the sufficiency
of (5.2.1.3) has been proved.

Let us prove sufficiency of (5.2.1.2). Let
$\eps_j\downarrow 0.$ By Theorem 5.2.2.2 we choose
a $C^0_0$-set $C_j$ outside
which $|u(z)||z|^{-\r(|z|)}<\eps_j.$

We construct the set $C$ by  5.2.2.7. Outside $C$
we have (5.2.1.3). And by Theorem 5.2.2.3 it is a
$C^0_0$-set.\qed
\edm
\newpage
\centerline {\bf 5.3.Asymptotic approximation of subharmonic
function} \subheading {5.3.1} One of the widely applied methods of
constructing entire functions with a prescribed asymptotic
behavior is the following: First construct a subharmonic function
behaving asymptotically as he logarithm of modulus of the entire
function,and then approximate it in some sense by the logarithm of
modulus of entire function such that the asymptotic is preserved.

Various requests a precision of the approximation and on the
metric in which such approximation was implemented generated a
specter of theorems of such kind we will demonstrate.

Historically the first theorems of such kind were
proved for concrete functions the masses of which were concentrated
on sufficiently smooth curves (in particularly, on lines, see ,e.g.
\cite {BM},\cite {Ev},\cite {Ki},\cite {Ar}, ...)

In such cases the approximation was very precise and
exceptional sets where the approximation failed were
small and determined.

The first general case was proved in \cite {Az(1969)}. Next
it was developed in \cite {Yu(1982)}, and made excellent in
\cite {Yu(1985)}. It is the following \proclaim {Theorem
5.3.1.1.(Yulmukhametov)} Let $u\in SH(\r)$. Then there exists an
entire function $f$ such that for every $\a\geq\r$
$$|u(z)-\log |f(z)||<C_{\a}\log |z|$$
for $z\notin E_{\a},$ where $E_{\a}$ is an exceptional set that can be covered by discs
$D_{z_j,r_j}:=\{z:|z-z_j|<r_j\}$ satisfying
the condition
$$\sum\limits_{|z_j|>R}r_j=o(R^{\r-\a}),\ R\ri\iy.$$
\ep
This theorem is precise in the following sense:
If
$$||z|-\log |f(z)||=o(\log |z|)
,\ z\ri\iy, \ z\notin E,$$
then for every covering of $E$ by discs $D_{z_j,r_j}$ and every $\eps>0$
$$\sum\limits_{|z_j|<R}r_j\geq R^{1-\eps},
\ R\ri\iy,$$
i.e. in any case this sum is not even bounded.

However it is necessary to remark that the construction from
\cite {Yulmukh.(1985)} `` rigidly'' fastens zeros of the entire
function, whereas the construction of \cite {Az(1969)}
and \cite {Yu(1982)} gives some possibilities
to move them ,which is needed in some constructions.

Let us also mention that such
approximation generates an approximation of a
plurisubharmonic function by logarithm of modulus of entire function in $\BC^p$ (see
\cite {Yu(1996)})

It is also useful to approximate subharmonic functions in an integral metric, for example $L^p,$ as was done in \cite {GG}.

Set
$$\|g\|_p:=\left (\int_0^{2\pi}|g(t)|^pdt
\right )^{1/p}.$$
Denote by $Q(r,u)$ a function that satisfies
the conditions:

1)if $u$ is of finite order, then $Q(r,u)=
O(\log r)$;

2)if $u$ is of infinite order, then
$Q(r,u)=O(\log r +\log \mu_u (r)).$
\proclaim {Theorem 5.3.1.2.(Girnyk,Gol$'$dberg)}
For every subharmonic in $\BC$ function $u$
there exists an entire function $f$ such that
$\|u(re^{i\cdot})-\log |f|(re^{i\cdot})\|_p=Q(r,u).$
\ep

This theorem also considers functions of infinite
order. In this case, it is possible replace $\mu_u (r)$ by
$T(r,u)$ or $M(r,u)$ in $Q(r,u)$ outside an exceptional set
$E\sbt \BR^+$ of finite measure. This theorem
is also unimprovable for subharmonic function
of finite order, because, for example,
$u=\frac {1}{2}\log |z|$ gives, as it is possible to prove:
$$\liminf\limits_{r\ri\iy}
\frac{\|u(re^{i\bul})-\log |f|(re^{i\bul})\|_p}
{\log r}>0$$

However it was find out \cite {LS},
\cite {LM} that  the reminder term
$O(\log|z|)$ that was regarded the best possible is not precise and in some ``regular'' cases can
be replaced with $O(1)$ outside a bigger (but still ``small'')
set .

Set for $E\sbt \BC$:
$$\Dl(E):=\limsup\limits_{r\ri\iy}
\frac{\text {mes} E\cap D_{0,R}}{R^2}.$$
\proclaim {Theorem 5.3.1.3.
(Lyubarskii,Malinnikova)}Let $u$ be a subharmonic
function in $\BC$ with $\mu_u$ satisfying the
conditions: $\mu_u(\BC)=\iy$ and there exists
$\a>0,\ q>1,\ R_0>0$ such that
$$\mu_u(D_{0,qR}\setminus D_{0,R})>\a$$
for all $R>R_0.$

Then there exists an entire function $f$ such that for every $\eps>0$
$$|u(z)-\log |f(z)||<C_\eps$$
for $z\in \BC\setminus E_\eps$ with
$\Dl(E_\eps)<\eps$
\ep
So if $\mu_u$ has no ``Hadamard's gaps'' such approximation is possible.

In this book we restrict ourself to weaker
and  simply proved theorem that is sufficient for
our aim
\proclaim {Theorem 5.3.1.4. (Approximation Theorem)} For every
$u\in SH(\r,\r(r))$
there exists an entire function $f$ such that
$$D'-\lim_{t\ri\iy}(u -\log |f|)_t=0.$$
\ep
Nevertheless this theorem has an important
\proclaim {Corollary 5.3.1.5} For every $u\in SH(\r,\r(r))$
there exists an entire function $f$ such that
 $$\Fr [u]=\Fr [f].$$
\ep

\subheading {5.3.2} Now we  prove Theorem 5.3.1.4. We can suppose, because of Theorem 3.1.6.1
(Dependence $\Fr$ on $V_\bullet$), that in the definition of
$(\bullet)_t$ (see 3.1.2.1) $V_t\equiv I$

 We prove this theorem for the case non-integer $\r.$
For proving this theorem we need
\proclaim {Lemma 5.3.2.1} Let $u\in \delta
SH(\r,\r(r)),\ for\ \r $
non-integer, and $\nu$ is its charge.

Then $u_t\ri 0$ iff $\nu_t\ri 0$ in $\Di'$ as $t\ri\iy.$
\ep
\demo {Proof} Sufficiency.
Suppose $u_t:=(u_1)_t-(u_2)_t\not\ri 0.$ There
exists a subsequence $t_j\ri\iy$ and subharmonic functions
$v_1$ and $v_2$ such that
$$u_{t_j}=(u_1)_{t_j}-(u_2)_{t_j}\ri v_1-v_2 :=v
\neq 0\tag 5.3.2.1$$
Applying to (5.3.2.1) the
continuity of $\Delta$ in $\Di'$ and using the conditions of the theorem,
we obtain
$$\nu_{t_j}\ri \frac {1}{2\pi}\Delta v=0.$$
Hence $v$ is harmonic. Since $v_1,v_2\in U[\r,]$ also
$v\in U[\r]$
(see Th. 2.8.2.1, t3),t4) and Th.2.8.2.3).

{\bf Exercise 5.3.2.1.}Prove this in details.

By Th.3.1.4.3 we obtain $v=0.$  Contradiction.

Necessity. Since the Laplace operator is continuous in
$\Di'$-topology, the assertion $u_t\ri 0$ implies
$\nu _t :=\frac {1}{2\pi}\Delta u_t \ri 0.$
\qed
\edm

Now we describe a construction of the zero distribution of
the future entire function.

Let $u\in SH(\r)$ and $\mu$ be its mass distribution.
Set
$$R_{j+1}:=R_j(j+1)^{4/\k}\tag 5.3.2.2$$
where
$$\k:=\min (\r -[\r],[\r]+1-\r).$$

Let us divide all the plane by circles of the form
$S_{R_j}:=\{|z|=R_j\}$ such that
$R_{j+1}/ R_j\ri\iy$ and $\mu (S_{R_j})=0.$

{\bf Exercise 5.3.2.2.} Prove that it is possible.

Chose a sequence $\dl_j\downarrow 0.$ Divide every annulus
$K_j:=\{z:R_j\leq|z|<R_{j+1}\}$ by circles $S_{R_{j,n}}$
for
$$R_{j,n}:=\left (\frac {1+\dl_j}{1-\dl_j}\right )^n R_j,
\  n=0,1,2,...,n_j,$$
where
$$n_j:=\left [\frac {\log \frac {R_{j+1}}{R_{j}}}
{\log \frac {1+\dl_{j}}{1-\dl_{j}}}\right ],$$
and by rays
$$L_k:={z:\arg z=k\dl_j},\ k=0,1,...,[2\pi/\dl_j].$$
They divide all the plane into sectors $K_{j,n,k}.$ We can choose $\dl_j$ in such way that $\mu(\p K_{j,n,k})=0$ because $ \mu( K_{j,n,k})$ is
monotonic function of $\dl_j$ and have only countable set of jumps.

{\bf Exercise 5.3.2.3.}Explain this in details.

Chose a point $z_{j,n,k}$ in every sector $K_{j,n,k}$ and
concentrate all the mass of the sector at this point.
In other words we consider a new mass distribution
$\hat \mu$ that
has masses concentrated in the points  $z_{j,n,k}$ and
$\hat \mu ({z_{j,n,k}})=\mu (K_{j,n,k}).$

The next lemma shows that $\hat\mu$ is close to $\mu.$
\proclaim {Lemma 5.3.2.3}One has
$$\hat\mu_t -\mu_t\ri 0$$
in $\Di'$ as $t\ri\iy.$
\ep
\demo {Proof}Assume the contrary, i.e.,
$\hat\mu_t -\mu_t\not\ri 0. $ Chose a sequence
$t_l\ri\iy$ such that
$\hat\mu_{t_l}\ri\hat\nu$ and
$\mu_{t_l}\ri\nu,\ \nu,\hat\nu\in \CM[\r],\ \nu\neq\hat\nu.$
Then there exists a disc
$K_{z_0,r_0}:=\{z:|z-z_0|<r_0\}$
such that
$\nu(K_{z_0,r_0})\neq\hat\nu(K_{z_0,r_0}).$
We can assume that this disc does not contain zero
 since for all the $\nu\in \CM[\r]$ the
condition $\nu (K_r)\leq \Delta r^\r, \forall r>0$ is fulfilled.

Suppose, for example,
$$\nu(K_{z_0,r_0})>\hat\nu(K_{z_0,r_0}).\tag 5.3.2.3 $$
Set
$a:=\nu(K_{z_0,r_0})-\hat\nu(K_{z_0,r_0})>0.$
Chose $\eps$ such that
$$\nu(K_{z_0,r_0})<\nu(\overline {K_{z_0,r_0-\eps}})+a/3.
\tag 5.3.2.4$$
This is possible because the countable additivity of $ \hat\nu$
implies
$\liml_{r'\uparrow r}\nu(K_{z_0,r'})=\nu(K_{z_0,r}).$

Consider now the sets
$t_lK_{z_0,r_0}, t_lK_{z_0,r_0-\eps}.$
For sufficiently large $t_l$ they are contained in the union
of the annuluses $K_{j_l}\cup K_{j_l+1}.$

As $j_l\ri\iy$ the diameters of all the sectors $K_{j_l,n,k}$
are $o(R_{j_l})$ uniformly. Thus they are $o(t_l).$
Hence for such $t_l$'s we can find a union $\Gamma_l$ of
sectors covering $t_lK_{z_0,r_0-\eps}$ that does not intersect  the circle of $t_lK_{z_0,r_0}.$

We have $\hat\mu (\Gamma_l)=\mu (\Gamma_l)$ by definition
of $\hat\mu.$ Using the monotonicity of measures, we obtain
$\mu (t_lK_{z_0,r_0-\eps})\leq \hat\mu (t_lK_{z_0,r_0})$
whence
$$\mu _{t_l}(K_{z_0,r_0-\eps})\leq
\hat\mu _{t_l}(K_{z_0,r_0}).$$
Passing to limit as $l\ri\iy$ and using Theorems 2.2.3.1 and
2.3.4.4, we obtain
$\nu(\overline {K_{z_0,r_0-\eps}})\leq
\hat\nu(K_{z_0,r_0}).$
Using (5.3.2.4), we obtain
$\nu(K_{z_0,r_0})-1/3[\nu(K_{z_0,r_0})-\hat\nu(K_{z_0,r_0})]
\leq \hat\nu(K_{z_0,r_0})$
and hence
$\nu(K_{z_0,r_0})\leq\hat\nu(K_{z_0,r_0}),$
that contradicts (5.3.2.3).
Since $\nu$ and $\hat\nu$ are symmetric in this reasoning
the lemma is proved.
\qed
\edm
Let us finish the proof of the Theorem 5.3.1.4
for non-integer $\r.$

We construct a distribution $n$ with  integer masses concentrated at
points $z_{j,k,n}.$ Set
$$n(z_{j,k,n}):=[\hat\mu (z_{j,k,n})]$$
and estimate the growth of the difference
$$\dl\mu:=\hat\mu -n$$
that is also a mass
distribution concentrated at the same points.

Since
$$\dl\mu(z_{j,k,n})\leq 1$$
it is sufficient to count the number of points in the disc $K_R.$

The number of points in the annulus $\{R_j\leq |z|<R\}$ is found from
(5.3.2.2)
$$\dl\mu(\{R_j\leq |z|<R\})\leq \left [\log
\left (\frac {1+\dl_j}{1-\dl_j}\right )\right ]^{-1}\frac {2\pi}{\dl_j}
\log\frac {R}{R_j}\leq $$
$$const\times\frac {\log (j+1)}{\dl_j^2}= const\times (j+1)^4\log (j+1).$$
The mass of the disc $K_R$ is estimated by the inequality
$$\dl\mu (K_R)\leq const\times \sum\limits_{k=0}^{n-1}(k+1)^4\log (k+1)=
o(n^6)=o(R^\eps)\tag 5.3.2.5 $$
for any $\eps>0$ because $R>R_{n-1}=((n-1)!)^{4/\kappa}.$

{\bf Exercise 5.3.2.4.}Check this in details.

The estimate (5.3.2.5) shows that
$$\dl\mu _t\ri 0\tag 5.3.2.6$$
 as $t\ri\iy$.

Lemma 5.3.2.2 and (5.3.2.6) implies that
$$\mu_t-n_t\ri 0\tag 5.3.2.7$$
Set
$$u_1(z):=\Pi (z,n,p)$$
(see (2.9.2.1)) where $\Pi$ is a canonical
potential. This is a subharmonic function in
the plane with integral masses. Thus it is
the logarithm of the modulus of the entire function
$$f(z)=\prod E(z/z_{j,k,n}).$$
 (5.3.2.7) implies by Lemma 5.3.2.1 that
$u_t-(u_1)_t \ri 0$ and this is the assertion of Theorem 5.3.1.4 for non-integer $\r.$
\qed
\newpage

\centerline {\bf 5.4. Lower indicator of A.A.Gol$'$dberg.Description of lower
indicator.}
\centerline {\bf Description of the pair:indicator-lower indicator }
\subheading {5.4.1} Now we consider the {\it lower indicator}. For an entire function of finite order $\r$ and normal type it can be defined in one of the following ways:
$${\underline h}_1(\phi,f):=\sup\limits_{C\in \Cal C}\{\liminf\limits_{
re^{i\phi}\ri\iy,re^{i\phi}\notin C}\log |f(re^{i\phi})|
r^{-\r(r)}\},\tag 5.4.1.1$$
where $\Cal C$ is the set of $C_0$-sets (see \cite {L(1980,Ch.II,\S1}),i.e. the sets that can be covered by union of discs $K_{\dl_j}(z_j):=\{z:|z-z_j|< \dl_j\}$ such that
$$\lim\limits_{R\ri\iy}\frac {1}{R}\sum\limits_{|z_j|<R}\dl_j=0.$$

The exclusion   of $C_0$-sets is necessary because we must exclude from our
consideration some neighborhoods of roots of $f(z)$ where $\log|f(z)|$ is near
$-\iy.$

Similarly, define

$${\underline h}_2(\phi,f):=\sup\limits_{E(\phi)\in \Cal E}\{\liminf\limits_{
 r\ri\iy,r\notin E(\phi)}\log |f(re^{i\phi})|
r^{-\r(r)}\},\tag 5.4.1.2$$
where $\Cal E$ is the set of $E_0$-sets (see \cite {L(1980),Ch.III}), i.e.  sets
$E\sbt (0,\iy]$ satisfying the condition
$$ \lim\limits_{R\ri\iy}\text {mes} \{E\cap(0,R)\}R^{-1}=0.$$
The definition (5.4.1.1) was introduced by A.A.Gol$'$dberg (see \cite {Go(1967}).
We will use the definition (3.2.1.2)
$$\underline h (\phi,f)=\inf \{v(e^{i\phi}):v\in \Fr [f]\}\tag 5.4.1.3$$
Proved in \cite {AP, Theorem 1} that  the definitions (5.4.1.1), (5.4.1.2) and (5.4.1.3) coincide.

Let us note that (5.4.1.3) uses the definition (3.2.1.2) only on the circle
$\{|z|=1\}.$ However, it is easy to check, by using Theorem 3.2.1.2 that for
$\underline h(z)=|z|^\r\underline h(\arg z)$ properties h1) and h2) preserved.

{\bf Exercise 5.4.1.1.} Check this.

We are going to prove
\proclaim {Theorem 5.4.1.1}Let  $g(\phi)$ be a  $2\pi$ -periodic function that is either semicontinuous from above
  or $\equiv -\iy$ and $\r(r)\ri\r$ be an arbitrary approximate order.Then there exists an entire
function $f\in A(\r,\r(r))$ such that
$$\underline h (\phi,f)= g(\phi)\tag 5.4.1.4$$
for all $\phi\in [0,2\pi).$
\ep

\subheading {5.4.2}We will use the following assertion that is a corollary
of Theorem 4.3.5.1 and Corollary 5.3.1.5:

\proclaim {Theorem 5.4.2.1}Let $\Lambda\sbt U[\r]$ be a compact, connected and $T_\bullet$--
invariant subset.

Then for any proximate order $\r (r)\rightarrow \r$ there exists
$f\in A(\r,\r (r))$ such that
$$h(\phi,f)=\sup \{v(e^{i\phi}):v\in \Lambda\}\tag 5.4.2.1$$
$$\underline h (\phi,f)=\inf \{v(e^{i\phi}):v\in \Lambda\}.\tag 5.4.2.2$$
\ep

{\bf Exercise 5.4.2.1.}  Prove Theorem 5.4.2.1.

For the sake of clarity let us restrict ourselves to non-integer
 $\r$.
We will construct a set $\Lambda$ such that
$$\inf \{v(e^{i\phi}):v\in \Lambda\}=g(\phi)$$
Denote
$$H(z,p):=\log |1-z|+\Re \sum\limits_{k=1}^{p}\frac {z^k}{k};\ p=[\r]$$

$$\g(z,K,\lm):=-\lm+K|z-1|,\ \lm,K\geq 0.$$
Note the following properties of these functions:

a)$\min\limits_{|z-1|\geq \dl}\dl H(z,p)|z|^{-\r}\ri 0$, as $\dl\ri 0;$

b)$\dl H(z,p)|z|^{-\r}\leq A\dl,$ for all $z\in \BC,$

where $A$ depends only on $p;$

c)$$\max\limits_{|z-1|\leq 0.5K}\g(z,K,\lm)\leq -\frac {1}{2}\tag 5.4.2.3$$

{\bf Exercise 5.4.2.2} Prove properties a),b),c).

Let us note that $H(1,p)=-\iy.$ Consider the family:
$$\Lambda_{\iy}= \{v_{\theta,\tau}(z):=H(ze^{-i\theta}\tau,p)\tau^{-\r}:\theta\in [0,2\pi),\
\tau\in (0,\iy)\}\cup {0}$$
This family is contained in $U[\r]$ because of b) and closed in $\Cal D$'-topology.
It is also $T_\bullet $-- invariant, hence, satisfies the conditions of Theorem
5.4.2.1.  For every $\phi\in [0,2\pi)$ there exists $\theta_0(=\phi)$, and
$\tau_0(=1)$ such that  $v_{\theta_0,\tau_0}(e^{i\phi})=H(1,p)=-\iy$
Hence
$$\inf \{v(e^{i\phi}):v\in \Lambda_{\iy}\}=-\iy \tag 5.4.2.4$$
For general case this construction will be improved, cutting the ``trunk'' of the function $H(ze^{-i\theta},p)$.

Take $\dl$ small enough so that the following conditions hold
$$ \dl H(z,p)|z|^{-\r}\geq -\frac {1}{4},\ \text {for}\ |z-1|\geq \dl
\tag 5.4.2.5$$
$$\dl H(z,p)\geq  -\frac {1}{4},\ \text{for}\ |z-1|=\dl,\tag 5.4.2.6$$
$$\dl\leq \frac {1}{2K}\tag 5.4.2.7$$
Then
$$\dl H(z,p)> \g (z,K,\lm),\text { for} |z-1|=\dl.\tag 5.4.2.8$$

Denote
$$W(z,K,\dl,\lm):=\cases \max \{\dl H(z,p),\g (z,K,\lm)\}, &\text {for}
\ |z-1|<\dl\\ \dl H(z,p), &\text {for}\ |z-1|\geq \dl\endcases \tag 5.4.2.9$$
\proclaim {Lemma 5.4.2.2}The following holds

aw) The function $W(z,K,\dl,\lm)$ is subharmonic in $\BC;$

bw)supp $ \mu_W \Subset \{|z-1|<\dl\};$

cw) $$\sup\limits_{z\in \BC}W(z,\bullet,\dl,\lm)|z|^{-\r}\leq A\dl, \tag 5.4.2.10$$
where $A$ depends only on $p.$
\ep
\demo {Proof of Lemma} For $|z-1|<\dl$ $W$ is subharmonic as the maximum of two
subharmonic functions. For $|z-1|\geq\dl$ it is harmonic even in the neighborhood of  the circle $|z-1|=\dl,$ because of inequality (5.4.2.8). So aw) and bw)
hold. The assertion cw) follows from b) and c) (5.4.2.3) above.
\edm

Now we get to the proof of (5.4.1.4).  Let $g_n\downarrow g$ a sequence of
continuously differentiable functions that converges to $g$ monotonically. This
is possible,
because $g$ is semicontinuous from above.

{\bf Exercise 5.4.2.3} Prove that  Theorem 2.1.2.9 and the
Weierstrass theorem of approximation of every periodic function by
trigonometrical polioses imply the last assertion.

Denote
$$M_n:=\max\limits_{\phi}g_n^{+}(\phi)$$
where as usual $a^+=\max (a,0)$ .

Set
$$ v_{\theta,n}(z):=W(ze^{-i\theta},K_n,\dl_n,M_n+1-g_n(\theta))+
(M_n+1)|z|^{\r},$$
where $\dl_n$ is chosen small and $K_n$ is chosen large.
Set $z=\tau e^{i\phi}.$ It is clear that
$$ v_{\phi,n}(e^{i\phi})=g_n(\phi)\tag 5.4.2.11 $$
for all $K_n,\dl_n.$

We can choose $K_n$ so large and $\dl_n$ so small that
$$\g(z,K_n,M_n+1-g(\theta))|z|^{-\r}\geq g_n(\phi)$$
for $|z-1|\leq \dl_n, $ because $g_n$ has bounded derivative.

After that we can make $\dl_n$ smaller so that for $|z-1|\geq \dl_n$ the
inequality (5.4.2.8) would hold.

{\bf Exercise 5.4.2.2.} Estimate exactly $K_n$ and $\dl_n$ via the
derivative of $g_n.$

Then
  $$v_{\theta,n}(z)|z|^{-\r}\geq g_n(\phi)$$
for all $z=re^{i\phi}.$
Thus
$$\min\limits_{\theta,\tau}v_{\theta,n}(\tau e^{i\phi})\tau^{-\r}=g_n
(e^{i\phi}),$$
and the minimum is attained for $\tau=1,\theta=\phi.$

Let us note that from (5.4.2.10) we have
$$\sup\limits_{\theta}\sup\limits_{z\in\BC}v_{\theta,n}(z)|z|^{-\r}\leq
A\dl_n+M_n+1\leq A+M_1+1. $$
Consider now the family of functions
$$\Lambda_0:=\{v_{\theta,n}(z\tau)|\tau|^{-\r}:\theta\in [0;2\pi),\ n=1,2,...,
\tau\in (0;\iy)\}$$
It is contained in $U[\r,\s]$ for $\s=A+M_1+1$ and is $T_\bullet$-invariant.
Let $\Lambda$ be its closure in $\Cal D'.$
Let us show that
$$g(\phi)=\inf\{v(e^{i\phi}):v\in \Lambda\}.\tag 5.4.2.12$$
Indeed, for every sequence $v_j\in \Lambda_1$
$$v_j(e^{i\phi})\geq\inf\limits_{n}g_n(\phi)=g(\phi)$$
Let $v\in \Lambda.$ By Theorem 2.7.4.1 ($\Cal D'$ and Quasi-everywhere Convergence)
$$v(z):=(\Cal D '-\lim\limits _{j\ri\iy} v_j)(z)=(\limsup\limits_{j\ri\iy} v_j)^*(z)$$
Hence
$$v(e^{i\phi})\geq g(\phi).$$
However, the infimum is attained for every $\phi$ on the sequence
$v_{\phi,n}(z)$ because of (5.4.2.11). Hence   (5.4.2.12) holds and
Theorem 5.4.1.4 is proved.\qed
\subheading {5.4.3}Now we describe of the pair: indicator-lower
indicator. Let $h$ be a $2\pi$-periodic, $\r$-trigonometrically convex function ($\r$-t.c.f) and let $g$ be a $2\pi$-periodic upper semicontinuous function. Further they are  indicator and lower indicator of entire function, and hence must satisfy
the condition
$$ h(\phi)\geq g(\phi),\ \phi\in [0,2\pi).\tag 5.4.3.1$$
An interval $(a,b)\sbt [0,2\pi)$ is called a {\it maximal interval of $\r$-trigonometricity}
 of the function $h$ if
$$h(\phi)=A\cos\r\phi+B\sin\r\phi,\ \phi\in (a,b)\tag 5.4.3.2$$
 for some  constants $A,B$, and  $h$ has no such representation on any larger
interval $(a',b')\spt(a,b).$

A function $h$ is said to be {\it strictly } $\r$-t.c.f. if it is a $\r$-t.c.f.
and is not $\r$ -trigonometrical on any interval.

If the function $h$ is  a strictly $\r$-t.c.f., then $h$ and $g$ (satisfying other previous bounds) could be an indicator and lower indicator of an entire function $f\in A(\r(r)).$ However this is not so if the function  $h$ has an intervals of trigonometricity.

Recall, for example, the famous M.Cartwright Theorem
\cite {L(1980), Ch.IV,\S2,Th.6} : if an indicator of an
entire function is trigonometrical on an interval $(a,b)$ with $b-a>\pi/\r$ then
the function is a CRG -function on this interval, i.e.,
$$h(\phi)=g(\phi),\ \phi\in (a,b).\tag 5.4.3.3$$

Let us formulate all the necessary conditions of such kind. Let $(a,b)$ be
a maximal interval of $\r$-trigonometricity of the function $h.$ The M.
Cartwright theorem can be formulated as the implication:
$$ (b-a>\pi/\r)\Rightarrow (5.4.3.3).\tag 5.4.3.4$$
The following implications are also necessary:
$$(\exists \phi_0\in (a,b): h(\phi_0)=g(\phi_0))\Rightarrow (5.4.3.3)\tag 5.4.3.5$$
$$(h(a)=g(a)\wedge h'_+(a)=h'_-(a)) \Rightarrow (5.4.3.3)\tag 5.4.3.6a$$
 $$(h(b)=g(b)\wedge h'_+(b)=h'_-(b)) \Rightarrow (5.4.3.3)\tag 5.4.3.6b$$
where $h'_\pm (a)$ and $h'_\pm (b)$ are the right and left derivatives of the
function $h$ at the points $a$ and $b.$
$$(b-a=\pi/\r\wedge h'_+(a)=h'_-(a))\Rightarrow (5.4.3.3)\tag 5.4.3.7a$$
$$(b-a=\pi/\r\wedge h'_+(b)=h'_-(b)\Rightarrow (5.4.3.3)\tag 5.4.3.7b$$
$$\left ( \liminf\limits_{\phi\ri a+0}\frac {h(\phi)-g(\phi)}{\phi -a}=0
\right )\Rightarrow (5.4.3.3)\tag 5.4.3.8a$$
$$\left ( \liminf\limits_{\phi\ri b-0}\frac {h(\phi)-g(\phi)}{b-\phi}=0
\right )\Rightarrow (5.4.3.3)\tag 5.4.3.8b$$

Now we shall  give an exact formulation.
The functions $h$ and $g$ are said to be {\it concordant} if at least one of the following conditions holds:

1.$h$ is strictly $\r$-t.c.;

2.for each $(a,b)$ that is a maximal interval of $\r$-trigonometricity of the
function $h$ the implications  (5.4.3.4)-(5.4.3.8b) are satisfied.

\proclaim {Theorem 5.4.3.1} Let $0<\r<\iy,$ $h(\phi)$ be a $2\pi$-periodic,
$\r$-t.c.f., $g(\phi)$ be an upper semicontinuous, $2\pi$-periodic function,
$h(\phi)\geq g(\phi)$ for all $\phi,$ and $h\not\equiv g.$

A function $f\in A(\r(r))$ which simultaneously satisfies the identity
$h_f\equiv h, \underline h_f\equiv g$ with an arbitrary proximate order
$\r(r)\ri\r$ exists if and only if the functions $h$ and $g$ are concordant.
\ep
\subheading {5.4.4} \demo {Proof of necessity} Note that implication (5.4.3.4) is a corollary of (5.4.3.6a) or (5.4.3.6b), because every $\r$-trigonometrical function is continuous and has continuous derivative in $(a,b).$ Recall that
$(\bul)_{[t]}$ was defined by (3.1.2.4a).

From properties of the limit set $\Fr [f]$ (Theorem 3.1.2.2, fr2),fr3)) and
the definition of indicators ( (3.1.2.1),(3.1.2.2)) we can obtain for every
function $v\in\Fr[f]$ the inequality
$$v(\tau e^{i\phi})\leq \tau^\r h(\phi),\ \phi\in[0,2\pi), \tau >0.\tag 5.4.4.1$$
Since $h(\phi)$ is $\r$-trigonometrical for $\phi\in (a,b)$, the
function $$H(re^{i\phi}):=r^\r h(\phi)$$ is harmonic in the angle
$$\Gamma (a,b):=\{re^{i\phi}:\phi\in (a,b), \ r\in (0,\iy)\},$$
whence the function $v-H$ is subharmonic and nonpositive in $\Gamma (a,b).$
By virtue of the maximum principle, either $v<H$ in $\Gamma (a,b)$ or
$v\equiv H$ in $\Gamma (a,b)$ for each $v\in\Fr[f].$ Note that the condition
$v\equiv H$ in $\Gamma (a,b)$ implies $v\equiv H$ in $\Gamma [a,b]$ for the closed
interval because of the upper semicontinuity of $v.$

Let us  prove (5.4.3.5).For every $v\in \Fr [f]$ we have
$v(re^{i\phi_0})-H(re^{i\phi_0})=0$ whence by the maximum principle
$v= H$ in $\Gamma (a,b).$ Hence (5.4.3.3) holds.

Let us prove (5.4.3.6a). Assume the contrary:$h(a)=g(a)\wedge h'_+(a)=h'_-(a)$ holds, but there exists
$\phi_0\in (a,b)$ such that $h(\phi_0)>g(\phi_0).$ Then there exists
$v\in\Fr [f]$ such that
$$g(\phi_0)\leq v(e^{i\phi_0})<h(\phi_0)$$
whence
$$v(\tau e^{i\phi})<\tau ^\r h(\phi)\in \Gamma (a,b).\tag 5.4.4.2$$
Without loss of generality ,we can assume that $v(z)>-\infty $ otherwise we
can replace $v$ with $\max (v,-C)$ for a large positive constant
$C>0.$

We choose $0<\tau_1<\tau_2$ and to the every function
$$W_j(re^{i\phi}):= v_{[\tau_j]}(re^{i\phi+a})-r^\r h(\phi +a),\ j=1,2, \g=b-a,\ re^{i\phi}\in
\Gamma (0,\g) $$
we apply the following lemma due to A.E.Eremenko and M.L.Sodin \cite {So}(see also \cite {PW},\cite {Ho}):
\proclaim {Lemma 5.4.4.1. (E.S.)} Let $W$ be a subharmonic nonpositive function inside the angle $\Gamma (0,\g)$, $\g>0.$ Then the following implication is valid
$$\left (\limsup\limits_{\phi\ri 0}\frac {W(e^{i\phi})}{\phi}=0\right )\Rightarrow W\equiv 0.$$
\ep
If the condition of this theorem is not satisfied for
$$ W^*(re^{i\phi})=\max\limits_{\tau\in [\tau_1,\tau_2]}v_[\tau](re^{i\phi})$$
it would be possible to insert a $\r$-t.c.function between $h(\phi)-\eps(\phi- a)$ (for a small $\eps$ ) and $v(e^{i\phi})$ .However, such
function does not exist, because of negative jump of derivative. So it will be a contradiction.
See further for details.

From Lemma 5.4.4.1 we get
$$ \liminf\limits_{\phi\ri a+0}\frac {h(\phi)- v_{[\tau_1]}(e^{i\phi})}{\phi -a}:=\a_1>0$$
and likewise
$$ \liminf\limits_{\phi\ri a+0}\frac {h(\phi)- v_{[\tau_2]}(e^{i\phi})}{\phi -a}:=\a_2>0.$$
 So a $\Delta >0$ can be chosen such that $a+\Delta<b$ and the inequalities
 $$H(\tau_je^{i\phi})-v_{[\tau_j]}(e^{i\phi})>\a\tau_j^\r(\phi -a),\ j=1,2,\tag 5.4.4.3$$
  where $\a:=1/2\min (\a_1,\a_2)$, hold for all $\phi\in [a,a+\Delta]$.

  We denote
  $$\beta:= \min\limits_{\tau\in [\tau_1,\tau_2]}(H(\tau e^{i(a+\Delta)})-
  v(\tau e^{i(a+\Delta)}))$$
  which is positive because of (5.4.4.2).

  Let us choose $\eps >0$  small enough to
  $$\eps< \min (\a,\be (\tau _2)^{ -\r } \Delta ^{-1})\tag 5.4.4.4$$
  and let us consider the $\r$-trigonometrical function
$$h_{\eps}(\phi):=\r^{-1}(h'(a)-\eps)\sin\r(\phi -a)+h(a)\cos\r(\phi -a), \ \phi\in (a,b)$$
that coincides with
$$h(\phi)=\r^{-1}h'(a)\sin\r(\phi -a)+h(a)\cos\r(\phi -a), \ \phi\in (a,b)$$
in the point $\phi=a$ but has a tangent that is lower then the tangent of $h.$

Further
$$ h(\phi)-h_\eps(\phi)=\r^{-1}\eps\sin\r(\phi -a)\leq \eps(\phi -a),\ \phi\in [a,a+\Delta]\tag 5.4.4.5$$
Combining (5.4.4.3)-(5.4.4.5) we obtain
$$ v_{[\tau_j]}(e^{i\phi})<{\tau_j}^{\r}h(\phi)-\a(\phi-a)<{\tau_j}^{\r}h(\phi)-\eps(\phi-a)\tag 5.4.4.6$$
$$\leq
{\tau_j}^{\r}h_\eps (\phi),\ \phi\in [a,a+\Delta],\ j=1,2 $$
$$v(\tau e^{i\phi})\leq {\tau}^{\r}h_\eps (a+\Delta)+\tau^{\r}\eps\Delta-\be<\tag 5.4.4.7$$
$${\tau}^{\r}h_\eps (a+\Delta),\ \tau\in[\tau_1,\tau_2]$$
We denote
$$G:=\{re^{i\phi}:\phi\in [a,a+\Delta],\ \tau\in[\tau_1,\tau_2]\}\tag 5.4.4.8$$
It follows from (5.4.4.6),(5.4.4.7) that
$$v(re^{i\phi})<r^{\r}h_\eps(\phi),\  re^{i\phi}\in\partial G,$$
where $\partial G$ is the boundary of the domain $G$. Since the functions $ v(re^{i\phi})$ and
$r^{\r}h_\eps(\phi)$ are subharmonic in $G$, by virtue of the maximum principle we have
$$v(re^{i\phi})<r^{\r}h_\eps(\phi),\ re^{i\phi}\in G.\tag 5.4.4.9$$
Let us consider the function
$$H_1(re^{i\phi}):=r^{\r}h_1(\phi),\  re^{i\phi}\in\Gamma (a-\Delta,a+\Delta)$$
where
$$h_1(\phi):=\cases h(\phi),\ \ \phi\in(a-\Delta,a],\\ h_\eps (\phi),\ \phi\in [a,a+\Delta).\endcases$$
The function $H_1$ is continuous in $\Gamma (a-\Delta,a+\Delta)$ and subharmonic in the angles $\Gamma ( a-\Delta,a)$ and $\Gamma (a, a+\Delta).$ Let us prove that it is subharmonic at
the point $z=e^{ia}.$ Let $\Cal M (z,R,v)$ be the mean value of $v$ over the circle $\{\z:|\z-z|=R\}$
(see  (2.6.1.1)). Taking into consideration (5.4.4.9) and subharmonicity of $v$ (see (2.6.1.1)), for all small $R$ we have
$$\Cal M(e^{ia},R,H_1)\geq M(e^{ia},R,v)\geq v(e^{ia})=H_1(e^{ia}).$$
Hence $H_1$ is subharmonic for $z=e^{ia}.$
Since $H_1$ is homogeneous, i.e. $H_1(kz)=k^\r H_1(z)$,
$$\Cal M(ke^{ia},kR,H_1)=k^\r\Cal M(e^{ia},R,H_1)\geq k^\r H_1(e^{ia})=H_1(ke^{ia})$$
So $H_1$ is subharmonic on the ray $\{z=ke^{ia}:k\in (0,\iy)\}$ and hence in the angle
$\Gamma (a-\Delta,a+\Delta).$ Thus $h_1(\phi)$ is a $\r$-t.c.f. for $\phi\in (a-\Delta,a+\Delta).$ However, by construction
$$(h_1)'_-(a)=h'_-(a)=h'_+(a)=(h_\eps)'_++\eps=(h_1)'(a)+\eps$$
and this contradicts the fact that $h_1$ is $\r$-t.c.f.

Concordance of the implication (5.4.3.6a) is proved.
\edm
\subheading {5.4.5} Here we continue proof of necessity. Pass to the proof of necessity of the condition (5.4.3.7a). Assume the contrary. Then there exists  $v\in \Fr [f]$ and
$\phi _0 \in [a,b]$ such that $g(\phi _0)\leq v(e^{i\phi _0})<h(\phi _0)$, whence by virtue of the maximum principle, $v(\tau e^{i\phi })<\tau ^\r h(\phi)$ for $\tau e^{i\phi}\in \Gamma (a,b).$ Actually $v(\tau e^{i\phi})\leq\tau ^\r h(\phi)$ everywhere and on the circle we have strict inequality. If $v(\tau e^{ia})= H(\tau e^{ia})$ for a $\tau >0$, then $v_{[\tau]} (e^{ia})=h(a)$, and it will suffice to repeat the arguments used in proving (5.4.3.6a) with $v_{[\tau]}$ instead of $v.$

{\bf Exercise 5.4.5.1} Do that.

So it is sufficient to examine the case $v(\tau e^{ia})< H(\tau e^{ia}),\ \tau>0.$
Denote
$$T(\phi):=h'(a)\r^{-1}\sin\r(\phi-a)+h(a)\cos\r(\phi-a).$$
This is a $\r$-trigonometrical function the graph of which is
tangent to the graph of $h(\phi)$ at the point $a.$

There are two possibilities for $T(\phi)$ on some small interval $\phi\in (a-\g,a),\ \g>0:$ either $T(\phi)<h(\phi)$ or  $T(\phi)=h(\phi).$

Inequality $T(\phi)>h(\phi)$ contradicts to $\r$ -t.convexity at the point $a.$
The equality on the sequence of points $\phi_j\ri a-0$ contradicts the maximum principle for $\r$-t.c.functions.

{\bf Exercise 5.4.5.2} Why is it?

If  $T(\phi)=h(\phi),\phi\in (a-\g,a),$ then $h$ is $\r$ -trigonometrical on the
interval $(a-\g,b)\spt (a,b)$ that was already considered in the case (5.4.3.4)
(M.Cartwright's Theorem).

So we assume $T(\phi)<h(\phi),\phi\in (a-\g,a).$ We set
$$h_1(\phi):=h(\phi)-T(\phi), \phi\in (a-\g,a)$$
$$v_1(re^{i\phi}):=v(re^{i\phi})-r^\r T(\phi), re^{i\phi}\in\Gamma (a-\g,b).$$
Then $h_1(\phi)=0$ for $\phi\in [a,b]$, $h_1(\phi)>0$ for $\phi\in (a-\g,a)$ and
$h'(a)=0.$

The function $v_1(e^{i\phi})<0,\phi\in [a,b) .$ Let us analyze the
behavior of the function $v_1(e^{i\phi})$ at the point $b.$ Either
$v_1(e^{ib})<0$ or $v_1(e^{ib})=0$ but $$\limsup_{\phi\ri
b-0}v_1(e^{i\phi})(b-\phi)^{-1}\leq -C$$ for some $C>0$ by Lemma
5.4.4.1 (E.S.)

From the other side $v_1(e^{i\phi})$ is strictly negative also in some left (say, $(a-\Delta,a))$ neighborhood of $a$ because of upper semicontinuity.
In any case $v_1(e^{i\phi})$ can be majorated on the interval $(a-\Delta,b)$ by the function
$$h_\eps:=-A\sin(\r-\eps)(b-\phi)$$
with sufficiently small $A.$

A point of intersection  of the graph of $h_\eps $  with the axis
$0,\phi$ can be regulated by $\eps$ and can be chosen so close to the
point $a$ that the graph of $h_\eps$ also intersect the graph of
$h_1(\phi)$, at some point $\theta_0<a$ because
$h_1(a)=h'_1(a)=0.$

 {\bf Exercise 5.4.4.1.} Make the precise proof with all the
estimates.

Let the parameters $A,\eps,\theta_0$ be fixed as above.
Denote
$$S:=\{re^{i\phi}:\phi\in (\theta_0,b),0<r<1\}.$$
 Then  $H_\eps (re^{i\phi}):=r^{\r-\eps}h_\eps(\phi)$ is harmonic in the
 sector $S$ and satisfies the inequality $H_\eps (re^{i\phi})\geq v_1(re^{i\phi})$
 on $\partial S.$ Hence $H_\eps (re^{i\phi})\geq v_1(re^{i\phi})$ on $S.$ Thus
 $$v(re^{i\phi})\leq H(re^{i\phi})+H_\eps (re^{i\phi}),\  re^{i\phi}\in S.\tag 5.4.5.1$$

 Let $\Cal M(r,v)$ be the mean value of the function on the circle $\{\z:|\z|=r\}$  (see 2.6.1.1). Using (5.4.5.1) we have
 $$\Cal M(r,v)\leq\int\limits_{\theta_0}^b [H(re^{i\phi})+H_\eps (re^{i\phi})]d\phi+
 \int\limits_{[0,2\pi)\backslash (\theta_0,b)}H(re^{i\phi})d\phi\leq d_1r^\r -d_2r^{\r-\eps},\ d_1,d_2 >0$$
 So we get $\Cal M(r,v)<0=v(0)$ for sufficiently small $r>0$ which contradicts  the
subharmonicity of the function $v$ at zero.
\subheading {5.4.6} Now we complete proof of necessity, proving
(5.4.3.8a,b). Assume the contrary:suppose
$$ \liminf\limits_{\phi\ri a+0}\frac {h(\phi)-g(\phi)}{\phi -a}=0\tag 5.4.6.1$$
but there exists a $\phi_0\in (a,b)$ such that $h(\phi_0)>g(\phi_0).$
Then there exists a function $v\in \Fr [f]$ such that
$$ g(\phi_0)\leq v(e^{i\phi_0})<h(\phi_0)\tag 5.4.6.2$$
Then the function $v_1(re^{i\phi}):=v(re^{i\phi})-H(re^{i\phi})$ is subharmonic and nonpositive in $\Gamma (a,b).$ By virtue of the maximum
principle $v_1(re^{i\phi})<0,\ re^{i\phi}\in\Gamma (a,b).$

From (5.4.6.1) we obtain
$$0= \liminf\limits_{\phi\ri a+0}\frac {h(\phi)-g(\phi)}{\phi -a}\geq \liminf\limits_{\phi\ri a+0}\frac {h(\phi)-v(e^{i\phi})}{\phi -a}=-\liminf\limits_{\phi\ri a+0}\frac {v_1(e^{i\phi})}{\phi -a}$$
whence, recollecting that $v_1(e^{i\phi})<0,$ we get
 $$\limsup\limits_{\phi\ri a+0}\frac {v_1(e^{i\phi})}{\phi -a}=0.$$
 Applying Lemma 5.4.4.1 (E.S.) to the function
 $$W(re^{i\phi})=v_1(re^{i\phi + a}),\ re^{i\phi}\in \Gamma (0,\gamma),\ \gamma=b-a$$
 we get $v_1\equiv 0$ in $\Gamma (a,b)$ which leads to a contradiction.
 The implication (5.4.3.8b) is proved in the same way. So the proof of  necessity in  Theorem 5.4.3.1 is completed.
\qed

We do not include here the proof of sufficiency and refer the readers to
the original paper \cite {Po(1992)}.
\newpage

\centerline {\bf 5.5. Asymptotic extremal problems.Semiadditive
integral.} \subheading {5.5.1 } Suppose  some class of entire
functions is determined by asymptotic behavior of their zeros, and
we want to know what is the restriction on  asymptotic behavior of
functions: for example, to estimate indicator of such function.
The first example of such problem was considered by B.Ya. Levin in
\cite {L(1980),Ch.IV,\S1,Example}. A developed theory of such
estimates was constructed in the papers of A.A.Gol$'$dberg \cite
{Go(1962)} and his pupils \cite {Kon(1967)}
 , \cite {KF(1972)}.We consider this theory from the point of view of limit sets.

Let $\Cal M \Subset \Cal M(\r)$ (see (3.1.3.4)) be a convex set of measures which is closed in $\Cal D '$ and is invariant with respect to the transformation $(\bullet)_t$ (see (3.1.3.1),(3.1.3.2))) and let $A(\Cal M)$ be a class of entire functions $f$ for which $\Fr [n_f]\sbt \Cal M.$ We suppose $\r$ is non-integer. Recall that canonical potential $\Pi (z,\nu,p)$ is defined by: (see (2.9.2.1))
$$\Pi (z,\nu,p):=\int\limits_\BC G_p(z/\z)\nu(d\z),$$
where $\nu$ is a measure and
$$G_p(z):=\log |1-z|+\Re \sum\limits_{k=1}^p \frac {z^k}{k}.$$

\proclaim {Theorem 5.5.1.1} \text { \cite {AP(1984)}}
The relation
$$h(\phi,f)=\sup\{ \Pi (e^{i\phi},\nu,p): \nu\in \Cal M\}\tag 5.5.1.1$$
 is valid. There exists $f\in A(\Cal M)$ for which the equality holds in (5.5.1.1) for all $\phi.$
\ep
\demo {Proof}We should only prove that there exists an entire function with
such indicator. Consider the set
$$\Lambda:=\{ \Pi (e^{i\phi},\nu,p): \nu\in \Cal M\}$$
It is a convex set contained in $U[\r].$ Thus there exists a subharmonic (see
Corollary 4.1.4.2) and hence entire (see Corollary 5.3.1.5) function $f$ such
that $\Fr [f]=\Lambda.$ By Theorem 5.4.2.1 (5.5.1.1) holds .
\qed\edm
For some $\Cal M$ it is possible to compute the supremum in (5.5.1.1)
and thus to obtain explicit precise estimates of indicators in the respective class $A(\Cal M).$
As an example, we shall present an estimate given by A.A.Gol$'$dberg.

We recall that the {\it upper density of zeros} of an entire function $f\in A(\r)$ is defined by the equality
$$\overline\Delta[n_f] :=\limsup\limits_{r\ri\iy}\frac {n_f(r)}{r^\r}$$
where $n_f$ is the distribution of zeros of the function $f,$ and denote

$$K(t,\phi):= -[\frac{d}{dt}G^+_p(e^{i\phi}/t)]^- \tag 5.5.1.2$$
where $a^+:=\max (a,0),\ a^-:=\min (a,0).$
This function is piece-wise continuous.
 \proclaim {Corollary  5.5.1.2}\text {\cite {Go(1962)}}Let the distribution of zeros $n_f$ of a function $f$ be concentrated on the positive ray, and
let $\overline\Delta[n_f]\leq \Delta<\iy .$ Then
$$h(\phi,f)\leq \Delta \intl_{0}^\iy t^\r K(t,\phi)dt , \
\phi\in[0,2\pi)\tag 5.5.1.3$$
and there exists a function from the same class for which  equality is attained for all $\phi.$
\ep

\demo {Proof of Corollary} We exploit Theorem 5.5.1.1. The class
of the functions $f$ satisfying the assumption  of the Corollary
coincides with the class of $f$ for which

$$\Fr [n_f]\sbt \Cal M=\{\nu\in \Cal M(\r): \text {supp} \nu \sbt [0,\iy]\wedge \nu(r)\leq \Delta r^\r\}.\tag 5.5.1.4 $$
 C
 {\bf Exercise 5.5.1.1} Show this by using Corollary 3.3.2.6.

Thus
$$ \Pi (e^{i\phi},\nu,p)=\int\limits_0^\iy G_p( e^{i\phi}/t)\nu(dt)\leq
\int\limits_0^\iy G^+_p( e^{i\phi}/t)\nu(dt)$$
Integrating by parts we obtain
$$\Pi (e^{i\phi},\nu,p)\leq
-\int\limits_0^\iy \nu (t)[\frac {d}{dt}G^+_p( e^{i\phi}/t)]^-dt$$
By (5.5.1.4) we get (5.5.1.3).
\qed
\edm

Denote
$$M_p(r):=\max \{G_p(re^{i\phi}):\phi\in [0,2\pi)\}$$
In the same way one can prove
\proclaim {Corollary 5.5.1.3}\text {\cite {Go(1962),Th.4.1}}.Let distribution of zeros of the function $f\in A(\r)$ satisfy only the condition
$\overline\Delta[n_f]\leq \Delta<\iy.$ Then
 $$h(\phi,f)\leq \Delta \r\intl_{0}^\iy t^{\r-1} M_p(1/t)dt , \
\phi\in[0,2\pi)\tag 5.5.1.5$$
and there exists a function from the same class for which  equality is attained for all $\phi.$
\ep
 {\bf Exercise 5.5.1.2.}Prove this Corollary exploiting
$$\Cal M:=\{\nu \in \Cal M(\r):\nu (r)\leq \Delta r^\r ,\ \forall r>0\}.$$
\subheading {5.5.2} To be able to obtain explicit estimates for
more diverse classes of entire functions defined by a restriction
on the density of zeros, Gol$'$dberg introduced an integral with
respect to a non-additive measure and obtained estimates for
indicators in terms of one-dimensional integral (along a
circumference) with respect to such a measure (\cite {Go(1962)}.Gol$'$dberg
initially constructed integral sum of a special
form.The construction presented here is based on the
Levin--Matsaev--Ostrovskii theorem (see (see \cite {Go(1962),Th.2.10}).
Fainberg (1983) developed this approach using a two-dimensional integral.
This made it possible to extend significantly the set of classes of
entire functions for which the estimate expressed by a nonadditive integral
is precise. We shall present these results after  the necessary definitions.

Let $\delta (X)$ be a non-negative monotonic function of
$X\sbt\BC$, the function being finite on bounded sets and $\delta
(\varnothing)=0.$ For a given family of sets $\Cal X:=\{X\}$ we
denote by $N(\delta,\Cal X)$ the class of countable-additive
measures $\mu$ defined by the relation
$$N(\delta,\Cal X):= \{\mu:\mu (X)\leq \delta (X),\ X\in \Cal X\}.$$
For a Borel function $f\geq 0$ we define the quantity

$$(\Cal X)\int fd\delta:=\sup \left\{ \int f d\mu:\mu\in N(\delta,\Cal X)\right\},$$
called an $(\Cal X)-integral$ with respect to a {\it nonnegative measure} $\delta.$ For a Borel set $E\sbt\BC$ we set
      $$(\Cal X)\int_E fd\delta:=(\Cal X)\int fI_Ed\delta,$$
where $I_E$ is an indicator of the set $E,$ i.e.,
$$I_E(z):=\cases 1, &\text{if }z\in E;\\ 0 &\text{if }z\notin E.\endcases$$
This integral possesses a number of natural properties: it is
monotonic with respect to $f$ and $\delta$ and the family $\Cal
X$,positively homogeneous and semi-additive with respect to the
function $f$ and $\delta.$ If $\delta$ is a measure, if $\Cal X$
is a Borel ring, and if $f$ is a measurable function, then ($\Cal
X$)-integral coincides with the Lebesgue -Stieltjes integral.

 {\bf Exercise 5.5.2.1.} Check these properties.

Let $\delta(\Theta)$ be a nonadditive measure on  the unit circle $\BT$, defined initially on the family of all open sets $\Theta\sbt \BT.$ It can be naturally extended to all closed sets $\Theta^F$ using the equality
$$\delta (\Theta^F):=\inf\{\delta(\Theta):\Theta \spt \Theta^F\}.$$
Let $\chi_\Theta$ be a set of open sets containing the set $\BT.$ Denote
$$D_{r,\Theta}:=\{z=te^{i\theta}: 0<t<r,e^{i\theta}\in \Theta\},\
\chi_z:=\{D_{r,\Theta}:r>0,\ \Theta\in \chi_\Theta\}$$
The subscripts $\Theta$ and $z$ at $\chi$ indicate that the families under consideration are located either on $\BT$ or on the plane, respectively.

Let us define a non-additive measure $\delta_z$ on $\chi _z$ by the equalities
$$\delta_z(D_{r,\Theta}):=r^\r\delta(\Theta), \ D_{r,\Theta}\in \chi _z.$$
Now the integral $(\chi _z)\int G^+_p(e^{i\theta}/\z)d\delta_z$ is defined.

Recall that the classical angular upper density of zeros of an entire function $f\in A(\r)$
is defined by the equality (compare (3.3.2.7))
$$\overline \Delta ^{cl}[n_f,\Theta]:=\limsupl_{r\ri\iy}n_f(D_{r,\Theta})r^{-\r}.$$
Consider the class of entire functions $A^{cl}(\delta,\chi _\Theta)$ defined by the equality
$$A^{cl}(\delta,\chi _\Theta):=\{ f: \overline \Delta ^{cl}[n_f,\Theta]\leq
\delta (\theta), \ \forall \Theta \in \chi _\Theta \}\tag 5.5.2.6a$$
for a given non-additive measure $\delta (\Theta)$ and a family $\chi _\Theta.$
\proclaim {Theorem 5.5.2.2}\cite {Fa} Let $\delta (\Theta)$ satisfy the condition

$$\delta (\Theta)=\delta (\overline \Theta) , \ \forall \Theta\in \chi _\Theta\tag 5.5.2.6 $$
(the dash means the closure of a set). Then

$$h(\phi,f)\leq (\chi _z)\int G^+_p(e^{i\theta}/\z)d\delta_z \tag 5.5.2.7$$
There exists a function $f\in A^{cl}(\delta,\chi _\Theta)$ such that  equality in (5.5.2.7) is attained for all $\phi\in [0,2\pi)$ simultaneously.
\ep
\demo {Proof} Let us note the following: If we replace in this theorem

$\overline \Delta ^{cl}[n_f,\Theta ]$ with its $\Cal D'$
counterpart $\vdelt (Co_{\Theta}(I_1))$ (see Theorem 3.3.1.2) and
consider the corresponding class of entire function $
A(\delta,\chi _\Theta)$  the assertion of the theorem holds
without conditions(5.1.5.6). You should only apply Theorem 5.5.1.1
with the corresponding $\Cal M.$ The condition (5.5.2.7) is
exploited only for replacing ``$\Cal D'$'' quantities by the
classic ones using results of \S 3.3.2. \qed \edm {\bf Exercise
5.5.2.2} Prove this theorem in details.

It is also worth  to note that every family $\chi _\Theta$ can be replaced for
a family $\chi '_\Theta$ that is dense in $\chi _\Theta$ (see 3.2.2) and such that for
$\chi '_\Theta$ (5.5.2.6) already holds (see Theorem 3.3.2.3).
\newpage

\centerline {\bf 5.6.Entire functions of completely regular growth. }
\centerline {\bf Levin-Pfluger Theorem. Balashov's theory.}
\subheading {5.6.1}The most famous definition of a function of  {\it completely regular growth }  (CRG-function) is the following:

A function $f\in A(\r(r))$ is a function of completely regular growth , if
the limit
$$\liml_{z\ri\iy} r^{-\r(r)}\log|f(z)|,\  r:=|z|$$
 exists  when $z\ri\iy$ uniformly outside some $C^1_0$ -set (see \S 5.2.1.)

Actually, it is equivalent to all other definitions of the functions
of completely regular growth in the plane (comp.\cite {Le, Ch.III},
\cite {Pf(1938},\cite {Pf(1939}).

By A.A.Gol$'$dberg (\cite {Go(1967)}) this definition was reduced to the following :

A function $f\in A(\r(r))$ is a function of completely regular growth , if
$$\underline h_f(\phi)=h_f(\phi),\ \forall \phi\in [0,2\pi).$$

Because of the formulae (3.2.1.1), (3.2.1.2) (see also \S 3.2.7) we have the following

\proclaim {Theorem 5.6.1.1} A function $f\in A(\r(r))$ is a function of completely regular growth   (CRG -function) iff $\Fr[f]$ consists of only
one subharmonic function $h(z).$
\ep

Because of  (3.2.1.11) the function $h(z)$ has the form
$$h(z)=r^\r h(e^{i\phi}) \tag 5.6.1.1$$
The function $h(\phi):=h(e^{i\phi})$ is $\r$ -trigonometrically convex and it was studied in \S \S 3.2.3, 3.2.4, 3.2.5.

\subheading {5.6.2} The initial definition of {\it regular zero distribution}
\cite {L(1980),Ch.II,\S1} is the following :

Let $n$ be a zero distribution (divisor,or mass distribution) of convergence exponent $\r_1:=\r [n]$ (see \S 2.8.3), and let $\r_1>[\r_1].$ Let $\r_1(r)\ri\r_1$ be a proper proximate order of $n(r)$ (see Th.2.8.1.2). It means that
$n\in \Cal M (\r(r)),\ \r (r)\ri\r_1$ (see \S 3.1.3).

The initial definition of {\it regular zero distribution} for  $\r_1$ being non-integer is:

A zero distribution $n$ is regular if the limit

$$\liml_{r\ri\iy}\frac {n(Co_{(\a,\be) }(I_t))}{t^{\r_1(t)}}:=\Delta ((\a,\be ))$$
\footnote {for definition $Co_{(\a,\be) }(I_t)$ see Exercise 3.3.1.5.}
exists for all $\a>\be$ except may be for a countable set on  the circle.

By using results of \S 3.3, one can prove
\proclaim {Theorem 5.6.2.1}The zero distribution $n$ is regular iff $Fr [n]$
consists only one measure $\nu_{reg}.$
\ep

{\bf Exercise 5.6.2.1.} Prove this exploiting Th's.3.3.3.1 and 3.3.2.4.

Recall that for $f\in A(\r(r)), \r(r)\ri\r$ we have $n_f\in \Cal M(\r(r)),\r(r)\ri\r.$ (see Th.2.9.3.2).  Now we can formulate

\proclaim {Theorem 5.6.2.2 (Levin-Pfluger)}\text{ \cite {L(1980),Ch.II,Ch.III}}
An entire function $f\in A(\r(r), \r(r)\ri\r$ of
non-integer order $\r$ is of completely regular growth function iff its zero distribution is regular.
\ep
After Theorems 5.6.1.1, 5.6.2.1 this theorem is a direct corollary of Th.3.1.5.1.

\subheading {5.6.3} Consider now the case of  integer $\r.$ In general, this case is  differs from the case of non-integer $\r.$ For example,  Th.
2.9.4.2  (Brelot-Lindel\"of) implies that
$$(f\in A(\r(r)), \r(r)\ri\r)\Longleftrightarrow  n_f\in \Cal M(\r(r)),\r(r)\ri\r$$
iff the family of polynomials (2.9.4.4a) is compact.

To describe the regularity of zero distribution for the case of  integer
$\r$ we assume that the limit
$$\liml_{R\ri\iy} \delta _R (z,\nu,\rho):=\Re [\delta_\iy z^\r] \tag 5.6.3.1$$
exists, where
$$\delta_\iy:=\liml_{R\ri\iy}\left[\intl _{|\z|<R}  |\z|^{-\r}\cos \arg \z n (d\z) +
i\intl _{|\z|<R}  |\z|^{-\r}\sin \arg \z n (d\z)\right]$$

Now a zero distribution $n\in \Cal M (\r (r)),\ \r(r)\ri\r$ with {\it integer} $\r$  is called  {\it regular } if  $\Fr[n]$ consists of only one measure
$\nu_{reg}$ as in Theorem 5.6.2.1 and  the limit (5.6.3.1) exists.

Under this definition the Theorem 5.6.2.2 still holds, because the set
$(\Cal H, \Fr)[\log |f|]$ from Theorem 3.1.5.2 consists of only one element
$(\Re [\dl_\iy z^\r], \nu_{reg}).$

Note also
\proclaim {Proposition 5.6.3.1 }The measure $\nu_{reg}$ has the following form
$$\nu_{reg}(drd\phi)=\r r^\r dr\otimes \Delta (d\phi)$$
where $\Delta$ is a measure of bounded variation on the unit circle.
 \ep
This assertion is a corollary of invariance of $\Fr[n]$, Th.3.1.3.3, frm3).
\subheading {5.6.4} In the papers \cite {Bal(1973)}\cite {Bal(1976)} functions of
{\it completely regular growth  along curves of regular rotation} were considered.
A {\it curve of regular rotation} is a curve that is described by the equation
$$z=te^{i(\g(t)\log t +\phi)}, 0<t<\iy$$
for a fixed $\phi.$

If $\g(t)\equiv \g$ then this curve is a logarithmic spiral. In general case
$\g (t)$  is a differentiable function such that
$$\g (t)\ri \g, t\g'(t)\ri 0,\  t\ri\iy.$$
To describe this theory in terms of limit sets we consider the transformation    $$P_t=te^{i\g(t)\log t}$$
$$u_t(z)=u(P_t z)t^{-\r(t)}$$
The following Theorem is similar to Th.3.1.2.1 \proclaim
{Proposition 5.6.4.1(Existence of spiral Limit Set)}The following
holds

esls 1) $u_t\in SH(\r(r))$ for all $t\in (0,\iy);$

esls 2) the family $\{u_t\}$ is precompact at infinity.
\ep
The set of all limits $\Di'-\liml_{j\ri\iy}u_{t_j}$ does not depend on $\g(t)$
but only on the constant $\g$ since
$$\liml_{t\ri\iy}(\g(t)-\g)\log t=\liml_{t\ri\iy}t\g'(t)=0$$
So it is the same that for
$$P_t=te^{i\g\log t},$$
i.e., the case that was already considered in the general theory.

In particular (3.2.1.8) for this case has the form
$$z^0(z)=e^{i(-\g\log r +\phi)}\tag 5.6.5.1$$
Hence , from Theorem 3.2.1.2 the indicator  (see (3.2.1.1)) has the form
$$h(re^{i\phi})=r^\r h(-\g\log r +\phi),\ z=re^{i\phi}, $$
 where $h(\phi)$ is a $\r$-trigonometrically convex $2\pi$-periodic function (see \S 3.2.3).

All other assertion of Levin-Pfluger theory can be obtained analogously from other general assertions as it was done in the previous \S\S.

Theorem 3.1.6.1 connects limit sets for every $\g.$

 {\bf Exercise 5.6.1.1.} Formulate and prove Balashov's analogy of
the Levin-Pluger Theorem 5.6.2.2 and Theorem 3.1.6.1 for $m=2.$

For other generalization of the Levin-Pfluger theory see \cite {AD}
\newpage

\centerline {\bf 5.7.General characteristics of growth of entire functions.}

\subheading {5.7.1} A  functional $\Cal F(u)$ acting in the unit circle and defined on
subharmonic functions $u\in SH(\r (r))$ is called a {\it growth characteristic}
if the following conditions are fulfilled:

1. {\it continuity}:
$$\Cal F(u_j)\ri \Cal F(u),\tag 5.7.1.1$$
if $u_j\ri u$ uniformly on
compacts (of course, for continuous functions $u$) or if $u_j\downarrow u$;

2. {\it positive homogeneity}:
$$\Cal F(cu)=c\Cal F(r,u); \tag 5.7.1.2$$
for every constant $c>0.$

Here we shall list some widely used functionals that satisfy these conditions:

$$H_\phi (u):=u(e^{i\phi});\tag 5.7.1.4$$
$$T(u)=\frac {1}{2\pi}\intl_0^{2\pi}u^+(e^{i\phi})d\phi;\tag 5.7.1.5$$
$$M_\a (u):=\max \{u(e^{i\phi}):|\phi|\leq \a\}\tag 5.7.1.6$$
$$ M(u):=M_\pi (u);\tag 5.7.1.7$$
$$ I_{\a\be}(u):=\intl_\a^\be u(e^{i\phi})d\phi \tag 5.7.1.8$$
$$I(u,g):=\intl_0^{2\pi} u(e^{i\phi})g(\phi)d\phi,\ g\in L^1[0,2\pi].\tag 5.7.1.9$$

 {\bf Exercise 5.7.1.1.} Check properties 1. and 2. for these functionals.

Let $\a (t)$ and $\a_\eps(\z)$ be the ``hats'' defined by the equalities
 (2.3.1.1)-(2.3.1.3) and let $R_\eps u$ be defined by (2.3.1.4).

This averaging has the following properties.

\proclaim {Proposition 5.7.1.1}

1. if $u$ is subharmonic, then $R_\eps u$ is subharmonic;

2. $R_\eps u\downarrow u $ as $ \eps \downarrow 0$ for every subharmonic function;

3. if $u_j\ri u$ in $\Di '$ and $u_j,u$ are locally summable functions,
$R_\eps u_j\ri R_\eps u$ uniformly on every compact set.
\ep

 {\bf Exercise 5.7.1.2.} Prove this  using Th.2.3.4.5, 2.6.2.3.

Now we can define  the {\it asymptotic characteristics of growth} of entire function $f\in A(\r (r)):$
 $$\overline {\Cal F} [f]:=\liml_{\eps \ri 0}\limsupl_{t\ri\iy}
\Cal F(R_\eps u_t(\bullet)),\tag 5.7.1.12$$
 $$\underline {\Cal F} [f]:=\liml_{\eps \ri 0}\liminfl_{t\ri\iy}
\Cal F(R_\eps u_t(\bullet)),\tag 5.7.1.13$$
where $u=\log |f|$ and $(\bullet)_t$ is defined by (3.1.2.1).
\proclaim {Proposition 5.7.1.2} For $\Cal F (u)$ defined by (5.7.1.4)
$$\overline {\Cal F}[f]=h_f(\phi);\ \
\underline {\Cal F}[f]=\underline h_f (\phi).$$
For other functionals from the list (5.7.1.5)-(5.7.1.9) one may replace $R_\eps u$ by $u$ and omit $\liml_{\eps\ri 0}$.
\ep
{\bf Exercise 5.7.1.3.} Prove this.

The following assertion connects the asymptotic growth characteristics with
limit sets.
\proclaim {Theorem 5.7.1.4} The relations
$$\overline {\Cal F}[f]=\sup\{\Cal F(v): v\in \Fr [f]\},$$
$$\underline {\Cal F}[f]=\inf\{\Cal F(v): v\in \Fr [f]\}$$
are true.
\ep
\demo {Proof} Let $v\in \Fr [f]$ and  $u_{t_j}\ri v$ in $\Di '.$ Then
$R_\eps u_{t_j}\ri R_\eps v$ uniformly on every compact set.Hence
$$\liml_{t_j\ri \iy} \Cal F(R_\eps u_{t_j})=\Cal F( R_\eps v).$$
Passing to the limit as $\eps\ri 0$ we obtain
$$\liml_{\eps \ri 0}\liml_{t_j\ri \iy} \Cal F(R_\eps u_{t_j})=\Cal F(v).$$
Choosing a sequence that corresponds to $\limsup$ or $\liminf$ we obtain the assertion of the Theorem.
\qed
\edm
Applying this theorem to the functional (5.7.1.4) we obtain the RHS's of
(3.2.1.1), (3.2.1.2) and hence another definition for the indicator and
lower indicator.

\subheading {5.7.2} A {\it family of growth characteristics}
$\chi_A:=\{ \Cal F_\a (r,\bullet): \a\in A\}$ is called {\it total} if the equation
$$ \Cal F_\a (v_1)=\Cal F_\a (v_2),\ \forall r>0,\ \a\in A \tag 5.7.2.1$$
implies $v_1\equiv v_2$ for $v_1,v_2\in U[\r]$ (see 3.1.2.4).

Here are some examples of the total families:
$$\chi_H:=\{H_\phi(u(e^{i\phi}):\phi\in [0,2\pi)\}; \tag 5.7.2.2$$
$$\chi_I:=\{I_{\a ,\be}(u):\a,\be \in [0,2\pi)\};\tag 5.7.2.3$$
$$\chi_{Fo}:=\{ c_k(u)=I(u,g_k):k\in \BZ\};\tag 5.7.2.4$$
where
$$g_0:=1,\ g_k:=\cos k\phi;\ g_{-k}=\sin k\phi,\ k\in \BN\tag 5.7.2.5$$
It is easy to deduce from Th.5.6.1.1
\proclaim {Theorem 5.7.2.1} Let a family $\{ \Cal F_\a (\bullet): \a\in A\}$
be a total family of characteristics. An entire function $f$ is a CRG -function
iff
$$\overline {\Cal F}_\a[f]=\underline {\Cal F}_\a[f]\tag 5.7.2.6$$
\ep
{\bf Exercise 5.7.2.1.} Check this.
\subheading {5.7.3} Let us consider a total family of characteristics of the form
$$\chi_\Psi:=\{I(u,\psi): \psi\in \Psi\},\tag 5.7.2.7$$
where $\Psi $ is a set which is complete in $L^1[0,2\pi].$ For instance, such are the families $\chi_I$ and $\chi_{Fo}.$
\proclaim {Theorem 5.7.3.1}\text {\cite {Po(1985}}Let  $f\in A(\r(r))$ is a CRG-function if and only if at least one of the following assertions are equivalent:

a)$\overline {\Cal F}[fg] =\overline {\Cal F}[f]+\overline {\Cal F}[g],\forall \Cal F\in\chi_\Psi$,

b)$\underline {\Cal F}[fg] =\underline {\Cal F}[f]+\underline {\Cal F}[g],\forall \Cal F\in\chi_\Psi$,

for all entire functions $g\in A(\r(r)).$

c)$f$ is a GRG-function.
\ep
Let us prove $c)\Longrightarrow a)$ and $c) \Longrightarrow b).$

Using the Theorem 5.7.1.4 we obtain for every characteristics $\Cal F$
$$\overline {\Cal F} [fg]= \sup \{ \Cal F(w): w\in \Fr [fg]\}.\tag 5.7.3.1$$

Because of Th.3.1.2.4  fru1)
$$ \Fr [fg]\sbt \Fr [f]+\Fr [g]$$
Since $f$ is CRG-function $Fr [f]$ consists of only one subharmonic function
$v_{reg}$ (see Th 5.6.1.1) and it is easy to check that in this case we have equality
$$\Fr [fg]=v_{reg} + \Fr [g].$$

{\bf Exercise 5.7.3.1.} Check this.

Since $\Cal F (v_{reg} +v_g)=\Cal F(v_{reg}) +\Cal F(v_g)$ We
We obtains
 $$\overline {\Cal F} [fg]= \Cal F(v_{reg}) +
\sup \{ \Cal F(v_g): v_g\in \Fr [g]\}=\overline {\Cal F}[f]+
\overline {\Cal F} [g]$$
So $c)\Longrightarrow  a)$  was proved. In the same way one can  prove $c)\Longrightarrow  b).$

{\bf Exercise 5.7.3.2.} Prove this.

In  the proof of sufficiency of the conditions of this theorem we can suppose that
$\psi$ belong to the space $\Di (\BT)$ of infinitely differentiable functions on the unit
circle $\BT$ because $\Di (\BT)$ is complete in $ L^1[0,2\pi].$ We prove now sufficiency of
b) in the Theorem 5.7.3.1.

We recall that (see (3.1.2.4a))
$$v_{[t]} (z)=v(tz)t ^{-\rho}, v\in U[\r] $$
to distinguish it from $(\bullet)_t$ that defined by
$$u_t(z)=u(tz)t^{-\r(t)}, u\in SH(\r(r))$$

The main constructive element for the proof of Theorem 5.7.3.1 is
\proclaim {Lemma 5.7.3.2} Let $\psi^0\in \Di (S).$ There exists $v\in U[\r]$ with the following properties:
$$\Di' -\liml_{t\ri 0}v_{[t]} =\Di' -\liml_{t\ri \iy}v_{[t]} =\tilde v,
\tag 5.7.3.2$$

$$\langle v_{[t]} (e^{i\bullet}),\psi^0\rangle\ >
\langle v(e^{i\bullet}),\psi^0\rangle\  \text{for} \ t \in (0,\iy), t    \neq 1,
\tag 5.7.3.3$$
$$\langle\tilde v(e^{i\bullet}),\psi^0\rangle\ >\ \langle v(e^{i\bullet}),\psi^0\rangle  \tag 5.7.3.4$$
\ep
\demo {Proof}Let $\psi^0$ be represented by  Fourier series
$$\psi^0 (\phi)=\frac {a_0}{2}+\suml_{n=1}^{\iy}(a_n \cos n\th + b_n \sin n\th)$$
Since $\psi^0\not\equiv 0$ there exists $a_k\neq 0$ or $b_k\neq 0 .$ Suppose there exists $a_k\neq 0.$ In the proof  we will consider three cases: 1.$k=0;$ 2.
$k\neq 0 \wedge k\leq p;$ 3. $k\geq p+1.$ The number $\r$ is supposed noninteger and $p=[\r].$

Consider the case $a_0\neq 0,\ a_0>0$. Set
$$\psi(x):=\log (-e^{-\a|x|}+C), \a>0, C>1$$
$$v(z):=|z|^\r e^{\psi (\log |z|)}=\exp (\r\log r+\psi(\log r))\tag 5.7.3.4a$$
Applying the Laplace operator, we obtain:
$$\Delta v=\frac {1}{r^2}r\frac {\p}{\p r}r\frac {\p}{\p r}v(r)=
e^{-2x}\frac{\p^2}{\p x^2}e^{\r x+\psi(x)}=$$
$$ =\exp ((\r-2)x+\psi(x))[(\r+\psi '(x))^2+ \psi''(x)], \ x=\log r
\tag 5.7.3.5$$
Since
$$\psi'(x)=\a \sgn x\exp (-\a|x|) \ri 0 ,\ \psi''(x)=-\a^2 \sgn x \exp (-\a|x|) \ri 0$$
as $x\ri \pm\iy,$ it is possible to chose $\a$ such that the expression (5.7.3.5)
 be positive. So $v(z)$ is subharmonic.

It is easy to check that all the assertions of the lemma are satisfied and
$\tilde v(z)=b|z|^\r$ where $b (>0)$ is a constant.
\edm
 {\bf Exercise 5.7.3.3.} Check this.
\demo {Continuation of the proof}
If $a_0=-|a_0|<0,$ consider the function
$$v^0(z):=\cases \log |z|, &|z|\geq 1,\\ 0, &|z|<1;\endcases$$
it is subharmonic and
$$<v^0_{[t]} (e^{i\bullet}),\psi^0> = a_0t^{-\r}\log^+ t \tag 5.7.3.6$$
Since the RHS of (5.7.3.6) is minimized for $t_0=e^{\r^{-1}},$ the function
$$v(z):=v^0_{[{t_0}^{-1}]}(z)$$
satisfies the assertions of the lemma with $\tilde v = \liml_{t\ri 0,\iy} v_{[t]}=0$

Now let $a_0=0, a_k\neq 0, 0<k<p.$
We will search for a  function $v$ of the form
$$v(re^{i\phi}):=\intl_0^{2\pi} G_p(re^{i(\phi -\th)})(1-\sgn a_k\cos k\th )
d\th\tag 5.7.3.7$$
This is the convolution $G_p(re^{i\bullet})*g$ of the primary kernel (see \S 2.9.1)
 $$G_p(z)=\log |1-z|+\Re \suml_{n=1}^p z^n/n$$
with a positive function $g(\th):=(1-\sgn a_k\cos k\th )$ on the circle. So
it is subharmonic. Recall that the cos -Fourier coefficients of the function
$G_p(re^{i\th})$ are (see Exercise 2.3.7.2)

$$\hat G_p(m,r)=\cases 0,\  &m=0,1,...,p\\ (1/m)r^m, \ &m=p+1,...\endcases \ \ \text {if}\ r\leq 1\tag 5.7.3.8$$
and
$$\hat {G}_p(m,r)=\cases\ \log r,\ &m=0\\ \frac {1}{m}(r^m- r^{-m}),\  &m=1,...,p\\ (1/m)r^m, \ &m=p+1,...\endcases \ \text{if}\  r\geq1\tag 5.7.3.9$$
All the sin-Fourier coefficients are equal to zero.
The Fourier coefficients of the function $g$ are $1$ and $-\sgn a_k.$

Using well known properties of Fourier coefficients, we obtain for $0<k\leq p$
$$\hat v_{[t]}(0)=\cases 0, &t   \leq 1 \\
\frac{\log t   }{t^\rho},\ &t   \geq 1\endcases $$
$$ \hat v_{[t]}(k)=\cases 0,&t   \leq 1 \\ -\frac {1}{k}
\frac{t^k-t^{-k}}{t^\rho}\sgn a_k  \ &t\geq 1\endcases$$
$$\langle v_{[t]} (e^{i\bullet}),\psi^0\rangle=\cases 0, &t\leq 1\\ -1/k(t^{k-\rho}-t^{-k-\rho})|a_k|\ &t\geq 1\endcases$$

The function $t\mapsto \langle v_t (e^{i\bullet}),\psi^0\rangle$ tends to zero when
$t\ri0,\iy$ and has its only minimum at the point
$$t_0=\left(\frac{\r+k}{\r-k}\right)^{1/k}$$
Thus $v_{(t_0)^{-1}}$ satisfies the conditions of the lemma with
$\tilde v=0.$

For $k\geq p+1$
we should take  the same $g$ and then
$$\langle v_t (e^{i\bullet}),\psi^0\rangle =\cases -(1/k)t^ {k-\rho}|a_k|,
&t\leq 1
\\ -(1/k)t^{-k-\rho}|a_k|\ &t\geq 1\endcases$$
So the corresponding function $t\mapsto
\langle v_t (e^{i\bullet}),\psi^0\rangle $ obtains
minimum at the point $t_0=1$ and the function $v$ satisfies the assertions of the lemma with $\tilde v=0.$
\qed
\edm
{\bf Exercise 5.7.3.4.} Prove the lemma for the case $b_k\neq 0$

\proclaim {Lemma 5.7.3.3} Let $v\in U[\r]$ with the following condition fulfilled
$$\Di ' -\liml_{t\ri 0}v_{[t]} =\Di '- \liml_{t\ri \iy}v_{[t]}= \tilde v$$
and let $u\in SH(\r(r))$ with some $v^0\in \Fr[u]$.

Then there exists $w^0 \in SH(\r(r))$ such that
$$\Fr [w^0]=\{v_{[t]}:t \in (0,\iy)\}\cup \tilde v\tag 5.7.3.10$$
and the following condition holds:

1.if the sequence   $\liml_{t_n\ri\iy}w^0_{t_n}=v_{[t]}$ for some
$t\in (0,\iy)$ and the sequence $u_{t_n}$ converges in $\Di'$ as $t_n\ri\iy$ then $\liml_{n\ri\iy}u_{t_n}=v^0_{[t]}.$
\ep
For proof see Corollary 4.4.1.3.
\proclaim {Lemma 5.7.3.4} Let $w\in SH(\r(r)),\ \psi\in \Di (S)$. Then the following holds:
$$\liminf\limits_{t\ri\iy}\langle w,\psi \rangle=\min\limits_{v\in \Fr\ w}
\langle v,\psi\rangle$$
$$ \limsup\limits_{t\ri\iy}\langle w,\psi \rangle=\max\limits_{v\in \Fr\ w}
\langle v,\psi\rangle$$
\ep
{\bf Exercise 5.7.3.5.}Prove this exploiting completeness of $\Fr.$
\demo {Proof of sufficiency in Theorem 5.7.3.1}In assumption b) we should prove that $f$ is a CRG-function, i.e.,by Th.5.6.1.1 its $\Fr[f]$ consists of  only one function. Since $\log |f|\in SH (\r(r))$ and   because of Th.5.3.1.4 (Approximation)
it is enough to prove the corresponding theorem for subharmonic functions.
Suppose

$$\underline{\Cal F} [u+w]=\underline{\Cal F}[u]+\underline{\Cal F}[w], \forall\Cal F\in \chi_\Psi\tag 5.7.3.11 $$
for all $w \in SH(\r(r).$
 We exploit Lemma 5.7.3.5 and write (5.7.3.10) in the form:
$$\min\limits_{v\in \Fr [u+v]}\langle v,\psi\rangle=\min\limits_{v\in \Fr u}
\langle v, \psi\rangle +\min\limits_{v\in \Fr w}
\langle v, \psi\rangle , \forall \psi\in\Psi.$$

  Suppose the contrary, i.e., $u$ is not a CRG-function and $\Fr u$ does not consists of only one $v_{\min}\in U[\r].$ Then there exists $v^0\neq v_{\min}.$  The family $\chi_\Psi$ is total; therefore
   there exists $\psi^0 \in \Psi$ such that
$\langle v^0,\psi^0\rangle\neq \langle v_{\min},\psi^0\rangle$ and hence
$$\langle v^0,\psi^0\rangle > \langle v_{\min},\psi^0\rangle.\tag 5.7.3.12$$
Using Lemma 5.7.3.2, construct for the function $\psi^0$ a function $v\in
U[\r]$ satisfying the conditions (5.7.3.2),(5.7.3.3) and (5.7.3.4). Apply
Lemma 5.7.3.3 to construct a function $w^0$ satisfying (5.7.3.10a) and
the condition 1.
Under conditions of the Theorem
$$\min\limits_{\om\in \Fr (u+w^0)}\langle \om,\psi^0\rangle=
\min\limits_{\om\in \Fr (u)}\langle \om,\psi^0\rangle +
\min\limits_{\om\in \Fr (w^0)}\langle \om,\psi^0\rangle\tag 5.7.3.13$$
Let $\g\in \Fr(u+w^0)$ be the function on which the minimum of LRH in
(5.7.3.11) is attained. Using  (5.7.3.3), (5.7.3.4) and (5.7.3.10), we can
rewrite (5.7.3.11) in the form
$$\langle \g,\psi^0\rangle =\min\limits_{\om\in \Fr u}\langle \om,\psi^0\rangle +\langle v,\psi^0\rangle\tag 5.7.3.14$$
Since $\g\in \Fr (u+w^0), \g=\Di'-\lim\limits_{n\ri\iy}(u+w^0)_{t_n}$.
Passing to subsequences, we can suppose that the sequences $\{u_{t_n}\}$ and
$\{w^0_{t_n}\}$ have  limits. Since $\Fr w^0$ has the form (5.7.3.10), there
are two possible cases : $w^0_{t_n}\ri v_{[t]},\ t\in (0,\iy)$ and
$w^0_{t_n}\ri\tilde v.$

Consider the first case. Because of condition 1. from Lemma 5.7.3.3
$ u_{t_n}\ri v^0_{[t]}$ and $\g= v^0_{[t]}+v_{[t]}.$ Substituting this in (5.7.3.13), we obtain
$$\langle v^0_{[t]},\psi^0\rangle - \min\limits_{\om\in \Fr u}\langle \om,\psi^0\rangle=\langle v,\psi^0\rangle
- \langle v_{[t]},\psi^0\rangle$$
This equality leads to contradiction because for  $t=1$ it contradicts
(5.7.3.12) and for $t\neq 1$ it contradicts (5.7.3.3).

Consider the second case, when $ w^0_{t_n}\ri\tilde v.$ Denote
$v^2=\liml_{n\ri\iy}u_{t_n}$ and rewrite (5.7.3.11) in the form
$$ \langle v^2,\psi^0\rangle -\min\limits_{\om\in \Fr u}\langle \om,\psi^0\rangle=\langle v,\psi^0\rangle
- \langle \tilde v,\psi^0\rangle$$
 The last equality contradicts (5.7.3.4).
\qed\edm

 Sufficiency of condition a) of Theorem 5.7.3.1 can be proved
 using the Lemmas 5.7.3.3 ,5.7.3.4 and the following  lemma.
\proclaim {Lemma 5.7.3.2'} Let $\psi^0\in \Di (S).$ There exists $v\in U[\r]$ with the following properties:
$$\Di' -\liml_{t\ri 0}v_{[t]} =\Di' -\liml_{t\ri \iy}v_{[t]} =\tilde v,
\tag 5.7.3.2'$$

$$\langle v_{[t]} (e^{i\bullet}),\psi^0\rangle\ <
\langle v(e^{i\bullet}),\psi^0\rangle\  \text{for} \ t \in (0,\iy), t \neq 1,
\tag 5.7.3.3'$$
$$\langle\tilde v(e^{i\bullet}),\psi^0\rangle\ <\ \langle v(e^{i\bullet}),\psi^0\rangle  \tag 5.7.3.4'$$
\ep
 {\bf Exercise 5.7.3.6.} Prove this lemma and sufficiency of a) in Theorem 5.7.3.1.
\subheading {5.7.4}
In this \S \ we consider the question of summing the asymptotic characteristics connected with the functional (5.7.1.4), i .e. indicator  and  lower indicator. Recall that $f\in A(\r(r))$ is completely regular on the ray $\{\arg z=\phi\}$ ( $f\in A_{reg,\phi}$) if
$$h_f(\phi)=\underline h_f (\phi)\tag 5.7.4.0$$

We are going to prove the following assertions:
\proclaim {Theorem 5.7.4.1}Let $f\in A_{reg,\phi}.$ Then for every $g\in A(\r(r))$
$$h_{fg}(\phi)=h_f(\phi)+h_g(\phi)\tag 5.7.4.1$$
$$\underline h_{fg}(\phi)=\underline h_f(\phi)+\underline h_g(\phi)\tag 5.7.4.2$$
   \ep
\proclaim {Theorem 5.7.4.2}Suppose the equality (5.7.4.1) holds for every $g\in A(\r(r)).$ Then $f\in A_{reg,\phi}.$
\ep
Let us note that the assertion of the Theorem 5.7.4.2 holds also if the equality
(5.7.4.1) fulfilled for some sequence $\phi_n\ri\phi,$ because indicator is
continuous function (see \S 3.2.5). So if the equality (5.7.4.1) holds for
the set $\Phi$ of $\phi$ that is dense in $[0,2\pi)$ (or the set
$$e^{i\Phi}:=\{e^{i\phi}:\phi\in \Phi\}\tag 5.7.4.3$$ is dense on the unit circle), then
$f\in A_{reg,\phi}$ for all $\phi,$ i.e. $f$ is a CRG -function.

On the other hand, the following assertion holds
\proclaim {Theorem 5.7.4.3}If the set $\Theta$ of $\theta$ is not dense in $[0,2\pi),$  there exists $f\in  A_{reg,\theta}, \theta\in \Theta$ that is not a CRG-function.
\ep
The situation  with lower indicator is analogous, but in another
 topology.

A set $E$ is called {\it non -rarefied} at a point $z_0$ if for every function $v$ subharmonic in a neighborhood of $z_0$ the following holds:
$$v(z_0)=\limsup\limits_{ z\in E,z\ri z_0,z\neq z_0}v(z)=
\limsup\limits_{ z\in E,z\ri z_0}v(z).$$
A set is {\it rarefied} if it is not non-rarefied.

Note  that if $\underline h_f(\phi)=-\iy$ then $\underline h_{fg}(\phi)=-\iy$ for every $g\in A(\r(r)).$ It is obvious that $f\notin A_{reg,\phi}.$

The next theorems was proved in \cite {GPS}.
\proclaim {Theorem 5.7.4.4}   Let (5.7.4.2) be fulfilled for $\psi\in E$ for all $g\in A(\r(r))$ and $e^{iE}$ is non rarefied at the point $e^{i\phi}$. Then
$f\in A_{reg,\phi}.$
\ep
\proclaim {Theorem 5.7.4.5}Let $E_0$ be a set such that $e^{iE_0}$ is rarefied at all the points of unit circle. Then there exists $f\in A(\r(r))$ for which  (5.7.4.2) fulfilled for all $\phi \in E_0$ and all $g\in A(\r(r)),$ but
$f\notin A_{reg,\phi}$ for all $\phi$ and $\underline h_f(\phi)>-\iy,\ \forall \phi.$
\ep

Let us note that $E_0$ can be dense in $[0,2\pi)$ and $E$ from Theorem 5.7.4.4
can even be of zero measure.

The proof of Theorems 5.7.4.4 and 5.7.4.5 is based on the following assertion
that gives a criterion for (5.7.4.2) in terms of limit sets $\Fr[f].$
\proclaim {Theorem 5.7.4.6}Let $f\in A(\r(r))$ and $\underline h_f(\phi)>-\iy.$ The condition (5.7.4.2)holds for every  $g\in A(\r(r)),$ such that $h_g(\phi)>-\iy$ iff
$$\liminf\limits_{t\ri 1}v(t e^{i\phi})=\underline h_f(\phi)\tag 5.7.4.4$$
 for all $v\in \Fr[f].$
 \ep
An analogous criterion holds for (5.7.4.1).
\proclaim {Theorem 5.7.4.7} Let  $f\in A(\r(r)).$ (5.7.4.2) holds for every  $g\in A(\r(r)),$ iff
$$\limsup\limits_{t\ri 1}v(t e^{i\phi})=h_f(\phi),\tag 5.7.4.5$$
for all $v\in \Fr[f].$
\ep
\proclaim {Corollary 5.7.4.8} The equality  (5.7.4.5) implies
$f\in A_{reg,\phi}$.
\ep
Actually , for every $v\in \Fr[f]$ we have, using semicontinuity of subharmonic functions and the definition (3.2.1.1) of the indicator,
$$h_f(\phi)=\limsup \limits_{t\ri 1}v(t e^{i\phi})\leq v(e^{i\phi})\leq
h_f(\phi)$$
for all $v\in \Fr[f].$ So $\Fr[f]$ consists of functions $v$ that coincide at
the point $e^{i\phi}$ and hence on the ray $\{re^{i\phi}:r\in (0,\iy)\}.$

Note also that the set $e^{iE}$ for which (5.7.4.1) holds is closed and Theorem 5.7.4.4 means that the set where (5.7.4.2) holds is {\it thinly closed}, i.e.,closed in {\it thin topology} (see
\cite {Br\S 6}.)

Therefore if $e^{i\phi_0}$ is a limit point of $e^{iE}$ in the euclidian (respectively, thin) topology, then (5.7.4.1)((5.7.4.2), respectively)is also a sufficient condition
for completely regular growth at $\phi_0.$

 \subheading {5.7.5} The main constructive element for proving Theorem
5.7.4.6 is
\proclaim {Lemma 5.7.5.1}Let $\eps>0,t_0>0$ and $\phi_0\in [0,2\pi)$ be fixed. Then there exists $v\in U[\r]$ with the following properties:

$$\Di '-\liml_{t\ri 0}v_{[t]}=\Di '-\liml_{t\ri \iy}v_{[t]}=0 \tag 5.7.5.1 $$
$$v(e^{i\phi_0})<v_{[t]}(e^{i\phi_0}),\ t \in (0,1)\cup (1,\iy)\tag 5.7.5.2$$
$$-\iy<v(e^{i\phi_0})<-\eps,\tag 5.7.5.3$$
and the inequality
$$v_{[t]}(e^{i\phi_0})-v(e^{i\phi_0})\leq \eps/2 \tag 5.7.5.4$$
implies
$$t\in [1/t_0,t_0]\tag 5.7.5.5$$
\ep
The last condition means that the function $\psi(t):=v_{[t]}(e^{i\phi_0})$ can be less than $\psi(1)+\eps/2$ only in  a neighborhood of $t=1.$

\demo {Proof}Set
$$w(z):=\max (\log |1-ze^{-i\phi_0}|,-N)+\Re \suml_{n=1}^p\frac {1}{n}
(ze^{-i\phi_0})^n,$$ $$\ N>0, p=[\r].\tag 5.7.5.6$$ It is obvious
that $w$ is subharmonic, with masses $\nu_w$ concentrated in a
neighborhood of the point $e^{i\phi_0}.$ Thus $\nu_w\in \Cal M
[\r]$ (see (3.1.3.4)) and
$$\Di'-\liml_{t\ri 0}(\nu_w)_{[t]}= \Di'-\liml_{t\ri \iy}(\nu_w)_{[t]}=0$$
Hence (see Th.3.1.4.2) $w\in U[\r],$ and (see (3.1.5.0))
 $$\Di'-\liml_{t\ri 0}w_{[t]}=\Di'-\liml_{t\ri \iy}w_{[t]}=0\tag 5.7.5.7  $$
Let us capitalize on the  behavior of $w_{[t]}$ on the ray $\{\arg
z=\phi_0\}.$
$$w_{[t]} (e^{i\phi_0}):=\psi(t)=(\max (\log |1-t|,-N)+\Re \suml_{n=1}^p\frac {1}{n})t^{-\r}
t^n,\tag 5.7.5.8 $$
It is possible to prove directly the following properties of $\psi(t).$

i) outside interval $(1-e^{-N},1+e^{-N}), \ \psi(t)=G_p(t)t^{-\r};$
where $G_p$ is the Primary Kernel (see \S 2.9.1) and inside this interval
the first summand is $-N;$

ii) $\psi(t)>0$ for $t>t_1$ where $t_1$ is a zero of the equation
$G_p(t)=0,\ \psi(t)$ decreases monotonically on the interval $(0, 1-e^{-N})$ and increases monotonically  on the interval $(1-e^{-N},t_1).$

{\bf Exercise 5.7.5.1.} Prove this.

Now set $t_2:=1-e^{-N}$ and $v(z):=Dw_{t_2}(z),$ where $D$ is a constant.
This function satisfies the conditions (5.7.5.1) and (5.7.5.2) of the lemma
and $v_{[t]}(e^{i\phi_0})$ has the only one negative minimum for $t=1$. Thus
it is possible to take $D$ sufficiently large  to satisfy the conditions
(5.7.5.3) and (5.7.5.4) for fixed $\eps$ and $t_0.$
\qed
\edm
 {\bf Exercise 5.7.5.2.} Prove this in details.

In the proof of Theorem 5.7.4.6 we also use Lemma 5.7.3.3. We can prove all
the assertions for subharmonic functions from $SH(\r(r)).$
\demo {Proof of Theorem 5.7.4.6}Necessity.
We should prove that if the equality
$$\underline h(e^{i\phi_0},u+w)=\underline h(e^{i\phi_0},u)+\underline h(e^{i\phi_0},w)\tag 5.7.5.9$$
holds for a fixed $u\in SH(\r(r)),\ \phi_0$ and every $w\in SH(\r(r)),$  then
$$ \liminf\limits_{t\ri 1}v(t e^{i\phi})=\underline h(e^{i\phi},u)\tag 5.7.5.10$$
 for all $v\in \Fr u.$
 Assume that $\underline h(e^{i\phi},u)>-\iy$ and $\underline h(e^{i\phi_0},w)> -\iy.$
Suppose the contrary, i.e. there exists $v^0\in \Fr u$ such that
$$\liminf\limits_{t\ri 1}v^0(t e^{i\phi_0})>\underline h(e^{i\phi_0},u) \tag 5.7.5.11$$
The inequality (5.7.5.11) implies that there exists $\eps >0$ and $t_0 >0$
such that for every $t\in [1/t_0,t_0]$ the inequality
$$  v^0(t e^{i\phi_0})> \underline h(e^{i\phi_0},u)+\eps\tag 5.7.5.11a$$
holds.
Let us construct by Lemma 5.7.5.1 for these $\eps, t_0, \phi_0$ a function
$v$ and by Lemma  5.7.3.3 for the functions $u, v^0 $ and the already found $v$ a function $w^0.$ Let us show that for $w^0$ the equality (5.7.5.9) does not hold.

Compute $\underline h(e^{i\phi},w^0).$ From (3.2.1.2)
$$\underline h(e^{i\phi_0},w^0)=\min \{0, \inf \{v_{[t]}(e^{i\phi_0}:t \in (0,\iy)\}\}.$$
The inequalities (5.7.5.3) imply that $0$ can be omitted and (5.7.5.2) implies
that the infimum is attained at $t =1,$ i.e.,
$$\underline h(e^{i\phi_0},w^0)=v(e^{i\phi_0}).\tag 5.7.5.12$$
Find $v^{\eps}\in \Fr (u+w^0)$ such that $\underline h(e^{i\phi_0},u+w^0)>
v^{\eps}(e^{i\phi_0})-\eps/3.$ Let $t_n\ri\iy$ and $(u+w^0)_{t_n}\ri v^{\eps}$ in $\Di'.$ Passing to subsequences we can assume that $u_{t_n}$ and $w^0_{t_n}$ also converge. Consider two cases. The first, when
$$\Di'-\lim w^0_{t_n}=v_{[t]}, \ t\in(0,\iy)\tag 5.7.5.13$$
By Lemma 5.7.3.3 $\lim u_{t_n}=v^0_{[t]}$ and hence
$v^{\eps}=\lim (w^0+u)_{t_n}=v_{[t]} +v^0_{[t]}.$ If $t\notin [1/t_0,t_0]$ then by (5.7.5.4)
$$v_{[t]} ( e^{i\phi_0})>v(e^{i\phi_0})+\eps/2=\underline h(e^{i\phi_0},w^0)+\eps/2\tag 5.7.5.14$$
In this case we have
$$\underline h(e^{i\phi_0},u+w^0)\geq v^{\eps}(e^{i\phi_0})-\eps/3\geq
v_{[t]}+v^0_{[t]} -\eps/3\tag 5.7.5.14a$$
Using (5.7.5.14a), we obtain
$$\underline h(e^{i\phi_0},u+w^0) \geq\underline h(e^{i\phi_0},w^0)+\underline h(e^{i\phi_0},u)+\eps/6\tag 5.7.5.15$$
If $t\in [1/t_0,t_0]$ then from (5.7.5.11) we have
$$\underline h(e^{i\phi_0},u+w^0)\geq \underline h(e^{i\phi_0},w^0)+\underline h(e^{i\phi_0},u)+2\eps/3 \tag 5.7.5.16$$
So the case (5.7.5.13) is settled.

Let $\Di'-\lim w^0_{t_n}=0.$ In this case   we have
$$\underline h(e^{i\phi_0},u+w^0)\geq v^{\eps}(e^{i\phi_0})-\eps/3\geq
\underline h(e^{i\phi_0},u)-\eps +\eps -\eps/3=
\underline h(e^{i\phi_0},u)-\eps+2\eps/3.$$
Using (5.7.5.12) and (5.7.5.3) we obtain
$$\underline h(e^{i\phi_0},u+w^0)\geq\underline h(e^{i\phi_0},u)+\underline h(e^{i\phi_0},w^0)+2\eps/3.$$
So we proved in any case that (5.7.5.9) does not hold if (5.7.5.10) does not
hold .

Let us prove sufficiency in Theorem 5.7.4.6. We prove it for subharmonic functions. Let $u\in SH(\r(r))$ and for every $v\in \Fr u $  (5.7.5.10) holds. Let us
show that for all $w\in SH(\r(r))$ (5.7.5.9) holds. It is sufficient to prove that
$$\underline h(e^{i\phi_0},u+w)\leq\underline h(e^{i\phi_0},u)+\underline h(e^{i\phi_0},w)\tag 5.7.5.17$$
holds since the inverse inequality holds for every $w\in  SH(\r(r))$ (see (3.2.1.5)).
Let us note for beginning that for every $v^2\in\Fr w$ there exist $v\in
\Fr (u+w)$ and $v^1\in\Fr u$ such that
$$v=v^1+v^2\tag 5.7.5.19$$
Indeed, let $t_n\ri\iy$ be a sequence such that $w_{t_n}\ri v^2.$ We can suppose , choosing subsequence, that $u_{t_n}\ri v^1$ and $(u+w)_{t_n}\ri v.$
Then (5.7.5.19) holds.

Let $\eps$ be arbitrarily small.Chose $v^2\in \Fr w$ such that $v^2(e^{i\phi})<h(e^{i\phi},w)+ \eps$ holds. From upper semicontinuity of $v^2$ we have
$$\limsup\limits_{t\ri 1}v^2 (e^{i\phi})\leq \underline h(e^{i\phi_0},w) +\eps.\tag 5.7.5.20$$
Let $v^1\in \Fr u$ and $v\in\Fr (u+w)$ satisfy (5.7.5.19).Then we have
$$ \underline h(e^{i\phi_0},u+w)\leq (v^1+v^2)_{[t]} ( e^{i\phi_0})=
v^1_{[t]} ( e^{i\phi_0})+v^2_{[t]} ( e^{i\phi_0}),\ \forall t.$$
Hence
$$\underline h(e^{i\phi_0},u+w)\leq \liminf\limits_{t\ri 1}v^1_{[t]} ( e^{i\phi_0})+
\limsup\limits_{t\ri 1}v^2_{[t]} ( e^{i\phi_0}).$$
Using (5.7.5.10) and (5.7.5.20) we obtain
$$\underline h(e^{i\phi_0},u+w)\leq \underline h(e^{i\phi_0},u)+
\underline h(e^{i\phi_0},w)+\eps.$$
This proves the inverse inequality and hence the equality (5.7.5.9), because $\eps$ is arbitrarily small.
\qed
\edm
\subheading {5.7.6} Now we are going to prove Theorem 5.7.4.4.We need the following assertion from Potential Theory.
\proclaim {Lemma 5.7.6.1}Let $E$ be a set that is non-rarefied at the point
$ e^{i\phi_0}.$  Let $E'$ be a set in $\BC , $ such that
$\forall e^{i\phi}\in E$ and $\forall \delta >0$ there exists a point
$z'\in E'$ on the ray $\{\arg z=\phi\}$ such that $|z'-  e^{i\phi}|<\delta$. Then $E'$ is also non-rarefied at the point $ e^{i\phi_0}.$
\ep
\demo {Proof} We can suppose without loss of generality that $E'$ have no intersection with some neighborhood of zero. Denote by $P(z)$ the map
$z\mapsto e^{i\arg z}.$ It is easy to see that for all pairs $z_1',z_2'\in E'$the inequality $|P(z_1')-P(z_2')|<A|z_1'-z_2'|$ hold for some constant $A.$ Thus
the logarithmic capacity (2.5.2.5) satisfies (\cite {La,Ch.II,\S 4,it.11,15})
$$\text {\bf cap}_l(M)<A\text{\bf cap}_l(M')\tag 5.7.6.2$$

where $M'\sbt E',M=P(M').$
Now we exploit the following properties of non-rarefied  sets. First, if $E$ is non-rarefied at a point $z_0$, then there exists
a compact set that is non-rarefied at $z_0$ (\cite {La,Ch.V,\S 1,it.5,\S 3,it.9} Second, for
a compact set $K$ that is non-rarefied at $z_0$
$$\sum\limits_{n=1}^{\iy}\frac {n}{\log (\text{\bf cap}_lK_n)^{-1}} =\iy\tag 5.7.6.3$$
where $K_n:=K\cap\{z:q^{n+1}\leq |z-z_0|\leq q^n\},\ 0<q<1.$

Using the inequality (5.7.6.2),we obtain that  divergence of the series
(5.7.6.3) for a compact $K\sbt E$ implies divergence for $K'\sbt E'$ where
$K=P(K'),$ i.e., $E'$ is non-rarefied at the point $P(e^{i\phi_0})=e^{i\phi_0}.$\qed
\edm
 \demo {Proof Theorem 5.7.4.4} Let $\eps(\phi)\ri 0$ as $\phi\ri\phi_0$ and let
$v\in\Fr u.$ Suppose (5.7.5.9)holds for $e^{i\phi}\in E.$
By Theorem 5.7.4.6 the equality (5.7.5.10) holds. Thus $\forall \Delta>0,\
\exists z'=z'(e^{i\phi},\Delta)$ such that
$$|z'-e^{i\phi}|<\Delta , \arg z'=\phi ,v(z')<\underline h(e^{i\phi})+\eps (\phi)\tag 5.7.6.4$$
Set
$$E':=\bigcup\limits_{\phi\in E}\bigcup\limits_{n=1}^{\iy}z'(e^{i\phi},1/n)$$
By (5.7.6.4) and upper semicontinuity of $h(e^{i\phi})$ we obtain
$$\limsup\limits_{z'\ri e^{i\phi_0},\ z'\in E'}v(z')\leq \underline h(e^{i\phi_0})\tag
5.7.6.5$$
Since $E'$ is non-rarefied by Lemma 5.7.6.1 the upper limit of $v$ coincides with
$v(e^{i\phi_0})$ and hence
$v(e^{i\phi_0})\leq\underline h(e^{i\phi_0}).$ The  inverse inequality holds always. Thus
$v(e^{i\phi_0})=\underline h(e^{i\phi_0}), \forall v\in \Fr u.$
Hence $h(e^{i\phi_0})= \underline h(e^{i\phi_0}).$
\qed
\edm
\subheading {5.7.7} Now we are going to prove Theorem 5.7.4.5. Before this we need to  describe a construction and prove  some auxiliary assertions.

Let $B_j:=\{z:T^j<|z|<T^{j+1}\},\ j=0,\pm 1,\pm 2,...$ where $T>1$ is a fixed number. Denote
$L_{E_0}:=\{z:e^{i\arg z}\in e^{iE_0}\}.$ Recall that $e^{iE_0}$ is a set rarefied at every point of the unit circle.
Let $Q$ be the set of rational numbers on the interval $(1,T).$ Set
$$S_Q:=\{z:|z|\in Q\},\ T^jS_Q:=\{zT^j:z\in S_Q\},\ A_j:=L_{E_0}\cap T^jS_Q,\
j=0,\pm 1,\pm 2,... $$
\proclaim {Lemma 5.7.7.1} There exists $v\in U[\r]$ such that
$$v(z)=-\iy \tag 5.7.7.1$$
for $z\in A_0$ and
$$\mu_v(e)=0 ,\ \forall e\sbt \BC\setminus B_0. \tag 5.7.7.2$$
\ep
\demo {Proof} The set $E$ is rarefied at  its every point, hence it is polar
(\cite {Br,Ch.7, \S 4}).Thus the set $\{z:|z|=r\}\cap L_{E_0}$ is polar (see
 \cite {Br,Ch.3,\S2}). A countable union of polar sets is polar (\cite {Br,Ch.3,\S 2}). Thus $A_0$ is polar.Hence there exists a positive measure $\mu$ concentrated on $B_0$ for which the potential $v(z):=\int G_p(z/\z)d\mu$ is equal to $-\iy$ on $A_0$ (see \cite {Br,Ch.4,\S 6, Applications}). It is easy to see that $\mu\in \Cal M (\r)$ and hence
$v\in U[\r]$ (see Th.3.1.4.2).
\qed
\edm
\proclaim {Lemma 5.7.7.2} There exists $\om \in U[\r]$ such that the following
conditions are fulfilled:
$$\om (z)=-\iy, \ z\in A:=\cup_{j=-\iy}^{+\iy}A_j;\
\om (Tz)=T^\r\om(z)\tag 5.7.7.3$$
\ep
\demo {Proof} Set for every $E\Subset \BC\setminus 0$
$$\nu (E):=\sum\limits_{j=-\iy}^{j=+\iy}T^{j\r}\mu_v(T^{-j}E\cup B_0)\tag 5.7.7.4$$
(compare Th.4.1.7.1). We have $\nu\in \Cal M (\r). $ Set
$$\om(z):=\int G_p(z/\z)\nu (d\xi d\eta),\ \z=\xi +i\eta .$$
This $\om$ satisfies (5.7.7.3).
\qed
\edm
{\bf Exercise 5.7.7.1.} Prove this using Th.4.1.7.1.
\proclaim {Lemma 5.7.7.3}Let $\om $ be a subharmonic function in $\BC.$
Denote
$$m(\phi):=\max\{\om (re^{i\phi}:r\in [1,T]\}.$$
Then there exists a
constant $C>-\iy$ such that $m(\phi)>C \ \forall \phi$.
\ep
\demo {Proof} If not, there exists a sequence $\phi_n$ that we can assume
to converge to $\phi_{\iy}$ such that $m(\phi_n)\ri -\iy.$ By upper semicontinuity of $\om$ we have $\om (z)=-\iy, \ ze^{-i\phi_\iy}\in [1,T].$ Thus
$\om (z)\equiv -\iy$ because the capacity of the segment in the plane is positive and hence it is not polar for some subharmonic function.
\qed
\edm
Recall that for $v\in U[\r]$ (see (4.1.3.1))
$$\Bbb C (v):=\Di'-\text {clos}\{v_{[t]}:0<t<\infty\},\tag 5.7.7.5$$
$$\Om (v):=\{v'\in U[\r] :(\exists t_k\rightarrow\iy)
(v'=\lim\limits_{k\rightarrow\infty} v_{[t_k]}\}\tag 5.7.7.6$$
$$A (v):=\{v'\in U[\r] :(\exists \tau_k\rightarrow 0)
(v'=\lim\limits_{k\rightarrow\infty} v_{[t_k]}\}\tag 5.7.7.7$$
By Theorems 4.1.3.4 and 4.2.1.2 if

$$ A (v)\cap\Om (v)\neq \varnothing ,\tag 5.7.7.8$$
there exists $u\in SH(\r(r))$ such that
$$\Fr u= \Bbb C (v).\tag 5.7.7.9$$
\proclaim {Lemma 5.7.7.4} There exists $v^1\in U[\r]$ such that the following holds:
$$A (v^1)=\Om (v^1)\tag 5.7.7.10$$
$$\inf \{v(e^{i\phi}):v\in \Bbb C (v^1)\}=\liminf\limits_{t\ri 1}
v (t e^{i\phi}) =0,\ \forall v\in \Bbb C (v^1),\ \forall e^{i\phi}\in e^{iE_0},\tag 5.7.7.11$$
$$\sup \{v( e^{i\phi}):v\in \Bbb C (v^1)\}\neq \inf \{v( e^{i\phi}):v\in \Bbb C (v^1)\}.\tag 5.7.7.12.$$
\ep
\demo {Proof} Let $\om (z)$ be constructed by Lemma 5.7.7.2. Set
$$v(z):=\om (z)+ D\log^+2|z|.$$
The condition (5.7.7.3) implies
$$A (\om)=\Om (\om)=\{\om_{[t]}:t\in [1,T]\}$$
because it is a Periodic Limit Set (see Th.4.1.7.1).

Since $(\log^+2|z|)_{[t]}\ri 0,\ t\ri 0, \ t\ri \iy$ the function $v$ satisfies the condition

 $$A (v)=\Om (v)=\{\om_{[t]}:t\in [1,T]\}$$

By Theorem 2.1.7.4 for the function $v^1:=v^+$ we have

$$A (v^1)=\Om (v^1)=\{\om^+_{[t]}:t\in [1,T]\}.$$
Note that $v^1(z)=0$ for $z\in A$ and since $A$ is dense in $L_{E_0}$ (5.7.7.11) holds. Choosing $D$ sufficiently large it is possible (using Lemma 5.7.7.3) to find on every ray $\{\arg z=\phi\}$ a point
$z_\phi$ where $v^1(z_\phi)>0.$ Hence $\sup \{v( e^{i\phi}):v\in \Bbb C (v^1)\}>0.$
Because of (5.7.7.11) and upper semicontinuity  of $\inf \{v(e^{i\phi}):v\in
\Bbb C (v^1)\}$ it is zero for every $e^{i\phi}.$ Thus (5.7.7.12) holds.
\qed
\edm
\demo {Proof of Theorem 5.7.4.5} Let us construct by Theorems 4.1.3.4 and 4.2.1.2 a function $u\in SH(\r(r)$ such that $\Fr u=\BC(v^1)$ where $v^1$ is taken from Lemma 5.7.7.4.
It does not belong to $A_{reg,\phi}$ for any $\phi$. The equality (5.7.5.9) holds for every $\phi\in E_0$ because of (5.7.7.11) by Theorem 5.7.4.6.
\qed
\edm
\subheading {5.7.8}The proof of Theorem 5.7.4.1 is a copy of the proof of Theorem 5.7.3.1.

 {\bf Exercise 5.7.8.1.} Prove Theorem 5.7.4.1.

Now we are going to prove Theorem 5.7.4.7 which implies (as it was shown in
Corollary 5.7.4.8) Theorem 5.7.4.2. The main constructive element of the proof of necessity is
\proclaim {Lemma 5.7.8.1} Let $\eps>0,t_0>0$ and $\phi_0\in [0,2\pi)$ be fixed. Then there exists $v\in U[\r]$ with the following properties:
$$\Di'-\liml_{t\ri 0}v_{[t}]=\Di'-\liml_{t\ri \iy}v_{[t]}=0\tag 5.7.8.1$$
$$v(e^{i\phi_0})>v_{[t]}(e^{i\phi_0}),\ t \in (0,1)\cap (1,\iy)\tag 5.7.8.2$$
and the inequality
$$v_{[t]}(e^{i\phi_0})-v(e^{i\phi_0})\geq -\eps/2\tag 5.7.8.3$$
implies
$$t\in [t_0,1/t_0]\tag 5.7.8.4$$
\ep The last condition means that the function $\psi(t):=v_{[t]}
(e^{i\phi_0})$ can be more than $\psi(1)-\eps/2$ only in a
neighborhood of $t=1.$ \demo {Proof} Consider the function
$$w(z):=\log ^+|z| \tag 5.7.8.5$$
It is subharmonic and satisfies (5.7.8.1). Since the function
$$\psi (t):=w_{[t]} (e^{i\phi_0})=t ^{-\r}\log ^+t$$ has its only
strict  maximum in  the point $t_{max}>1,$  the function
$$v(z):=w(z/t_{max})$$
has  all the properties (5.7.8.1)-(5.7.8.4).
\qed
\edm

After this lemma all the proof of Theorem 5.7.4.6 can be repeated with minimal
changes.

 {\bf Exercise 5.7.8.2.} Prove Theorem 5.7.4.7.

\subheading {5.7.9} Now we are going to prove Theorem 5.7.4.3.
Let us prove the following
\proclaim {Lemma 5.7.9.1} Let $\overline\Theta$ be a closed subset of
$[0,2\pi)$.  Then for every $\s >0$ there exists a $2\pi$- periodic $\r$ -trigonometrically convex function $h(\phi)$ such that
$$h(\phi)=\s\tag 5.7.9.1$$
for $\phi\in \overline \Theta$ and
$$h(\phi)>\s\tag 5.7.9.2$$
for $\phi\notin \overline\Theta.$
\ep
\demo {Proof} We can suppose that $0\in \overline \Theta$ otherwise we can shift it a little. The set $[0,2\pi)\setminus \overline \Theta$ is open and it can
be represented as union of non-intersecting open intervals. If length of an interval is $\leq \pi/\r$ we can construct a $\r$-trigonometrical function that equals to $\s$ on the ends of the interval. It is greater than $\s$ in the all inner points of the interval because $f(\phi)\equiv \s$ is strictly $\r$-trigonometrical function. If the length of the interval is greater than $\pi/\r$, for example this is $(-l/2,l/2)$ with $l>\pi/2\r,$  we cover it by intersecting intervals of length less then $\pi/\r,$ construct $h_I(\phi)$ as before for every interval $I$
and   set $h(\phi)=\max\limits_{I}h_I(\phi).$ It is obvious that $h(\phi)$ is
greater than $\s$ and it is $\r$- trigonometrically convex.
\qed
\edm
Theorem 5.7.4.3 is a corollary of Lemma 5.7.9.1 and
the following
\proclaim {Theorem 5.7.9.2} Let $h_1$ and $h_2$ be two $\r$ -trigonometrically convex functions. Then there exists a function $f\in A(\r(r))$ such that
$$h_f(\phi)=\max (h_1(\phi),h_2(\phi)),\  \underline h_f(\phi)=\min (h_1(\phi),h_2(\phi))$$
\ep
\demo {Proof} Consider the set
$$U:=\{v(z)=cr^\r h_1(\phi)+(1-c)r^\r h_2(\phi):0\leq c\leq 1\}\tag 5.7.9.3$$
It consists of invariant subharmonic functions,it is contained in $U[\r]$ and  satisfies the condition
of Theorem 4.1.4.1. Hence (Theorems 4.2.1.2, Corollary 5.3.1.5) there exists
a function $f\in A(\r(r))$ such that
$$\Fr f=U\tag 5.7.9.4$$
By formulae (3.2.1.1),(3.2.1.2) we obtain the assertion of the theorem, using
(5.7.9.3).
\qed
\edm
{\bf Exercise 5.7.9.1.} Prove Theorem 5.7.4.3.
\subheading {5.7.10} The family of characteristics $\{\Cal F_\a,\ \a\in A\}$
is called {\it independent} if for every subset $A'\sbt A$ (or subset in some class of subsets , for example, measurable or closed) there exists a function
$f=f_{A'}\in A(\r(r))$ such that
$$\underline {\Cal F}_\a[f]=\overline{\Cal F}_\a[f],\ \a\in A'$$
$$\underline {\Cal F}_\a[f]\neq\overline{\Cal F}_\a[f],\ \a\in A\setminus A'$$
It means that for every pointed subset of characteristics there exists a function that has  regular growth with respect to this subset of characteristics and is not of regular growth with respect to all other characteristics.

Theorem 5.7.4.3 can be considered as an assertion of independence of the family (5.7.2.2).
\proclaim {Theorem 5.7.10.1} The family $\chi_{Fo}$ (5.7.2.4) is independent,\ep
 i.e., for every $A\sbt \BZ$ there exists $f\in A(\r(r)$ such that
$$\liml_{r\ri\iy} r^{-\r(r)}\int_0^{2\pi} \log |f(re^{i\phi})g_k(\phi)d\phi
$$
exists for all $k\in A$ and does not exists for $k\in \BZ\setminus A.$
For beginning we prove
\proclaim {Lemma 5.7.10.2} There exist two $\r$-trigonometrically convex functions $h_1$ and $h_2$ for which
$$\int_0^{2\pi}h_1(\phi)g_k(\phi)d\phi=\int_0^{2\pi}h_2(\phi)g_k(\phi)d\phi,\
k\in A\tag 5.7.10.11$$
 $$\int_0^{2\pi}h_1(\phi)g_k(\phi)d\phi\neq\int_0^{2\pi}h_2(\phi)g_k(\phi)d\phi,\ k\in \BZ\setminus A\tag 5.7.10.12$$
 \ep
 \demo {Proof} Let $g(\phi)\in C^2$ be a function, the Fourier coefficients of which with indices $k\in A$ are equal to zero . We can represent it as a difference of $\r$-trigonometrically convex functions in the following way. Suppose for simplicity  that $\r$ is non-integer. Then take $T_\r g=g''+\r^2g$ and consider
$$h_1(\phi):=\frac {1}{2\r\sin\pi\r}\int_0^{2\pi}\widetilde {\cos\rho}(\phi-\psi-\pi)(T_\r g)^+(\phi)d\phi;$$
$$h_2(\phi):=\frac {1}{2\r\sin\pi\r}\int_0^{2\pi}\widetilde {\cos\rho}(\phi-\psi-\pi)(T_\r g)^-(\phi)d\phi.$$
By Th.3.2.3.3 these functions are $\r$ trigonometrically convex  and
$h_1-h_2=g.$ Hence (5.7.10.11) , (5.7.10.12) holds.
\qed
\edm
\demo {Proof of Th.5.7.10.1} We consider a function $f\in A(\r(r))$ with the limit set $U:=\{v(z)=cr^\r h_1(\phi)+(1-c)r^\r h_2(\phi):0\leq c\leq 1\ \}$ with
$h_1,h_2$ from Lemma and exploit Theorem 5.7.1.4.
\qed
\edm
{\bf Exercise 5.7.10.1.} Do this in details.
\newpage

\centerline {\bf 5.8. A generalization of Valiron- Titchmarsh theorem.}
\subheading {5.8.1}The point of depart on this topic is the following
\proclaim {Theorem VT}\text{ \cite {Va},\cite {Tich}}Let $f\in A(\r),\
\r<1$ have its zeros on the negative ray. If the limit
$$\liml_{r\ri\iy} r^{-\r}\log |f(r)|$$
exists, then the limit
$$ \liml_{r\ri\iy} r^{-\r}n(r)$$
exists.
\ep
The latter means that $f$ is a CRG-function.

The general problem is the following. Let $\r$ be any non-integer number,
$f\in A(\r(r))$ and
suppose all zeros of $f$ lie on a finite system of rays
$$K_{S_1}:=\{z=re^{i\phi}:0<r<\iy, \phi\in S_1\}\tag 5.8.1.1$$
where
$$S_1:=\{e^{i\th_j}:j=1,2,...,m\}\tag 5.8.1.2$$
We write $n_f\in \Cal M_{S_1}.$

Let $n_j$ be a zero distribution on the ray $\{\arg z=\th_j\}$ and all the limits
$$\liml_{r\ri\iy} r^{-\r}n_j(r):=\Delta_j\tag 5.8.1.3$$
exist. In such case we  write $n_f\in \Cal M_{reg,S_1}$.

Let $K_S$ be one more system of rays
$$S=\{e^{i\psi_k}:k=1,2,...,n\}\tag 5.8.1.4$$
Some $\psi_k$ can coincide with some $\th_j.$
Suppose that $f$ has regular growth on this system, i.e.,
$$ h_f(\phi)=\underline h_f(\phi), \ e^{i\phi}\in S\tag 5.8.1.5 $$
In such case we write $f\in A_{reg, S}.$

The problem is, what is the connection between $S$ and $S_1$ so that the implication
$(f\in A_{reg, S})\Longrightarrow (n_f\in \Cal M_{reg,S_1})$ holds.

This problem can be reformulated in another way. For $n_f\in \Cal M_{reg,S_1}$
if $n_f\in \Cal M_{S_1}$ it is necessary and sufficient that $f$ is a CRG -function, because existence the angle density is equivalent to  the existence of all the
limits. So the problem can be reformulated in the form:
what is connection between $S$ and $S_1,$ so that the implication
$(f\in A_{reg, S})\Longrightarrow (f\  is\ CRG-function)$  holds.

Denote
$$G(t,\g,\r):=G_p(e^{t-i\g})e^{-\r t},\ p=[\r]\tag 5.8.1.6$$
where $G_p$ is the Primary Kernel:
$$G_p(z)= \log |1-z|+\Re \sum_{k=1}^p z^k/k $$
Set
$$\hat G(s,\g,\r):=\intl_{-\iy}^{\iy}G(t,\g,\r)e^{-ist}dt$$
This is the Fourier transformation of $G(t,\g,\r)$. It can be computed
(see, e.g,\cite {Oz,Lemma3})
  $$\hat G(s,\g,\r)=\frac {\pi\cos (\pi +\g)(\r +is)}{(\r + is)\sin\pi(\r + is)} $$
Consider the matrix
$$\hat \BG(s,S_1-S):=\|\hat G(s,\th_j -\psi_k,\r)\|.\tag 5.8.1.7$$
We are going to prove
\proclaim {Theorem 5.8.1.1}The implication
$$\{f\in A_{reg,S}\}\wedge \{n_f\in \Cal M_{S_1}\}\Longrightarrow \{f\ is\ a \ CRG-function\}$$ holds iff
$$\text{rank}\ \hat \BG(s,S_1-S)=m,\ \forall s\in (-\iy,\iy)\tag 5.8.1.5 $$

\ep
As a corollary we obtain the following (\cite {De})
\proclaim {Theorem 5.8.1.2.(Delange)} Suppose that $S_1$ and $S$
consist of  one ray,i.e.,
$$S_1=\{e^{i\th_1}\},\ S=\{e^{i\psi_1}\}.$$
The implication (5.8.1.5) holds iff
$$\th_1-\psi_1\neq(1-(2k+1)/2\r)\pi,\ k=1,2,...\tag 5.8.1.6$$
\ep
\subheading {5.8.2} A Fourier transformation for distribution $\nu$ on the real axes is a distribution in  the standard space $\Cal S'$ (see \cite {H\"o, v.1,Ch.7,\S 7.1} For locally bounded measure  which variation  ``not very quickly''
growing it can be defined by
$$(\Cal F\nu)(s):=\liml_{\eps\ri 0}\intl_{-\iy}^{\iy}e^{its}e^
{-\frac{\eps t^2}{2}}\nu (dt).$$
where the right side is understood in the sense of distributions.

For example, if $\nu(dt):=e^{is_0t}dt,$ we have $\Cal F \nu (s)=\delta
(s-s_0)$
where $\delta$ is the Dirac function.

{\bf Exercise 5.8.2.1.} Check this.

For  distribution and a summable function one can b define the convolution  , for which the property
$\Cal F(f*\nu)(s)=\Cal F f(s)\Cal F\nu (s)
$ holds.

\demo {Proof of Theorem 5.8.1.1.}
Since $f\in \Cal M_{S_1}$ the limit set $\Fr n_f$ is concentrated on $K_{S_1}.$
So every $v\in \Fr[f] $ can be represented in the form (see Th.3.1.5.1):
$$v(z)=\sum\limits_{j=1}^{j=m}\intl_{0}^{\iy}G_p(z/re^{i\th_j})\mu_j(dr)\tag 5.8.2.1$$
where $\mu_j$ is concentrated on the ray $\{\arg \z=\th_j\}$ and belong to
$U[\r].$ After changing variables:
$$ r=e^t, |z|=e^t$$ we obtain from (5.8.2.1)
$$v^1(te^{i\phi})=\sum\limits_{j=1}^{j=m}\intl_{-\iy}^{\iy}
G(t-\tau,\phi-\th_j,\r) \mu^1_j(d\tau)\tag 5.8.2.1$$
where
$$e^{\r \tau}\mu^1(d\tau):=\mu_j(dr),\ v^1(te^{i\phi}):=v(|z|e^{i\phi})
e^{-\r |z|}.\tag 5.8.2.2$$
The equality (5.8.2.1) can be written as
$$ v^1(te^{i\phi})=\sum\limits_{j=1}^{j=m}[G(\bullet,\phi-\th_j,\r)*\mu^1_j](t)\tag 5.8.2.3$$
where * stands for convolution.
Then $f\in A_{reg,S}$ with $n_f\in \Cal M_{S_1},$ iff every pair
$v_1,v_2 \in \Fr[f]$ satisfies the condition
$$v_1(z)=v_2(z), z\in K_S\tag 5.8.2.4$$
Denote by $\mu_{1,j},\mu_{2,j}$ the restriction of $ \mu_{v_1} \mu_{v_1},$
to the ray $\{\arg z=\th_j\}. $ Set  $\nu_j:=\mu_{1,j}-\mu_{2,j}$
Using (5.8.2.3) we can rewrite (5.8.2.4) in the form

$$\sum\limits_{j=1}^{j=m}[G(\bullet,\phi_k-\th_j,\r)*\nu^1_j](t)\equiv 0 \ , k=1,2,...,n.\tag 5.8.2.5$$
Applying Fourier transforms we obtain a system of linear equations:
$$\sum\limits_{j=1}^{j=m}[\hat G(\bullet,\psi_k-\th_j,\r)\cdot\hat \nu^1_j](t)
\equiv 0 \ , k=1,2,...,n.\tag 5.8.2.6$$

Suppose now that rank $\hat \BG(s,S-S_1)=m$ for every $s\in \BR.$ The system (5.8.2.6) has only the trivial solution for every $s$. Thus
$\hat \nu^1_j(s)\equiv 0,$ for $j=1,2,..., m.$ This implies
$\nu^1_j (t)\equiv 0$ for $j=1,2,...,m$ and $\nu_j\equiv 0$ for $j=1,2,...,m. $
Thus $\mu_{v_1}=\mu_{v_2}$, i.e, (by (5.8.2.3) )$\Fr [f]$ consists of one function $v\in U[\r].$ Thus $f$ is a CRG-function.

Conversely, suppose that rank $\hat \BG(s,S-S_1)<m$ for some $s_0$

Then there exists a nontrivial solution $(b_,...b_m)$ that satisfies the
corresponding system. We obtain that $\{\hat \nu^1_j b_j\dl(s-s_0),\
j=1,2...m\}$
is a solution of (5.8.2.6) for all $s\in\BR$ and hence the
$$\nu^1_j(dt)=b_je^{its_0}dt,\ j=1,2,...m$$
Since $\nu^1_j$ have  bounded densities $d\nu^1_j/dt,$ we can find a constant $C$ such that $\sup\{|d\nu^1_j/dt|:0<t<\iy,\ j=1,2... m\}\leq C.$

Set $$\mu^1_{1,j}(dt)=Cdt+\nu^1_j(dt);\ \mu^1_{2,j}=Cdt .\tag 5.8.2.7$$
Both of these are measures. Now we pass to $m_{1,j},\ m_{2 ,j}$ via (5.8.2.2).It is easy to check that $m_{1,j},\ m_{2 ,j}\in \Cal M(\r).$

{\bf Exercise 5.8.2.2.} Check this.

Consider $\mu _1,\mu_2 \in \Cal M (\r)$ which are defined uniquely by their restrictions $\mu_{1,j},\mu_{2,j}$ respectively on $K_{S_1}.$ Set
$$v_1(z):=\intl_{\BC}G_p(z/\z)\mu_1(d\xi d\eta);
v_2(z):=\intl_{\BC}G_p(z/\z)\mu_2(d\xi d\eta);\ \z=\xi+i\eta.$$
It is easy to check that the equality
$$v_1(z)=v_2(z), \ z\in K_S\tag 5.8.2.8$$
holds.

{\bf Exercise 5.8.2.3.} Check this.

Since $\mu_1$ and $\mu_2$ is a finite sum of trigonometrical functions, for
$v_1$ and $v_2$  the condition (4.1.3.3) is satisfied. Thus by Theorem 4.3.6.1
there exists a function $f\in A(\r(r))$ for which
$$\Fr[f]=\bigcup\limits_{0\leq c\leq 1} \BC (cv_1+(1-c)v_2)$$
Since for $v\in \BC (cv_1+(1-c)v_2) $ (5.8.2.8) also holds, the same
holds for $v\in\Fr [f]$ and this function is not a CRG-function.
\qed
\edm
\newpage

\centerline {\bf 6. Application to the completeness of exponential
systems} \centerline {\bf in convex domains and the multiplicator
problem} \vskip .10in The completeness of exponential systems in
convex domains is intimately connected to the multiplicator
problem. Considering a special form of exponent system is related
to the study of special subharmonic functions, that determine the
periodic limit set, so called automorphic subharmonic functions.
The next \S\S 6.1, 6.2 are devoted to these problems.

\centerline {\bf 6.1. Problem of multiplicator.}

\subheading {6.1.1} Let $\Phi\in A(\r(r))$ and let $H(\phi)$ be a $\r$-trigonometrically convex function. A function $g\in A(\r(r)) $ is called an $H$-multiplicator of $\Phi$ if the indicator $h_{g\Phi}$ of the product $g\Phi$ satisfies the inequality
$$h_{g\Phi}(\phi)\leq H(\phi),\forall \phi.$$
In some questions (see \S 6.3) we need to determine whether a given  function $\Phi$ has a multiplicator.We shall study this problem in terms of the limit set of $\Phi.$
Define $H(z):=r^\r H(\phi),\ z=re^{i\phi}.$ Let $v\in U[\r]$ (see (3.1.2.4)). Consider the function
$$m(z,v,H):=H(z)-v(z).$$
As  will be proved in Corollary 6.1.9.3, the maximal subharmonic minorant of
$m(z,v,H)$ exists and is continuous.
The maximal subharmonic minorant of $m$ (m.s.m.) belonging to $U[\r]$ will be denoted
by $\Cal G_Hv,$ while the domain of definition of the operator $\Cal G_H$ will
be denoted by $D_H.$ Though $m(0,\bullet,\bullet)=0$, the m.s.m. of $m$ can differ from zero (as was remarked by A.E.Eremenko and M.L.Sodin), but if the m.s.m. of $m$ equals zero at zero , then it belongs to $U[\r].$

{\bf Exercise 6.1.1.1.} Consider the function
$$w(z)=\cases |z|\log |z|, &\text {if}\  |z|\leq 1;\\
|z|-1, &\text {if}\  |z|\geq 1\endcases.$$ It is subharmonic and belongs to $U[1].$ Show that the  maximal subharmonic minorant of $K|z|-w(z)$ is different from zero in 0 for every $K>0.$

\proclaim {Theorem 6.1.1.1} $\Phi\in A(\r(r))$ has an $H$-multiplicator iff
$$\Fr[\Phi]\sbt D_H.\tag 6.1.1.1$$
\ep
\demo {Proof of necessity} Let $g$ be a multiplicator of $\Phi$, i.e.,
$$h_{g\Phi}(\phi)\leq H(\phi)\tag 6.1.1.2$$
and let $v\in \Fr[\Phi].$ We can choose $v_{g\Phi}\in \Fr [g\Phi]$ and
$v_g\in \Fr[g]$ such that $v_{g\Phi}=v+v_g$ (see Th.3.1.2.4, fru1)).

{\bf Exercise 6.1.1.2.} Prove this directly.

By definition of indicator (3.2.1.1) and (6.1.1.2) we have
$v_{g\Phi}(z)\leq H(z)$ or $v_g(z)\leq m(z,v,H).$
Since $v_g\in U[\r],$ $v\in D_H.$
\qed
\edm
For proving sufficiency we need the following
\proclaim {Theorem 6.1.1.2} The operator $\Cal G_H$ is

1. upper semicontinuous in the $\Di'$ -topology, 6.1.1.5i.e.,
$$(v_j\ri v)\wedge (\Cal G_Hv_j\ri w)\Longrightarrow (w\in U[\r])\wedge
(w(z)\leq \Cal G_H (z), z\in \BC);$$

2.invariant: $(\Cal G_H v)_{[t]}= \Cal G_H v_{[t]};$ (see (3.1.2.4a) for
$P\equiv I;$)

3.concave:
$$(\forall v_1,v_2\in D_H ,\  c\in [0;1])\Longrightarrow (v_c:=cv_1+(1-c)v_2\in D_H )$$
and
$$ \Cal G_H (v_c)\geq c\Cal G_H (v_1)+(1-c)\Cal G_H (v_2).$$
\ep
\demo {Proof} Let us prove 1). Suppose $v_j\in U[\r]\ri v$ and $\Cal G_H v_j\ri w.$
Then $$ \Cal G_H v_j\leq H(z)-v_j(z),\ z\in \BC, \tag 6.1.1.3$$
Applying $(\bullet)_\eps$ from (2.6.2.2) and Th.2.3.4.5 reg 3), we obtain
$$w_\eps\leq (H)_\eps (z)-(v)_\eps (z),\  z\in \BC\ w_\eps(0)\geq 0\ $$
Passing to the limit as $\eps\downarrow 0$ we obtain by Th.2.6.2.3, ap2)
$$w(z)\leq H(z)-v(z)=m(z,H,v),z\in \BC.$$
Since  $0\leq w(0)\leq m(0,H,v)=0$ we have $w(0)=0$ and , hence, $w\in U[\r].$
Thus $v\in D_H$ and $w(z)\leq\Cal G_H v(z)$.

Let us prove 2). Since $H(z)$ is invariant with respect to $(\bullet)_{[t]},$
$$ (\Cal Gv)_{[t]}\leq H(z)-v_{[t]}.$$
Hence,
$$(\Cal Gv)_{[t]}(z)\leq( \Cal G(v_{[t]}))(z),\tag 6.1.1.4$$
because $\Cal G(v_{[t]})$ is   m a x i m a l subharmonic minorant.
We can replace $v$ with $v_{[1/t]}$  and obtain $(\Cal
Gv_{[1/t]})_{[t]}(z)\leq \Cal Gv(z).$ Applying $(\bullet)_{[t]}$
to the two sides of the inequality, we obtain $\Cal
Gv_{[1/t]}(z)\leq (\Cal Gv(z))_{[1/t]}.$ Now we can replace $1/t$
with $t$ and obtain the reverse inequality  to (6.1.1.4), which,
 together  with
(6.1.1.4),  proves 2).

 3). Let $v_1,v_2 \in D_H$ and $c\in [0;1].$ One has
$$\Cal Gv_i (z)\leq H(z)-v_i(z),\ i=1,2,\ \forall z.$$
Then
$$[c\Cal Gv_1 +(1-c)\Cal Gv_2](z)\leq H(z)-[cv_1+(1-c)v_2](z).$$
Thus
$$[c\Cal Gv_1 +(1-c)\Cal Gv_2](z)\leq \Cal G[cv_1+(1-c)v_2](z).$$
\qed  \edm

\demo {Proof of sufficiency in Theorem 6.1.1.1} Assume that $\Fr[\Phi]\sbt
D_H$ and consider the set
$$\bold U:=\{(v',v''):v''\leq \Cal Gv', v'\in \Fr[\Phi]\}.\tag 6.1.1.5 $$
Then $\bold U$ is non-empty, because of (6.1.1.1), closed, because of
Th. 6.1.1.2, 1), and invariant, because of Th.6.1.1.2, 2).

 Every fiber
$\bold U''=\{v'': v''\leq \Cal Gv'\}$ is convex because of Th.6.1.1.2, 3).
By Th. 4.4.1.2 there exists $u''\in U(\r(r))$ such that for the curve
$\bold u:= (u',u'')$
 $$\Fr[\bold u]=\bold U.\tag 6.1.1.6$$

By the Th.5.3.1.4 (Approximation Theorem) the function $u''$ can be replaced with $\log|g|,$ where $g\in A(\r(r)),$ retaining the property (6.1.1.6).

Let us prove that $g$ is an $H$-multiplicator of $\Phi.$ Indeed, set
$\Pi:=g\Phi.$ It is enough to prove that for every $v_\Pi\in \Fr [\Pi]$
$$v_\Pi(z)\leq H(z)\tag 6.1.1.7$$
Note that every $v_\Pi$ has the form
$v_\Pi=v_g+v,$ where $(v,v_g)\in \bold U.$ Thus, because of the definition (6.1.1.5)
$v_\Pi$ satisfies (6.1.1.7).
\qed
\edm

Let us note that the pair $(v,\CG_Hv)\in \bold U$ because of closeness of
$\bold U.$ Hence the following assertion holds
\proclaim {Proposition 6.1.1.3} Every $\Phi\in A(\r)$ that satisfies (6.1.1.1)
has a multiplicator $g\in A(\r)$ such that
$$v+\CG_Hv\in \Fr [g\Phi].\tag 6.1.1.8$$\ep

{\bf Exercise 6.1.1.3.} Check this in details.

Although $v\in U[\r]$ is in general an upper semicontinuous function,

\proclaim {Theorem 6.1.1.4} The function $\CG_Hv(z),\ v\in U[\r], $ is a continuous function that is harmonic outside the set $E= \{z:\CG_Hv(z)=m(z,v,H)\}.$\ep

\demo {Proof}$\CG_Gv(z)$ is continuous because of Corollary 6.1.9.3. If $\CG_Hv(z_0)<v(z_0)$ and if $\CG_Hv(z)$ is not harmonic in a neighborhood of $z_0,$ we can make sweeping of masses from a small disc $\{|z-z_0|<\eps\}$
(see Th.2.7.2.1). The obtained subharmonic function will be grater than
$\CG_Hv(z),$ contradicting maximality.
\qed\edm

\subheading {6.1.2}Suppose that some $H$-multiplicator $g=g(z,\Phi,H)$ of the function $\Phi$ is found. We examine the function $\Pi=g\Phi$.The structure of its limit set is described by the following statement:
\proclaim {Proposition 6.1.2.1} Every $v_\Pi\in \Fr[g\Phi]$ can be written as
$v_\Pi=v+w_1,$ where $v\in\Fr[\Phi]$ and $w_1\in U[\r]$ with the condition
$$w_1(z)\leq \CG _H(z), \forall z\in\BC, \tag 6.1.2.1$$
and, conversely, for every $v\in \Fr[\Phi]$ there exists a
$v_g,\ v_g(z)\leq \CG_Hv(z),$ such that
$$v+v_g\in \Fr[g\Phi].$$.
\ep

{\bf Exercise 6.1.2.1.} Prove this, like in Exercise 6.1.1.1.

An $H$-multiplicator $G$ of the function $\Phi$ will be called {\it ideally
complementing} if it satisfies the condition
$$\Fr[G\Phi]=\{v_\Pi=v+\CG_Hv:v\in \Fr[\Phi]\}.$$
If a multiplicator is ideally complementing then equality is achieved in
(6.1.2.1) for all $v\in\Fr[\Phi].$ This make the multiplicator optimal in
another respect. Recall that an entire function $f$ is of {\it minimal type}
with respect to a proximate order $\r(r),\r(r)\ri \r$ if (see (2.8.1.6))
$$\s_f:=\limsupl_{r\ri\iy} \log M(r,f)r^{-\r(r)}=0.$$
\proclaim {Proposition 6.1.2.2} Let $G=G(\bullet,\Phi,H)$ be an ideally
complementing $H$-multiplicator of a function $\Phi$.Then each $H$-multiplicator of  the function $\Pi=G\Phi$ is of minimal type.
\ep
This proposition is proved in \S 6.1.3.

A function $\Phi$ is said to be ideally complementable if for each $H$ the
condition (6.1.1.1) implies that $\Phi$ has an ideally complementing multiplicator.
For instance, if $\Phi$ is a function of completely regular growth (see \S 5.6)
then it is ideally complementable.

{\bf Exercise 6.1.2.2.} Prove this.

\proclaim {Theorem 6.1.2.3} Every function with periodic limit set is ideally
complementable.
\ep
This theorem is proved in \S 6.1.6.

Let $C\sbt \BR^l$ be an $l$-dimensional connected compact and  let $\{h(\phi,c):
c\in C\}$ be a set of $\r$-t.c. functions that is continuous   with respect to $c\in C.$ For example, $c\in [0,1]$ and $h(\phi,c)=ch_1(\phi)+(1-c)h_2(\phi).$
The set
$$U_{ind}:=\{v(re^{i\phi}s)=r^\r h(\phi,c):c\in C\}\tag 6.1.2.2 $$
is the limit set of an entire function.

{\bf Exercise 6.1.2.3} Prove this using  Theorem 4.3.6.1.

Such a set is called a {\it set of indicators}. Entire functions with such limit sets can be also considered as a generalization of CRG-functions.
\proclaim {Theorem 6.1.2.4} Every function $\Phi$ whose limit set is a set of indicators is ideally complementable.
\ep
This theorem is proved in \S 6.1.7.

The existence  of an ideally complementing $H$- multiplicator depends, of course , both on $\Phi\in A(\r(r))$ (or, more precisely, on its limit
set $\Fr [\Phi]$) and on $H.$
\proclaim {Theorem 6.1.2.5} Let $\Phi$ and $H$ be such that the condition        (6.1.1.1) is satisfied.The function $\Phi$ has an ideally complementing
$H$-multiplicator if and only if the operator $\CG_H$ is continuous on
$\Fr \Phi.$
\ep
This theorem is proved in \S 6.1.6.

Now we formulate a sufficient condition for continuity of the operator
$\CG_H.$ We shall say that {\it the maximum principle for $U[\r]$ is valid
in the domain} $G, $ (which is, generally speaking,  unbounded), if the conditions $w\in U[\r],\ w(z)=0$ for $z\notin G$ imply $w(z)\equiv 0.$

Let us denote by $\CH_w$ a region of harmonicity of $w\in U[\r]$, i.e.,
a region where the conditions ``$w$ is harmonic in $G$'' and ``$G\spt \CH_w$'' imply
 $G=\CH_w.$

We remark that $\CH_w$ is a connected component of the open set on which
$w$ is harmonic. Generally it is not unique.

The image of $U\in U[\r]$ will be denoted by  $\CG_H U$, while its closure in
the $\Di'$ -topology will be denoted by clos $\CG_H U.$
\proclaim {Theorem 6.1.2.6} Suppose for every $w\in$clos $\CG_H U$ and every
$\CH_w$ the maximum principle for $U[\r]$ holds. Then $\CG_H$ is continuous on $U.$
\ep
This theorem is proved in \S 6.1.5.

In \S 6.1.8 we will construct an example of $\Phi$ and $H$ such that the operator
$\CG_H$ is not continuous on $\Fr[\Phi].$ This is also an example of the function that has not ideally complementing multiplicator.
\subheading {6.1.3} \demo {Proof of Proposition 6.1.2.2} Let $g$ be an ideally
complementing multiplicator of the function $\Pi=G\Phi.$ We denote
$$\th(z):=(gG\Phi)(z).\tag 6.1.3.1$$
Let $v_g\in \Fr[g].$ Let us choose $t_j\ri\iy$ such that:
$$(\log |g|)_{t_j}\ri v_g;\ (\log|\Pi|)_{t_j}\ri v_\Pi\in \Fr[\Pi];\
(\log |\th|)_{t_j}\ri v_\th\in \Fr[\th].$$
It follows from (6.1.3.1) that $v_\th=v_g+v_\Pi$. Since $g$ is a multiplicator of $\Pi,$ we have
$$v_\th(z)=v_g(z)+v_\Pi(z)\leq H(z).\tag 6.1.3.2$$
Since $G$ is an ideally complementing multiplicator, $v_\Pi=v+\CG_Hv$. So for
all $z\in\BC$ (6.1.3.2) implies
$$(v_g+\CG_Hv)(z)\leq (H-v)(z).$$
Since $\CG_Hv$ is the maximal subharmonic minorant, $v_g(z)\leq 0$ and hence
$v_g(z)\equiv 0.$ Thus (see (3.2.1.1)) we have $h_g(\phi)\equiv 0$ and
therefore
$$\s_g=\max\limits_{0\leq \phi\leq 2\pi} h_g(\phi)=0$$
\qed
\edm
\subheading {6.1.4} In order to prove Theorem 6.1.2.6 we need a number of auxiliary statements.
\proclaim {Lemma 6.1.4.1}Let the maximum principle be valid in $G$ for $U[\r]$
and for some continuous functions $w_1,w\in U[\r]$
satisfy:

a)w is harmonic in $G;$

b)$w_1(z)=w(z)$ outside of $G.$

 Then
$$w_1(z)\leq w(z),\ z\in G.\tag 6.1.4.1$$
\ep
\demo {Proof} We set
$$w_0(z):=\cases (w_1-w)^+(z),&\ z\in G\\0,&\ z\notin G.\endcases$$
This function is continuous in $\BC$ and, evidently, subharmonic both in $G$
and in $\BC\setminus \overline G.$ Since $w_0(z)\geq 0,$ the inequality for the mean values
$$0=w_0(z)\leq\frac{1}{2\pi}\intl_0^{2\pi}w_0(z+\eps re^{i\phi})d\phi,\ z\in \partial G,$$
implies the subharmonicity on $\p G.$ Since $G$ satisfies the maximum
principle for $U[\r]$, we have $w_0\equiv 0$, which is equivalent to (6.1.4.1).
\qed
\edm
Now we shall dwell on some properties of maximal subharmonic minorants and,
in particular, of $w= \CG_H v$. Let
$$E_v:=\{z\in\BC:\CG_Hv(z)=m(z,v,H)\}\tag 6.1.4.2$$
We remark that $m(z,v,H)$ is a $\dl$ -subharmonic function in $\BC$ whose
charge will be denoted by $\nu(\bullet,v)$, its positive and negative parts will be denoted by $\nu^+$ and $\nu^-$.

Let us denote by $\mu_H$ the measure of $H(z),$ It is decomposed into the product of measures (see \S 3.2)
$$\mu_h=\Delta_H\otimes\r r^{\r-1}dr,\tag 6.1.4.3$$
where $\Delta_H$ is the measure on the unit circle and  $\r r^{\r-1}dr$ is the measure on the ray. It is obvious that
$$\nu^+(\bullet,v)\leq\mu_H(\bullet).\tag 6.1.4.4$$
We shall denote the mass distribution of $w\in U[\r]$ by $\mu_w.$

The modulus of continuity of $w$ (if $w$ is continuous) will be denoted
$\om _w(z,h),\ z\in\BC,\ h>0$.

The following lemma lists various properties of $w\in \CG_H U,\  U\sbt U[\r]$ which will be useful in the sequel:
 \proclaim {Lemma 6.1.4.2}Let $w\in\CG_HU$.Then

1.$w\in U[\r,\s]$ where
$$\s=4\cdot 2^\r[\max\{H(e^{i\phi}):\phi\in [0,2\pi]\}+2\s_1]$$
$$\s_1=\max\{v(z)|z|^{-\r}:z\in\BC, v\in U\};$$

2.the charge restriction $\nu(\bullet,v)|_{E_v}$ to $E_v$ is nonnegative, i.e.,
$$\nu(\bullet,v)|_{E_v}= \nu^+(\bullet,v)|_{E_v};$$

3. outside ${E_v}$ the function $w$ is harmonic, i.e.,
$$\mu_w|_{\BC\setminus E_v}=0;$$

4.the measure $\mu_w$ is bounded from above by $\nu^+(\bul,v),$ i.e.,
$$\mu_w\leq\nu^+(\bul,v);$$

5. $\CG_H U$ is equicontinuous on each compact set, i.e.,
$$\om_w(z,h)\leq C(R,\s,\r)\sqrt h\log(1/h),\ |z|\leq R,$$
where $C(R,\s,\r)$ is independent of $w\in \CG_H U.$
\ep
\demo {Proof} Let us prove property 1. We have
$$T(r,w):=\frac {1}{2\pi}\intl_0^{2\pi}w^+(re^{i\phi})d\phi$$
$$\leq \frac{1}{2\pi}\left [r^\r\intl_0^{2\pi}H^+(e^{i\phi}) d\phi +
\intl_0^{2\pi}v^+(re^{i\phi})d\phi + \intl_0^{2\pi}v^-(re^{i\phi})d\phi\right].$$
Since $v(0)=0$, we have
    $$\intl_0^{2\pi}v^-(re^{i\phi})d\phi \leq \intl_0^{2\pi}
v^+(re^{i\phi})d\phi.$$
Therefore
    $$T(r,w)\leq[\max\{H(e^{i\phi}):\phi\in[0,2\phi]\}+2\s_1]r^\r.$$
It is known (see Theorem 2.8.2.3, (2.8.2.5)) that $M(r)\leq 4T(2r)$. So
we conclude that
$$w(z)\leq4\cdot 2^\r[\max \{H(e^{i\phi};\phi \in[0,2\pi]\} +2\s_1]|z|^\r=\s |z|^\r.$$
Let us prove property 2. To this end we shall use the following theorem
(Grishin's Lemma)\cite {Gr}
\proclaim {Theorem A.F.G} Let $g$ be a nonnegative  $\dl$-subharmonic
function, and let $\nu_g$ be its charge. Then the restriction $\nu_g|_E$ to the
set $E=\{z:g(z)=0\}$ is a measure.
\ep
Applying this theorem to the function $g:=m(z,v,H)-\CG_Hv(z)$, we get
$$\nu(\bul,v)|_{E_v}\geq\mu_w|_{E_v},\tag 6.1.4.5$$
hence, we obtain property 2.

Let us prove property 3. Since $w$ and $v$ are upper semicontinuous, and $H$
is continuous (see Th.3.2.5.5) , the set $\{z:(w+v)(z)-H(z)<0\}$ is open.

Let us take a neighborhood of an arbitrary point of this set and replace the function $w$ within it with the Poisson integral constructed using this function, i.e., let
us sweep out the mass from this neighborhood. The subharmonic function obtained would be strictly greater than the initial one, if the latter were not harmonic. This means that the initial $w$ was not the maximal minorant. We have arrived at a contradiction, which proves property 3.

Property 4. immediately follows from property 3. and (6.1.4.5).

In order to prove property 5. we shall need auxiliary statement which will be stated as lemmas. Let
$$P(z,\phi,R):=\frac {1}{2\pi}\frac {R^2-|z|^2}{|z-Re^{i\phi}|}$$
be the Poisson kernel in the disc $K_R=\{z:|z|<R\}$.

Below  $C'$s with indices will denote constants.
\proclaim {Lemma 6.1.4.3} In the disc $K_{R/2},$ we have
$$|grad_zP(z,\phi,R)|\leq C_1(R),$$
where $C_1(R)$ depends only on $R.$
\ep

{\bf Exercise 6.1.4.1.} Prove this.

We shall introduce the notation for the Green function for the Laplace operator in the disc $K_R:$
$$G(z,\z,R):=\log\left |\frac {R^2-\z\overline z}{R(z-\z)}\right |.$$
The disc $\{\z:|\z-z|<t\}$ will be denoted by $K_{z,t}.$
\proclaim {Lemma 6.1.4.4} Let $z\in K_{R/2}\setminus K_{\z,\sqrt h}.$ Then for a small $h$
$$|grad_zG(z,\z,r)|\leq C_2(R)/\sqrt h$$
\ep
{\bf Exercise 6.1.4.2.} Prove this.

Let us denote $\mu(z,t):=\mu(K_{z,t}).$
\proclaim {Lemma 6.1.4.5} For $z\in K_{R/2},0<t<R/10,$ we have
$$\mu_H(z,t)\leq C_3(\s,R)t.$$
 \ep
  \demo {Proof} Applying Th.2.6.5.1 (Jensen-Privalov) to the function $H(z),$
we obtain
$$M_H=\max\{H(e^{i\phi}:\phi\in [0;2\pi]\}=\Delta_H(\BT)/\r$$
where $\BT$ is the unit circle.

Now
$$\mu(z,t)\leq \Delta_H(\BT)\intl_{|z|-t}^{|z|+t} r^{\r-1}dr \leq
\r^2M_H R^{\r-1}t \leq \s C(\r)R^{\r-1}t$$
where $C(\r)$ is a constant depending only on $\r$.
This proves the lemma.\qed
\edm
\proclaim {Lemma 6.1.4.6} Let $h<1$ and suppose that a monotonic function
$\mu(t)$ satisfies the condition
$$\mu(t)<ct\tag 6.1.4.6$$
Then
$$\intl_{0}^{\sqrt h}\log(1/t)\mu(dt)\leq (3/2)c\sqrt h\log h$$
\ep
{\bf Exercise 6.1.4.3.} Prove this  integrating by parts and using (6.1.4.6).

\proclaim {Lemma 6.1.4.7}Let $z\in K_{R/2}$ and $\z\in K_R.$ Then
$$|\log|(R^2-\z\overline z/R||\leq C_4(R)$$
\ep

{\bf Exercise 6.1.4.4.} Prove this.

Now we pass to the proof of assertion 5. from Lemma 6.1.4.2. According to the F.Riesz theorem (Th.2.6.4.3) we represent $w$ in the circle  as
$$w(z)=H(z,w)-\intl_{K_R}G(z,\z,R)\mu_w(d\xi d\eta), \ \z=\xi+i\eta,\tag 6.1.4.7$$
where
$$H(z,w)=\frac{1}{2\pi}\intl_0^{2\pi}P(z,\phi,R)w(Re^{i\phi})d\phi.$$
It follows from Lemma 6.1.4.3 and  1. of Lemma 6.1.4.2 that
$$|grad_z H(z,w)|\leq C_1(R)\frac{1}{2\pi}\intl_0^{2\pi}|w|(Re^{i\phi})d\phi\leq C_1(R)2\s R^\r . \tag 6.1.4.8$$
We split the integral in (6.1.4.7) into three terms
$$\psi_1(z,h):=\intl_{K_r\setminus K_{z_0,\sqrt h}}G(z,\z,R)\mu_w(d\xi d\eta),$$
$$\psi_2(z,h)=\intl_{ K_{z_0,\sqrt h}}\log|(R^2-\overline z\z)/R|\mu(d\xi d\eta)$$
$$\psi_3(z,h)\intl_{ K_{z_0,\sqrt h}}\log|\z-z|\mu(d\xi d\eta),$$
where $z_0$ is an arbitrary fixed point in $K_{R/2}.$

Combining property 4. and inequality (6.1.4.4) we have
$$\mu_w(E)\leq \mu_H(E),\ \forall E\sbt K_R.\tag 6.1.4.9$$
For all $z\in K_{z_0,{\sqrt h}/2}$ Lemma 6.1.4.4 yields
$$|\operatorname {grad} \psi_1(z,h)|\leq C_2(r)\s R^\r/\sqrt h.\tag 6.1.4.10$$Combining Lemmas 6.1.4.5 and 6.1.4.7 with inequality (6.1.4.9), we get
$$|\psi_2(z,h)|\leq C_4(R)C_3(\s,R)\sqrt h.\tag 6.1.4.11$$
Further, from Lemmas 6.1.4.5 and 6.1.4.6, taking into account the fact \linebreak that
$\log|\z-z|<0,$ we obtain
$$|\psi_3(z,h)|\leq(3/2)C_3(\s,R)\sqrt h\log h.\tag 6.1.4.12$$
Now consider the difference
$$\Delta w:=w(z_0+\Delta z)-w(z_0),\ |\Delta z|<h<{\sqrt h}/2.$$
It can be represented as
$$\Delta w=\Delta\psi_1+\Delta\psi_2+\Delta\psi_3 +\Delta H(z,w).\tag 6.1.4.13$$
Choosing $h$ small enough, one may assume that $z_0+\Delta z\in K_{{\sqrt h}/2,z_0}.$Thus, according to (6.1.4.11),
$$|\Delta\psi_2(z_0,h)|\leq|\psi_2(z_0,h)|+|\psi_2(z_0+h,h)|\leq C_6(\s,R)\sqrt h.\tag 6.1.4.14$$
Likewise (6.1.4.12) yields
$$|\Delta\psi_3(z_0,h)|\leq|\psi_3(z_0,h)|+|\psi_3(z_0+h,h)|\leq  C_7(\s,R)\sqrt h\log h.\tag 6.1.4.15$$
Finally, from (6.1.4.10) and (6.1.4.8), respectively, we obtain
$$|\Delta \psi_1|\leq C_3(\s,R)\sqrt h,\ |\Delta H(z_0,w)|\leq C_8(\s,R)h.\tag 6.1.4.16$$
Substituting (6.1.4.14)-(6.1.4.16) into (6.1.4.13), we obtain relation 5. of Lemma 6.1.4.2.
\qed
\edm
The lemma is proved.
\subheading {6.1.5} In this item we are going to prove Theorem 6.1.2.6. However, before that, we prove
\proclaim {Lemma 6.1.5.1} Let $w_n=\CG_Hv_n,\ v_n\overset{\Di'}\to{\ri}v$ and
$w_n\overset{\Di'}\to{\ri}w_\iy$. Set
$$E_\iy:=\{z:w_\iy(z)=H(z)-v(z)\}.$$
Then $w_\iy$ is harmonic in $\BC\setminus E_\iy.$
 \ep
Let us note that $w_\iy,$ in general, is not the maximal subharmonic minorant
because the operator $\CG_H$ can  be only upper semicontinuous, as  will be demonstrated by
example in \S 6.1.8. However, $w_\iy$ is a minorant of $H-v$ because of Th. 6.1.1.2, 1.
\demo{Proof} Let $z_0\notin E_\iy.$ Then there exists a $\dl>0$ such that
$$w_\iy(z_0)+v(z_0)\leq H(z_0)-2\dl$$
Since the function $b(z):=w_\iy+v(z)-H(z)$ is upper semicontinuous, there
exists an $\eps=\eps (\dl)$ such that $b(z)<-\dl$ for all
$z\in \{|z-z_0|<2\eps\}.$

Let $(\bul)_\eps$ be a smoothing operator from (2.6.2.3). If $w_n\overset{\Di'}\to{\ri}w$ then $(w_n)_\eps\ri w$ uniformly on every compact set (Th.2.3.4.5,
reg3) and for every subharmonic function $v$ the sequence $v_\eps (z)\downarrow
v(z)$,when $\eps\downarrow 0$ (Th.2.6.2.3, ap2).

Then $(b)_\eps (z)<-\dl,$ for $|z-z_0|<\eps$ or $(w_\iy)_\eps(z)+(v)_\eps(z)\leq(H)_\eps (z)-\dl.$ The function $H$ is continuous,  hence uniformly continuous on the circle $\{z:|z-z_0|\leq \eps\}$. Thus we can replace $(H)_\eps$ in the last inequality with $H$ and  $ \dl$ with $\dl/2.$ So, we have
$$(w_\iy)_\eps(z)+v_\eps(z)\leq H(z)-\dl/2, \ |z-z_0|<\eps.\tag 6.1.5.1$$
Since $(\bul)_\eps$ is monotonic on subharmonic functions, we can replace $\eps$
in (6.1.5.1) with any $\eps_1<\eps.$ So we obtain
$$ (w_\iy)_{\eps_1}(z)+v_{\eps_1}(z)\leq H(z)-\dl/2, \ |z-z_0|<\eps\tag 6.1.5.2$$
Since $(w_n)_{\eps_1}\ri (w_\iy)_{\eps_1}$ uniformly in the disc $|z-z_0|\leq \eps$ we can replace in (6.1.5.2) $w_\iy$ with $w_n$ and respectively $v$ with
$v_n,$ changing $\dl/2$ with $\dl/4.$ After that we can pass to the limit as
$\eps_1\downarrow 0$ for every sufficiently large $n.$
So we obtain
$$w_n(z)+v_n(z)\leq H(z)-\dl/4,\ |z-z_0|<\eps$$
It means that the disc $\{|z-z_0|<\eps\}\sbt \BC\setminus E_{v_n} .$ Because of
Lemma 6.1.4.2, 3. $w_n $ is harmonic in this disc for all large $n.$ Thus
$w_\iy$ is also harmonic, as the $\Di'$-limit of $w_n.$
\qed
\edm

\demo {Proof of Theorem 6.1.2.6} Let $v_n\overset{\Di'}\to{\ri}v$. Then the set
$w_n=\CG_Hv_n$ is equicontinuous by Lemma 6.1.4.2,  5., and we can choose
from it a subsequence uniformly converging to a continuous function
$w_\iy.$ Let
$w=\CG_Hv,E=E_w, E_\iy$ being defined in Lemma 6.1.5.1.

Since
$$(w_\iy+v)(z)\leq H(z), \ (w+v)(z)\leq H(z), \ \forall z\in \BC$$
and $v$ is upper semicontinuous, whereas $w$ and $H$ are continuous, the sets
 $E$ and $E_\iy$ are closed.

Since $w_\iy(z)\leq H(z)-v(z)$ we have
$$w_\iy(z)\leq w(z),\ \forall z\in\BC,\tag 6.1.5.3$$
and therefore
$E_\iy\sbt E.$

The function $w$ is subharmonic in $\BC\setminus E_\iy$, and $w_\iy$ is
harmonic in $\BC\setminus E_\iy$ by Lemma 6.1.5.1. They take the same values on $E_\iy$. As the maximum principle holds in $\CH_{w_\iy}$ by assumption we have, according to Lemma 6.1.4.1 the inequality
$$w(z)\leq w_\iy(z),\  \forall z\in\BC.\tag 6.1.5.4$$
The inequalities (6.1.5.4) and (6.1.5.3) imply that $w(z)=w_\iy(z),$ i,e.,
$\CG_H$ is continuous.
\qed
\edm
 \subheading {6.1.6} \demo {Proof of Theorem 6.1.2.5} Sufficiency. We exploit
the following criterion for existence of a limit set that follows from Theorems
4.2.1.1, 4.2.1.2, 4.3.1.2 and Corollary 5.3.1.5:
\proclaim {Proposition 6.1.6.1}
In order that $U\sbt U[\r]$ be a limit set
of an entire function $f\in A(\r(r))$ it is necessary and sufficient that there
exists a piecewise continuous, $\om$-dense in $U$ asymptotically dynamic pseudo-trajectory (a.d.p.t) $v(\bul|t).$\ep

{\bf Exercise 6.1.6.1.} Check this.

Let $v_\Phi(\bul|t)$ be an a.d.p.t. corresponding to $\Fr \Phi.$ Consider the
pseudo-trajectory $v_g(\bul|t):=\CG_Hv_\Phi(\bul|t).$ It exists because of (6.1.1.1). Prove that this pseudo-trajectory is asymptotically dynamical, i.e.,
(4.3.1.1) is fulfilled. Recall that $T_\tau\bul =(\bul)_{[e^\tau]}.$

Using the property of invariance of $\CG_H$ (Th 6.1.1.2, 2.) we
have
$$T_\tau v_g(\bul|e^t)- v_g(\bul|e^{t+\tau})=\CG_H [T_\tau v_\Phi(\bul|e^t)- v_\Phi(\bul|e^{t+\tau})].$$
Thus (4.3.1.1) is fulfilled because of continuity of $\CG_H$. Also the condition
of $\om$-denseness (4.3.1.4) is fulfilled and
 $$\{w\in U[\r]:(\exists t_j\ri\iy)\  w=\Di'-\lim v_g(\bul|e^{t_j})\}
=\CG_H(\Fr\Phi)$$
The corresponding entire function $g\in A(\r(r))$ with the limit set $U_g=\CG_H(\Fr\Phi)$ is an ideally complementing multiplicator, because
$$\Fr[g\Phi]=\{v+\CG_Hv:v\in\Fr[\Phi]\}.$$

{\bf Exercise 6.1.6.2.} Check this.

Necessity. Let $G$ be an ideally complementing multiplicator  of
$\Phi.$ Let us show that $\CG_H$ is continuous on $\Fr[\Phi].$
Assume this is not the case, i.e., there exists a sequence $v_j\ri
v$ such that $\CG_Hv_j\ri  W$ and $W\neq \CG_Hv.$ Since the limit
set $\Fr[G\Phi]$ is closed, we have $v_j+\CG_Hv_j\ri v+\CG_Hv, \
v_j\in\Fr[\Phi]$. On the other hand , $v_j+\CG_Hv_j\ri v+W.$ Thus,
$W=\CG_Hv$ that is a contradiction. \qed \edm

 \demo{Proof of Theorem 6.1.2.3} Let $\Fr [\Phi]$ be a periodic limit set, that is
$$\Fr [\Phi]=\BC(v)=\{v_{[t]}:1\leq t\leq e^P\},$$
where $v\in U[\r].$ We shall show that $\CG_H$ is continuous on $U[\r].$
 By Theorem 6.1.1.2, 2) the equality $(\CG_Hv)_{[t]}=\CG_Hv_{[t]}$ holds. Since
the operation $(\bul)_[t]$ is continuous for all $t,$ $\CG_H$ is continuous on
$\BC(v).$\qed\edm
\subheading {6.1.7} Now we are going to prove Theorem 6.1.2.4. However, we need some preparations.

Let $h(\phi),\ \phi\in [0,2\pi)$ be a $2\pi$-periodic $\r$-t.c.function, satisfying the condition
$$\max\limits_{\phi \in [0,2\pi]}h(\phi)=\s$$
We denote this class as $TC[\r,\s]$ and denote
$$TC[\r]:=\bigcup\limits_{\s>0}TC[\r,\s].$$
The class of function $w=h_1-h_2$ where $h_1,h_2\in TC[\r,\s]$ will be denoted as
$\dl TC[\r,\s]$ and denote also
$$\dl TC[\r]:=\bigcup\limits_{\s>0}\dl TC[\r,\s].$$
From properties of $\r$-t.c.function (see \S\S 3.2.3-3.2.5) we can obtain the following properties of $\dl-\r$ t.c.functions:
\proclaim {Proposition 6.1.7.1} For $w\in \dl TC[\r]$ the following holds:

1.$w'(\phi-0)$ and $w'(\phi+0)$  exist at each point and are bounded in $[0;2\pi];$

2.$w'(\phi-0)= w'(\phi+0)$ for all $\phi\in [0;2\pi]$ , except, perhaps, a countable set;

3.the charge $\Delta_w$ generated by the function
$$\Dl_w:=w'(\phi)+\r^2\int^\phi w(\th)d\th$$
has  bounded variation $|\Dl_w|;$ the variation $|\Dl_w|(\a,\be)$ of the charge on the interval $(\a;\be)$ and the variation of the charge generated by derivative $|w'|(\a;\be)$on the same interval satisfy the relation:
$$|\Dl_w|(\a,\be)\geq |w'|(\a;\be)+\r^2(\be-\a).$$

4.For all $w\in \dl TC[\r,\s]$
$$\max(|w'(\phi-0)|,|w'(\phi+0)|)\leq C(\r,\s),\ \phi\in [0;2\pi];$$

5.if $r^\r w_n{\Di'}\overset\to{\ri} r^\r w$ and $w_n \in \dl TC[\r,\s],$ then
$w_n\ri w$ uniformly on $[0;2\pi].$
\ep

{\bf Exercise 6.1.7.1.} Prove this using properties of $\r$ -t.c.functions.

We also need a technical
\proclaim {Lemma 6.1.7.2} Let $M_n(\phi)$ be a sequence of functions which satisfies the conditions:

1. $M_n\geq 0;\ M_n(0)=0;$

2. $M_n$ converges uniformly to $M_\iy(\phi)\geq A\sin \r\phi,\ A>0;$

3. $M'_n(\phi -0),\ M'_n(\phi +0)$ exist at every point, and they coincide
almost everywhere;

4.there exists a sequence $\phi_n\downarrow 0$ such that for each arbitrarily small $\eps >0$ and arbitrarily large $n_0\in \BN$ there exists $n>n_0$ for which
the inequality $M_n(\phi_n)<\eps\phi_n$ holds.

Then there exists a sequence $(\z_n,\eta_n)$ of disjoint intervals and a subsequence $M_{k_n}$ such that
$$M'_{k_n}(\z_n)-M'_{k_n}(\eta_n)\geq A\r/2.\tag 6.1.7.1$$
\ep
\demo {Proof} Set $\eps_0=1/2,\ \eta_0=\pi/4$ and choose the required sequence
recurrently. Let $\eps_n,\eta_n,\z_n$ be already chosen. Set $\eps_{n+1}=
\eps_n/2,$  find $\phi_{n+1}<\eta_n$ and choose $k_0=k_0(n)$ so that for $k>k_0$
$$M_k(\phi_{n+1}) -A\r\phi_{n+1}>-\eps_{n+1}\phi_{n+1}.$$
This is possible because of 2. and $\sin\r\phi\sim \phi, \ \phi\ri 0.$
So we have
$$\frac {M_k(\phi_{n+1})}{\phi_{n+1}}>A\r -\eps_{n+1}.\tag 6.1.7.2$$

Now, choose $\psi_{n+1}<\phi_{n+1}$ and $k_{n+1}>k_0$ so that
$$M_{k_{n+1}}(\psi_{n+1})<\eps_{n+1}\psi_{n+1}\tag 6.1.7.3$$
This  is possible by 4.
Thus for small $\eps_{n+1}$ from (6.1.7.2) and (6.1.7.3) we obtain
$$\frac{M_{k_{n+1}}(\phi_{n+1})-M_{k_{n+1}}(\psi_{n+1})}{\phi_{n+1}-\psi_{n+1}}
>(2/3)A\r \tag 6.1.7.4$$
On the  other hand
$$\frac{M_{k_{n+1}}(\psi_{n+1})-M_{k_{n+1}}(0)}{\psi_{n+1}-0}<\eps_{n+1}\tag 6.1.7.5$$
On the interval $(\psi_{n+1},\phi_{n+1})$ there is a point $\eta_{n+1}$ where the derivative exists and the inequality
  $$M'_{k_{n+1}}(\eta_{n+1})\geq \frac{M_{k_{n+1}}(\phi_{n+1})-M_{k_{n+1}}(\psi_{n+1})}{\phi_{n+1}-\psi_{n+1}}\tag 6.1.7.6$$
is valid.
Also there is a point $\z_{n+1}\in (0,\psi_{n+1})$ where derivative exists and
the inequality
$$M'_{k_{n+1}}(\z_{n+1})\leq\frac{M_{k_{n+1}}(\psi_{n+1})-M_{k_{n+1}}(0)}
{\psi_{n+1}-0}\tag 6.1.7.7$$
is valid.

From the inequalities (6.1.7.4)-(6.1.7.7) we obtain (6.1.7.1)
\qed
\edm
\demo {Proof of Theorem 6.1.2.4 } Denote by $\hat {\CG}_Hh$ the maximal $\r$-t.c.minorante of $H(e^{i\phi})-h(\phi).$ It follows from Th.6.1.1.2, 2. that
$$\CG_H(r^\r h(\phi))(re^{i\phi})=r^\r\hat {\CG}_H h(\phi)$$

{\bf Exercise 6.1.7.2.} Prove this.

So taking in consideration Proposition 6.1.7.1, 5., one must prove
\proclaim {Proposition 6.1.7.3}The operator $\hat {\CG}_H$ is continuous on the set
$$\hat {U}_{ind}:=\{h(\phi,c):c\in C\}$$
in the uniform topology.
\ep
\demo {Proof} Let $h_n\ri h,\ h_n,h\in \hat {U}_{ind}.$ Denote $\hat w_n=\hat\CG_Hh_n,\
\hat w=\hat\CG_Hh,\ \hat w_\iy=\lim_{n\ri\iy}\hat w_n.$ We set also
$\hat M_n=H-h_n-\hat w_n,\ \hat M\iy=H-h-\hat w_\iy,\ \hat M=H-h-\hat w.$ Let $(\a_n;\be_n)$
be a maximum interval where $\hat M_n(\phi)>0.$ We shall show that
$\be_n-\a_n\leq \pi/\r.$ Indeed, for a fixed $n$ let us consider the function
$$\hat W_n:=\hat w_n+\eps_n L(\phi-(\a_n+\be_n)/2)$$
where
$$L(\phi)=\cases \cos|\phi|,\ &\phi\in (-\pi/2\rho;\pi/2\rho);\\0, &\phi\in
[-\pi;\pi]\setminus (-\pi/2\rho;\pi/2\rho),\endcases$$
 and $\eps_n$ is small enough. If $\be_n-\a_n>\pi/\r$, then $\hat W_n$ is also a
 $\r$-t.c. minorant of $H-h_n,$ i.e., $\hat w_n$ is not maximal.

If $\be_n-\a_n=\pi/\r,$ then, to ensure that $\hat w_n$ is a
maximal minorant, at least one of the conditions
$$\liminf\limits_{\phi\ri\a_n+0}\frac {\hat M_n (\phi)}{\phi-\a_n}=0,\
\liminf\limits_{\phi\ri\be_n-0}\frac {\hat M_n (\phi)}{\be_n -\phi}=0.\tag 6.1.7.8$$
must be satisfied.

Let us choose (and preserve the previous notation) a subsequence
$ \hat M_n (\phi)$ for which $\a_n\ri\a,$ and $\be_n\ri\be$.

If $\be -\a<\pi/\r,$ then the maximum principle for $\r$-t.c.functions is valid. Repeating arguments of proof of Theorem 6.1.2.6, we obtain $\hat w_\iy=\hat w$ for all $\phi$, which proves Proposition 6.1.7.3 for the case considered.

{\bf Exercise 6.1.7.3.} Repeat them.

Consider the case when $\be -\a=\pi/\r.$

Set $q(\phi)=(w-w_\iy)(\phi).$ The function $q$ is $\r$ -trigonometric on the
interval $(\a;\be)$ since $\hat w$ and $\hat w_\iy$ are $\r$-trigonometric, i.e.,
have  the form $A\sin\r\phi+B\cos\r\phi.$

{\bf Exercise 6.1.7.4.} Explain this.

Besides, we have $q\geq 0$ and $q(\a)=q(\be)=0.$ It is easy to see that $q$ has the form
$$q(\phi)=A\sin\r(\phi-\a),\ A>0. \tag 6.1.7.9$$
{\bf Exercise 6.1.7.5.} Prove this.

Since $\hat M(\phi)\geq 0,$ we have $\hat M_\iy(\phi)=\hat
M(\phi)+(\hat w-\hat w_\iy)(\phi)\geq (\hat w-\hat w_\iy)(\phi), \ \forall \phi,$ whence
$$\hat M_\iy(\phi)\geq A\sin\r(\phi-\a),\ A>0. \tag 6.1.7.10$$
Since $\be_n-\a_n\leq \pi/\r,$ the segment $[\a,\be]$ contains the infinite
sequence $\a_n$ or $\be_n.$ Let us single out a subsequence, let it be, for
example, $\a_n\ri\a+0, \ \a_n\in[\a;\be]$.

Consider the sequence $M_n(\phi)=\hat M_n(\phi -a).$ From the
definition of $M_n$ and from relation 6.1.7.10 it follows that
conditions 1. and 2. of Lemma 6.1.7.1 are fulfilled. Condition 3.
is fulfilled because of property 1. of Proposition 6.1.7.1.
Further, if $\a_n\not\equiv\a$, then condition 4. of Lemma 6.1.7.1
is trivially true,since $M_n(\a_n-\a)=0$; otherwise, if
$\a_n\equiv \a$, condition 4. follows from (6.1.7.8).

Applying Lemma 6.1.7.1, we obtain the union of intervals
satisfying (6.1.7.7). The equality $H(\phi)=\hat M_n-h_n-\hat w_n$
yields the following inequality for the measure
$\Dl_H$:$\Dl_H((\eta_n;\z_n))\geq A\r/2.$ Summing this inequality
and taking into account the fact that the intervals do not
intersect, we obtain $\Dl_H(\cup_n (\eta_n,\z_n))=\iy,$ which is
impossible. So, Proposition 6.1.7.3 is proved. \qed \edm Hence,
Theorem 6.1.2.4 is proved. \qed \edm \subheading {6.1.8} In this
item we show an example of $H$ and   an entire function without an
ideally complementing $H$- multiplicator.

According to Th.6.1.2.5 , to construct such an example it is sufficient to construct a limit set on which $\CG_H$ is not continuous.

We set
$$L(\eta)=\cases \cos|\eta|,\ &\eta\in (-\pi/2\rho;\pi/2\rho);\\0, &\eta\in
[-\pi;\pi]\setminus (-\pi/2\rho;\pi/2\rho).\endcases$$
 Let us define $X\in C^\iy$    so that $X(\xi)=1$ for $\xi<0$ and $X=0$ for $\xi>\a.$

We set
$$\k:=(1/\r^2)\max\limits_{(-\iy;+\iy)}[2\r X'+X''](\xi),\ \
H_0(\eta):=L(\eta)+\k.\tag 6.1.8.1$$
We also set
$$v(\z,c):=[H_0-X(\xi-c)L(\eta)]e^{\r\xi},\ \ \z=\xi+i\eta,$$
where $H_0$ and $L$ have been periodically extended from the interval $[-\pi;\pi]$ to $(-\iy,+\iy).$

As $H(z)$ we take
$$H(z):=H_0(\phi)r^\r$$
\proclaim {Lemma 6.1.8.1} We have
$$v(\log z,c)\in U[\r,\s],\ \   \s=1+\k,\tag 6.1.8.2$$
$$\CG_{H}v(\bul,c)\equiv 0,\tag 6.1.8.3$$
$$\liml_{c\ri\iy}v(\log z,c)=\k r^\r\tag 6.1.8.4$$
uniformly with respect to $z\in K\Subset \BC ,$ and
   $$\CG_H(\k r^\r)=L(\phi)r^\r. \tag 6.1.8.5$$
\ep
\demo  {Proof} For the Laplace operators in $\z$ and $z$ it is true that
$\Dl_\z=\Dl_z/|\z|^2$.
Let us check that $v(\z,c)$ is subharmonic in $\z.$ We have
$$\Dl_\z v(\z,c)=\{[1-X(\xi-c)](L''+\r^2 L)(\eta)+[
\r^2\k -L(\xi)[X''(x-c)+2\r X'(x-c)]\}e^{\r\xi}$$

{\bf Exercise 6.1.8.1.} Check this computation.

Since $X(\xi)\leq 1$ and $L(\eta)$ is $\r$-t.c.
$$[1-X(\xi-c)](L''+\r^2 L)(\eta)\geq 0$$
Since $L(\xi)\leq 1 $ and $[X''(x-c)+2\r X'(x-c)]\leq \k\r^2$ we have
$$[\r^2\k -L(\xi)[X''(x-c)+2\r X'(x-c)]\geq 0.$$
Thus $v(\log z,c)$ is subharmonic.

{\bf Exercise 6.1.8.2.} Prove that $v(\log z, c)\in U[\r,\s]$ for $\s=1+\k.$

Let us prove (6.1.8.3). We have
$$ H(z)-v(\log z,c)=X(\log r -c)L(\phi)r^\r.$$
Since $X=0$ for $r>e^{c+\a}$, the maximal subharmonic minorant of $H-v$ is
zero by the maximum principle.

Relation (6.1.8.4) is obvious, since $X(\log r -c)$ converges to 1 uniformly
on every disc $\{|z|\leq R\}.$ Relation 6.1.8.5 follows from the equality
$H(z)-\k r^\r=L(\phi)r^\r$, since $L(\phi)r^\r\in U[\r].$
\qed
\edm

Now we pass to the construction of the example.
Examine the set
$$U_1:=clos \{v(\log z,c):c\in [0;\iy)\}.$$
It contains the function
$$\Di'-\lim_{c\ri\iy}v(\log z,c)=\k r^\r.$$
Let us consider the minimal convex $(\bul)_{[t]}$ -invariant set $U$
containing $U_1$.The set is contained in $U[\r,1+\k].$ It is a limit set for
a certain entire function $\Phi.$ Let us show that $\CG_H$ is not continuous on
$\Fr[\Phi]$. We take an arbitrary sequence $c_j\ri\iy$ and set
$v_j(z):=v(\log z, c_j)\in U$. Now $\Di'-\lim_{j\ri\iy}v_j=\k r^\r$ by (6.1.8.4) and $\CG_Hv_j(z)=0,$ so $\Di'-\lim_{j\ri\iy} \CG_Hv_j=0$ but
$$\CG_H(\lim v_j)=\CG_H(\k r^\r)=L(\phi)r^\r\not\equiv 0$$
which shows the lack of continuity.

By virtue of Th.6.1.2.5 $\Phi$ it is not ideally complementable.
\subheading {6.1.9}Here we prove  existence and continuity of
maximal subharmonic minorant for some classes of functions $m(z)$
. \proclaim {Theorem 6.1.9.1} Let $m(z)$ be a continuous function
such that the set of subharmonic minorants is nonempty. Then the
maximal subharmonic minorant of $m$ exists and is continuous. \ep
\demo {Proof} The set of subharmonic minorants is not empty and
partially ordered. Indeed, for every subset $\{u_\a,\a\in A\}$ of
subharmonic minorants there exists $u_A=(\sup \{u_\a:\a\in A\})^*$
which is subharmonic and is a minorant of $m$, because $m$ is
continuous.

{\bf Exercise 6.1.9.1.} Explain this in details.

Thus there exists a uniquely maximal element m.s.m.$(z,m)$, which is a subharmonic minorant of $m.$

Let us prove that it is continuous at every point $z_0.$ Since m.s.m.$(z,m)$
is upper semicontinuous,
$$m.s.m (z,m)>m.s.m (z_0,m)-\eps $$
for $|z-z_0|<\dl$ for arbitrary small $\eps$ and corresponding $\dl=\dl(\eps).$
So we need to prove the inequality
$$m.s.m (z,m)<m.s.m (z_0,m)+\eps $$
for arbitrary small $\eps$ and corresponding $\dl=\dl(\eps).$

Perform sweeping m.s.m $(z,m)$ from the disc $|z-z_0|<\dl$ such that the result
$u(z,\dl)$ satisfies the inequality
$$m.s.m(z,m)<u(z,\dl)<m.s.m(z,m)+\eps<m(z)+\eps.$$
Thus $u(z,\dl)-\eps<m(z).$ Hence $m.s.m(z,m)>u(z,\dl)-\eps $ for all $z.$ Since $u(z,\dl)$
is continuous $u(z,\dl)>u(z_0,\dl)-\eps$ in the disc $\{|z-z_0|<\dl_1\}.$ So
m.s.m.$(z,m)>u(z_0,\dl)-\eps>$m.s.m $(z_0,m)-\eps$\qed
\edm
\proclaim {Theorem 6.1.9.2} Let $m=m_1-m_2,$ where $m_1,m_2$ are subharmonic
functions. Then the maximal subharmonic minorant of $m$ exists. If $m_1$
is continuous, then the maximal subharmonic minorant is continuous.
\ep
\demo {Proof} Set $\CM_\eps (z,m):=\CM_\eps (z,m_1)-\CM_\eps (z,m_2),$ where
$\CM_\eps (z,m_i), i=1,2$ is defined by  2.6.1.1.
Since $\CM_\eps (z,m)$ is continuous (see Theorem 2.6.2.3 (Smooth approximation)), there exists m.s.m.$(z,\CM_\eps (z,m)).$
We have
$$u(z,m):=\limsup\limits_{\eps\ri 0}m.s.m (z, \CM_\eps (\bul, m)\leq \lim\limits_{\eps\ri 0}\CM_\eps (z,m)=m_1-m_2(z)=m(z)$$
 Now we prove that the upper semicontinuous regularization $u^*(z,m)$ also satisfies the inequality $u^*(z,m)\leq m(z).$
Indeed, $m_2 +u(z,m)\leq m_1(z).$ Hence,
$$ \CM_\eps (z,m_2)+\CM_\eps (z,u(\bul,m))\leq\CM_\eps (z,m_1).$$
Passing to limit we obtain three subharmonic function and inequality
$$m_2(z)+u^*(z,m)\leq m_1(z).$$
We prove that $u^*(z,m)$ is the m.s.m.$(z,m).$ If not, there would exist
a subharmonic function $u_1$ which exceeds  $u^*(z,m)$ on a set of positive measure (otherwise they coincide);thus we would have for some $z$ and $\eps$
$$ u^*(z,m)<\CM_\eps (z,u_1)\leq  m.s.m.(z,\CM_\eps (z,m))$$
This contradicts the definition of $ u^*(z,m).$

Now suppose that $m_1$ is continuous at a point $z_0$. From Th.2.6.5.1
(Jensen- Privalov) we
obtain that it is equivalent to
$$\int _0^\eps \frac {\mu_{m_1}(\{z:|z-z_0|<t\})}{t}dt=o(1),\ \eps\ri 0$$

Similarly to the  proof of Th.6.1.4.2, 5, we obtain
$$\mu_{m.s.m(z,m)}\leq \mu_{m_1}.$$
Hence $ m.s.m.(z,m)$ is also continuous.
\qed
\edm

{\bf Exercise 6.1.9.2.} Prove continuity in details.

\proclaim {Corollary 6.1.9.3} For $m=m(z,v,H)$ the function
$\CG_Hv(z):=m.s.m. (z,m)$ exists and is continuous.
\ep

{\bf Exercise 6.1.9.3.} Prove Corollary 6.1.9.3.
\newpage

\centerline {\bf 6.2 A generalization of $\r$-trigonometric convexity. }
\subheading {6.2.1}One of the important and useful kinds of limit sets is periodic limit sets.
They are determined by one subharmonic function $v\in U[\r]$ that satisfies the condition
$$v(Tz)=T^\r v(z),\ z\in\BC\tag 6.2.1.1$$
Such a function is called {\it automorphic }. They generate the class of so called
$L_\r$-subfunctions, that is a generalization of $\r$ -trigonometrically convex functions. In this part we are going to review of properties of   such functions    from different point of view, that will be useful for applications (see \cite {ADP}).

In connection with property (6.2.1.1) it is natural to consider so called $T$ -{\it homogeneous domains} in $\BC,$ i.e., such domains $G$ that  satisfy the condition
$\{Tz:z\in G\}=G$ or shortly $TG=G.$ As we can see they are invariant with respect dilation by $T.$
For example, every component of an open set of harmonicity of an automorphic function is a $T$-homogeneous domain.

 Let $v$ satisfy (6.2.1.1).Then the function
$$ q(z):=v(e^z)e^{-\r x}\tag 6.2.1.2 $$
is $2\pi$ periodic function in $y$ and $P$-periodic in $x,$ where $P=\log T.$

The function $q$ can be considered as a function on a torus $\BT^2_P,$
 obtained by identifying the opposite sides of the rectangle
$\Pi = (0,T) \times (-\pi, \pi).$

The homology
group of $\BT^2_P $ is nontrivial, and generated by  the cycles
$\gamma_x,\gamma_y$, where $\gamma_x=\BT^2_P\cap \{y=0\},\
\gamma_y=\BT^2_P\cap \{x=0\}$.

Let $\pi$ be the covering map of $\Bbb C$ onto $\Bbb T^2_P$, then
$\phi =\pi \circ \log$ is a well-defined covering map of $\BC
\setminus\{0 \}$ onto
$\Bbb T^{2}_{P}$, where the group of deck transformations is given by
the dilations by $T^{m}$ for $m\in \BZ$.
So if $G$ is a given $T$-homogeneous domain, then
$$D = \pi \circ \log G=\phi (G)\tag 6.2.1.3$$
is a domain in $\Bbb T^2_P$.
On the other hand, not every domain in $\BT^2_P$ has a
$T$-homogeneous domain as its preimage under $\phi$.
The preimage
$\phi^{-1} (D)$ under $\phi$ is a possibly disconnected set which is invariant
under dilations by $T^{m}$ for $m\in \BZ$.
An intrinsic description is given by the
next proposition.

 \proclaim {Proposition 6.2.1.1} Let $\g$ be a closed curve in a
domain $D\sbt\BT^2_P$ that is homologous in $\BT^2_P$ to a cycle $\gamma =
n_x\gamma_x+n_y\gamma_y,\ n_x,n_y\in \BZ.$ Then

1. If $n_x=0$ for every  such $\g$ in $D,$  , then
$$\phi^{-1}(D)=\cup_{j=-\iy}^{\iy}G_j,$$ where $G_j=T^jG_0,$ $G_0$ is
an arbitrary connected component of $\phi^{-1}(D),$ and $G_j\cap
G_l=\varnothing$ for $j\neq l.$

2. If there exists a curve $\g$ as
above with $n_x\neq 0,$ then
$$\phi^{-1}(D)=\cup_{q=0}^{k-1}G_q,$$ where $k =\min |n_x|$ with the
minimum taken over all such curves $\g$; $G_0$ is an arbitrary
component of
$\phi^{-1}(D);$ $G_j, \ j=0,1,...,k-1$,
are disjoint $T^k$-homogeneous domains, and for every $m\in\BZ,\ T^{
m} G_0=G_q$,
provided $ \ m=lk+q,$ for some
$q\in \BZ,\ 0\leq q\leq k-1,\ l\in \BZ.$ \ep

We call domains as in part 2 of Proposition
6.2.1.1 {\it connected on spirals}. In particular, this proposition
shows that  for every $D$ connected on spirals, we can  find a connected $T^k$-
homogeneous domain that relates to $D$ by (6.2.1.3).

Let us give some examples. The domain $D^\prime=\BT^2_P\cap \{|x-P/2|<P/4\}$ is not
connected on spirals, whereas $D^{\prime\prime}=\BT^2_P\cap \{|y| <
\pi/4\}$ is.  It follows that $D^\prime\cap D^{\prime\prime}$ is not
connected on spirals whereas $D^\prime\cup D^{\prime\prime}$ is.

The situation can be more complicated.
Set
$$x'(x,y,\a):=x\cos\a+y\sin\a;\ y'(x,y,\a):=-x\sin\a+y\cos\a;\
0\leq \a<\pi/4 ;$$
$$P_1:=(1/2)\ |x'(P,2\pi,-\a)|;\ P_2:=(1/2)\
|y'(P,2\pi,-\a)|.$$
Then
$R'=\{z'=x'+i'y':-P_1 < x'<P_1;\ -P_2 < y'<P_2 \}$
is a fundamental rectangle for
$\BT^2_P$ in the corresponding coordinates.
Set
$f(y'):=(P_2 -y')^{-1}-
(y'+P_2)^{-1}$
and  $D_{0,0}:=\{z':-P_2 < y'<P_2;
\ f(y')<x'<f(y')+d\}$
where $0<d<P_1.$ Then the domains
$D_{l,m}:=D_{0,0}+2P_1 l+2P_2mi',\ l,m\in \BZ$ are
disjoint, and their union $D$ determines a domain
$\hat D\sbt\BT^2_P.$ This $\hat D$ is determined completely by the intersection of $D$ with the rectangle $R=(0,P)\times (-\pi,\pi).$ The domain $\hat D$ is not connected on spirals.

One more example. Consider the family of lines
$L_l:=\{z=x+iy:y=\pi/(kP)x+l\pi/k,\ x\in \BR\}, \ l\in \BZ.$ It
determines a closed curve (spiral) $\g$ on $\BT^2_P$ with $n_1=k.$
The open set $D_k=\{z:|z-\z|<\eps, \z\in L_l,\ l\in\BZ\},\
0<\eps<P/2\sqrt {\pi^2+k^2},$ determines a domain $\hat D_k$ on
$\BT^2_P$ that is connected on spirals, and such that $\phi^{-1}(\hat
D_k)$ consists of $k$ components, every one of them $T^k$-homogeneous.

 Since the function $v$ in (6.2.1.1) is subharmonic, the function
$q$ of (6.2.1.2) is upper semicontinuous and in the $D'$ topology on $\Bbb
T^2_P$ satisfies the inequality $L_\r q\geq 0$, where
$$
L_\rho:=\Delta+2\rho {\partial\over \partial x}+\rho^2.\tag 6.2.1.4
$$
Such functions $q$ are called {\sl subfunctions with respect to }
$L_\r,$ or $L_\r$-{\sl subfunctions}.

 $L_\r q$ is a positive measure on $\BT^2_P.$

The operator $L_\r$ arises naturally by changing variables $z\mapsto \log  z$
in the Laplace operator $\Dl_\z .$

{\bf  Exercise 6.2.1.1} Check this. Set $\z=e^z.$

Let us note that if $q$ depends only on the variable $y,$ it is a $2\pi$-periodic
$\r$- trigonometric convex function because $L_\r$ turns into $T_\r=(\bul)''+\r^2(\bul)$
(c.f. \S3.2.3).

\subheading {6.2.2} Consider the solution of the homogeneous boundary problem
$$\aligned
&L_\r q = 0\,  \quad \text{in }D;\\
&q\bigm|_{\p D} = 0,\endaligned \tag 6.2.2.1 $$
where $D$ is a domain in ${\BT}^2_P$
and $q$ is bounded in a neighborhood of $\p D$ with
boundary value zero quasi-everywhere.
This is a spectral problem for a {\sl pencil
} of differential operators ( \cite {Ma}).

A solution of this problem can be
defined for an arbitrary domain $D\subset{\BT}^2_P$ with a boundary of  positive capacity.

The {\it spectrum} of the problem (6.2.2.1) consists of those (complex) $\r$ for which (6.2.2.1) holds for some function $q\not\equiv 0$.
The minimal positive point of the spectrum $\r(D)$ exists iff there exists the
spectrum. The spectrum exists iff the domain $D$ is connected on spirals.
In this case $\r (D)$ is the order of the minimal harmonic function in everyone of the
domains $G_i$ that corresponds to $D$ by Prop.6.2.1.1.

The quantity $\r(D)$ is {\it strictly monotonic}. It means that if two domains
$D_1,D_2\in {\BT}^2_P$ are such that $D_1\sbt D_2$ and the capacity of
$D_2\setminus D_1$ is positive, then $\r(D_2)<\r(D_1).$ For example, thus is the case of
 $D_2=\{|y|<d, d<2\pi \}$ and $D_1$ is the same strip without the segment
$\{it:0\leq t\leq d\}$.

In connection with the problem of multiplicator we considered the maximal subharmonic minorant of a function $m=H-v$ where $v$ is a $T$- automorphic function.
From Th.6.1.1.2, 2. we can  obtain that if $v$ is a $T$- automorphic function,
  then $\CG_Hv$ is also $T$- automorphic.

{\bf Exercise 6.2.2.1} Check this.

Thus for this case finding $D_H$ in Th. 6.1.1.1 is reduced to finding a maximal
$L_\r$-subfunction $q$ that satisfies the inequality
$$ q(z)\leq m(z):=[H(e^z)-v(e^z)]e^{-\r x},\ z\in {\BT}^2_P.\tag 6.2.2.2$$
We say that $m(z)$ has an $L_\r-subminorante.$

The idea of $\r(D)$ give possibility for the \proclaim {Theorem
6.2.2.1} If  $m$ has a non-zero $L_\r$-subminorant, then
$\r(D)\leq \r$ for some component $D$ of the open set
$\CM_+:=\{z:m(z)>0\}.$

Conversely, if $\r(D)<\r$ (strict inequality) for some component
$D$ of the set $\CM_+$, and $m(z)\geq 0$ for all $z\in {\BT}^2_P,
$ then $m$ has a non-zero $L_\r$-subminorant.\ep

{\bf Exercise 6.2.2.2} Prove that  $\CM_+$ is open.
\subheading {6.2.3} If $\r\notin \BZ,$ the operator $L_\r$ has a fundamental solution $E_\r(\bul-\z)$ in ${\BT}^2_P,$ where $\z$ is a shift by the torus,
i.e., by the modulus $P+i2\pi.$ It means that
$$L_\r E_\r(\bul-\z)=\dl_\z, $$
in $\Di'({\BT}^2_P),$ where $\dl_\z$ is the Dirac function, concentrated at
$\z.$

If $\r\in \BZ,$ there exists, like for operator $T_\r$ and spherical operator
(see Th.3.2.4.2, Th.3.2.6.3), a generalized fundamental solution $E'_\r$ that
satisfies the equation
$$L_\r E'_\r(\bul-\z)=\dl_\z -\cos \r (y-\eta),\ \z=\xi +i\eta$$
in $\Di'({\BT}^2_P).$
\proclaim {Theorem 6.2.3.1}Let $\r>0,\ \r\notin \BZ.$ Then every $L_\r$-subfunction on  ${\BT}^2_P$ can be represented in the form
$$q(z)=\int_{{\BT}^2_P}E_\r(z-\z)\nu(d\z),\tag 6.2.3.1$$
where $\nu=L_\r q.$\ep This theorem is the counterpart of Th.'s
3.2.3.3,\ 3.2.6.2. \proclaim {Theorem 6.2.3.2}Let $\r>0,\ \r\in
\BZ.$ Then the mass distribution $\nu=L_\r v$ satisfies the
condition
$$\int_{\Bbb T^2_P} e^{\pm i \r y} \nu (dz)=0,
\tag 6.2.3.2$$
and the representation
$$q(z)=\Re(Ce^{i\r y})+\int_{\Bbb T^2_P}E_\r' (z-\zeta)\nu(d\z)\tag 6.2.3.3$$
holds with $C$ that is a complex scalar.
\ep
This theorem is the counterpart of Th.'s 3.2.4.2, 3.2.6.2.

Let $D\sbt {\BT}^2_P$ and $\r(D)>\r.$ Then the operator $L_\r$ has in $D$ the
Green function $-G_\r(z,\z,D).$ Thus for every $q$ that is a $L_\r$ -subfunction
in $D$ and bounded from above in $\overline D$ we have the representation
$$q(z)=g(z) -\int_D G_\r(z,\z,D)\nu(d\z),\tag 6.2.3.4$$
in which $\nu=L_\r q$ and $g$ is the minimal majorant on $\p D$ of
the function $q$, satisfying $L_\r g=0$  in $D.$

This is the counterpart of Th.2.6.4.3 (F.Riesz representation) and Th.3.2.5.1.

From (6.2.3.4) one can easily obtain
\proclaim {Theorem 6.2.3.3. (Maximum principle)} If $\r (D)>\r$ and $q(z)$ is an
$L_\r$ -subfunction such that $q(z)\leq 0,\ z\in\p D,$ then $q(z)\leq 0,\ z\in
D.$ \ep

{\bf Exercise 6.2.3.1} Prove this.
\proclaim {Theorem 6.2.3.4} An $L_\r$-subfunction in ${\BT}^2_P$ can not attain
zero maximum if it is not zero identically.\ep

{\bf Exercise 6.2.3.1} Prove this exploiting (6.2.1.2) and properties of subharmonic functions.

\proclaim {Theorem 6.2.3.5} Let $q$ be an $L_\r$ -subfunction in $\BT_P^2.$ If
$q(z)\leq 0$ for $z\in \BT_P^2$ then $q(z)\equiv 0.$\ep

{\bf Exercise 6.2.3.2} Prove this using Th.3.1.4.7 (**Liouville).

\proclaim {Proposition 6.2.3.6} Let $q_D$ be the solution of the problem (6.2.2.1) in a domain $D$ with a smooth boundary, corresponding to $ \r=\r(D).$ Suppose that $q_D(z_0)=1$ for some $z_0\in D.$ Then
$$\frac {\p q_D}{\p n}>0, \forall z\in \p D.$$
\ep {\bf Exercise 6.2.3.3} Prove this, using properties of
positive harmonic functions. \subheading { 6.2.4 } In the part
devoted to completeness of exponential system (\S 6.3) we will
need the notion of minimality of a subharmonic function from
$U[\r].$ A function $v\in U[\r]$ is called {\it minimal} if the
function $v-\eps r^\r$ has no subharmonic minorant for arbitrary
small $\eps>0.$ If $v$ is $T$- automorphic, the corresponding
$L_\r$ -subfunction $q$ is called minimal if the function $q-\eps$
has no $L_\r$-subminorant in ${\BT}^2_P.$ We formulate  one
sufficient condition for minimality and one sufficient condition
for nonminimality. \proclaim {Theorem 6.2.4.1}   Let $\Cal
H_\rho(q)$ be the maximal open set on which $L_\rho q=0$.  If
there exists a connected component $M\subset \Cal H_\rho(q)$ such
that $\rho(M)<\rho,$ then $q$ is a minimal $L_\rho$-subfunction.
\endproclaim For example, $q\equiv 0$ is minimal.
\proclaim{Proposition 6.2.4.2} The function $q$ is nonminimal if
$q(z)\geq c$ or $L_\rho q-c > 0$ for some positive $c$ for all
$z\in \Bbb T^2_P$.  \endproclaim For example, $q\equiv c>0$ is
nonminimal.

\newpage
\centerline {\bf 6.3. Completeness of exponential systems in convex domains.}

\subheading {6.3.1} Let $\Lm :=\{\lm_k\},\ k=1,2,...$ be a set of points in the complex plane $\BC,$ satisfying the condition $\lm_k\neq 0$ and $\lm_j\neq \lm_k ,$ if $k\neq j.$

Consider the canonical product
$$\Phi_\Lm (\lm):=\prod_k(1-\lm/\lm_k)\exp \lm/\lm_k\tag 6.3.1.1$$
We suppose in this \ \S\   that $\Phi_\Lm (\lm)$ is an entire function of order
one and normal type, i.e., a function of {\it exponential type}(see
\cite {Levin, Ch.1,\S 20}.

This fact can be expressed in terms of $\Lm$ by using the Brelot-Lindel\"of Theorem
2.9.4.2.

{\bf Exercise 6.3.1.1} Formulate this theorem for entire functions of order one
and normal type under assumption that $\r (r)\equiv 1.$

We will suppose that the upper density of zeros (see \S 2.8,\S
5.1). $\overline\Dl_\Lm>0.$ \subheading {6.3.2} Let $G\sbt \BC$ be
a convex bounded domain containing zero . This last request does
not restrict any of the further considerations connected to
completeness, because $\exp\Lm:=\{e^{\lm_j z}:\lm_j\in \Lm\}$ can
be replaced by the system $ \{e^{\lm_j(z-z_0)}:\lm_j\in \Lm\}$ and
$e^{\lm_j(z-z_0)}=C_je^{\lm_jz}.$ Let  $A(G)$ be the space of
holomorphic functions in $G$ with the topology of
 uniform convergence on compacts. We will study the completeness of the exponential
systems
$$\exp \Lm:=\{e^{\lm_jz}:\lm_j\in \Lm\}\tag 6.3.2.1$$
in $A(G).$

We  will be interested in the following questions:

1. {\it completeness} of $\exp \Lm$ in $A(G).$

2. {\it maximality} of $G$ for $\exp \Lm,$ which is complete in  $A(G);$

3.{\it extremal overcompleteness} of $\exp \Lm$ in $A(G)$ for a maximal $G.$

Let us give a precise definitions  of maximality and extremal  overcompleteness.
The completeness  means that every function
$f\in A(G)$ can be approximated on every compact set $K\Subset G$ with  arbitrary precision by linear combinations
of functions from $\exp \Lm.$

A convex domain $G$ is called {\it maximal} for a system  $\exp \Lm,$ which is complete in $A(G)$ if for every domain $G_1$ such that $G\Subset G_1$ $\exp \Lm$ is not complete in $A(G_1).$

A system  $\exp \Lm$ is called {\it extremely overcomplete } in $A(G)$ for
a maximal $G,$ if for every sequence $\Lm_1:=\{\lm^1_j\}$ such that
$\Lm_1\cap \Lm =\varnothing$ and $\overline\Dl_{\Lm_1}>0$ the domain $G$ is not
maximal for the system $\exp \Lm\cup \Lm_1.$

Another words, every essential enlargement of an extremely overcomplete system
enlarges also the maximal domain of completeness.
\subheading {6.3.3} Let
$$h_\Lm (\phi):=\limsupl_{r\ri\iy}\log|\Phi_\Lm(re^{i\phi})|r^{-1}$$
be the indicator of $\Phi_\Lm$. It is 1-trigonometrically convex function or
simply {\it trigonometrically convex function} (t.c.f). Let $G_\Lm$ be the
{\it conjugate indicator diagram} of $\Phi_\Lm,$ i.e., a convex domain of the form
$$G_\Lm:=\{z:\max\limits_{z\in G_\Lm }\Re (ze^{i\phi})\leq h_\Lm(\phi)\}.$$

Let us describe conditions for completeness, maximality and extremal overcompleteness when $\Lm$ is a {\it regular set} (see \S 5.6) and $\Phi_\Lm$ is a CRG -
function (see \S 5.6).

We say that $G_\Lm$ is {\it enclosed}  in $G$ if it can be enclosed in $G$ by
parallel translation, it is {\it enclosed with  sliding}, if it can be moved after enclosing only in one direction, {\it enclosed rigidly} if it is impossible to move after enclosing, {\it freely enclosed} in every other case of enclosing.
\proclaim {Theorem 6.3.3.1} Let $\Lm$ be a regular set.Then the following
holds:

1.\{$\exp\Lm$ is not complete in $A(G)\}\Longleftrightarrow \{G_\Lm$ is freely enclosed to $G$\};

2.\{$G$ is maximal for $\exp\Lm\}\Longleftrightarrow \{G_\Lm$ is not freely enclosed to $G$\};

3.\{{$\exp\Lm$ is  extremely overcomplete in $A(G)\}\Longleftrightarrow
\{G_\Lm$
is enclosed  rigidly in $G\}.$
\ep

Let us note that $G$ is maximal for {$\exp\Lm$ but not extremely overcomplete if and only if $G_\Lm$ is enclosed with sliding in $G.$

This theorem is a corollary of the more general Theorem 6.3.4.1, but will be proved independently in \S 6.3.10.
\subheading {6.3.4} If $\Lm$ is not  regular, it is natural to exploit the notion of limit set (see \S 3.1) to characterize of $\exp \Lm.$

Suppose the limit set of $\Phi_\Lm$ has the form
$$\Fr [\Phi_\Lm]:=\{v(\lm)=|\lm|(ch_1+(1-c)h_2)(\arg \lm):c\in [0;1]\}$$
where $h_1,h_2$ are t.c.f.

Such limit set is a particular case of  $U_{ind}$ (6.1.2.2). It is called {\it indicator limit set}  and it is indeed a limit
set of an entire function (see Exercise 6.1.2.4)).

The asymptotic behavior of the  set $\Lm$ (i.e., the limit set of the corresponding mass distribution) can be described completely  using Th.3.1.5.2.

{\bf Exercise 6.3.4.2.} Do that.

We will call such $\Lm$ an {\it indicator set}. Denote by $G_1,G_2$ the conjugate diagram of $h_1,h_2.$ Since $G_1,G_2$ are convex, the set
$$\a G_1+\be G_2:=\{\a z_1+\be z_2:z_1\in G_1,z_2\in G_2\},\  \a,\be>0$$
is also convex and is a conjugate diagram of the t.c.f.
$h:=\a h_1+\be h_2.$
\proclaim {Theorem 6.3.4.1} Let a set $\Lm$ be an indicator set.Then the following holds:

1.\{$\exp\Lm$ is not complete in $A(G)\}\Longleftrightarrow \{G_1$ and $G_2$ are freely enclosed to $G$\};

2.\{$G$ is maximal for $\exp\Lm\}\Longleftrightarrow \{G_1$ and $G_2$ are enclosed in $G$ and at least one of them is not freely enclosed to $G$\};

3.\{{$\exp\Lm$ is extremely overcomplete in $A(G)\}\Longleftrightarrow
\{cG_1+(1-c)G_2$
is enclosed  rigidly in $G \ \forall c\in [0;1]\}.$
\ep
 This theorem is proved in \S 6.3.10.

The equality holds:
$$h_\Lm=\max (h_1,h_2)\tag 6.3.4.1$$
Thus the conjugate diagram $G_\Lm$ of the function $h_\Lm$ is the convex hull
of $G_1$ and $G_2$.

Let us note that the indicator $h_\Lm$ does not determine the  completeness of
the system $\exp \Lm$ if $\Lm$ is not regular set, as the following example
shows

{\bf Example 6.3.4.1}. Let
$$G_1:=\{z=x+iy:x=1;-1\leq y\leq 1\},$$
$$G_2:=\{z=x+iy:x=-1;-1\leq y\leq 1\},$$
and
$$G=\{z:|z|<1+\eps\}$$
with a small $\eps.$

{\bf Exercise 6.3.4.3} Prove that $G_1$ and $G_2$ are freely enclosed
to $G$ and their convex hull is not enclosed.

Let $\Lm$ be a set such that the interior of $G_\Lm$ coincides with $ G.$ If $\Lm$ is regular set  then
$\exp \Lm$ is complete in $A(G)$, $G$ is maximal for $\exp \Lm$ and
$\exp \Lm$ is extremely overcomplete in $A(G).$

If $\Lm$ is an indicator set, then the first two assertions hold but
$\exp \Lm$ can
be not extremely overcomplete:

{\bf Example 6.3.4.2} Set
$$G_1:=\{z=x+iy:-1\leq x\leq0;y=0\};\ G_2:=\{z=x+iy:x=1;-1\leq y\leq 1\}.$$
Here $G_\Lm$ is triangle in which $G_1$ is freely enclosed and $G_2$ is  rigidly
enclosed, but $cG_1+(1-c)G_2$ is free enclosed for all $c:0<c<1.$

{\bf Exercise 6.3.4.4} Check this.

{\bf Example 6.3.4.3} Set
$$G_1:=\{z=x+iy:x=-1;y\in [-1;1]\};$$
$$G_1:=\{z=x+iy:x=1;y\in [-1;1]\}.$$

{\bf Exercise 6.3.4.5} Check that $G_1$ and $G_2$ are inclosed with
sliding in $G_\Lm.$

If $G_1$ and $G_2$ are  rigidly enclosed  in $G$ it does not imply in general that $cG_1+(1-c)G_2$ are  rigidly enclosed for all $c\in [0;1].$

{\bf Example 6.3.4.4} Let $G_1$ be  an equilateral  triangle
inscribed in the circle $|z|=1$ , let $G_2$ be the same triangle rotated by
the angle $\pi/6,$ and let $G$ be the unit disc.

{\bf Exercise 6.3.4.6} Show that $\frac {1}{2}(G_1+G_2)$ is freely enclosed in
$G.$

If $G_1\cap G_2$ is  rigidly enclosed in $G$ then  $cG_1+(1-c)G_2$ is  rigidly
enclosed for $c\in [0;1].$

 {\bf Exercise 6.3.4.7} Check this.

However this is not a necessary condition.

{\bf Example 6.3.4.5} Set
$$G:=\{z=x+iy:|x|<1;|y|<1\};$$
$$G_1:=\{z=x+iy:x\in (-1,1);-x>y>-1\};$$
$$G_2:=\{z=x+iy:x\in(-1,1); -1<y<x\}.$$

 {\bf Exercise 6.3.4.8} Check that every triangle $cG_1+(1-c)G_2$ is
 rigidly enclosed in $G$ and $G_1\cap G_2$ is freely enclosed.
\subheading {6.3.5} Consider in more details the conditions for extremal
overcompleteness in the case when $\Lm$ is indicator set and  $G_\Lm =G$ or, in other words, if
$$h_\Lm=h_G.\tag 6.3.5.1$$
We can suppose that $h_1$ and $h_2$ are linearly independent,
otherwise we exploit Theorem 6.3.3.1. If, for example,  the
inequality $h_1(\phi)\leq h_2(\phi),\forall \phi,$ holds, the
extremal  overcompleteness is in the case when $G_1$ is  rigidly
enclosed in $G_2$ because $G_1\cap G_2=G_1$  , and this case  was
mentioned above (Exercise 6.3.4.7).

Consider the general case. Denote $g(\phi):=|h_1-h_2|(\phi)$, and
set $\Theta_\Lm:=\{\phi:g(\phi)>0\}$. This is an open set on the
unit circle. Denote as $I_\Lm:=(\a_1,\a_2)$ the maximal interval
contained in $\Theta_\Lm$ and denote by $d_\Lm$ its length. Since
$g(\phi)$ is continuous
$$g(\a_j)=0,\ j=1,2.\tag 6.3.5.1$$
If also   at least one of the conditions:
$$\liminf\limits_{\phi\in I_\Lm, \phi\ri\a_j}\frac {g(\phi)}{\phi-\a_j}=0,
\ j=1,2,$$
is fulfilled we say $g$ is {\it  zero with tangency} on $\partial I_\Lm.$
\proclaim {Theorem 6.3.5.1} Suppose $\Lm$ is an indicator set that satisfies
(6.3.5.1). In order that $exp \Lm$ be extremely overcomplete in $A(G)$  it is necessary and sufficient that at least one of the following condition holds:

1.$d_\Lm <\pi;$

2.$d_\Lm =\pi$ and $g$ is zero  with tangency on $\partial I_\Lm.$
\ep
This theorem is proved in \S 6.3.11.
\subheading {6.3.6} We call $\Lm$ {\it periodic} if $\Fr [\Phi_\Lm]$ is
a periodic limit set (see Th.4.1.7.1). In such case all the limit set is determined by one subharmonic function $v\in U[1]$ (see (4.1.3.1)). Let us characterize the system $\exp \Lm$ for periodic $\Lm.$

Set
$$h_G(\phi):=\max\{\Re(ze^{i\phi}):z\in G\}\tag 6.3.6.1$$
 $$m(\lm,G,v):=|\lm|h_G(\arg \lm)-v(\lm)\tag 6.3.6.2$$

Denote as $\CG_Gv$ -the maximal subharmonic minorant of the function $m(\lm,G,v)$
A function $w\in U[1]$ is called {\it minimal} if the function $w-\eps|\lm|$ has  no subharmonic minorant in $U[1]$ for every small $\eps>0$.
The harmonic function of the form
$$H(\lm):=|\lm|(A\cos(\arg\lm)+B\sin (\arg\lm)),\tag 6.3.6.3$$
for example, is minimal.

We will denote as HARM the set of the functions of the form (6.3.6.3) .
\proclaim {Theorem 6.3.6.1} Let $\Lm$ be a periodic set. The following holds:

 1.\{$\exp\Lm$ is not complete in $A(G)\}\Longleftrightarrow \{\CG_Gv$ exists and is non minimal\};

2.\{$G$ is maximal for $\exp\Lm\}\Longleftrightarrow \{\CG_Gv$ exists and is minimal\};

3.\{{$\exp\Lm$ is extremely overcomplete in $A(G)\}\Longleftrightarrow
\{\CG_Gv \in HARM\}.$
  \ep

This theorem is proved in \S 6.3.12.
\subheading {6.3.7} Let us characterize the completeness of $\exp_\Lm$ for
periodic $\Lm$ in other terms. For this we need the information that was presented in \S 6.2. We will take $\r=1.$
Denote
$$q_\Lm(z):=v_\Lm(e^z)e^{-x}\tag 6.3.7.1$$
(compare with (6.2.1.2)).
As it was explained in \S 6.2 this function is an $L_1$-subfunction on the torus
$\BT_P^2.$ Set
$$m(z,G,q_\Lm)=h_G(y)-q_\Lm,\ D(G,\Lm):=\{z:m(z,G,q_\Lm)>0\}\sbt \BT_P^2$$
The set $D(G,\Lm)$ is open because $-m$ is an upper semicontinuous function (see
Th.2.1.2.4), denote
$$ \r (\Lm, G):=\min \r(M)$$
where the minimum is taken over all components $M$ of $D(G,\Lm),$ and it is attained on one of the components because they are not intersecting and $\BT_P^2$ is compact.

{\bf Exercise 6.3.7.1} Explain this in details, using properties of $\r(D)$
(\S 6.2).
\proclaim {Theorem 6.3.7.1} If
$$\r(\Lm,G)\geq 1\tag 6.3.7.2$$
then $\exp \Lm$ is complete in $G.$
\ep
This theorem is proved in \S 6.3.12.

Let $w:=g_Gq_\Lm(z)$ be the maximal $L_1$-subminorant of $q_\Lm.$
Denote by $\CH_\Lm$ the open set in $\BT_P^2$ where $L_1w =0.$
\proclaim {Theorem 6.3.7.2} If there exists a component $M$ of
$\CH_\Lm$ such that $\r(M)< 1$ then $w$ is minimal , \ep and,
hence, $G$ is maximal for  $\exp\Lm.$

This theorem follows directly from Th.6.2.3.1.

It is not known if the condition (6.3.7.2) is necessary.

Consider in details the situation, in which the domain $G$ coincides with
$G_\Lm$ , the conjugated indicator diagram of $h_\Lm,$ i.e., we suppose that

$$h_G(\phi)=h_\Lm (\phi),\ \forall \phi. \tag 6.3.7.3$$
In this case $m(z,G,q_\Lm)\geq 0$ and we obtain the following criterion:
\proclaim {Theorem 6.3.7.3} In order that  $\exp \Lm$ be complete in $G_\Lm$ it
is necessary and sufficient that
$$\r(\Lm,G_\Lm)\geq 1.\tag 6.3.7.4$$
\ep
This theorem is proved in \S 6.3.12.
The condition (6.3.7.3) automatically  implies maximality if there is completeness.

Since
$$h_\Lm(y)=\max \{q_\Lm(x+iy):x\in[0;P\}, \tag 6.3.7.5$$
the function $m(z,G_\Lm,q_\Lm)$ has a zero in $x$ for every fixed $y$.

Thus the set $D(G,\Lm)$ does not contain any curve $y=const$ on the torus.
\proclaim {Theorem 6.3.7.4} Let $G_0$ be a strictly convex domain and let
$D_0\sbt \BT_P^2$ be such that $\BT_P^2\setminus D_0$ intersect every line
$\{y=y_0\},\ y_0\in [0,2\pi].$

Then there exists a periodic $\Lm$ such that
$$G_\Lm=G_0,\ D(G_\Lm,\Lm)=D_0.\tag 6.3.7.6$$
\ep
This theorem is proved in \S 6.3.13.

{\bf Example 6.3.7.1} Let $D_0$ be the complement in $\BT_P^2$ to the set
$$M:=\{z=x+iy:x=f(y),y\in [0;2\pi]\}\tag 6.3.7.7$$
where $f(y)$ is a continuous $2\pi$-periodic function satisfying the condition
$$0<f(y)<P.$$ Then $\r(D_0)=\iy$ , because this domain is not connected on
spirals (see \S 6.2.). It means that for every strictly convex $G_0$ there exists
a periodic $\Lm$ such that $G_\Lm=G_0$ and $\exp\Lm$ is  extremely overcomplete in $G_0$.

{\bf Example 6.3.7.2.} Let $D_0$ be the complement to the set
$$M:=\{z=x+iy:x=\frac{P}{2\pi}y,\ 0\leq y\leq 2\pi\}$$
Then
$$\r(D_0)=\frac{1}{2}\left (1+(2\pi/P)^2\right )\tag 6.3.7.8$$
(see \S 6.3.13)

Thus, choosing $P$, and using Th.6.3.7.4,it is possible make $\exp \Lm$ complete or non-complete in $G_0$($=G_\Lm$) for every strictly
convex domain $G_0.$
\subheading {6.3.8} Now pass to generalizations. Denote by $D_G$ the natural
domain of definition of the operation $\CG_G,$ i.e. the set of $v\in U[1]$ for which $m(\lm,G,v)$ (see (6.3.6.2)) has a subharmonic minorant belonging to
$U[1].$

Let $\Phi_\Lm$ be defined by the equality (6.3.1.1).The condition that for every
$v\in\Fr [\Phi_\Lm]$ the function $m(\lm,G,v)$ has a subharmonic minorant belonging to $U[1]$ is possible to express by the relation
$$\Fr [\Phi_\Lm]\sbt D_G\tag 6.3.8.1$$
(compare with (6.1.1.1)).

We call the set $U\sbt U[1]$  {\it minimal} ($U\in MIN$) if for
arbitrary small $\eps>0$ there exists $w=w_\eps\in U$ such that
the function $w_\eps-\eps |\lm|$ has no subharmonic minorant,
belonging to $U[1].$

Let us note that if $U$ contains a minimal function in the sense of (\S 6.3.6),
then $U\in MIN.$

We denote the image of $\Fr [\Phi_\Lm]$ under the mapping by the operator $\CG_G$  as
$J_G(\Lm).$
\proclaim {Theorem 6.3.8.1}The follow holds:

1.\{$\exp_\Lm$ is not complete in $A(G)$\}$\Longleftrightarrow$
\{(6.3.8.1)\ holds $\wedge J_G\notin MIN;$\}

2.\{ $G$ is maximal for $\exp_\Lm$\}$\Longleftrightarrow$
\{(6.3.8.1)\ holds $\wedge J_G\in MIN;$\}

3.\{$\exp \Lm$ is  extremely overcomplete for maximal $G$\}
$\Longleftrightarrow$ \{(6.3.8.1)\ holds $\wedge J_G \in HARM;$\}
\ep
\subheading {6.3.9} In the proof of Theorem 6.3.8.1 that we are going to prove now we exploit
\proclaim {Theorem 6.3.9.1. (A.I.Markushevich) \cite {see,Lev,,Ch.4,\S7}}Let
$A(\BC\setminus \overline G)$ be a class of functions $\psi$ which are holomorphic in
$\BC\setminus \overline G$ and equal to zero in infinity. In order that the system
$\exp \Lm$ be complete in $A(G),$ it is necessary and sufficient that the function
$$\Phi (\lm):=\intl_{L_\psi}e^{\lm z}\psi (z)dz,\tag 6.3.9.1$$
where $\psi\in A(\BC\setminus \overline G),$ and $ L_\psi\Subset G$ is a rectifiable closed curve,has the following property: the condition
$$\Phi(\lm_k)=0,\ \forall \lm_k\in \Lm\tag 6.3.9.2$$
implies $\Phi(\lm)\equiv 0.$
\ep
\demo {Proof Theorem 6.3.8.1, 1.}Necessity. Let $\exp \Lm$ is not complete .
By Theorem 6.3.9.1 $\Phi(\lm_k)=0,$ but $\Phi (\lm)\not\equiv 0.$ The function
$g(\lm):=\Phi(\lm)/\Phi_\Lm(\lm),$ where $\Phi_\Lm$ is from (6.3.1.1), is
an entire function and it has  order one and normal or minimal type by
Th.2.9.3.1. Set
$$u^g:=\log|g(\lm)|;\ u^\Phi(\lm):=\log |\Phi(\lm)|;\ u^\Lm(\lm):=\log |\Phi(\lm)|.$$
We have from (6.3.9.1)
$u^\Phi(\lm)\leq \max\{\Re (\lm z):z\in L_\psi\}+C_\psi\},$
where $C_\psi$ is a constant, depending possibly on  $\psi.$

This implies that
$$u^\Phi(\lm)\leq h_{G_1}(\phi)r +C_\psi,\ \lm=re^{i\phi},\tag 6.3.9.3$$
for some convex domain $G_1\Subset G.$

Let $v\in \Fr [\Phi_\Lm].$ Choose a sequence $t_j\ri\iy$ for which
$(u^\Lm)_{t_j}\ri v,$ and the sequences $(u^\Phi)_{t_j}$ and $(u^g)_{t_j}$ also
converges to $v^\Phi$ and $v^g$ respectively. From the equality
$u^g(\lm)=u^\Phi(\lm)-u^\Lm(\lm)$ we obtain
$v^g(\lm)=v^\Phi(\lm)-v(\lm)$
where $v^g\in \Fr [g],\ v^\Phi\in \Fr [\Phi].$

Since (6.3.9.3) implies $v^\Phi(\lm)\leq h_{G_1}(\phi)r$
$$v^g(\lm)\leq  h_{G_1}(\phi)r-v(\lm)\tag 6.3.9.4$$
and it means that for every $v\in\Fr[\Phi_\Lm]$ $\CG_{G_1}v$ and hence $\CG_Gv$  exist, i.e., the condition (6.3.8.1) holds.

Let us show that the condition $J_G(\Lm)\notin MIN$ is satisfied. We have for some $\dl>0$ the relation
$$h_{G_1}(\phi)-h_G(\phi)\leq -\dl.$$
From (6.3.9.4) we obtain
$$v^g(\lm)+\dl r\leq m(\lm,G,v)\tag 6.3.9.5$$
The left hand side of the inequality (6.3.9.5) belongs to $U[1].$ Thus
$w_v:=\CG_{G_1}v$ satisfies the condition $v^g(\lm)+\dl r\leq w_v(\lm)$ for every $v\in\Fr[\Phi_\Lm].$ It means that $J_G(\Lm)\notin MIN.$

Necessity is proved.\qed\edm

For proving sufficiency we exploit the following assertion
\proclaim {Theorem 6.3.9.2. (I.F.Krasichkov-Ternovskii)} Suppose there exists an entire function $g$ such that
$$h_{g\Phi_\Lm}(\phi)<h_G(\phi),\ \forall \phi .\tag 6.3.9.6$$
Then the system $\exp \Lm$ is not complete for some convex domain
$G_1\Subset G.$
\ep
This theorem connects the problem of completeness to the multiplicator problem.
  \demo {Proof of Th.6.3.9.2} Let $g(\lm)$ satisfy (6.3.9.6). Denote by
$\psi(z)$ the Borel transformation for $\Phi(\lm):=g(\lm)\Phi_\Lm(\lm).$ By
P\'olya Theorem (see, for example, \cite {Lev.,Ch.1,\S 20}) all the singularities of $\psi$ are contained in a convex domain $G_\Phi$ which is the conjugate diagram
of the indicator $h_\Phi(\phi).$ Thus the representation (6.3.9.1) holds with $L_{\psi}$ that embraces $G_\Phi.$ It follows from (6.3.9.6) that $G_\Phi\Subset G.$ Thus it is possible to choose $L_\psi$ between $\p G_\Phi$ and $\p G.$
Since (6.3.9.2)for $\Phi$ is fulfilled  and $\Phi (\lm)\not\equiv 0,$  $\exp \Lm$ is non-complete in some convex $G_1\Subset G$ such that
$L_\psi \Subset G_1$ by Th.
6.3.9.1.\qed\edm

Now we can prove sufficiency in Th.6.3.8.1, 1.  From the condition
$J_G(\Lm)\notin MIN$ it follows that one can choose $\dl>0$ such that
$\forall v\in\Fr[\Phi_\Lm]$ the functions $w_v-\dl r$ where $w_v:=\CG_G$, have subharmonic minorants. As we already said in \S 6.3.2 completeness does not depend
of shift by any fixed $z_0.$ Thus we can suppose that $0\in G$ and, hence,
$h_G(\phi)>0$ for all $\phi.$  Let $\g<2\dl$ be such that
$h_G(\phi)-\g>0$ and $G_1\Subset G$ satisfy
$$h_{G_1}(\phi) -\g/3>0,\ h_G(\phi)-h_{G_1}(\phi)\leq \g/2\tag 6.3.9.7$$
Let us check that
$$D_{G_1}\spt \Fr[\Phi_\Lm],\tag 6.3.9.8$$
Indeed, for $v\in \Fr [\Phi_\Lm]$ we have
$$m(\lm,G_1,v):=h_{G_1}(\phi)r -v(\lm)\geq h_{G}(\phi)r -(\g/2) r -v(\lm)\geq$$
         $$h_{G}(\phi)r  -v(\lm)-\dl r\geq w_v -\dl r \tag 6.3.9.9$$
Since the right hand side of (6.3.9.9) has a subharmonic minorant from $U[1],$
then (6.3.9.8) is proved. By Theorem 6.1.1.1 there exists a multiplicator $g(z)\in A(1)$ such that
$$h_{g\Phi}(\phi)\leq h_{G_1}(\phi) <h_G(\phi).\tag 6.3.9.10$$
From Th.6.3.9.2 we obtain that $\exp \Lm$ is non-complete in $G.$

\demo {Proof Theorem 6.3.8.1, 2.}Necessity. Let $G_j,\ j=1,2,...$ be a sequence of convex domains, satisfying the conditions $G_j\Supset G,\ G_j\downarrow G.$ Since $\exp \Lm$ is not complete in every $A(G_j),$
$D_{G_j}\spt \Fr [\Phi_\Lm]$ by Th.6.3.8.1, 1.

The sequence $w_j:=\CG_{G_j} v$ satisfies
$$w_j(\lm)\leq h_{G_j}(\phi)r-v(\lm),\ \lm\in \BC.$$
Since $\{w_j\}$ is compact and $h_{G_j}\ri h_G,$ one can find a subsequence with the limit $w\in U[1].$ Then $w(\lm)\leq h_G(\phi)r-v(\lm).$ Hence $\CG_G v$
exists.

If $J_G\in MIN$ would not hold, then, by Th.6.3.8.1, 1., $\exp \Lm$ were non-complete in $A(G),$ which contradicts maximality.

Necessity is proved. Let us prove sufficiency.

Completeness of $\exp \Lm$ in $A(G)$ follows from Th.6.3.8.1, 1.. We will prove
that $\exp_\Lm$ is non-complete in $A(G_1)$ for every $G_1\Supset G$ under the condition $D_G\spt\Fr[\Phi_\Lm]$. Set
$$\dl:=\min\limits_\phi[h_{G_1}(\phi)-h_G(\phi)]>0$$
Then $\forall v\in\Fr[\Phi_\Lm]$
$$\CG_Gv +\dl r\leq h_{G_1}(\phi)r-v(\lm),\ \lm\in\BC.$$
This means that $\CG_{G_1}v\geq \CG_Gv+\dl r.$ Hence $J_{G_1}(\Lm)\notin
MIN$ and, by Th.6.3.8.1, 1., $\exp \Lm$ is non-complete in $A(G_1).$
\qed
\edm

\demo {Proof of Th.6.3.8.1, 3.} Necessity. By  Th.6.3.8.1, 2. from maximality $G$  (6.3.8.1) follows. We will prove that $\CG_Gv\in HARM \ \forall
v\in\Fr [\Phi_\Lm].$ Suppose it is not fulfilled, i.e., there exists
$v_0\in\Fr[\Phi_\Lm]$ such that the mass distribution $\nu _0$ of the function $w_0=\CG_Gv_0$ is not zero. By Proposition 6.1.1.3 there exists a multiplicator
$g$ such that $v_0+w_0\in \Fr[g\Phi_\Lm].$ Let $\Lm_0$ be the set of zeros of $g.$
Since $\nu_0\in \Fr\Lm_0,$ $\overline\Dl (\Lm_0)>0,$ because $\nu_0\neq 0$ and
by the definitions in \S 3.3.1.

We can shift a little zeros of $g$ and suppose without lack of generality  that they are simple and $\Lm_0\cap \Lm=\varnothing.$

The condition for a multiplicator gives the inequality:
$$h_{g\Phi_\Lm}(\phi)\leq h_G(\phi),\ \forall \phi.$$
It implies
$$m(\lm,G,v_\Pi)=rh_G(\phi)-v_\Pi\geq 0$$
for all $v_\Pi\in \Fr[g\Phi_\Lm].$
It means that $m(\lm,G,v_\Pi)$ has zero as a minorant
$\forall v_\Pi\in \Fr[g\Phi_\Lm]$, i.e.,
$D_G\spt \Fr[g\Phi_\Lm].$
So the domain $G$  maximal although the system $\exp\Lm$ is replaced with the system $\exp (\Lm\cup \Lm_0).$ This contradicts to  extremal
overcompleteness. Hence, $\nu_0\equiv 0$ and $w_0=\CG_Gv_0\in HARM$.

Necessity is proved. Let us prove sufficiency.

Let the condition $\CG_Gv\in HARM \ \forall v\in\Fr[\Phi_\Lm]$ hold. Suppose that there exists $\Lm_0$ such that $\overline \Dl_{\Lm_0}>0$ and $G$ is maximal for the system $\exp(\Lm\cup\Lm_0).$

Theorem 6.3.8.1,, 2. implies
$$D_G\spt\Fr[\Phi_{\Lm_1}],\tag 6.3.9.11$$
where $\Lm_1=\Lm\cup \Lm_0.$

For every $v_0\in\Fr[\Phi_{\Lm_0}]$ one can find $v\in\Fr[\Phi_\Lm]$ such that
$$v_1:=v_0+v\in \Fr[\Phi_{\Lm_1}]$$
The condition $\overline\Dl_{\Lm_0}>0$ implies that one can choose $v_0$ for which the Riesz measure $\nu_0\not\equiv 0$. For $w_1=\CG_Gv_1$ one has the inequality
 $w_1\leq rh_G-v_1$ by (6.3.9.11), so $w_1+v_0\leq rh_G-v$ holds.
Hence $w_v:=\CG_Gv$ satisfies the inequality
$$(w_1+v_0)(\lm)\leq w_v(\lm),\ \forall \lm\in \BC.\tag 6.3.9.12$$
Let us show that (6.3.9.12) is impossible.Indeed, since $w_v\in HARM$
$w:=w_1+v_0-w_v\leq 0$ and $ w\in U[\r]$. Thus $w\equiv 0.$ However the Riesz
measure $\nu_w\geq \nu_0\not\equiv 0,$ hence $w\not\equiv 0.$ This contradiction proves sufficiency.
\qed
\edm
\subheading {6.3.10}Now we prove Theorems 6.3.3.1, 6.3.4.1 and 6.3.5.1.
We need some auxiliary assertions.
\proclaim {Lemma 6.3.10.1} Let $v:=rh_1(\phi)$ and $G_1$ be the conjugated diagram of $h_1.$ Then the following holds:

1.\{ $G_1$ is freely enclosed in $G\}\Longleftrightarrow
\{\CG_Gv$ is non-minimal\};

2.\{$G_1$ is enclosed to $G$ but not free enclosed\}$\Longleftrightarrow
\{\CG_Gv$ is minimal\};

3.\{ $G_1$ is  rigidly enclosed in $G\}\Longleftrightarrow
\{\CG_Gv\in HARM$ \};

4.\{ $G_1$ is not enclosed in $G\}\Longleftrightarrow
\{\CG_Gv$ does not exist\};
\ep

To prove this lemma we need the following two:
\proclaim {Lemma 6.3.10.2} Let $v:=rh(\phi).$ Then $\CG_Gv=rh_1(\phi)$ where
$h_1$ is the maximal trigonometrically convex minorant of the function
$$m(\phi,G,h):=h_G(\phi)-h(\phi)$$.
\ep
\demo {Proof} Let $v_1=\CG_Gv.$ Since $v_{[t]}=v$ for all $t>0$
$$(v_1)_{[t]}=\CG_G v_{[t]}=\CG_Gv$$
by Th.6.1.1.2, 2.

Thus the function
$$\hat v_1:=\left(\sup\limits_t (v_1)_{[t]}\right )^*(\lm)\geq v_1(\lm)$$
and is also a subharmonic minorant belonging to $U[1]$. Thus $v_1=\hat v_1$
. However, the function $\hat v_1$ is invariant with respect to the transformation $(\bul)_{[t]}.$ Hence it has the form $rh_1(\phi).$The maximality of $h_1(\phi)$
follows from the maximality $v_1.$
\qed
\edm
\proclaim {Lemma 6.3.10.3}In order that $v:=rh_1$ be a minimal function it is necessary and sufficient that $G_1,$ the conjugate diagram of $h_1,$ be a segment (in particular, a point).
\ep
\demo {Proof}Let $v=rh_1$ be minimal and let $G_1$ be the conjugate diagram of
$h_1.$ If $G_1$ is not segment, then it contains some disc of radius
$\dl>0.$ Hence there exists a trigonometric function
$A \cos (\phi-\phi_0)$ such that
$$\dl+A\cos(\phi-\phi_0)\leq h_1(\phi).$$
Multiplying this inequality by $r,$ we obtain that $v-\dl r$ has a
harmonic (and  hence subharmonic) minorant. This contradicts to
minimality.

Inversely, Suppose $v$ is not minimal. Then there exists $\dl>0$ and  t.c.f. $h_2(\phi)$ such
that
$$h_2(\phi)\leq h_1(\phi)-\dl.\tag 6.3.10.1$$
For every t.c.f. $h_2$ there exists a trigonometric function
$A\cos(\phi-\phi_0)$ such that
$$h_2(\phi)+A\cos(\phi-\phi_0)\geq 0\tag 6.3.10.2$$
This corresponds to a shift of the diagram  which contains zero. From
(6.3.10.1) and (6.3.10.2) we obtain
$$\dl-A\cos(\phi-\phi_0)\leq h_1(\phi),$$
which means that $G_1$ contains some disc of radius $\dl>0.$ So it
is not a segment. \qed \edm \demo {Proof of Lemma 6.3.10.1} $G$ is
freely enclosed iff the following assertion holds: there exists
$\dl>0$ and a trigonometrical function $A\cos(\phi-\phi_0)$ such
that the inequality
$$h_1(\phi)+\dl-A\cos(\phi-\phi_0)\leq h_G(\phi).\tag 6.3.10.3$$
holds.

{\bf Exercise 6.3.10.1.} Prove this.

Let $\CG_Gv$ be non-minimal. By Lemma 6.3.10.2 it has the form $w_2=rh_2,$
where $h_2$ is the maximal trigonometrically convex minorant of
$m(\phi,G,h_1).$ There exists $\dl>0$ such that the function $w_2-\dl r$ has
the maximal subharmonic minorant $v_3=rh_3(\phi).$ Let $A\cos(\phi-\phi_0)$ be a trigonometric function for which
$$h_3(\phi)+A\cos (\phi -\phi_0)\geq 0.$$
In addition,
$$h_3(\phi)\leq h_2(\phi)-\dl, \ h_2(\phi)\leq h_G-h_1(\phi)$$
From this we obtain (6.3.10.3) and hence that $G_1$ is free enclosed.

Inversely, let $G_1$ is freely enclosed in $G.$ From (6.3.10.3) it follows that
$$\dl-A\cos(\phi-\phi_0)\leq h_G(\phi)-h_1(\phi).\tag 6.3.10.4$$
Multiplying (6.3.10.4) by  $r,$ we obtain that $m(\lm,G,v)$ has a minorant
$v_0=r(\dl-A\cos (\phi-\phi_0))$ which obviously is non-minimal. Hence,
$\CG_Gv $ is non-minimal.

 $G_1$ is enclosed in $G$ with sliding, hence there does not exists
$\dl>0$ such that (6.3.10.3) is fulfilled, but there exists a segment with
support function $$
E(\phi)=B|\sin\phi|+A\cos(\phi-\phi_0),$$
such that the inequality
$$h_1(\phi)+E(\phi)\leq h_g(\phi)\tag 6.3.10.5$$
holds.

{\bf Exercise 6.3.10.2} Prove this.

Let $\CG_Gv$ be  minimal. By Lemma 6.3.10.2 it has the form $w_2=rh_2$ and by
Lemma 6.3.10.3 $h_2=E(\phi).$ Thus $E(\phi)\leq (h_G-h_1)(\phi),$ which is equivalent to (6.3.10.5).

Prove 2., suppose $G$ is not freely enclosed and hence only (6.3.10.5) is possible.
If $\CG_Gv$ were non-minimal, (6.3.10.3) would  follow, as it was proved above. This  contradicts the supposition.

The rigid enclosure is equivalent only to the inequality of the form
$$h(\phi)-A\cos(\phi-\phi_0)\leq h_G(\phi) \ \forall \phi ,$$
and impossibility of enclosure is equivalent to the impossibility of even such
an inequality.
Thus all other assertions of the lemma can be proved analogously.

{\bf Exercise 6.3.10.3} Do this in details.

\qed
\edm

\demo {Proof Theorem 6.3.3.1} Regularity of $\Lm$ means that
$Fr [\Phi_\Lm]=\{v_0\}$ where $v_0=rh_\Lm.$ Thus $J_G(\Lm)=\{\CG_Gv_0\}$ and all the assertions of Theorem 6.3.3.1 follows from Th. 6.3.8.1 and Lemma 6.3.10.1.
\qed
\edm

{\bf Exercise 6.3.10.4} Check this in details.

For proving Theorem 6.3.4.1 we need an additional
\proclaim {Lemma 6.3.10.4} Let $\Lm$ an indicator set, $v_1=rh_1,\ v_2=rh_2.$
Then
$$\{J_G(\Lm)\notin MIN\}\Longleftrightarrow \{ \CG_Gv_1\ and\ \CG_Gv_2 \ are\ non-minimal\}.$$
\ep
 \demo {Proof}Suppose $w_1:=\CG_Gv_1$ and $w_2:=\CG_Gv_2$ are not minimal,
i.e., $w_1-\dl r$ and $w_1-\dl r$ have  subharmonic minorants $g_1$ and $g_2.$

Then $cg_1+(1-c)g_2$ is a minorant of the function $cw_1+(1-c)w_2 -\dl r$,
 i.e., $J_G(\Lm)\notin MIN.$ The inverse implication is trivial.
\qed
\edm

{\bf Exercise 6.3.10.5} Prove this.

\demo {Proof of Theorem 6.3.4.1} Suppose $\exp \Lm$ is not complete. By theorem 6.3.8.1 $J_G\notin MIN.$ By Lemma 6.3.10.4 $\CG_Gv_1$ and $\CG_Gv_2$ are not minimal. Hence $G_1$ and $G_2$ are freely enclosed in $G$ by Lemma 6.3.10.1.
Since every one of these assertions is reversible, the inverse implication also holds. Analogously the other cases are proved.
\qed\edm

{\bf Exercise 6.3.10.6} Prove all this in details.

\subheading {6.3.11} To prove Theorem 6.3.5.1 we need some auxiliary assertions.
\proclaim {Lemma 6.3.11.1} Let $\phi_0$ be a maximum point of t.c.f. $h(\phi)$ and $h(\phi_0)\geq 0.$ Then
$$h(\phi)\geq h(\phi_0)\cos(\phi-\phi_0), \ |\phi-\phi_0|\leq \pi/2 .\tag 6.3.11.1$$
\ep
\demo {Proof} Denote $y(\phi):=h(\phi_0)\cos (\phi -\phi_0).$ We have
$y(\phi_0)=h(\phi_0)$ and $y(\phi)$ is a trigonometric function. If
$y(\phi_1)=h(\phi_1)$ for some $\phi_1$ such that $|\phi_1-\phi_0|<\pi/2$ this
contradict to Th.3.2.5.2. If $y(\phi)$ does not intersect $h(\phi),$ this contradicts to Prop.3.2.5.6.
\edm
\qed
\proclaim {Lemma 6.3.11.2}Let $H(\phi)$ be a trigonometric function on the interval $I=(\a,\be)$ of length $\leq \pi,$ such that $H(\phi)=0$ at one of the ends of $I$. Then every of the conditions

1.$H(\phi_0)=0,\ \phi_0\in (\a;\be);$

2.$H(\phi)$ is zero on $\p I$ with tangency;

implies $H(\phi)\equiv 0,\ \phi\in I.$
\ep

{\bf Exercise 6.3.11.1} Prove this.

\proclaim {Lemma 6.3.11.3} Let $g\geq 0$ be a continuous periodic
function, and let $\Theta_\Lm,I_\Lm,d_\Lm$ be defined as in
Theorem 6.3.5.1. In order that its maximal t.c. minorant be a
trigonometrical function, it is necessary and sufficient
satisfying at least one of the conditions:

1.$d_\Lm<\pi;$

2.$d_\Lm=\pi$ and $g(\phi)$ is zero with tangency on $\p I.$
\ep
\demo {Proof} Necessity. Suppose $d_\Lm>\pi.$ Without loss of generality we can suppose that $I_\Lm=(\a;-\a)$, where $\a>\pi/2$.

Set $\cos^+\phi :=\max (\cos\phi,0),$
$$a=\inf\left(\frac {g(\phi)}{\cos^+\phi}:\phi\in (-\a;\a)\right).\tag 6.3.11.2$$
We have $a>0.$

Set
$$\cases a_1 \cos \phi,&\ |\phi|\leq \pi/2\\ 0&\ |\phi|>\pi/2,\endcases\tag 6.3.11.3$$
where $a_1\leq a.$

The function $h(\phi)$ is a t.c.minorante of $g(\phi)$ and it is not a trigonometric function, which contradicts the supposition. Thus $d_\Lm\leq\pi.$

Suppose $d_\Lm =\pi$ and the condition to be zero with tangency on
$\p I$ does not hold. Then for $a$ defined by (6.3.11.2) the
condition $ minorant a>0$ holds and $h(\phi)$ defined by
(6.3.11.3) is a non-trigonometric of $g.$

Sufficiency. The  first condition holds and $I=(\a;\be)$ be an arbitrary interval belonging to $\Theta_\Lm ;$ let $h(\phi)$ be the maximal t.c.minorante of
$g(\phi).$

Set
$$H(\phi):=h(\phi_0)\cos(\phi-\phi_0),$$
where $\phi_0$ is the maximum point of $h(\phi)$ on $I.$ From inequality
 (6.3.11.1) and the conditions $g(\a)=g(\be)=0$ follows $H(\a)=H(\be)=0.$
Then, by Lemma 6.3.11.2, we obtain $H(\phi)\equiv 0.$ Thus $h(\phi_0)=0$
and $h(\phi)\equiv 0$ for $\phi\in (\a;\be)$, i.e., $h(\phi)$ is trigonometric.

Let the second condition be fulfilled.  Lemma 6.3.11.1 implies that $H(\phi)$
is zero with tangency on $\p I.$ By Lemma 6.3.11.2 we obtain that $h(\phi)\equiv 0.$
\edm
\qed

\demo {Proof of Theorem 6.3.5.1} Necessity. Let us note that if
$v\in \Fr[\Phi_\Lm]$ then for $c\in [0;1]$ we have the equality
$$m(\lm,G,v)=r\left (c(h_2-h_1)^++(1-c)(h_2-h_1)^-\right)(\phi):=
rm(\phi,c)\tag 6.3.11.4$$
Let $\exp\Lm$ be  extremely overcomplete in $A(G)$. By Th. 6.3.8.1
$J_G\sbt HARM,$ i.e., for every $c\in[0;1]$ the maximal t.c.minorant of the function $m(\phi,c)$ is trigonometric. Since
$$\forall c\in [0;1],\ \Theta_\Lm =\{\phi:m(\phi,c)>0\}$$
the necessity follows from Lemma 6.3.11.3.

Sufficiency. It follows directly from Lemma 6.3.11.3. \edm\qed

{\bf Exercise 6.3.11.2} Explain this.

\subheading {6.3.12} To prove Theorem 6.3.6.1 we need
\proclaim {Lemma 6.3.12.1} If $w\in U[1]$ is non-minimal then
$$\BC(w):=\{w_{[t]}:1\leq t\leq e^P\}\notin MIN$$
\ep
It follows from Th.6.1.1.2, 2.

{\bf Exercise 6.3.12.1} Explain this in details.

\demo {Proof of Theorem 6.3.6.1} By Theorem 6.1.1.2, 2.
$J_G (\Lm)=\BC(\CG_Gv).$ Thus  Lemma 6.3.12.1 implies that $J_G(\Lm)\notin MIN$
 if and only if $\CG_Gv$ is not minimal. Thus Theorem 6.3.8.1 implies
Th. 6.3.6.1, 1. and 2.

Suppose $J_G(\Lm)\sbt HARM$. Hence, $\CG_Gv=rH_0(\phi),$ where $H_0$ is trigonometric. Inversely, Lemma 6.3.12.1 implies $J_G(\Lm)=\{rH_0(\phi)\}.$
\qed\edm

\demo {Proof of Theorem 6.3.7.1} Let $\r(\Lm,G)>1.$ Suppose $w_q:=g_Gq$ exists.
By definition of $\r(\Lm,G)$ we have $w_q(z)\leq 0$ for $z\in \p D(\Lm,G).$ By
Th.6.2.3.3 (Maximum principle $w_q\leq 0$ for $z\in D(G,\Lm).$ Also
 $w_q\leq 0$ for $z\in \BT_P^2\setminus D(G,\Lm)$ by definition of $D(\Lm,G).$
By Theorem 6.2.3.4 $w_q\equiv 0$ and hence is minimal. So $exp\Lm$ is complete by Th.6.3.6.1.
If $\r(\Lm,G)=1,$ then the system $\exp \Lm$ is complete for every $G_n\Supset G,$ because of strict monotonicity of $\r(\bul)$ (see \S 6.2.2) so $G$ is the
maximal domain.
\qed\edm

For proof of Theorem 6.3.7.3 we need an auxiliary assertion. We suppose that
$D$ is an image on $\BT_P^2$ by the map (6.2.1.3) of the domain $G$ with a smooth boundary.
\proclaim {Theorem 6.3.12.2} Let $D\sbt \BT_P^2$ and $\r (D)\leq 1.$ Then
$\r (\BT_P^2\setminus \overline D)>1,$ if $D_\Lm\neq \{\Re z>0\}.$\ep

For the proof we need the following assertion which was proved  originally by
A.Eremenko and M.Sodin:

\proclaim {Theorem 6.3.12.3(Eremenko,Sodin)} Let $\Gamma$ be a Jordan curve, connecting $0$ and
$\iy,\  T\Gamma=\Gamma$ for some $T>1.$ Let $D_+,D_-$ be domains, into which $\Gamma$
divides the plane, and let $\r_1,\r_2$ be the orders of the minimal harmonic functions in
$D_+$ and $D_-$ respectively.

Then
$$\frac{1}{\r_1}+\frac{1}{\r_2}\leq 2,$$
and equality is attained only if $\Gamma$ consists of two rays.\ep

 We will give prove this theorem in \S 6.3.14 .
\demo {Proof of Theorem 6.3.12.2}  Let $\r_1=\r(D),$ and suppose $ q_1(z)$ is
a solution of boundary problem (6.2.2.1),
$\r_2=\r(\BT_P^2\setminus \overline D), \ q_2(z)$ is a solution of corresponding
boundary problem. Then the image of the boundary under the map $\lm=e^z$ (we denote
it as
$\Gamma$) satisfies the conditions of Theorem 6.3.12.3  and the functions $v_1(\lm):=q_1(\log \lm)|\lm|^\r_1$ and
$v_2(\lm):=q_2(\log \lm)|\lm|^\r_2$ are positive harmonic functions in
$D_+,D_-$ with orders $\r_1$ and $\r_2$ respectively. By Theorem 6.3.12.3
we obtain
$$1/\r (D)+1/\r(\BT_P^2\setminus \overline D)\leq 2.$$
and  equality holds only if $\Gamma$ is a pair of rays, i.e.,
$D_\Lm=\{\Re z>0\}.$ \qed \edm \demo {Proof Theorem 6.3.7.3}
Necessity. Suppose $\r(\Lm,G_\Lm)<1.$ Let us prove that $exp \Lm$
is not complete. To this end we construct an $L_1$ -minorant of
$m(z,G_\Lm,\Lm)$ and prove that it is not minimal.

Let $D_0\Subset D(G_\Lm,\Lm)$ be a domain with smooth boundary for which
$\r(D_0)=1.$ This is possible because of strict monotonicity $\r(D)$
(\S 6.2.2). Let $q_0$ be a solution of the problem (6.2.2.1) satisfying the condition
$$0<\max\{q_0(z):z\in D_0\}\leq\{m(z,G_\Lm,\Lm);z\in D_0\}-2\eps$$
for sufficiently small $\eps.$
By Theorem 6.3.12.2 $\r(\BT_P^2\setminus \overline D_0)>1.$ Thus the potential
$$\Pi(z)=-\int_D G_\r(z,\z,D)\nu(d\z)$$
exists and $\nu$ can be chosen in such way that supp $\nu\Subset
\BT_P^2 \setminus \overline D_0.$ By Proposition 6.2.3.6
$$\frac{\p q_0}{\p n}>0,\ z\in D_0$$
Thus $\nu$ can be chosen in such way that
$$-\frac {\p \Pi}{\p n}<\min \frac{\p q_0}{\p n}, \ z\in \p D_0.$$
Then the function
$$q(z)=\cases q_0(z), &\ z\in D_0\\ \Pi(z), &\ z\in \BT_P^2\setminus D_0,
\endcases$$
is an $L_1$-subfunction on $\BT_P^2.$

{\bf Exercise 6.3.12.2} Explain this in details, exploiting Th.2.7.2.1.

The function $q(z)$ satisfies the condition
$$q(z)\leq m(z,G_\Lm,\Lm)-2\eps,\  \ forall z\in \BT_P^2, $$
because of negative potential. Hence,
$$q_1(z):=q(z)+\eta$$
for some $\eta>0$ also is a minorant of $m(z,G_\Lm,\Lm)$ and it is not minimal.
Necessity is proved. Sufficiency follows from Th.6.3.7.1.
\qed
\edm
\subheading {6.3.13} Now we pass to the proof of Theorem 6.3.7.4 and construction of Example
6.3.7.2.
\demo {Proof of Theorem 6.3.7.4} The set $\BT_P^2\setminus D_0$ is closed.
Let $\phi (z)$ be an infinitely differentiable function equal to zero on
$\BT_P^2\setminus D_0$ and positive on $D_0.$
Set
$$q(z):=h_0(y)-\eps \phi(z),$$
where $h_0(y)$ is a t.c.f., corresponding to $G_0$ ,
and let $\eps $ be small enough to satisfy $L_1q(z)>0,\ z\in \BT_P^2.$
It is possible, because $L_1h_0(y)>0$ by the condition of the theorem.

Then we have
$$m(z,G_0,q)=\eps\phi(z).$$
hence
$$\{z:m(z,G_0,q)>0\}=D_0.$$
Take
$$v(\lm):=|\lm|q(\log\lm)$$
and construct an entire function $\Phi_\Lm$ for which
$$\Fr[\Phi_\Lm]=\{v_{[t]}:1\leq t\leq e^P\}$$
It is easy to check that the zero distribution of this function has all the properties demanded by Theorem 6.3.7.4.
\qed\edm

{\bf Exercise 6.3.13.1} Check this.

\demo {Proof of (6.3.7.8)} Consider the problem
$$L_1 q(z)=0 ,\  q|_{x=(2\pi/P)y}=0\tag 6.3.13.1$$

Let us pass in the equation to new coordinates
$$\cases \xi=x\cos\a+y\sin\a \\ \eta=-x\sin\a +y\cos\a,\ \tan\a =2\pi/P\endcases.$$
Then the equation takes the form:
$$\left [\frac {\p ^2}{\p \xi ^2}+\frac {\p ^2}{\p \eta  ^2}
+2\r\left (\cos\a \frac {\p }{\p \xi}-\sin\a \frac {\p }{\p \eta} \right )
r^2\right ]R(\xi,\eta)=0$$
The condition of being zero on $D_0$ is
$$R(\xi,2\pi l\cos\a)=0.\ l\in \BZ.$$
The condition of periodicity gives
$$R(\xi+(P/\cos\a)k,\eta)=R_1 (\xi,\eta),\ k\in\BZ.$$
We search for a solution that does not depend of $\xi.$
We have
$$R''(\eta)-2\r\sin\a R'(\eta)+\r^2R(\eta)=0,\ R(0)=R(2\pi\cos\a)=0$$
Further,
$$R(\eta)=C_1e^{(\r\sin\a)\eta}\cos((\r\cos\a)\eta)+
 C_1e^{(\r\sin\a)\eta}\sin((\r\cos\a)\eta)$$
Exploiting the boundary condition, we have
$$\r_{\min}=(2\cos^2\a)^{-1}=\frac {1}{2}
\left [1+\left (\frac {2\pi}{P}\right )^2\right ].$$
The corresponding  eigenfunction
$$R=\exp{(\r_{\min}\sin\a)\eta}  \sin((\r_{\min}\cos\a)\eta).$$
It is zero on $\BT_P^2\setminus D_0$ and positive in $D_0,$ so it determined
up to a constant multiple.
\qed
\edm

\subheading {6.3.14}We are going to prove Theorem 6.3.12.3. Actually we
prove
\proclaim {Theorem 6.3.14.1}Let $\Gamma_1,\Gamma_2,...,\Gamma_n$ be Jordan
curve, such that

1.$\Gamma_i,\ i=1,2,...,n$ connect $0$ and $\iy;$

2.there exists a number $T,|T|>1$ (not necessarily real) for which
$T\Gamma_i=\Gamma_i,\ i=1,2,...,n.$

Let $D_i,\ i=1,2,...,n$ be domains into which the plane is divided, and let $\r _i$ be
the order of the minimal harmonic function in  $D_i$. Then
$$\sum\limits_i 1/\r _i\leq 2\tag 6.3.14.1$$
and equality holds if and only if $\Gamma_i$ are a logarithmic spirals (or
rays, when $T\in \BR_+$).\ep

\demo {Proof}Denote by $H_i$ the minimal harmonic function in $D_i$. Then
$H_i=\Im \phi _i$ where $\phi _i:D _i\mapsto \Pi ^+$ is a conformal map of
$D_i$ to upper half plane, $\phi (0)=0.$ The maps
$g_i:=\phi_i(T\phi _i^{-1}):\Pi^+\mapsto \Pi^+$ does continued by isomorphism to
$\BC,$ and $g_i(0)=0.$  Thus $g_i(z)=\s _i z$, where $\s _i>1.$ Hence,
$\phi _i(Tz)=\s_i\phi _i(z)$ or
$$Th_i(z)=h_i(\s _iz),\ h_i:=\phi _i^{-1}:\Pi^+\mapsto D _i.$$
Now we exploit the following inequality from \cite {LevG}
$$\sum\limits_{i=1}^n \frac {1}{\log\s _i}\leq \frac {2\log T}{|\log T|^2}\leq
\frac{2}{\log T}.\tag 6.3.14.2$$ \edm

\newpage

\centerline {\bf Notation} \subheading {2.1} $\Bbb R^m ,\
M(f,x,\varepsilon),\ f^* (x),\ C^+(E), \ C^- (E),\chi _G,\
\chi_F,\ \G_A,\ F^A,\ K,\ M(f,K),$ $\ K_n,\ K_{\max}.\  $
 \subheading
{2.2} $\sigma (G),\ \mu,\ G_0 (\mu),\ \text {supp} \mu, \ \Cal M
(G),\ \mu _F (E),\ \nu,\ {\Cal M}^d,\ \nu ^+,\ \nu ^-,\ |\nu|,\
\roman{supp}\ \phi,$ $\ \overset{*}\to\ri,\ \overset \circ \to E,\
\overline E,\ \sigma ({\Bbb R}^{m_1} \times {\Bbb R}^{m_2}),\
\Phi_1 \otimes \Phi_2 ,\ \mu_1\otimes \mu_2.$ \subheading {2.3}
$\varphi_n \overset {\Cal D} \to {\rightarrow} \varphi,\ \alpha
(t),\ \alpha _\varepsilon (x),\ \psi_\varepsilon (x),\
<f,\varphi>,\ <\delta _x ,\varphi>,\ <\delta^{(n)} _x ,\varphi>,\
<\mu,\varphi>,\ <\alpha f,\varphi >,\ <f_1+f_2,\varphi>,\
<\frac{\partial}{\partial x_k} f,\varphi >,$ $\ f_\epsilon (x),\
f\mid _{G_1},\ \widetilde {\cos\r}(\phi).$ \subheading {2.4}$
\Delta _{\bold x ^0},\ \Cal E_m (x),\ \theta _m ,\ G(x,y,\Omega),\
G(x,y,K_{a,R})$ \subheading {2.5}$\Pi (x,\mu,D),\ G_{N}(x,y),\
\Pi_N (x,\mu,D),\ \Pi (x,\mu),\ \Pi(z,\mu),\ U[\r]\text {{\bf
cap}} _G (K,D),\ \text {\bf cap}_m (K) ,\ $ $\text {\bf cap}_m
(D),\ \overline {\text {\bf cap} }_m (E),\ \underline {\text {\bf
cap} }_m (E),\ \text {\bf cap}_l (K)$ \subheading {2.6}$\Cal M
(x,r,u),\ \Cal N (x,r,u),\ E^\epsilon,\ D^{-\epsilon},\ u_\epsilon
(x),\ K_R,\ M(r,u),\ \mu (r,u),\ \Cal M (r,u),\ N(r,u),$ $\ M(z).$
\subheading {2.7}$\tilde u(x),\ \mu_x(t),\
E(\alpha,\alpha',\epsilon,\mu ),\ E_{n,\delta_0}.$ \subheading
{2.8}$a(r),\ \rho[a],\ \sigma [a],\ \r(r),\ \sigma [a,\rho (r)],\
V(r),\ L(r),\ \delta SH(\Bbb R^m),\ T(r,u),\ \rho_T [u],\ \sigma_T
[u],$ $\ \sigma_T [u,\rho(r)],\ $ $ \rho_M[u],\ \sigma_M [u],\
\sigma_M [u,\rho(r)],\ \rho [\mu],\ \bar \Delta [\mu],\ \bar\Delta
[\mu,\rho(r)],\ N(r,\mu),\ \rho_N[\mu],\ \delta\Cal M (\Bbb R^m).$
\subheading {2.9}$H(z,\cos\gamma,m),\ G(x,y,\Bbb R^m),\ D_k
(x,y),\ H(z,\cos\gamma,m,p),\ G_p(x,y,m),\ G_p(z,\zeta,2),$ $\ \Pi
(x,\mu,p),$ $\ \delta SH (\rho),\ \Pi_<^R (x,\nu,\rho - 1),\ $ $
\Pi_>^R (x,\nu,\rho),\ \delta _R(x,\nu,\rho),\ \delta _R
(z,\nu,\rho),\ \delta_R (x,u,\rho),$
 $\ M(r,\delta),\ \bar \Delta _\delta [u,\rho],$
$ \ \Omega [u,\rho(r)],\ T(r,\lambda,>),\ T(r,\lambda,<).$
\subheading {3.1}$V_t,\ P_t,\ SH (\Bbb R^m,\rho,\rho(r)),\ SH
(\rho(r)),\  u_t (x),\ \bold {Fr}[u,\rho(r),V_\bullet,\Bbb R^m],\
U[\rho,\sigma],\ U[\r],$  $\ v_{[t]},\ \Cal M (\Bbb R^m,\rho (r)),\ $
$\mu\in \Cal M (\rho (r)),$$\ \bold
{Fr}[\mu,\rho(r),V_\bullet,\Bbb R^m],\ \Fr [\mu],\ \Cal
M[\rho,\Delta],\ \Cal M[\rho],\ \nu_{[t]}.$ \subheading
{3.2}$h(x,u),\ \lh (x,u),\ l_{\bold x ^0},\ x^0 (x),\ T_\rho,\ G_I
(\phi,\psi),\ \Pi_I(\phi,ds),\ TC_\rho,\ Co_\Om .$ \subheading
{3.3}$\overline \Delta (G,\mu),\ \overline \Delta (E,\mu),\
\underline \Delta (K,\mu),\ \underline \Delta (E,\mu),\
Co_{\Om}(I),\ \vdelt ^{cl}(E),\ \ndelt ^{cl}(E),\ \Om^G(\eps),\
\Om^K(\eps).$ \subheading {4.1}$T^t,\ (T^\bul,M),\ d(\bul,\bul),\
\Om (T^\bul),\ \Bbb C (m),\Om (m),\ A(m),\ T_t v.$ \subheading
{4.2}$U_0,\ \be(x),\ b_0,\ k(s),\ R_\eps v(x),\ Str(\dl),\
v(x|t),\ v(\bul,t)$ \subheading {4.3}$w(\bul|t),\ w(\bul|\bul).$
\subheading {4.4}$\bold u,\ (\bold u)_t,\ \Fr [\bold u],\ \bold
U[\r]. $ \subheading {5.1}$M(r,f),\ T(r,f),\ \r_T [f],\  \r_M
[f],\ \s_T[f,\r(r)],\ \s_M[f,\r(r)],\ n(K_r),\ n(r),\ \r[n],\
\vdelt [n],$
 $ \ N(r,n),\ \r_N [n],\ \vdelt_N [n],\ p[n],\ \Fr [f],\ \Fr [n],\  Mer(\r,\r(r)),\ T(r,f),\ \r_T[f],\
\s_T[f,\r(r)].  $ \subheading {5.2}$\a -\overline {\mes }C,\
C_0^\a,\ C_0^0.  $ \subheading {5.3}$\|g\|_p. $ \subheading
{5.4}${\underline h}_1(\phi,f),\ {\underline h}_2(\phi,f),\
\underline h (\phi,f). $ \subheading {5.5}$N(\delta,\Cal X),\
(\Cal X)\int fd\delta,\ (\Cal X)\int_E fd\delta,\ \delta
(\Theta^F),\ D_{r,\Theta},\ \delta_z(D_{r,\Theta}),\
A^{cl}(\delta,\chi _\Theta)$ \subheading {5.7}$\Cal F(u),\ H_\phi
(u),\ T(u),\ M_\a (u),\ M(u), \ I_{\a\be}(u),\ I(u,g),\ \overline
{\Cal F} [f], \ \underline {\Cal F} [f],\ \chi_H,\ \chi_I,\
\chi_{Fo}.$
 \subheading {5.8}$K_{S_1},\ S_1,\ G(t,\g,\r),\ \hat \BG(s,S_1-S),\ (\Cal F\nu)(s).$
\subheading {6.1} $H(z),\ m(z,v,H),\ \Cal G_H,\ D_H,\ U_{ind}, \
\hat {U}_{ind}.$ \subheading {6.2}${\BT}^2_P,\ \Di'({\BT}^2_P),\
q(z),\ L_\rho,\ E_\r(\bul-\z), \ E'_\r(\bul-\z),\ q_D,\
G_\r(z,\z,D),\ \Cal H_\rho(q).$ \subheading {6.3}$\Lm ,\ \Phi_\Lm
(\lm),\ \exp\Lm,\ A(G),\ h_\Lm (\phi), \ G_\Lm ,\ \a G_1+\be G_2,\
\Theta_\Lm,\ I_\Lm,,\ d_\Lm\ h_G(\phi), \ m(\lm,G,v),$ $ \ H(\lm),$
 $\ q_\Lm(z),\ D(G,\Lm),\ \r (\Lm, G),\ g_Gq,\ \CG_G,\ D_G,\ MIN,
\ J_G(\Lm),\ HARM,$  $\ m(\phi,G,h),\  E(\phi). $

\newpage

\centerline {\bf List of Terms}
\subheading {2.1}
2.1. upper semicontinuous regularization

2.1.  upper semicontinuous function

2.1. lower semicontinuous function

\subheading {2.2}
2.2. measure

2.2. mass distribution

2.2.  support of $\mu$

2.2. $\mu$ is {\it concentrated} on $E\in \sigma (G)$

2.2. {\it restriction} of $\mu$ onto $F\in \sigma (G)$.

2.2.  charge

2.2. {\it positive} and {\it negative}, respectively, {\it variations} of $\nu$

2.2. {\it full variation} of $\nu$

2.2. variation

2.2. {\it Borel function}

2.2. {\it restriction} of $\mu$ on the set $E$

2.2.  product of measures
\subheading {2.3}
2.3. linear space

2.3. {\it topological} space

2.3. {\it linear continuous functional} on $\Cal D$

2.3. L.Schwartz {\it distribution}

2.3.  Dirac delta-function

2.3. {\it the n-th derivative} of the Dirac delta-function

2.3. {\it regular}  distribution

2.3. {\it positive} distribution

2.3. {\it product } of a distribution $f$ by an {\it infinitely differentiable}
function $\alpha (x)$

2.3. {\it sum} of distributions $f_1$ and $f_2$

2.3. {\it partial derivative} of distribution

2.3. sequence of distributions $f_n$ {\it converges} to
a distribution $f$

2.3. {\it regularization} of the distribution $f$

2.3. {\it restriction } of distribution $f\in \Cal D' (G)$ to $G_1\sbt G$

2.3. {\it fundamental solution } of $L$ at the point $y$

2.3. {\it  sequence of distributions $f_n$ {\it converges} to
a distribution $f$

2.3. {\it regularization} of the distribution $f$

2.3. {\it restriction } of distribution $f\in \Cal D' (G)$ to $G_1\sbt G$

2.3 {\it spherical} operator

\subheading {2.4}
2.4. {\it harmonic} distribution

2.4. {\it Lipschitz} boundary,{\it Lipschitz} domain

2.4. harmonic measure

2.4. {\it spherical function } of a {\it degree} $\rho$

2.4.  {\it Green potential } of $\mu$ relative to $D$

2.4. {\it Newton} potential

2.4. {\it logarithmic} potential
\subheading {2.5}
2.5. {\it Green capacity} of the compact set $K$ relative to the domain
$D.$

2.5. {\it Wiener } capacity

2.5. {\it external} and {\it inner} capacity of any set $E$

2.5. {\it capacible} set

2.5. {\it logarithmic} capacity

2.5.  {\it irregular} point

2.5. {\it equilibrium} mass distribution

2.5. {\it h -Hausdorff} measure

2.5. {\it Carleson} measure
\subheading {2.6}
2.6. {\it mean value} of $u(x)$ on the sphere $S_{x,r}:=\{y:|y-x|=r\}$

2.6. {\it subharmonic} function

2.6. {\it the least harmonic majorant} of $u$ in $K$

2.6. {\it Riesz} measure of   the subharmonic function $u$
\subheading {2.7}
2.7. {\it precompact} family of functions

2.7. a sequence $f_n$ of locally summable functions
{\it converges in }$L_{loc}$

2.7. {\it quasi-everywhere} convergence

2.7. a sequence of  functions $u_n$ {\it converges} to a function $u$
{\it relative} to $\alpha$- Carleson measure

2.7.  a point $x\in \Bbb R^m$
$(\alpha,\alpha',\epsilon)$-{\it normal} with respect to the measure $\mu$
\subheading {2.8}
2.8. {\it order} of $a(r)$

2.8. {\it type number} of $a(r)$

2.8. $a(r)$ of {\it minimal type}

2.8. $a(r)$ of {\it normal type}

2.8. $a(r)$ of {\it maximal type}

2.8. {\it  convergence exponent} for the sequence $\{r_j\}$

2.8. a {\it proximate order} with respect to order $\rho$

2.8. {\it equivalent} proximate orders

2.8. {\it type number  with respect to a proximate order}

2.8.{\it proper} proximate order

2.8.{\it Nevanlinna characteristic}

2.8.{\it order of $u(x)$ with respect to $T(r).$}

2.8.{\it characteristics} $\rho_M[u],\ \sigma_M [u],\ \sigma_M [u,\rho(r)]$

2.8.{\it convergence exponent} of $\mu$

2.8.{\it upper density} of $\mu.$

2.8.{\it genus} of $\mu$

2.8. {\it N-order of $\mu$}

2.8.{\it N-type of  $\mu$}
\subheading {2.9}
2.9.{\it Gegenbauer } polynomials

2.9.{ \it Chebyshev} polynomials

2.9.{\it primary kernel}

2.9.{\it canonical potential}

2.9.{\it zero distribution}

2.9.{\it canonical Weierstrass product}.
\subheading {3.1}
3.1.{\it limit set } of the function $u(x)$

3.1.{\it limit set } of the mass distribution $\mu$
\subheading {3.2}
3.2.{\it indicator} of growth of $u$

3.2.{\it lower indicator}

3.2.$\rho$-{\it subspherical} function

3.2.$\rho$-{\it trigonometrically convex} ($\rho$-t.c.)

3.2.{\it fundamental relation of indicator.}
\subheading {3.3}
3.3.{\it upper} ( {\it lower) density} of $\mu$

3.3.{\it subadditivity } of $\overline \Delta (E,\bullet)$

3.3.{\it superadditivity} of $\underline \Delta(E,\bullet):$

3.3.{\it monotonic} function of $E\in\Rm.$

3.3.to be {\it dense in}

3.3.{\it angular } densities
\subheading {4.1}
4.1.{\it  dynamical system}

4.1.{\it $(\eps,s)$-chain from $m$ to $m'$}

4.1.{\it chain recurrent} dynamical system

4.1.{\it non-wandering} point

4.1.{\it attractor}

4.1.{\it completely regular growth}

4.1.{\it polygonally connected} set

4.1.{\it periodic} dynamical system
\subheading {4.2}
4.2.{\it partition of unit}
\subheading {4.3}
4.3.{\it pseudo-trajectory}

4.3.{\it asymptotically dynamical } pseudo-trajectory with
 {\it  dynamical asymptotics } $T_\bul$ (a.d.p.t.)

4.3.{\it piecewise continuous} pseudo-trajectory $w(\bul|\bul)$

4.3. $\om${\it --dense} pseudo-trajectory
\subheading {4.4}
4.4.{\it subharmonic curve }
\subheading {5.1}
5.1{\it is an entire function of order $\r$ and normal type with respect to
proximate order} $\r(r)$

5.1.{\it entire} function with prescribed limit set

5.1.{\it  meromorphic function of order  $\r$ and normal type with respect to
a proximate order} $\r(r)$
\subheading {5.2}
5.2.{\it relative} Carleson $\a$ -measure
\subheading {5.3}
5.3.{\it lower indicator} of entire function
\subheading {5.4}
5.4.{\it maximal interval of $\r$-trigonometricity}

5.4.{\it strictly } $\r$-t.c.f.

5.4. {\it concordant} $h$ and $g$
\subheading {5.5}
5.5.{\it upper density of zeros} of entire function

5.5. $(\Cal X)-integral$ with respect to a {\it nonnegative measure} $\delta.$
\subheading {5.6}
5.6.{\it completely regular growth } function

5.6.{\it regular zero distribution}

5.6.regular zero distribution with integer $\r$ 

5.6.{\it completely regular growth functions along curves of regular rotation}

5.6.{\it curve of regular rotation}
\subheading {5.7}
5.7.{\it growth characteristic}

5.7.{\it continuity},{\it positive homogeneity}

5.7.{\it asymptotic characteristics of growth}

5.7.{\it total}{\it family of growth characteristics}

5.7.{\it non -rarefied} set

5.7.{\it rarefied} set

5.7.{\it thinly closed} set

5.7.{\it independent}family of characteristic
\subheading {6.1}
6.1.{\it ideally  complementing} $H$-multiplicator

6.1.entire function is of {\it minimal type} with respect to a proximate
order $\r(r),\r(r)\ri \r$

6.1.{\it limit set of indicators}

6.1.{\it the maximum principle for $U[\r]$ is valid in the domain} $G$

\subheading {6.2}
6.2.{\it automorphic }

6.2.{\it connected on spirals}

6.2.{\it spectrum}

6.2.{\it strictly monotonic}

6.2.{\it minimal} $v\in U[\r]$ \subheading {6.3} 6.3.function of
{\it exponential type}

6.3.{\it completeness}

6.3.{\it maximality}

6.3.{\it extremal overcompleteness}

6.3.{\it maximal} domain of completeness

6.3.{\it extremely overcomplete}  system  $\exp \Lm$

6.3.{\it trigonometrically convex function} (t.c.f)

6.3.{\it conjugate indicator diagram}

6.3{\it regular set}

6.3.$G_\Lm$ is {\it enclosed}  in $G$

6.3.{\it enclosed with  sliding}

6.3.{\it enclosed hardly}

6.3.{\it enclosed freely}

6.3.{\it indicator limit set}

6.3.{\it indicator set}

6.3.{\it  zero with tangency}

6.3.$\Lm$ is {\it periodic}

6.3.$w\in U[1]$ is {\it minimal}

6.3.$U\sbt U[1]$ is {\it minimal}

\end